\documentclass{article}

\usepackage[final,nonatbib]{nips_2016}

\usepackage[utf8]{inputenc} 
\usepackage[T1]{fontenc}    
\usepackage{hyperref}       
\usepackage{url}            
\usepackage{booktabs}       
\usepackage{amsfonts}       
\usepackage{nicefrac}       
\usepackage{microtype}      

\usepackage{graphicx}
\usepackage{epstopdf}
\usepackage{graphics}
\usepackage{algorithm}
\usepackage{algpseudocode}
\usepackage{amsfonts}
\usepackage{amsmath}
\usepackage{amssymb}
\usepackage{amsthm}
\usepackage{xspace}
\usepackage{float}

\usepackage{todonotes} 

\newtheorem{theorem}{Theorem}[section]
\newtheorem{lemma}[theorem]{Lemma}

\usepackage{lipsum}
\graphicspath{{./}{./}}

\usepackage{enumitem}
\newcommand{\multi}{\textit{multi-batch}\xspace}

\newcommand{\defeq}{\stackrel{\rm def}{=}}

\newtheorem{innercustomthm}{Theorem}
\newenvironment{customthm}[1]
  {\renewcommand\theinnercustomthm{#1}\innercustomthm}
  {\endinnercustomthm}

\newtheorem{innercustomlemma}{Lemma}
\newenvironment{customlemma}[1]
  {\renewcommand\theinnercustomlemma{#1}\innercustomlemma}
  {\endinnercustomlemma}

\title{A Multi-Batch L-BFGS Method for Machine Learning}

\author{
  Albert S. Berahas\\
  Northwestern University\\
  Evanston, IL\\
  \texttt{albertberahas@u.northwestern.edu} \\
\AND
Jorge Nocedal \\
  Northwestern University\\
  Evanston, IL\\
  \texttt{j-nocedal@northwestern.edu} \\
\And
 Martin Tak\'{a}\v{c} \\
Lehigh University \\
Bethlehem, PA \\
\texttt{takac.mt@gmail.com} \\
}

\begin{document}

\maketitle
 \begin{abstract}
The question of how to parallelize the stochastic gradient descent (SGD) method has received much attention in the literature. In this paper, we focus instead on batch methods that use a sizeable fraction of the training set at each iteration to facilitate parallelism, and that employ second-order information. In order to improve the learning process, we follow a \emph{multi-batch} approach in which the batch changes at each iteration. This can cause  difficulties because 
L-BFGS employs gradient differences to update the Hessian approximations, and when these gradients are computed using different data points the process can be unstable. This paper shows how to perform stable quasi-Newton updating in the multi-batch setting, illustrates the behavior of the algorithm in a distributed computing platform, and studies its convergence properties for both the convex and nonconvex cases.

 \end{abstract}

\section{Introduction}
\label{sec:intro}

 It is common in machine learning to encounter optimization problems involving millions of parameters and very large datasets.  To deal with the computational demands imposed by such applications, high performance implementations of stochastic gradient and batch quasi-Newton methods have been developed \cite{agarwal2014reliable,dean2012large,chen2014large}. In this paper we study a  batch approach based on the L-BFGS method \cite{mybook}  that strives to reach the right balance between efficient learning and productive parallelism.

In supervised learning, one seeks to minimize  
 empirical risk,
 \begin{align*} 
 {F(w) :=  \frac{1}{n} \sum  _{i=1}^{n}f(w;x^{i},y^{i}) \defeq \frac{1}{n} \sum _{i=1}^{n}f_i(w)},
\end{align*}
where $ (x^i, y^i)_{i=1}^n$  denote the training examples and $f(\cdot;x,y) : \mathbb{R}^d \rightarrow \mathbb{R}$ is the composition of a prediction function (parametrized by $w$) and a loss function. The training problem consists of finding an optimal choice of the parameters $w \in \mathbb{R}^d$ with respect to $F$, i.e.,
{
\begin{align}  \label{eq:obj}
     \min _{w\in\mathbb{R}^d}F(w) =\frac{1}{n} \sum_{i=1}^{n}f_i(w).
\end{align}
}
At present, the preferred optimization method  is the stochastic gradient descent (SGD) method \cite{RobMon51, bottou-lecun-2004}, and its variants \cite{johnson2013accelerating,Schmidt2016,
defazio2014saga}, which are implemented either in an asynchronous manner (e.g. when using a parameter server in a distributed setting) or following a synchronous mini-batch approach that exploits parallelism in the gradient evaluation \cite{bertsekas1989parallel,recht2011hogwild,Goodfellow-et-al-2016-Book}.
A drawback of the asynchronous approach is that it cannot use large batches, as this would cause updates to become too dense
and  compromise the stability and scalability  of the method
\cite{mania2015perturbed,recht2011hogwild}.
 As a result, the algorithm spends more time in communication as compared to computation. On the other hand, using a synchronous mini-batch approach one can achieve a near-linear decrease in the number of SGD iterations as the mini-batch size is increased, up to a certain point after which the increase in computation is not offset by the faster convergence \cite{takavc2013mini}.
 
An alternative to  SGD  is  a batch method, such as L-BFGS, which is able to reach high training accuracy and allows one to perform more computation per node, so as to achieve a better balance with communication costs \cite{zhang2015disco}.  Batch methods are, however,  not as efficient learning algorithms as SGD in a sequential setting \cite{bousquet2008tradeoffs}. To benefit from the strength of both methods some high performance systems employ SGD at the start and later switch to a batch method \cite{agarwal2014reliable}. 

{\bf Multi-Batch Method.} In this paper, we follow a different approach consisting of a single method that  selects a \emph{sizeable} subset (batch) of the training data to compute a step, and changes this batch at each iteration to improve the learning abilities of the method. We call this a \emph{multi-batch} approach to differentiate it from the mini-batch approach used in conjunction with SGD, which employs a very small subset of the training data. When using large batches it is natural to employ a quasi-Newton method, as incorporating second-order information imposes little computational overhead and improves the stability and speed of the method. We focus here on the L-BFGS method, which employs  gradient information to update an estimate of the Hessian and computes a step in $O(d)$ flops, where $d$ is the number of variables.  The multi-batch approach can, however, cause difficulties to L-BFGS because this method employs gradient differences  to update Hessian approximations. When the gradients used in these differences are based on different data points, the updating procedure can be unstable.  Similar difficulties arise in a parallel implementation of the standard L-BFGS method, if some of the computational nodes devoted to the evaluation of the function and gradient are unable to return results on time --- as this again amounts to using different data points to evaluate the function and gradient at the beginning and the end of the iteration.
 The goal of this paper is to show that stable quasi-Newton updating can be achieved in both settings without incurring extra computational cost, or  special synchronization. The key is to perform quasi-Newton updating based on the overlap between consecutive batches. The only restriction is that this overlap should not be too small, something that can be achieved
in most situations.

{\bf Contributions.} 
We describe a novel implementation of the batch L-BFGS method that is robust  in the absence of sample consistency; i.e., when different samples are used to evaluate the objective function and its gradient at consecutive iterations.
The numerical experiments show that the method proposed in this paper --- which we call the \multi L-BFGS method --- achieves a good balance between computation and communication costs. 
We also analyze the convergence properties of the new method (using a fixed step length strategy) on both convex and nonconvex problems.
%




\section{The Multi-Batch Quasi-Newton Method}
\label{sec:method}

In a pure batch approach, one applies a gradient based method, such as L-BFGS \cite{mybook}, to the deterministic optimization problem
\eqref{eq:obj}. When the number $n$ of training examples is large, it is natural to parallelize the evaluation of $F$ and $\nabla F$ by assigning the computation of the component functions $f_i$ to different processors. If this is done on a distributed platform, it is possible for some of the computational nodes 
to be  slower than the rest. In this case, the contribution of the slow (or unresponsive) computational nodes could  be ignored given the stochastic nature of the objective function. This leads, however, to an inconsistency in the objective function and gradient at the beginning and at the end of the iteration, which can be detrimental to quasi-Newton methods
. Thus, we seek to find a \emph{fault-tolerant} variant of the batch L-BFGS method that is capable of dealing with slow or unresponsive computational nodes.

A similar challenge arises in a  \multi implementation of the L-BFGS method in which the entire training set $ T= \{ (x^i, y^i)_{i=1}^n\}$ is not employed at every iteration, but rather,  a subset of the data is used to compute the gradient. Specifically, we consider a method in which the dataset is randomly divided into  a number of batches --- say 10, 50, or 100 --- and  the minimization is performed with respect to a different batch at every iteration.  At the $k$-th iteration, the algorithm chooses a batch $S_k \subset \{1, \ldots, n\}$, computes 
\begin{align}   \label{eq:batch_fg}
{F}^{S_{k}}(w_k)=\frac{1}{\left|S_{k}\right|} \sum_{i\in S_{k}}f_i\left(w_{k}\right), \qquad \nabla{F}^{S_{k}}(w_k) = {g}_{k}^{S_{k}} = \frac{1}{\left|S_{k}\right|} \sum_{i\in S_{k}}\nabla f_i\left(w_{k}\right) ,
\end{align}
and takes a step along the direction $- H_k g_k^{S_k}$, where $H_k$ is an approximation to $\nabla^2 F(w_k)^{-1}$. Allowing the sample $S_k$ to change freely at every iteration gives this approach  flexibility of implementation and is beneficial to the learning process, as we show in Section~\ref{sec:num_res}. (We refer to $S_k$ as the sample of training points, even though $S_k$ only indexes those points.)

The case of unresponsive computational nodes and the  {multi-batch} method are similar. The main difference is that node failures create unpredictable changes to the samples $S_k$, whereas a multi-batch method has control over sample generation. In either case, the algorithm employs a stochastic approximation to the gradient and can no longer be considered deterministic.  We must, however, distinguish our setting from that of the classical SGD method, which employs small mini-batches  and noisy gradient approximations. Our algorithm operates with much larger batches so that distributing the function evaluation is beneficial and the compute time of $g_k^{S_k}$ is not overwhelmed by communication costs. This gives rise to gradients with relatively small variance and justifies the use of a second-order method such as L-BFGS.

{\bf Robust Quasi-Newton Updating.}
The  difficulties created by the use of a different sample $S_k$ at each iteration can be circumvented if consecutive samples $S_{k}$ and $S_{k+1}$ overlap,  so that
$      O_k= S_{k} \cap S_{k+1} \neq \emptyset.  
$  One can then perform stable quasi-Newton updating  by computing  gradient differences based on this overlap, i.e.,  by defining 
\begin{align}   \label{pairs}
       y_{k+1}=g_{k+1}^{O_{k}}-g_{k}^{O_{k}}, \qquad s_{k+1} = w_{k+1}-w_k,
\end{align}
in the notation given in \eqref{eq:batch_fg}. The correction pair $(y_k, s_k)$ can then be used in the BFGS update. When the overlap set $O_k$  is not too small, $y_k$ is a useful approximation of the curvature of the objective function $F$ along the most recent displacement, and will lead to a productive quasi-Newton step. This observation is based on an important  property of Newton-like methods, namely that there is much more freedom in choosing a Hessian approximation than in computing the gradient \cite{byrd2011use,bollapragada2016exact}. Thus, a smaller sample $O_k$ can be employed for updating the inverse Hessian approximation $H_k$ than for computing the batch gradient $g_k^{S_k}$ in the search direction $- H_k g_k^{S_k}$. In summary, by ensuring that unresponsive nodes do not constitute the vast majority of all working nodes in a {fault-tolerant} parallel implementation, or by exerting a small degree of control over the creation of the samples $S_k$ in the multi-batch method, one can design a robust method that naturally builds upon the fundamental properties of BFGS updating. 

We should mention in passing that a commonly used strategy for ensuring stability of quasi-Newton updating in machine learning is to enforce gradient consistency \cite{schraudolph2007stochastic
}, i.e., to use the same sample $S_k$ to compute gradient evaluations at the beginning and the end of the iteration. Another popular remedy is to use the same batch $S_k$ for multiple iterations \cite{ngiam2011optimization}, alleviating the gradient inconsistency problem at the price of slower convergence.  In this paper, we assume that achieving such \emph{sample consistency is not possible} (in the fault-tolerant case) or \emph{desirable} (in a multi-batch framework), and wish to design a new variant of L-BFGS that imposes minimal restrictions in the sample changes. 

\subsection{Specification of the Method}

At  the $k$-th iteration, the multi-batch BFGS algorithm chooses a set $S_k \subset \{1, \ldots, n\}$ and computes a new iterate
\begin{align}  \label{eq:update}
w_{k+1}=w_k-\alpha_{k}H_{k} g_{k}^{S_{k}} ,
\end{align}
where $\alpha_{k}$ is the step length, $g_k^{S_k}$ is the batch gradient \eqref{eq:batch_fg} and $H_{k}$ is the inverse BFGS
Hessian matrix approximation that is updated at every iteration by means
of the formula
\begin{align*}   
H_{k+1}=V_{k}^{T}H_{k}V_{k}+\rho_{k}s_{k}s_{k}^{T}, \qquad \rho_{k}=\tfrac{1}{y_{k}^{T}s_{k}}, \qquad V_{k}=I-\rho_{k}y_{k}s_{k}^{T} .
\end{align*}
 To compute the correction vectors $(s_k, y_k)$, we determine the overlap set $O_k = S_{k} \cap S_{k+1}$ consisting of the samples that are common at the $k$-th and $k+1$-st iterations. We define  
 \begin{align*}  
         {F}^{O_{k}}(w_k)=\frac{1}{\left|O_{k}\right|} \sum_{i\in O_{k}}f_i\left(w_k\right), 
         \qquad \nabla F^{O_k}(w_k)={g}_{k}^{O_{k}}=\frac{1}{\left|O_{k}\right|} \sum_{i\in O_{k}}\nabla f_i\left(w_k\right),
\end{align*}
and compute the correction vectors as in \eqref{pairs}. 
In this paper we assume that $\alpha_k$ is  constant.

In the limited memory version, the matrix $H_k$ is defined at each iteration as the result of applying $m$ BFGS updates  to a multiple of the identity matrix, using a set of $m$ correction pairs $\{s_i, y_i\}$ kept in storage. The memory parameter $m$ is typically in the range  2 to 20.  When computing the matrix-vector product in \eqref{eq:update} it is not necessary to form that matrix $H_k$ since one can obtain this product via the two-loop recursion \cite{mybook}, using the $m$ most recent correction pairs $\{s_i, y_i\}$. After the step has been computed, the oldest pair $(s_j, y_j)$ is discarded and the new curvature pair is stored. 

A pseudo-code of the proposed method is given below, and depends on several parameters.  The parameter $r$ denotes the fraction of samples in the dataset used to define the gradient, i.e., $r = \frac{\left| S\right|}{n}$. The parameter $o$ denotes the length of overlap between consecutive samples, and is defined  as a fraction of the number of samples in a given batch $S$, i.e., $o = \frac{\left| O\right|}{\left| S\right|}$.

\begin{algorithm} 
\caption{Multi-Batch L-BFGS}
  \label{alg:multi}
 {\bf Input:} $w_{0}$ (initial iterate),  $ T= \{ (x^i, y^i)$, for $i=1, \ldots, n\}$ (training set),
$m$ (memory parameter), $r$ (batch, fraction of $n$), $o$ (overlap, fraction of batch), $k\leftarrow0$ (iteration counter). 

  \begin{algorithmic}[1]
  \State Create initial batch $S_{0}$ \Comment{As shown in Firgure~\ref{fig:sample_creation}}
  \For {$k=0,1,2,...$}
\State Calculate the search direction $p_{k}=-H_{k}g_{k}^{S_{k}}$ \Comment{Using L-BFGS formula}
\State Choose the step length $\alpha_{k} >0$
\State Compute  $w_{k+1}=w_k+\alpha_{k}p_{k}$ 
\State Create the next batch $S_{k+1}$
\State Compute the curvature pairs $s_{k+1}=w_{k+1}-w_k$ and
$y_{k+1}=g_{k+1}^{O_{k}}-g_{k}^{O_{k}}$ 
\State Replace the oldest  pair $(s_i, y_i)$  by $s_{k+1}, y_{k+1}$
\EndFor
  \end{algorithmic}
\end{algorithm}

\subsection{Sample Generation}
\label{sec:sampling}
We now discuss how the sample $S_{k+1}$  is created at each iteration (Line 8 in Algorithm \ref{alg:multi}). 

\begin{figure}[h!]
\begin{centering}
\includegraphics[scale=0.75]{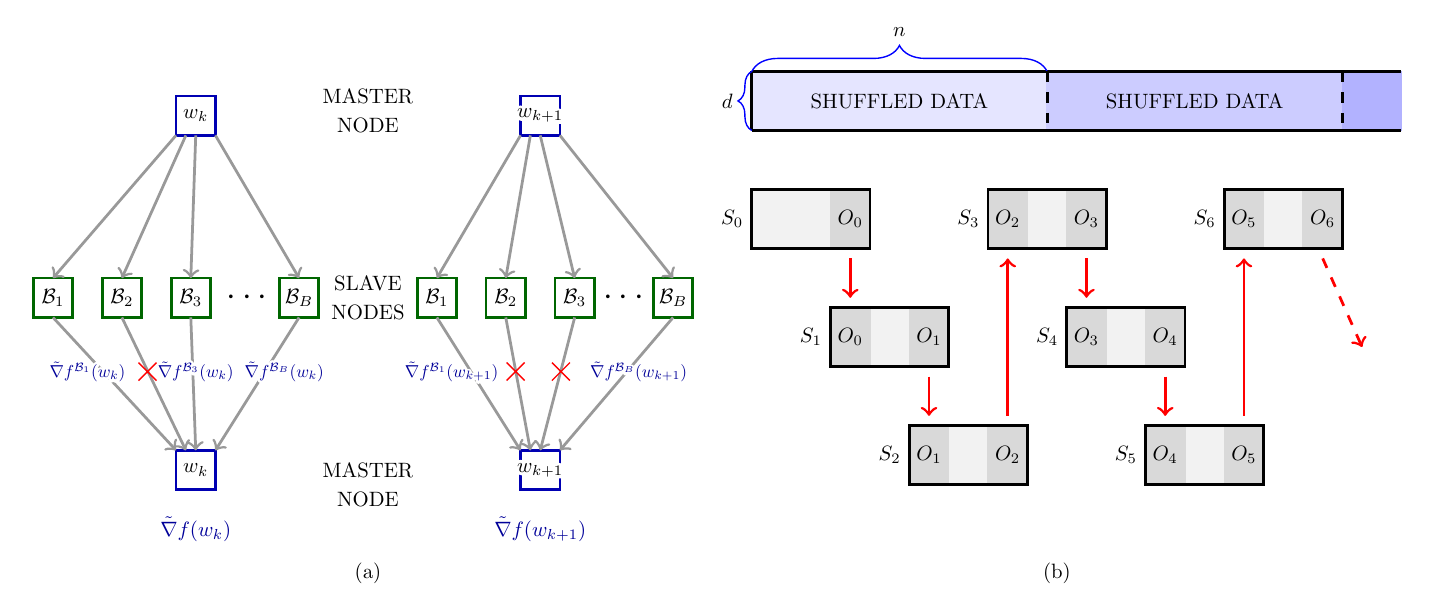}
\par\end{centering}
\caption{Sample and Overlap formation.}
\label{fig:sample_creation}
\end{figure}

{\bf Distributed Computing with Faults.}
Consider a distributed implementation in which  slave nodes read the current iterate $w_k$ from the master node, compute a local gradient on a subset of the dataset, and send it back to the master node for aggregation in the calculation \eqref{eq:batch_fg}. 
 Given a time (computational) budget, it is possible for some nodes to fail to return a result. The schematic in Figure~\ref{fig:sample_creation}a illustrates the gradient calculation across two iterations, $k$ and $k+1$, in the presence of faults.  Here $\mathcal{B}_i$,  $i=1,...,B$ denote the batches of data that each slave node $i$ receives (where $T = \cup_i \mathcal{B}_i$), and 
$\tilde{\nabla}f(w)$ is the gradient calculation using all nodes that responded within the preallocated time. 

Let $\mathcal{J}_k\subset \{1,2,...,B\}$ and $\mathcal{J}_{k+1}\subset \{1,2,...,B\}$ be the set of indices of all nodes that returned a gradient at the $k$-th and $k+1$-st iterations, respectively. Using this notation $S_k = \cup_{j\in \mathcal{J}_k} \mathcal {B}_j$ and $S_{k+1} = \cup_{j\in \mathcal{J}_{k+1}} \mathcal {B}_j$, and we define $O_k = \cup_{j \in \mathcal{J}_k\cap \mathcal{J}_{k+1}} \mathcal {B}_j$.
The simplest implementation in this setting preallocates the data on each compute node, requiring minimal data communication, i.e., only one data transfer. In this case the samples $S_k$ will be independent if node failures occur randomly. On the other hand, if the same set of nodes fail, then sample creation will be biased, which is harmful both in theory and practice.  One way to ensure independent sampling 
 is to shuffle and redistribute the data to all nodes after a certain number of iterations.

{\bf Multi-batch Sampling.}
We propose two strategies for the \textit{multi-batch} setting.

Figure~\ref{fig:sample_creation}b illustrates the sample creation process in the first strategy. 
The dataset is shuffled and batches are generated by collecting subsets of the training set, in order. Every set (except $S_0$) is of the form $S_k= \{ O_{k-1}, N_k, O_k\}$, where $O_{k-1}$ and $O_k$ are the overlapping samples with batches $S_{k-1}$ and $S_{k+1}$ respectively, and $N_k$ are the samples that are unique to batch $S_k$.
After each pass through the dataset, the samples are reshuffled, and the procedure described above is repeated. In our implementation samples are drawn without replacement, guaranteeing that after every pass (epoch) all samples are used. This strategy has the advantage that it requires no extra computation in the evaluation of $g_k^{O_k}$ and $g_{k+1}^{O_k}$, but the samples $\{S_k\}$ are not independent.

The second sampling  strategy is simpler and requires less control. At every iteration $k$, a batch $S_k$ is created by randomly selecting $\left| S_k \right|$ elements from $\{1,\ldots n\}$. The overlapping set $O_k$ is then formed by randomly selecting $\left| O_k \right|$ elements from $S_k$ (subsampling). This strategy is slightly more expensive since $g_{k+1}^{O_k}$ requires extra computation, but if the overlap is small this cost is not significant.

\section{Convergence Analysis}
\label{sec:conv}

In this section, we analyze the convergence properties of the multi-batch L-BFGS method (Algorithm \ref{alg:multi}) when applied to the minimization of \emph{strongly convex}  and \emph{nonconvex} objective functions, using a fixed step length strategy. We assume that the goal is to minimize the empirical risk $F$ given in \eqref{eq:obj}, but note that a  similar analysis could be used to study  the minimization of the expected risk.

\subsection{Strongly Convex case}

Due to the 
stochastic      nature of the multi-batch approach, every iteration of Algorithm~\ref{alg:multi} employs a gradient that contains errors that do not converge to zero. Therefore, by using a fixed step length strategy one cannot establish convergence to the optimal solution $w^{\star}$, but only convergence to a neighborhood of $w^{\star}$ \cite{nedic2001convergence}. Nevertheless, this result is of interest as it reflects the common practice of using a fixed step length and decreasing it only if the desired testing error has not been achieved. It also illustrates the tradeoffs that arise between the size of the batch and the step length.

In our analysis, we make the following 
assumptions about the objective function and the algorithm.

\textbf{Assumptions A.}
\emph{ 
\begin{enumerate} [leftmargin=0.5cm,topsep=0pt,itemsep=0ex,,partopsep=0ex,parsep=0ex]
\item $F$ is twice continuously differentiable.
\item There exist positive constants $\hat{\lambda}$ and $\hat{\Lambda}$ such that
$\hat{\lambda} I \preceq \nabla^2F^O(w) \preceq \hat{\Lambda} I$
for all $w \in \mathbb{R}^d$ and all sets $O \subset  \{1,2,\ldots,n\}$.
\item There is a constant $\gamma$ such that $\mathbb{E}_{S}\left[ \| \nabla  F^{S}(w) \| \right]^2 \leq \gamma^2
$ for all $w \in \mathbb{R}^d$ and all sets $S\subset  \{1,2,\ldots,n\}$.
\item The samples $S$ are drawn independently and $\nabla F^{S}(w)$ is an unbiased estimator of the true gradient $\nabla F(w)$ for all $w \in \mathbb{R}^d$, i.e.,
$
\mathbb{E}_{S}[ \nabla F^{S}(w)] = \nabla F(w).$
\end{enumerate}
}

Note that Assumption $A.2$ implies that the entire Hessian $\nabla^2F(w)$ also satisfies 
\begin{align*}  
 \lambda I \preceq \nabla^2F(w) \preceq  \Lambda I,  \quad \forall w \in \mathbb{R}^d,
\end{align*}
for some constants $ \lambda,  \Lambda>0$. Assuming that every sub-sampled function $F^O(w)$ is strongly convex is not unreasonable  as a regularization term is commonly added in practice when that is not the case.
 
We begin by showing that the inverse Hessian approximations $H_k$ generated by the multi-batch L-BFGS method have eigenvalues that are uniformly bounded above and away from zero. The proof technique used is an adaptation of that in \cite{Sammy_SQN}.

\begin{lemma}	\label{lemma1}
If Assumptions A.1-A.2 above hold, there exist constants $0<\mu_1\leq \mu_2$ such that the Hessian approximations $\{H_k\}$ generated by Algorithm~\ref{alg:multi} satisfy
\begin{align*}    
\mu_1 I \preceq H_k \preceq \mu_2 I,\qquad \text{for } k=0,1,2,\dots 
\end{align*}
\end{lemma}

Utilizing Lemma \ref{lemma1}, we show that the multi-batch L-BFGS method with a constant step length converges to a neighborhood of the optimal solution.

\begin{theorem}
\label{thm:const}
Suppose that Assumptions A.1-A.4  hold and let $F^{\star} = F(w^{\star})$, where $w^{\star}$ is the minimizer of $F$. Let $\{w_k\}$ be the iterates generated by Algorithm~\ref{alg:multi} with 
$\alpha_k = \alpha \in  (0,\frac{1}{2\mu_1 \lambda})$,
starting from $w_0$.
Then for all $k\geq 0$,
\begin{align*}   
\mathbb{E} [ F(w_k) - F^{\star} ] & \leq   ( 1-2\alpha \mu_1 \lambda  )^k  [ F(w_0) - F^{\star}  ] +  [ 1-(1-\alpha\mu_1 \lambda)^k ]\frac{\alpha \mu_2^2 \gamma ^2 \Lambda}{4 \mu_1 \lambda}  
\\    &
      \xrightarrow[]{k\rightarrow \infty} \frac{\alpha \mu_2^2 \gamma ^2 \Lambda}{4 \mu_1 \lambda}. 
\end{align*}
\end{theorem}


The bound provided by this theorem  has two components: (i) a term decaying linearly to zero, and (ii) a  term identifying the neighborhood of convergence. Note that a larger step length yields a more favorable constant in the linearly decaying term, at the cost of an increase in the size of the neighborhood of convergence.  We will consider again these tradeoffs in Section~\ref{sec:num_res}, where we also note that larger batches increase the opportunities for parallelism and improve the limiting accuracy in the solution, but slow down the learning abilities of the algorithm.

 One can establish convergence of the multi-batch L-BFGS method to the optimal solution $w^\star$ by employing a  sequence of step lengths $\{ \alpha_k \}$ that converge to zero according to the schedule  proposed by Robbins and Monro \cite{RobMon51}. However, that provides only a sublinear rate of convergence, which is of little interest in our context where large batches are employed and  some type of linear convergence is expected.  In this light, Theorem~\ref{thm:const} is more relevant to practice.

\subsection{Nonconvex case}

The BFGS method is known to fail on noconvex problems \cite{mascarenhas2004bfgs,dai2002convergence}. Even for L-BFGS, which makes only a finite number of updates at each iteration, one cannot guarantee that the Hessian approximations have eigenvalues that are uniformly bounded above and away from zero. To establish convergence of the BFGS method in the nonconvex case \emph{cautious} updating procedures have been proposed \cite{li2001global}. Here we employ a cautious strategy that is well suited to our particular algorithm; we skip the update, i.e., set $H_{k+1} = H_k$,  if the curvature condition  
\begin{align} \label{curv}
	y_k^Ts_k \geq {\epsilon} \| s_k \|^2
\end{align}
 is not satisfied, where $\epsilon>0$ is a predetermined constant.
Using said mechanism 
we show that the eigenvalues of the Hessian matrix approximations generated by the multi-batch L-BFGS method are bounded above and away from zero (Lemma \ref{lemma2}). The analysis presented in this section is based on the following assumptions.

\textbf{Assumptions B.}
\emph{
\begin{enumerate}[leftmargin=0.5cm,topsep=0pt,itemsep=0ex,,partopsep=0ex,parsep=0ex]
\item $F$ is twice continuously differentiable.
\item The gradients of $F$ are $\Lambda$-Lipschitz continuous, and the gradients of $F^{O}$ are $\Lambda_{O}$-Lipschitz continuous for all $w \in \mathbb{R}^d$ and all sets $O \subset  \{1,2,\ldots,n\}$. 
\item The function $ F(w)$ is bounded below by a scalar $\widehat F$ .
\item There exist constants $\gamma \geq 0$ and $\eta>0$ such that 
$\mathbb{E}_{S}\left[ \| \nabla  F^{S}(w) \| \right]^2 \leq \gamma^2 + \eta \| \nabla F(w)\|^2$
for all $w \in \mathbb{R}^d$ and all sets $S\subset  \{1,2,\ldots,n\}$. 
\item The samples $S$ are drawn independently and $\nabla F^{S}(w)$ is an unbiased estimator of the true gradient $\nabla F(w)$ for all $w \in \mathbb{R}^d$, i.e.,
$\mathbb{E} [ \nabla F^{S}(w) ] = \nabla F(w).
$
\end{enumerate}
}


\begin{lemma}		\label{lemma2}
Suppose that Assumptions B.1-B.2  hold and let $\epsilon >0$ be given. Let $\{H_k \}$ be the Hessian approximations generated by Algorithm~\ref{alg:multi}, with the modification that $H_{k+1} = H_k$ whenever \eqref{curv} is not satisfied. Then,
 there exist constants $0<\mu_1\leq \mu_2$ such that 
\begin{align*}    
\mu_1 I \preceq H_k \preceq \mu_2 I,\qquad \text{for } k=0,1,2,\dots 
\end{align*}
\end{lemma}
We can now follow the analysis in \cite[Chapter ~4]{bottou2016optimization} to establish the following result about the behavior of the gradient norm for the  multi-batch L-BFGS method with a cautious update strategy. 

\begin{theorem}
Suppose that Assumptions B.1-B.5 above hold, and let $\epsilon >0$ be given. Let $\{w_k\}$ be the iterates generated by Algorithm~\ref{alg:multi}, with 
$\alpha_k = \alpha \in  (0,\frac{\mu_1}{\mu_2^2\eta \Lambda} )$, 
starting from $w_0$, and  with the modification that $H_{k+1} = H_k$ whenever \eqref{curv} is not satisfied.
Then,
\begin{align*}	
\mathbb{E} \Big[\frac{1}{L}\sum_{k=0}^{L-1} \| \nabla F(w_k) \|^2 \Big] & \leq \frac{\alpha \mu_2^2 \gamma^2 \Lambda}{ \mu_1 } + \frac{2[ F(w_0) - \widehat{F} ]}{\alpha \mu_1 L }\\
& \xrightarrow[]{L\rightarrow \infty}\frac{\alpha \mu_2^2 \gamma^2 \Lambda}{ \mu_1 }.
\end{align*}
\end{theorem}
\vskip-5pt
This result bounds the average norm of the gradient of $F$  after the first $L-1$ iterations, and shows that the iterates spend increasingly more time in regions where the objective function has a small gradient.



\section{Numerical Results }   
\label{sec:num_res}
\vskip-5pt
 
\begin{figure}[ht]
\centering
\includegraphics[width=4.5cm]{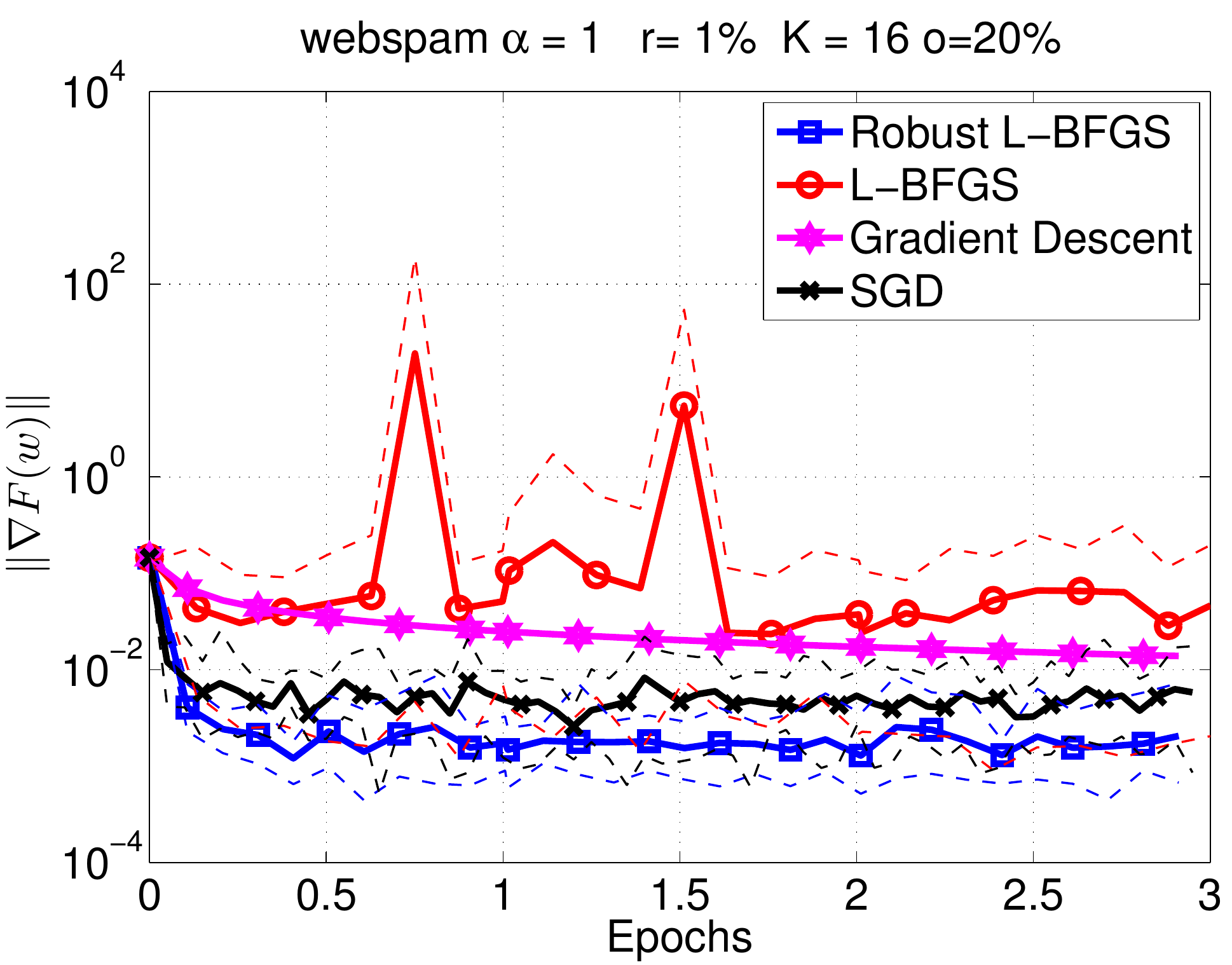}
\includegraphics[width=4.5cm]{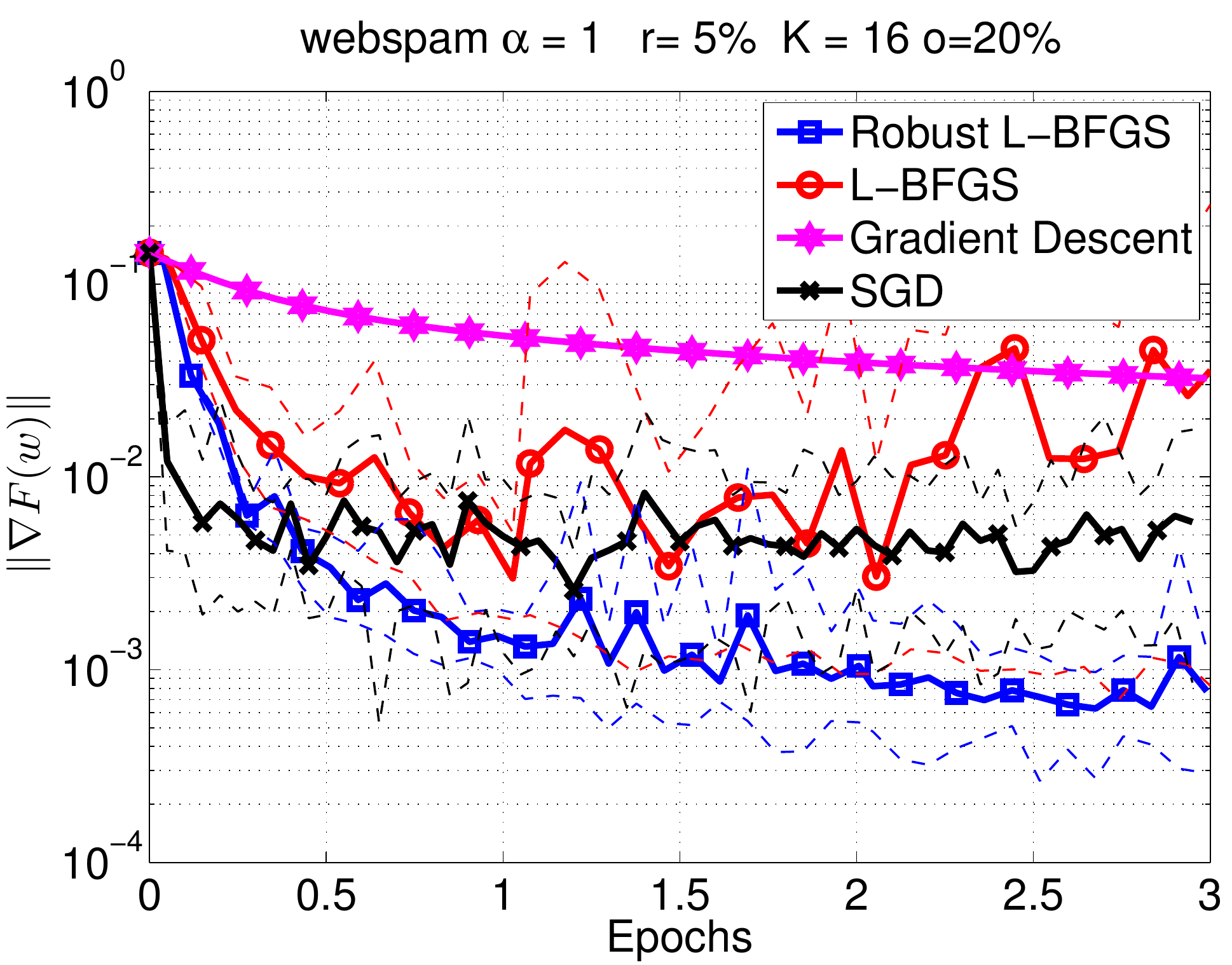}
\includegraphics[width=4.5cm]{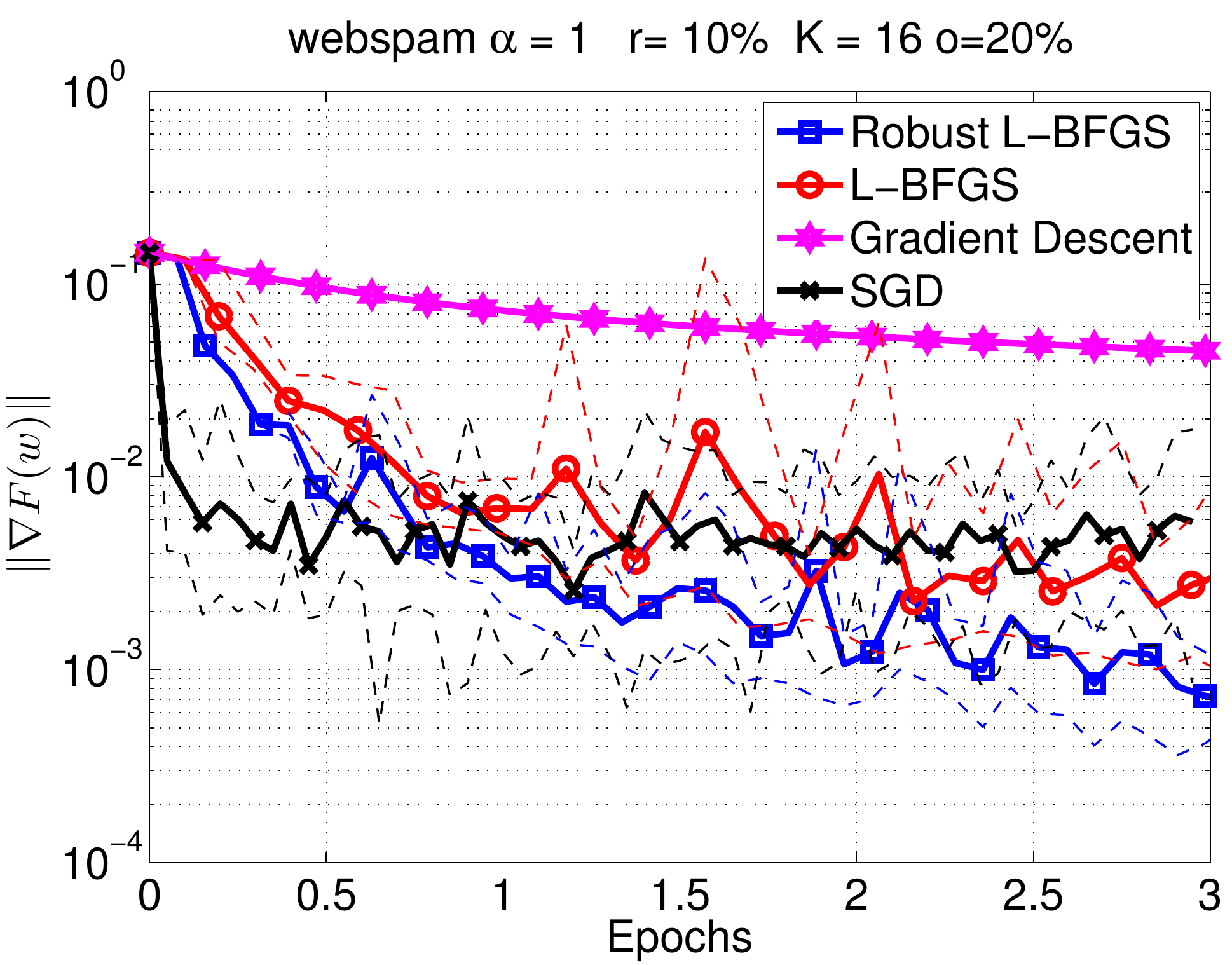}

\hrule 

\includegraphics[width=4.5cm]{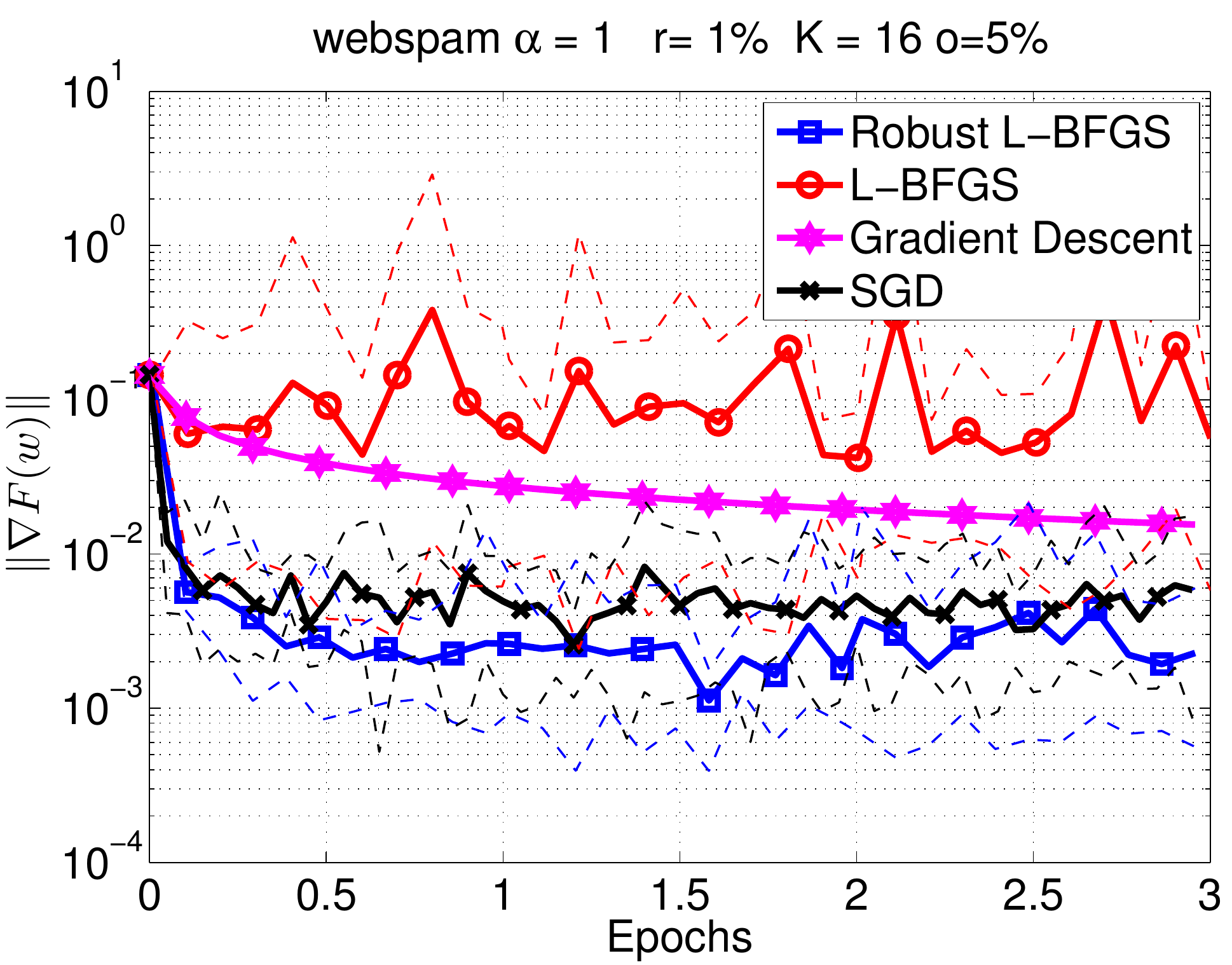}
\includegraphics[width=4.5cm]{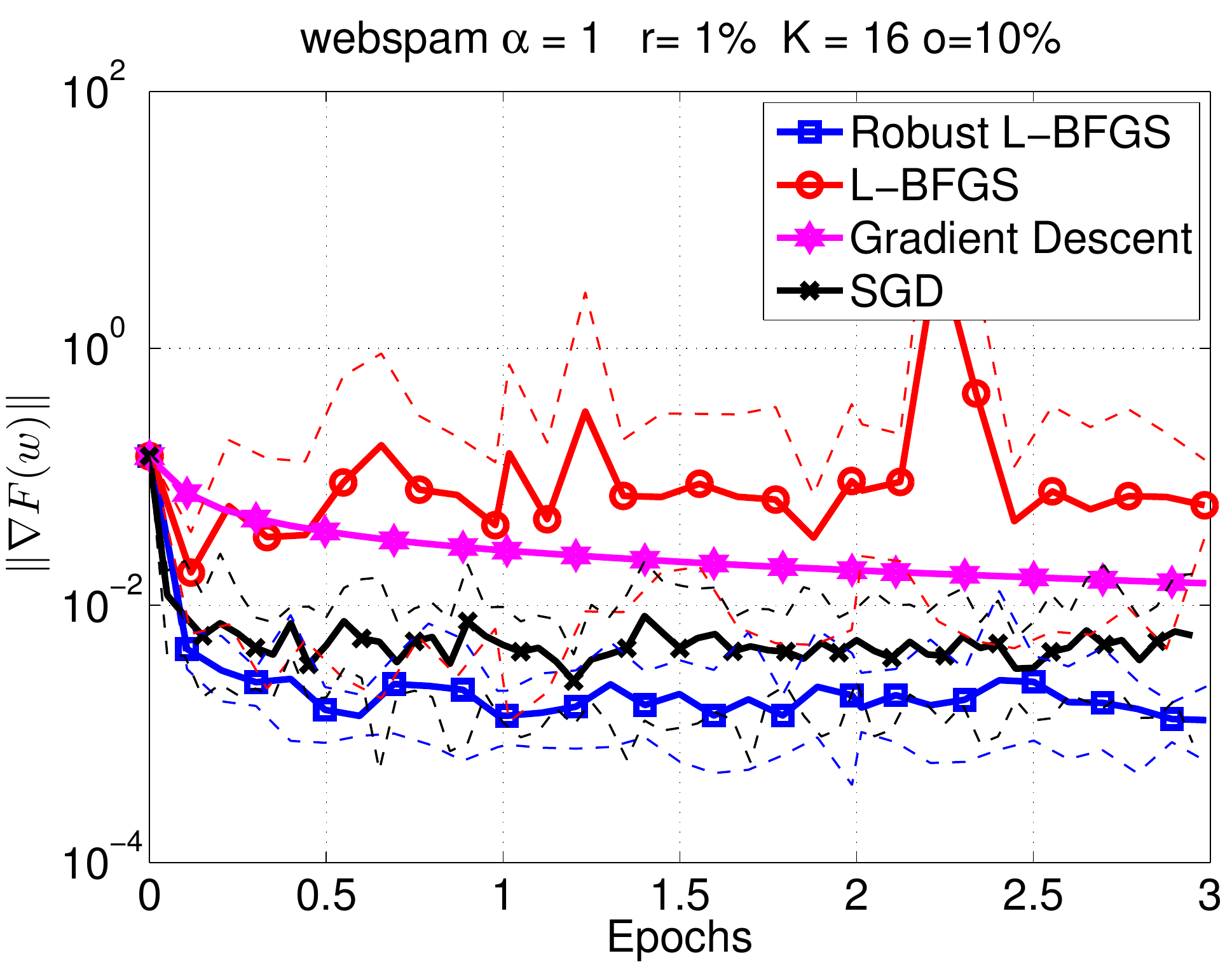}
\includegraphics[width=4.5cm]{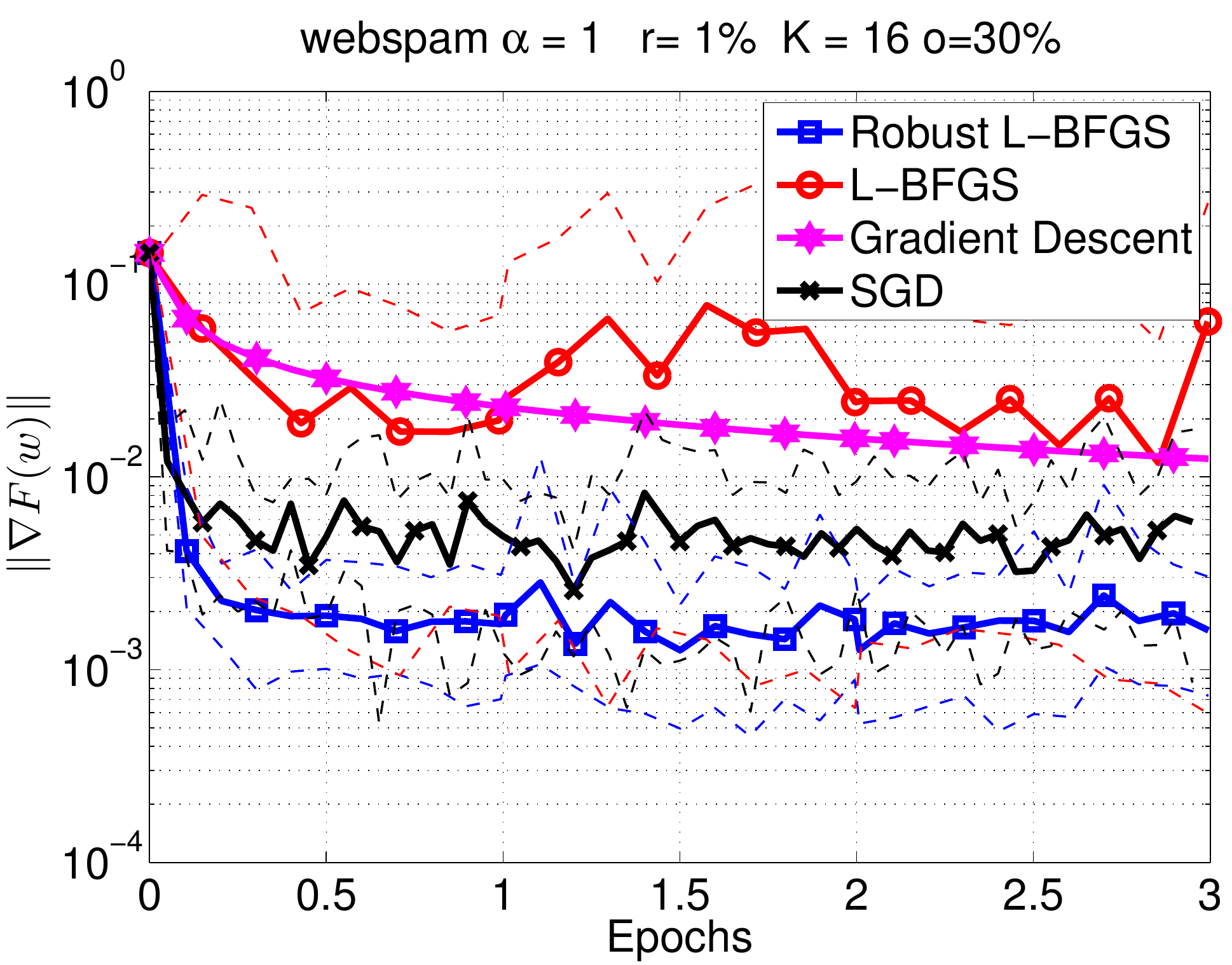}

\caption{\textbf{webspam dataset}. Comparison of Robust L-BFGS, L-BFGS (multi-batch L-BFGS without enforcing sample consistency), Gradient Descent (multi-batch Gradient method) and SGD for various batch ($r$) and overlap ($o$) sizes.  Solid lines show average performance, and dashed lines show worst and best performance, over 10 runs (per algorithm). $K=16$ MPI processes.}
\label{fig:demo:MB}
\end{figure}

In this Section, we present numerical results that evaluate the proposed robust multi-batch L-BFGS scheme (Algorithm \ref{alg:multi}) on logistic regression  problems. Figure \ref{fig:demo:MB} shows the performance on the  webspam dataset\footnote{LIBSVM: \url{https://www.csie.ntu.edu.tw/~cjlin/libsvmtools/datasets/binary.html}. 
}, where we compare it against three methods:
(i) multi-batch L-BFGS without enforcing sample consistency (L-BFGS),  where gradient differences are computed using different samples, i.e., $y_k = g_{k+1}^{S_{k+1}}-g_{k}^{S_{k}}$; (ii) multi-batch gradient descent (Gradient Descent), which is obtained by setting $H_k = I$ in Algorithm \ref{alg:multi}; and, (iii) serial SGD, where at every iteration one sample is used to compute the gradient. We run each method with 10 different random seeds, and, where applicable, report results for different batch ($r$) and overlap ($o$) sizes.  The proposed method is more stable than the standard L-BFGS method; this is especially noticeable when $r$ is small. On the other hand, serial SGD achieves similar accuracy as the robust L-BFGS method and at a similar rate (e.g., $r=1\%$), at the cost of $n$ communications per epochs versus $\frac{1}{r(1-o)}$ 
communications per epoch. Figure \ref{fig:demo:MB} also indicates that the robust L-BFGS method is not too sensitive to the size of overlap. Similar behavior was observed on other datasets, in regimes where $r\cdot o$ was not too small; see Appendix \ref{sec:ext_numerical_multi}. We mention in passing that  the L-BFGS step was computed using the a vector-free implementation proposed in \cite{chen2014large}.

\begin{figure}[ht]
\centering

\includegraphics[width=4.5cm]{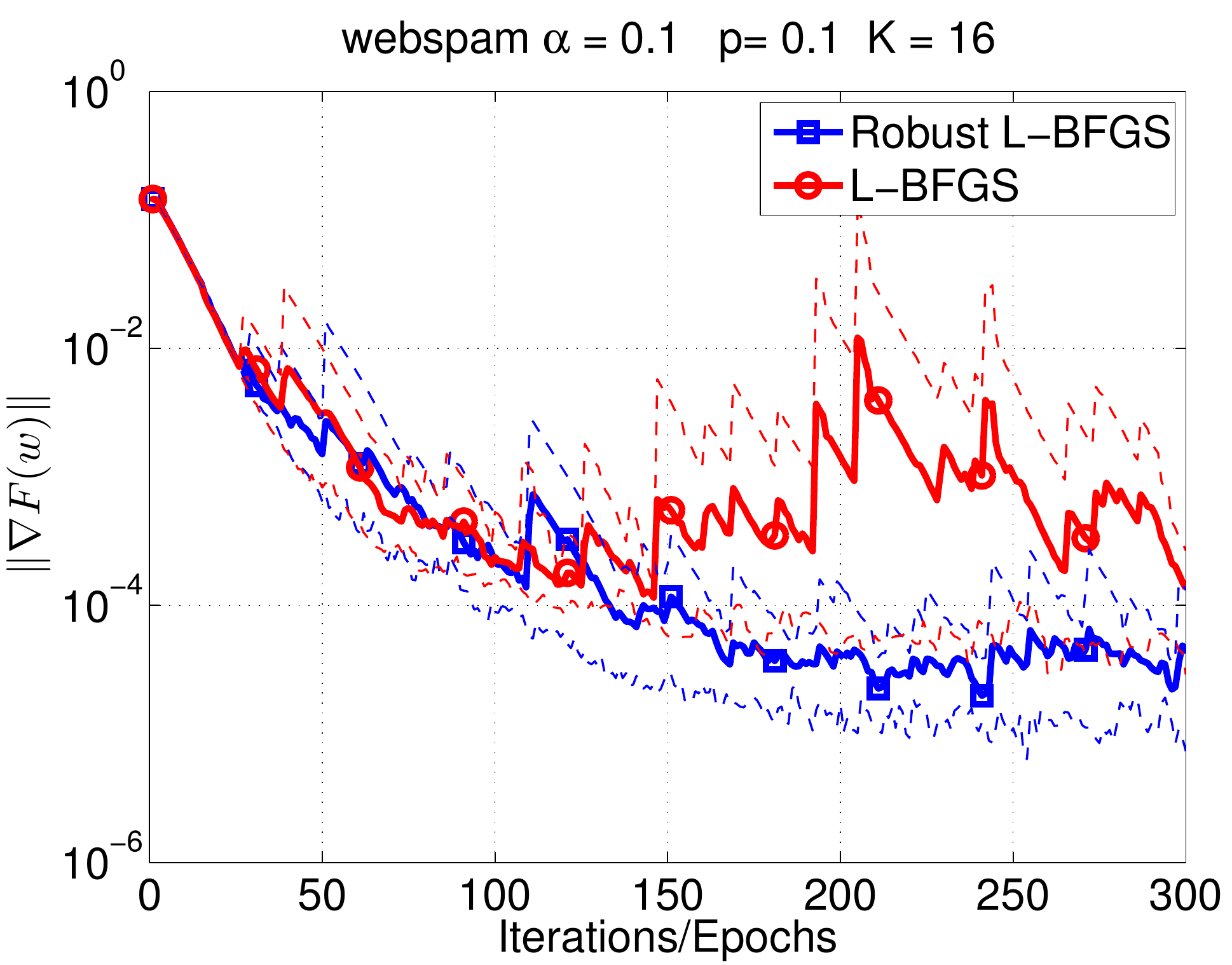}
\includegraphics[width=4.5cm]{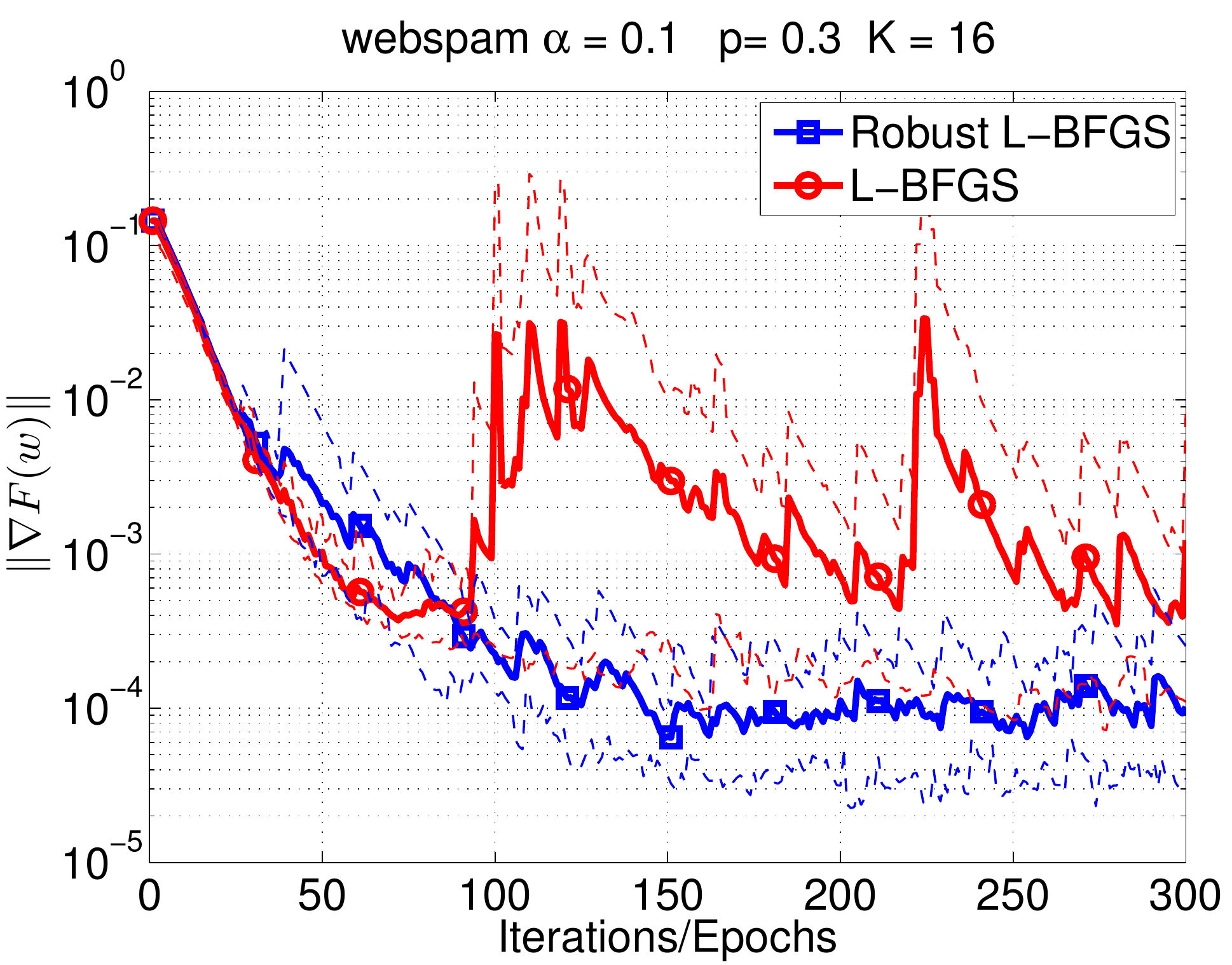}
\includegraphics[width=4.5cm]{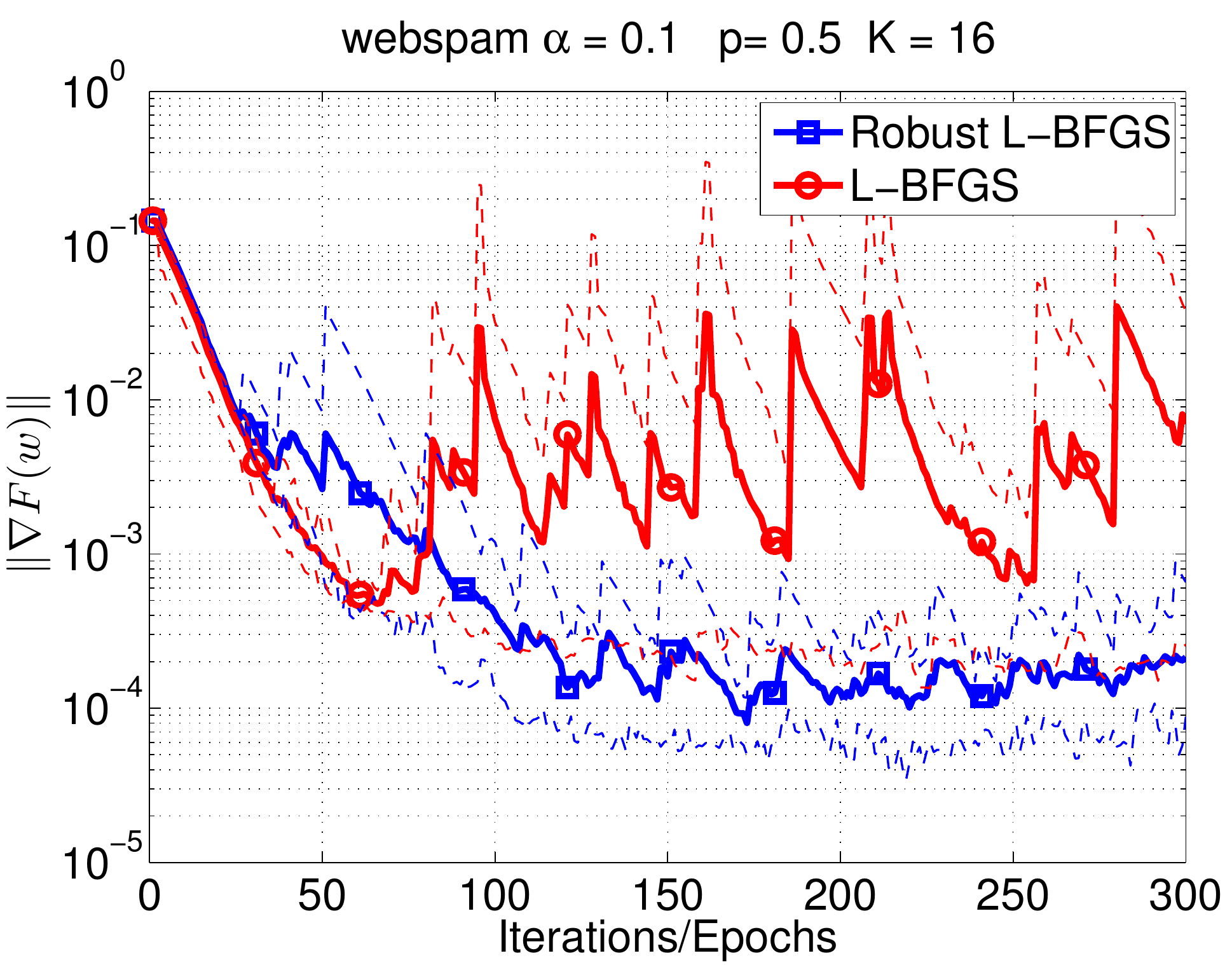}

\caption{\textbf{webspam dataset}. Comparison of Robust L-BFGS and L-BFGS (multi-batch L-BFGS without enforcing sample consistency), for various node failure probabilities $p$. Solid lines show average performance, and dashed lines show worst and best performance, over 10 runs (per algorithm). $K=16$ MPI processes.
}\label{fig:FT:DEMO}
\end{figure}

We also explore the performance of the robust multi-batch L-BFGS method in the presence of node failures (faults), and compare it to the  multi-batch variant that does not enforce sample consistency (L-BFGS). Figure~\ref{fig:FT:DEMO} illustrates the performance of the methods on the webspam dataset, for various probabilities of node failures $p \in \{0.1, 0.3, 0.5 \}$, and suggests that the robust L-BFGS variant is more stable; see Appendix \ref{sec:ext_numerical_fault} for further results.

\begin{figure}[ht]
  \centering

\includegraphics[width=4.5cm]{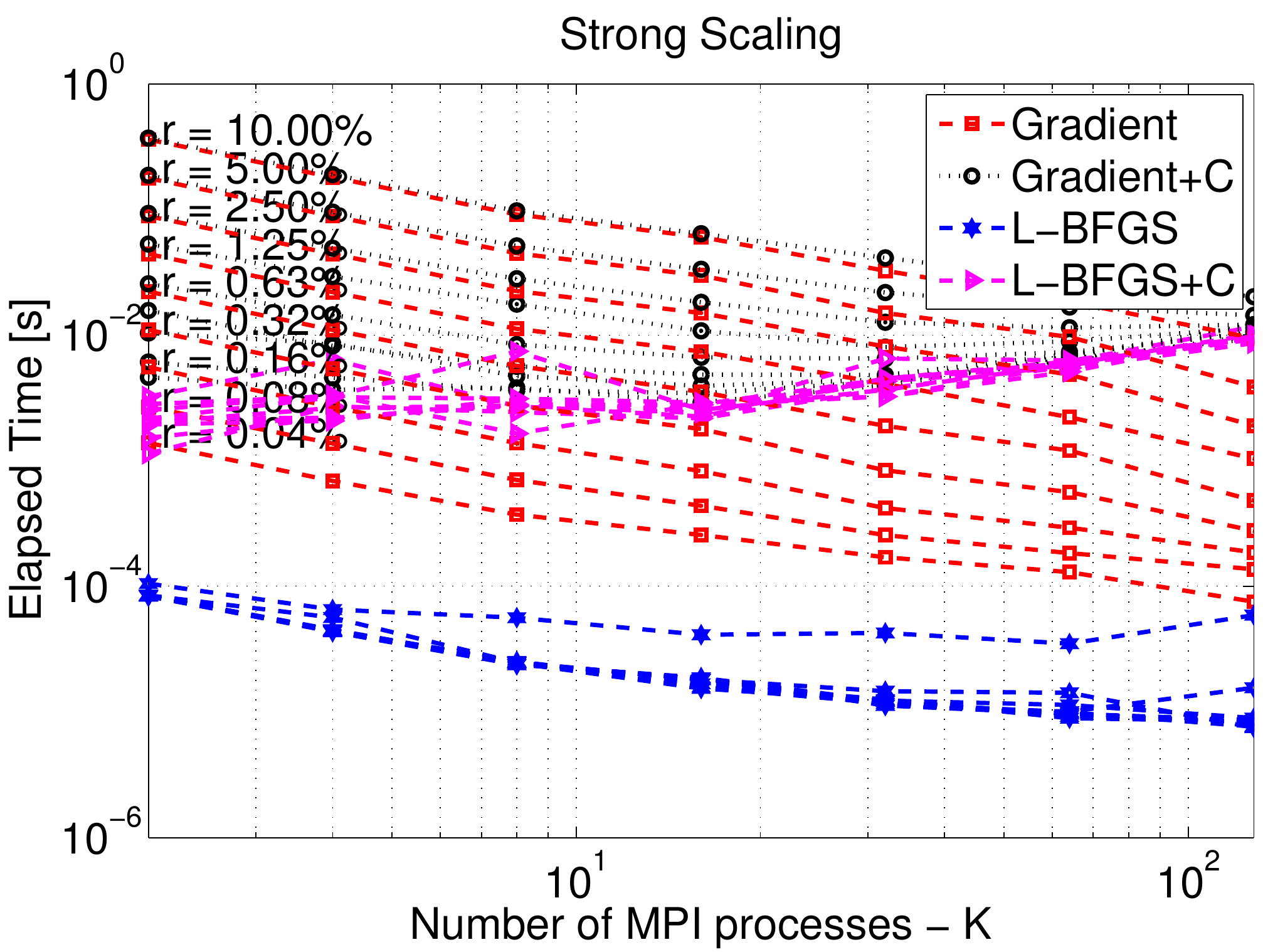}
\includegraphics[width=4.5cm]{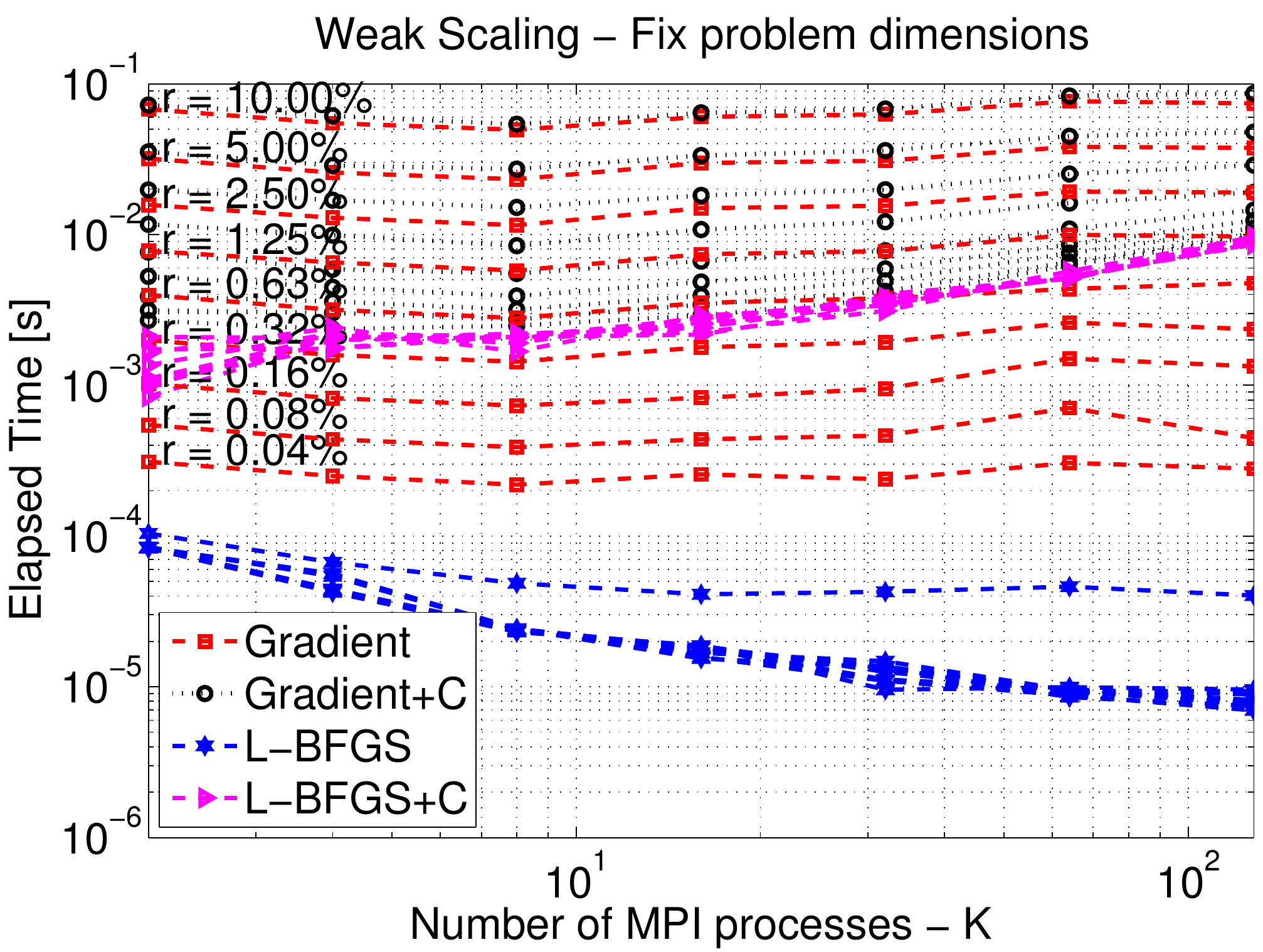}
  
 \caption{Strong and weak scaling of multi-batch L-BFGS method.}
 \label{fig:scaling}
 
   \end{figure}

 Lastly, we study the strong and weak scaling properties of the robust L-BFGS method on artificial data (Figure \ref{fig:scaling}). We measure the time needed to compute a gradient (Gradient) and the associated communication (Gradient+C), as well as, the time needed to compute the L-BFGS direction (L-BFGS) and the associated communication (L-BFGS+C), for various batch sizes ($r$).
The figure on the left shows strong scaling of multi-batch LBFGS on a $d=10^4$ dimensional problem with $n=10^7$ samples. The size of input data is 24GB, and we vary the number of MPI processes, $K \in \{1,2,\dots, 128\}$. The time it takes to compute the gradient decreases with $K$, however, for small values of $r$, the communication time exceeds the compute time. 
 The figure on the right shows weak scaling on a problem of similar size, but with varying sparsity. Each sample has $10\cdot K$ non-zero elements, thus for any $K$ the size of local problem is roughly $1.5$GB (for $K=128$ size of data 192GB). We observe almost constant time for the gradient computation while the cost of computing the L-BFGS direction decreases with $K$; however, if communication is considered, the overall time needed to compute the L-BFGS direction increases slightly. For more details see Appendix \ref{sec:scaling_multi}.

\section{Conclusion}
\label{sec:rem}

This paper describes a novel variant of the L-BFGS method that is robust and efficient in two settings. The first occurs in the presence of node failures in a distributed computing implementation; the second arises when one wishes to employ a different batch at each iteration in order to accelerate learning. The proposed method avoids the pitfalls of using inconsistent gradient differences by performing quasi-Newton updating based  on the overlap between consecutive samples. Numerical results show that the method is efficient in practice, and a convergence analysis illustrates its theoretical properties. 


\subsubsection*{Acknowledgements}

The first two authors were supported by the Office of Naval Research award N000141410313, the Department of Energy grant DE-FG02-87ER25047 and the National Science Foundation grant DMS-1620022.  Martin Tak\'{a}\v{c} was supported by National Science Foundation grant 
CCF-1618717.

\clearpage

\small
\bibliographystyle{abbrv}

\bibliography{multi_refs.bib}

\normalsize

\clearpage 
\appendix

\section{Proofs and Technical Results}
\label{sec:app}

\subsection{Assumptions}

We first restate the Assumptions that we use in the Convergence Analysis section (Section \ref{sec:conv}). Assumption $A$ and $B$ are used in the strongly convex and nonconvex cases, respectively.

\paragraph{Assumptions A}
\emph{
\begin{enumerate}
\item $F$ is twice continuously differentiable.
\item There exist positive constants $\hat{\lambda}$ and $\hat{\Lambda}$ such that
\begin{align} \label{eq:b_grad}
\hat{\lambda} I \preceq \nabla^2F^O(w) \preceq \hat{\Lambda} I,
\end{align}
for all $w \in \mathbb{R}^d$ and all sets $O \subset  \{1,2,\ldots,n\}$.
\item There is a constant $\gamma$ such that 
\begin{align}    \label{eq:b_grad}
\mathbb{E}_{S}\left[ \| \nabla  F^{S}(w) \| \right]^2 \leq \gamma^2,
\end{align}
for all $w \in \mathbb{R}^d$ and all batches $S\subset  \{1,2,\ldots,n\}$.
\item The samples $S$ are drawn independently and $\nabla F^{S}(w)$ is an unbiased estimator of the true gradient $\nabla F(w)$ for all $w \in \mathbb{R}^d$, i.e.,
\begin{align}		\label{eq:unbiased}
\mathbb{E}\Big[ \nabla F^{S}(w)\Big] = \nabla F(w).
\end{align}
\end{enumerate}
}

Note that Assumption $A.2$ implies that the entire Hessian $\nabla^2F(w)$ also satisfies 
\begin{align}   \label{eq:b_hess1}
 \lambda I \preceq \nabla^2F(w) \preceq  \Lambda I,   \forall w \in \mathbb{R}^d,
\end{align}
for some constants $ \lambda,  \Lambda>0$.

\paragraph{Assumptions B}
\emph{
\begin{enumerate}
\item $F$ is twice continuously differentiable.
\item The gradients of $F$ are $\Lambda$-Lipschitz continuous and the gradients of $F^{O}$ are $\Lambda_{O}$-Lipschitz continuous for all $w \in \mathbb{R}^d$ and all sets $O \subset  \{1,2,\ldots,n\}$.
\item The function $ F(w)$ is bounded below by a scalar $\widehat F$.
\item There exist constants $\gamma \geq 0$ and $\eta>0$ such that 
\begin{align}    \label{eq:b_grad1}
\mathbb{E}_{S}\left[ \| \nabla  F^{S}(w) \| \right]^2 \leq \gamma^2 + \eta \| \nabla F(w)\|^2,
\end{align}
for all $w \in \mathbb{R}^d$ and all batches $S\subset  \{1,2,\ldots,n\}$.
\item The samples $S$ are drawn independently and $\nabla F^{S}(w)$ is an unbiased estimator of the true gradient $\nabla F(w)$ for all $w \in \mathbb{R}^d$, i.e.,
\begin{align}		\label{eq:unbiased1}
\mathbb{E}\Big[ \nabla F^{S}(w)\Big] = \nabla F(w).
\end{align}
\end{enumerate}
}

\subsection{Proof of Lemma 3.1 (Strongly Convex Case)}
\label{sec:app1}

The following Lemma shows that the eigenvalues of the matrices generated by the multi-batch L-BFGS method are bounded above and away from zero if $F$ is strongly convex. 

\begin{customlemma}{3.1} 
If Assumptions A.1-A.2 above hold, there exist constants $0<\mu_1\leq \mu_2$ such that the Hessian approximations $\{H_k\}$ generated by the \multi L-BFGS method (Algorithm 1) satisfy
\begin{align*}    
\mu_1 I \preceq H_k \preceq \mu_2 I,\qquad \text{for } k=0,1,2,\dots 
\end{align*}
\end{customlemma}

\begin{proof}
Instead of analyzing the inverse Hessian approximation $H_k$, we study the direct Hessian approximation $B_k = H_k^{-1}$. In this case, the limited memory quasi-Newton updating formula is given as follows
\begin{enumerate}
\item Set $B_k^{(0)}=\frac{y_k^Ty_k}{s_k^Ty_k}I$ and $\tilde{m} = \min\{k,m\}$; where $m$ is the memory in L-BFGS.
\item For $i=0,...,\tilde{m}-1$ set $j=k-\tilde{m}+1+i$ and compute
\begin{align*}    
B_k^{(i+1)}=B_k^{(i)}-\frac{B_k^{(i)}s_js_j^TB_k^{(i)}}{s_j^TB_k^{(i)}s_j} + \frac{y_jy_j^T}{y_j^Ts_j}.
\end{align*}
\item Set $B_{k+1} = B_k^{(\tilde{m})}.$
\end{enumerate}

The curvature pairs $s_k$ and $y_k$ are updated via the following formulae
\begin{align}    \label{eq:curv_upd}
y_{k+1}=g_{k+1}^{O_{k}}-g_{k}^{O_{k}}, \qquad s_k = w_{k+1}-w_k.
\end{align}
A consequence of Assumption $A.2$ is that the eigenvalues of any sub-sampled Hessian ($\left| O \right|$ samples) are bounded above and away from zero. Utilizing this fact, the convexity of component functions and the definitions  \eqref{eq:curv_upd}, we have
\begin{align} \label{eq:expr1}
y_k^Ts_k \geq \frac{1}{\hat{\Lambda}}\|y_k\|^2  \quad & \Rightarrow \quad \frac{\|y_k\|^2 }{y_k^Ts_k} \leq \hat{\Lambda}.
\end{align}

On the other hand, strong convexity of the sub-sampled functions, the consequence of Assumption $A.2$ and definitions  \eqref{eq:curv_upd}, provide a lower bound,
\begin{align} \label{eq:expr2}
y_k^Ts_k \leq \frac{1}{\hat{\lambda}}\|y_k\|^2  \quad & \Rightarrow \quad \frac{\|y_k\|^2 }{y_k^Ts_k} \geq \hat{\lambda}.
\end{align}   

Combining the upper and lower bounds \eqref{eq:expr1} and \eqref{eq:expr2}  
\begin{align}    \label{eq:bound}
\hat{\lambda} \leq \frac{\|y_k\|^2}{y_k^Ts_k} \leq \hat{\Lambda}.
\end{align}

The above proves that the eigenvalues of the matrices $B_k^{(0)}=\frac{y_k^Ty_k}{s_k^Ty_k}I$ at the start of the L-BFGS update cycles are bounded above and away from zero, for all $k$. We now use a Trace-Determinant argument to show that the eigenvalues of $B_k$ are bounded above and away from zero. 

Let $Tr(B)$ and $\det(B)$ denote the trace and determinant of matrix $B$, respectively, and set $j_i = k-\tilde{m}+i$. The trace of the matrix $B_{k+1}$ can be expressed as, 
\begin{align} \label{eq:trace}
Tr(B_{k+1}) &= Tr(B_k^{(0)}) - Tr\sum_{i=1}^{\tilde{m}}\big(\frac{B_k^{(i)}s_{j_i}s_{j_i}^TB_k^{(i)}}{s_{j_i}^TB_k^{(i)}s_{j_i}}\big) + Tr\sum_{i=1}^{\tilde{m}} \frac{y_{j_i}y_{j_i}^T}{y_{j_i}^Ts_{j_i}}\nonumber\\
&\leq Tr(B_k^{(0)}) + \sum_{i=1}^{\tilde{m}} \frac{\|y_{j_i}\|^2}{y_{j_i}^Ts_{j_i}}\nonumber\\
&\leq Tr(B_k^{(0)}) + \tilde{m}\hat{\Lambda} \nonumber\\
& \leq C_1, 
\end{align}
for some positive constant $C_1$, where the inequalities above are due to \eqref{eq:bound}, and the fact that the eigenvalues of the initial L-BFGS matrix $B_k^{(0)}$ are bounded above and away from zero.

Using a result due to Powell \cite{powell1976some}, the determinant of the matrix  $B_{k+1}$ generated by the multi-batch L-BFGS method can be expressed as, 
\begin{align}	\label{eq:det}
\det (B_{k+1}) &= \det (B_{k}^{(0)}) \prod_{i=1}^{\tilde{m}} \frac{y_{j_i}^Ts_{j_i}}{s_{j_i}^TB_{k}^{(i-1)}s_{j_i}} \nonumber\\
& =  \det (B_{k}^{(0)}) \prod_{i=1}^{\tilde{m}} \frac{y_{j_i}^Ts_{j_i}}{s_{j_i}^Ts_{j_i}}  \frac{s_{j_i}^Ts_{j_i}}{s_{j_i}^TB_{k}^{(i-1)}s_{j_i}}\nonumber\\
& \geq  \det (B_{k}^{(0)}) \Big( \frac{\hat{\lambda}}{C_1} \Big)^{\tilde{m}} \nonumber\\
& \geq C_2,
\end{align}
for some positive constant $C_2$, where the above inequalities are due to the fact that the largest eigenvalue of $B_{k}^{(i)}$ is less than $C_1$ and Assumption $A.2$.

The trace \eqref{eq:trace} and determinant \eqref{eq:det} inequalities derived above imply that largest eigenvalues of all matrices $B_k$ are bounded above, uniformly, and that the smallest eigenvalues of all matrices $B_k$ are bounded away from zero, uniformly.
\end{proof}

\subsection{Proof of Theorem 3.2 (Strongly Convex Case)}
\label{sec:app2}

Utilizing the result from Lemma \ref{lemma1}, we now prove a linear convergence result to a neighborhood of the optimal solution, for the case where Assumptions $A$ hold.

\begin{customthm}{3.2}
Suppose that Assumptions A.1-A.4 above hold, and let $F^{\star} = F(w^{\star})$, where $w^{\star}$ is the minimizer of $F$. Let $\{w_k\}$ be the iterates generated by the \multi L-BFGS method (Algorithm 1) with 
\begin{align*}
\alpha_k = \alpha \in  (0,\frac{1}{2\mu_1 \lambda}),
\end{align*}
starting from $w_0$.
Then for all $k\geq 0$,
\begin{align*}   
\mathbb{E} [ F(w_k) - F^{\star} ] &\leq  ( 1-2\alpha \mu_1 \lambda  )^k  [ F(w_0) - F^{\star}  ] +  [ 1-(1-\alpha\mu_1 \lambda)^k ]\frac{\alpha \mu_2^2 \gamma ^2 \Lambda}{4 \mu_1 \lambda} \\
&\xrightarrow[]{k\rightarrow \infty} \frac{\alpha \mu_2^2 \gamma ^2 \Lambda}{4 \mu_1 \lambda}.
\end{align*}
\end{customthm}

\begin{proof}
We have that
\begin{align} \label{eq:proof}
F(w_{k+1}) & = F(w_k -\alpha H_k \nabla F^{S_k}(w_k)) \nonumber \\
 & \leq F(w_k) + \nabla F(w_k)^T (-\alpha H_k \nabla F^{S_k}(w_k)) + \frac{\Lambda}{2}\| \alpha H_k \nabla F^{S_k}(w_k)\|^2 \nonumber \\
 & \leq F(w_k) - \alpha \nabla F(w_k)^T  H_k \nabla F^{S_k}(w_k) + \frac{\alpha^2 \mu_2^2 \Lambda}{2} \| \nabla F^{S_k}(w_k)\|^2,
\end{align}
where the first inequality arises due to \eqref{eq:b_hess1}, and the second inequality arises as a consequence of Lemma \ref{lemma1}.

Taking the expectation (over $S_k$) of equation \eqref{eq:proof}
\begin{align} \label{eq:proof_2}
\mathbb{E}_{S_k}[ F(w_{k+1})] & \leq F(w_k) - \alpha \nabla F(w_k)^T  H_k \nabla F(w_k) + \frac{\alpha^2 \mu_2^2 \Lambda}{2} \mathbb{E}_{S_k} \Big[ \| \nabla F^{S_k}(w_k)\| \Big]^2  \nonumber\\
& \leq F(w_k) - \alpha \mu_1 \| \nabla F(w_k) \|^2  + \frac{\alpha^2 \mu_2^2 \gamma^2 \Lambda}{2},
\end{align}
where in the first inequality we make use of Assumption $A.5$, and the second inequality arises due to Lemma \ref{lemma1} and Assumption $A.4$.

Since $F$ is $\lambda$-strongly convex, we can use the following relationship between the norm of the gradient squared, and the distance of the $k$-th iterate from the optimal solution.
\begin{align*}   
2\lambda [F(w_k) - F^{\star}] \leq \| \nabla F(w_k)\|^2.
\end{align*}
Using the above with \eqref{eq:proof_2},
\begin{align} \label{eq:proof_3}
\mathbb{E}_{S_k} [F(w_{k+1})] & \leq F(w_k) - \alpha \mu_1 \| \nabla F(w_k) \|^2  + \frac{\alpha^2 \mu_2^2 \gamma^2 \Lambda}{2} \nonumber \\
& \leq F(w_k) - 2\alpha \mu_1 \lambda  [F(w_k) - F^{\star}]  + \frac{\alpha^2 \mu_2^2 \gamma^2 \Lambda}{2}.
\end{align}

Let 
\begin{align}			\label{eq:phi}
\phi_k = \mathbb{E}[F(w_k) - F^{\star}],
\end{align}
where the expectation is over all batches $S_0,S_1,...,S_{k-1}$ and all history starting with $w_0$. Thus \eqref{eq:proof_3} can be expressed as,
\begin{align}    \label{eq:proof_4}
\phi_{k+1} \leq (1 - 2\alpha \mu_1 \lambda ) \phi_k + \frac{\alpha^2 \mu_2^2 \gamma^2 \Lambda}{2},
\end{align}
from which we deduce that in order to reduce the value with respect to the previous function value, the step length needs to be in the range 
\begin{align*}
\alpha \in \Big(0, \frac{1}{2\mu_1 \lambda}\Big).
\end{align*}

Subtracting $\frac{\alpha\mu_2^2\gamma^2\Lambda}{4\mu_1 \lambda}$ from either side of \eqref{eq:proof_4} yields
\begin{align}    \label{eq:proof4.1}
\phi_{k+1} - \frac{\alpha\mu_2^2\gamma^2\Lambda}{4\mu_1 \lambda} & \leq (1 - 2\alpha \mu_1 \lambda ) \phi_k + \frac{\alpha^2 \mu_2^2 \gamma^2 \Lambda}{2} - \frac{\alpha\mu_2^2\gamma^2\Lambda}{4\mu_1 \lambda} \nonumber \\
& = (1 - 2\alpha \mu_1 \lambda ) \Big[ \phi_k  - \frac{\alpha\mu_2^2\gamma^2\Lambda}{4\mu_1 \lambda} \Big].
\end{align}

Recursive application of \eqref{eq:proof4.1} yields
\begin{align*}   
\phi_{k} - \frac{\alpha\mu_2^2\gamma^2\Lambda}{4\mu_1 \lambda}& \leq  (1 - 2\alpha \mu_1 \lambda )^k \Big[ \phi_0  - \frac{\alpha\mu_2^2\gamma^2\Lambda}{4\mu_1 \lambda} \Big],
\end{align*}

and thus, 
\begin{align}	\label{eq: blah}
	\phi_{k} \leq (1 - 2\alpha \mu_1 \lambda )^k \phi_0  + \Big[ 1-(1-\alpha\mu_1 \lambda)^k\Big]\frac{\alpha \mu_2^2 \gamma ^2 \Lambda}{4 \mu_1 \lambda}.
\end{align}

Finally using the definition of $\phi_k$ \eqref{eq:phi} with the above expression yields the desired result,
\begin{equation*}  
\mathbb{E}\Big[ F(w_k) - F^{\star}\Big] \leq \Big( 1-2\alpha \mu_1 \lambda \Big)^k \Big[ F(w_0) - F^{\star} \Big] + \Big[ 1-(1-\alpha\mu_1 \lambda)^k\Big]\frac{\alpha \mu_2^2 \gamma ^2 \Lambda}{4 \mu_1 \lambda}. \qedhere 
\end{equation*}
\end{proof}

\subsection{Proof of Lemma 3.3 (Nonconvex Case)}
\label{sec:app3}

The following Lemma shows that the eigenvalues of the matrices generated by the multi-batch L-BFGS method are bounded above and away from zero (nonconvex case). 

\begin{customlemma}{3.3}		
Suppose that Assumptions B.1-B.2  hold and let $\epsilon >0$ be given. Let $\{H_k \}$ be the Hessian approximations generated by  the multi-batch L-BFGS method (Algorithm~\ref{alg:multi}), with the modification that the Hessian approximation $H_k$ update is performed only when  
\begin{align} 	\label{eq:skip}		
	y_k^Ts_k \geq {\epsilon} \| s_k \|^2,
\end{align}
else $H_{k+1} = H_k$. Then, there exist constants $0<\mu_1\leq \mu_2$ such that 
\begin{align*}    
\mu_1 I \preceq H_k \preceq \mu_2 I,\qquad \text{for } k=0,1,2,\dots 
\end{align*}
\end{customlemma}

\begin{proof}
Similar to the proof of Lemma \ref{lemma1}, we study the direct Hessian approximation $B_k = H_k^{-1}$.

The curvature pairs $s_k$ and $y_k$ are updated via the following formulae
\begin{align}    \label{eq:curv_upd1}
y_{k+1}=g_{k+1}^{O_{k}}-g_{k}^{O_{k}}, \qquad s_k = w_{k+1}-w_k.
\end{align}

The skipping mechanism \eqref{eq:skip} provides both an upper and lower bound on the quantity $\frac{\|y_k\|^2 }{y_k^Ts_k}$, which in turn ensures that the initial L-BFGS Hessian approximation is bounded above and away from zero. The lower bound is attained by repeated application of Cauchy's inequality to condition \eqref{eq:skip}. We have from \eqref{eq:skip} that
\begin{align*}
	\epsilon \| s_k \|^2 &\leq y_k^Ts_k \leq  \| y_k \| \| s_k \|,
\end{align*} 
and therefore
\begin{align*}
	\| s_k \| \leq \frac{1}{\epsilon} \| y_k \|.
\end{align*}
It follows that 
\begin{align*}
	s_k^Ty_k \leq \| s_k \| \| y_k \| \leq \frac{1}{\epsilon} \| y_k \|^2
\end{align*}
and hence
\begin{align} \label{eq:lower}
	\frac{\| y_k \|^2}{s_k^Ty_k} \geq \epsilon.
\end{align}

The upper bound is attained by the Lipschitz continuity of sample gradients,
\begin{align*}
	 y_k^Ts_k & \geq \epsilon \| s_k \|^2\\
	&\geq  \epsilon \frac{ \| y_k \|^2}{\Lambda_{O_k}^2},
\end{align*} 
Re-arranging the above expression yields the desired upper bound,
\begin{align}	\label{eq:upper}
	  \frac{\| y_k \|^2}{s_k^Ty_k} \leq \frac{\Lambda_{O_k}^2}{\epsilon}.
\end{align}

Combining \eqref{eq:lower} and \eqref{eq:upper},
\begin{align*}
 \epsilon \leq \frac{\|y_k\|^2 }{y_k^Ts_k} \leq \frac{\Lambda_{O_k}^2}{\epsilon}.
\end{align*}

The above proves that the eigenvalues of the matrices $B_k^{(0)}=\frac{y_k^Ty_k}{s_k^Ty_k}I$ at the start of the L-BFGS update cycles are bounded above and away from zero, for all $k$. The rest of the proof follows the same trace-determinant argument as in the proof of Lemma \ref{lemma1}, the only difference being that the last inequality in \ref{eq:det} comes as a result of the cautious update strategy.
\end{proof}

\subsection{Proof of Theorem 3.4 (Nonconvex Case)}
\label{sec:app4}

Utilizing the result from Lemma \ref{lemma2}, we can now establish the following result about the behavior of the gradient norm for the  multi-batch L-BFGS method with a cautious update strategy. 

\begin{customthm}{3.4}
Suppose that Assumptions B.1-B.5 above hold. Let $\{w_k\}$ be the iterates generated by the \multi L-BFGS method (Algorithm~\ref{alg:multi}) with 
\begin{align*}
\alpha_k = \alpha \in  (0,\frac{\mu_1}{\mu_2^2\eta \Lambda} ),
\end{align*}
where $w_0$ is the starting point. Also, suppose that if 
\begin{align*}		
	y_k^Ts_k < {\epsilon} \| s_k \|^2,
\end{align*}
for some ${\epsilon}>0$, the inverse L-BFGS Hessian approximation is skipped, $H_{k+1}=H_k$. Then, for all $k\geq0$,
\begin{align*}	
\mathbb{E} \Big[\frac{1}{L}\sum_{k=0}^{L-1} \| \nabla F(w_k) \|^2 \Big] & \leq \frac{\alpha \mu_2^2 \gamma^2 \Lambda}{ \mu_1 } + \frac{2[ F(w_0) - \widehat{F}]}{\alpha \mu_1 L }\\
& \xrightarrow[]{L\rightarrow \infty}\frac{\alpha \mu_2^2 \gamma^2 \Lambda}{ \mu_1 }.
\end{align*}
\end{customthm}

\begin{proof}
Starting with \eqref{eq:proof_2}, 
\begin{align*}
\mathbb{E}_{S_k}[ F(w_{k+1})] &\leq F(w_k) - \alpha \mu_1 \| \nabla F(w_k) \|^2 + \frac{\alpha^2 \mu_2^2 \Lambda}{2} \mathbb{E}_{S_k} \Big[ \| \nabla F^{S_k}(w_k)\| \Big]^2  \nonumber\\
& \leq F(w_k) - \alpha \mu_1 \| \nabla F(w_k) \|^2  + \frac{\alpha^2 \mu_2^2  \Lambda}{2} (\gamma^2 + \eta \| \nabla F(w)\|^2)\\
& = F(w_k) - \alpha \big(\mu_1 - \frac{\alpha \mu_2^2  \eta\Lambda}{2}\big) \| \nabla F(w_k) \|^2 + \frac{\alpha^2 \mu_2^2 \gamma^2\Lambda}{2}\\
& \leq F(w_k) - \frac{\alpha \mu_1}{2} \| \nabla F(w_k) \|^2 + \frac{\alpha^2 \mu_2^2 \gamma^2\Lambda}{2},
\end{align*}
where the second inequality holds due to Assumption $B.4$, and the fourth inequality is obtained by using the upper bound on the step length. Taking an expectation over all  batches $S_0,S_1,...,S_{k-1}$ and all history starting with $w_0$ yields
\begin{align}	\label{eq:nc1}
	\phi_{k+1}-\phi_k \leq - \frac{\alpha \mu_1}{2}  \mathbb{E}\| \nabla F(w_k) \|^2  + \frac{\alpha^2 \mu_2^2 \gamma^2 \Lambda}{2},
\end{align}
where $\phi_k = \mathbb{E}[F(w_k)] $. Summing \eqref{eq:nc1} over the first $L-1$ iterations
\begin{align} \label{eq:nc2}
	\sum_{k=0}^{L-1} [\phi_{k+1}-\phi_k] &\leq - \frac{\alpha \mu_1}{2}  \sum_{k=0}^{L-1}  \mathbb{E}\| \nabla F(w_k) \|^2 + \sum_{k=0}^{L-1} \frac{\alpha^2 \mu_2^2 \gamma^2 \Lambda}{2} \nonumber \\
	&= - \frac{\alpha \mu_1}{2}   \mathbb{E} \Big[\sum_{k=0}^{L-1} \| \nabla F(w_k) \|^2 \Big] +  \frac{\alpha^2 \mu_2^2 \gamma^2 \Lambda L}{2}.
\end{align}
The left-hand-side of the above inequality is a telescoping sum
\begin{align*}
	\sum_{k=0}^{L-1} [\phi_{k+1}-\phi_k]  &= \phi_{L}-\phi_0 \nonumber \\
	&= \mathbb{E}[F(w_{L})] -F(w_0) \nonumber\\
	& \geq \widehat{F} -F(w_0). 
\end{align*}
Substituting the above expression into \eqref{eq:nc2} and re-arranging terms
\begin{align*}
	\mathbb{E} \Big[\sum_{k=0}^{L-1} \| \nabla F(w_k) \|^2 \Big] \leq \frac{\alpha \mu_2^2 \gamma^2 \Lambda L}{ \mu_1 } + \frac{2[ F(w_0) - \widehat{F}]}{\alpha \mu_1 }.
\end{align*}

Dividing the above equation by $L$ completes the proof.
\end{proof}

\clearpage 
\section{Extended Numerical Experiments - Real Datasets}
\label{sec:ext_numerical}
In this Section, we present further numerical results, on the datasets listed in Table \ref{tbl:alldatasets}, in both the multi-batch and fault-tolerant settings. Note, that some of the datasets are too small, and thus, there is no reason to run them on a distributed platform; however, we include them as they are part of the standard benchmarking datasets.

\textbf{Notation.} Let $n$ denote the number of training samples in a given dataset, $d$ the dimension of the parameter vector $w$, and $K$ the number of MPI processes used. The parameter $r$ denotes the fraction of samples in the dataset used to define the gradient, i.e., $r = \frac{\left| S\right|}{n} $. The parameter $o$ denotes the length of overlap between consecutive samples, and is defined  as a fraction of the number of samples in a given batch $S$, i.e., $o = \frac{\left| O\right|}{\left| S\right|}$.

\begin{table}[h!]
\caption{Datasets together with basic statistics. All datasets are available at \url{https://www.csie.ntu.edu.tw/~cjlin/libsvmtools/datasets/binary.html}.} 
\label{tbl:alldatasets}
\centering
\begin{tabular}{l|r|r|r|r}
\multicolumn{1}{c|}{Dataset} & \multicolumn{1}{|c|}{$n$} & \multicolumn{1}{|c}{$d$} 
&
\multicolumn{1}{|c}{Size (MB)}
&
\multicolumn{1}{|c}{ K}

\\ \hline \hline

ijcnn (test) &  91,701 & 22 & 14 & 4
\\
cov & 581,012 & 54 & 68 & 4
\\
news20 & 19,996 & 1,355,191 & 134 & 4

\\
rcvtest & 677,399 & 47,236 & 1,152 & 16 \\

url & 2,396,130 & 3,231,961& 2,108&16 \\
kdda  &	8,407,752	&	20,216,830 & 2,546 & 16\\
kddb  & 19,264,097	&	29,890,095 & 4,894 & 16 \\
webspam & 350,000 & 16,609,143 & 23,866 & 16 \\
splice-site & 50,000,000 &  11,725,480  & 260,705 & 16
\end{tabular}
\end{table}

We focus on logistic regression classification; the objective function is given by 
\begin{align*}
	\min_{w \in \mathbb{R}^d} F(w) =  \frac{1}{n}\sum_{i=1}^{n}\log(1+e^{-y^i(w^Tx^i)})
  + \frac{\sigma}{2} \|w\|^2,
\end{align*}
where $ (x^i, y^i)_{i=1}^n$  denote the training examples and $\sigma = \frac1n$
is the regularization parameter.


\subsection{Multi-batch L-BFGS Implementation}
\label{sec:ext_numerical_multi}

\begin{figure}[h!]
\centering
\includegraphics[width=4.6cm]{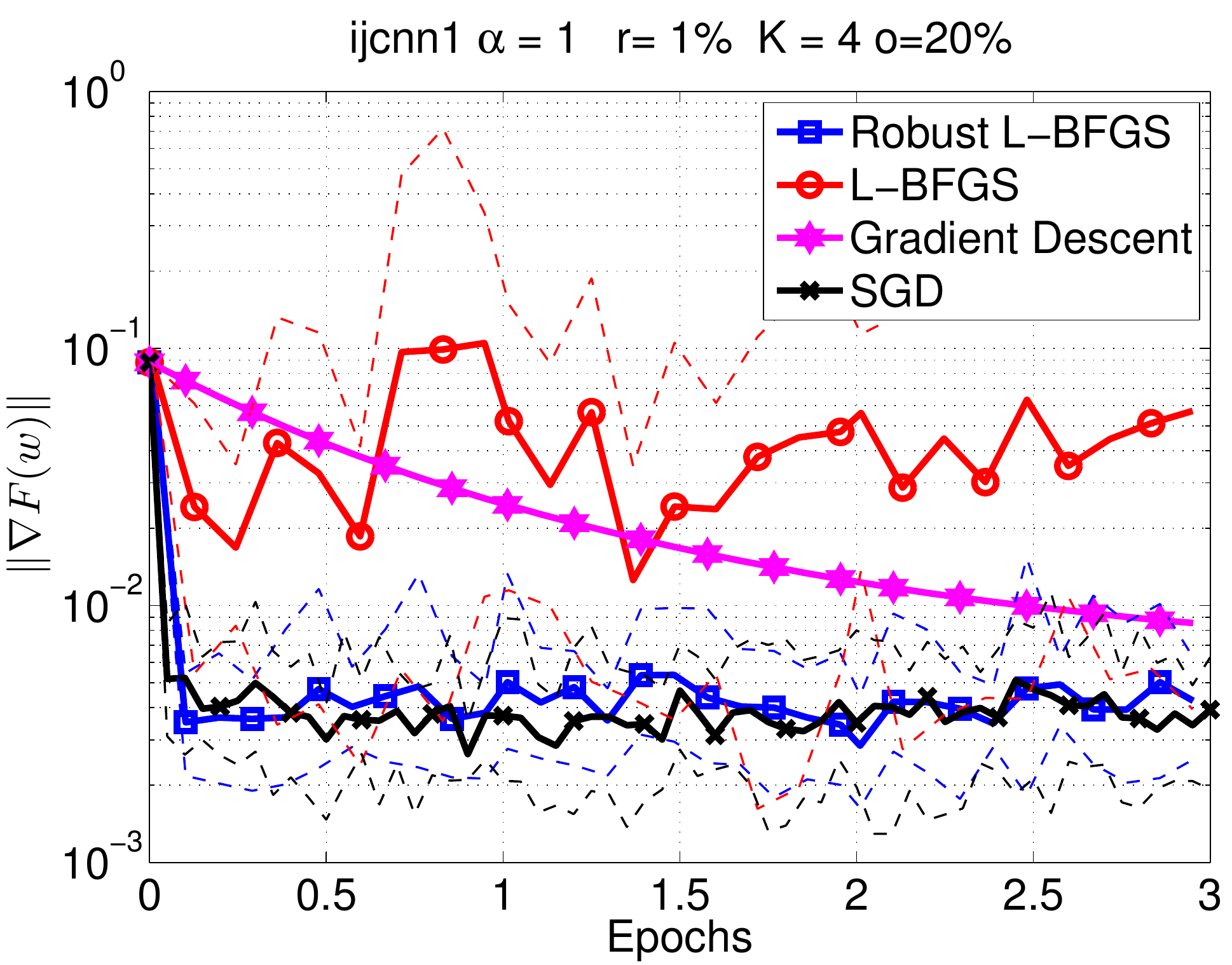}
\includegraphics[width=4.6cm]{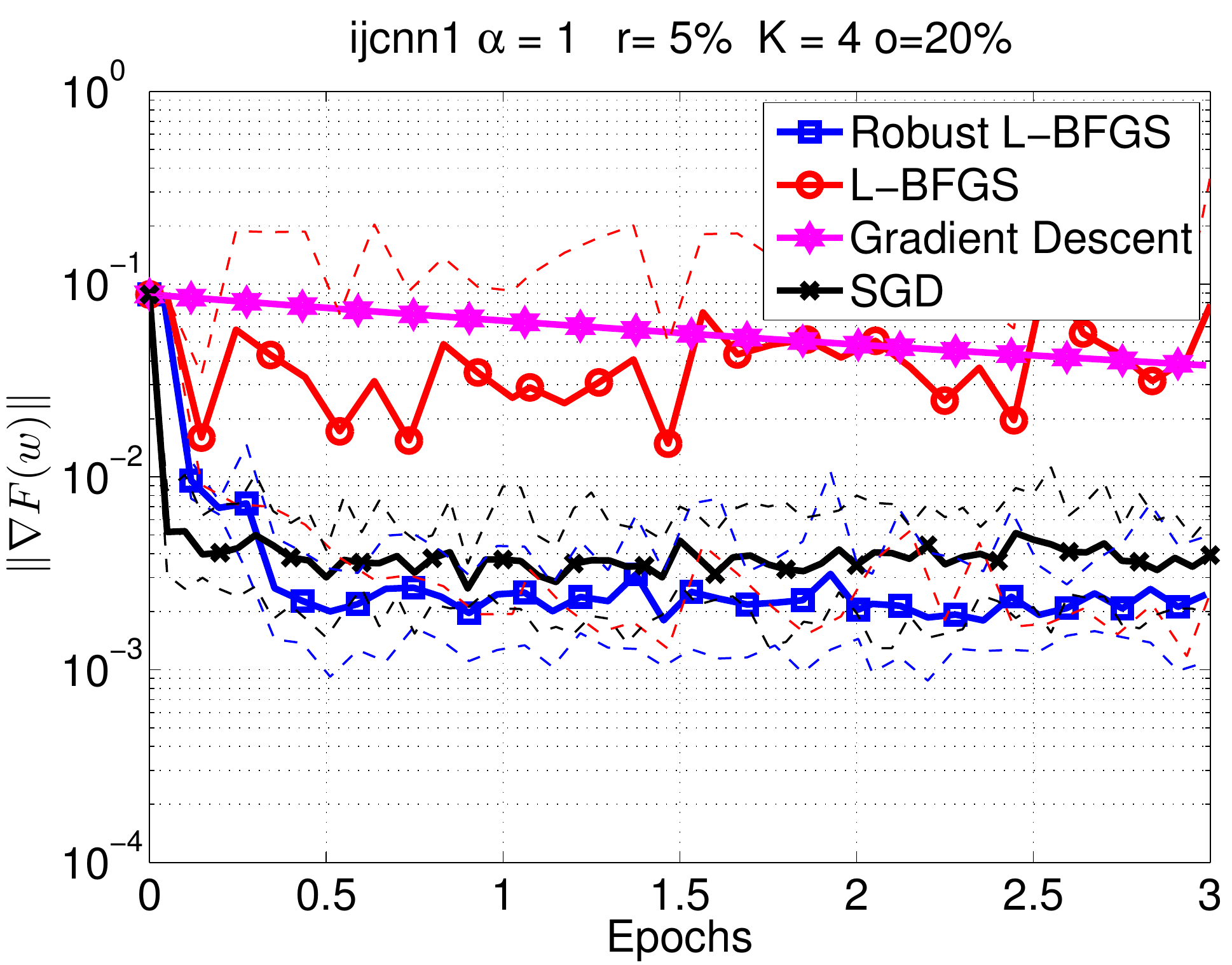}
\includegraphics[width=4.6cm]{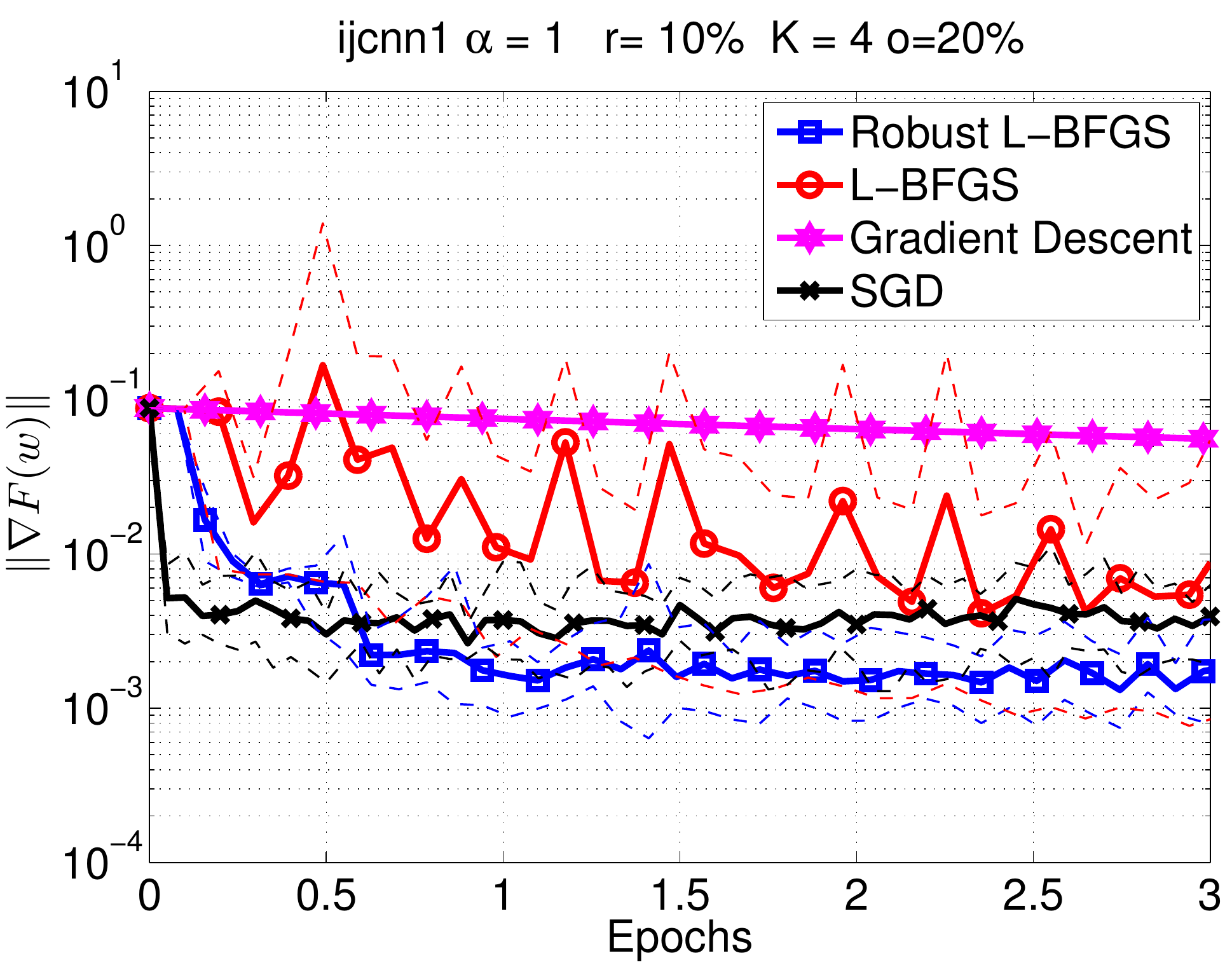}

\includegraphics[width=4.6cm]{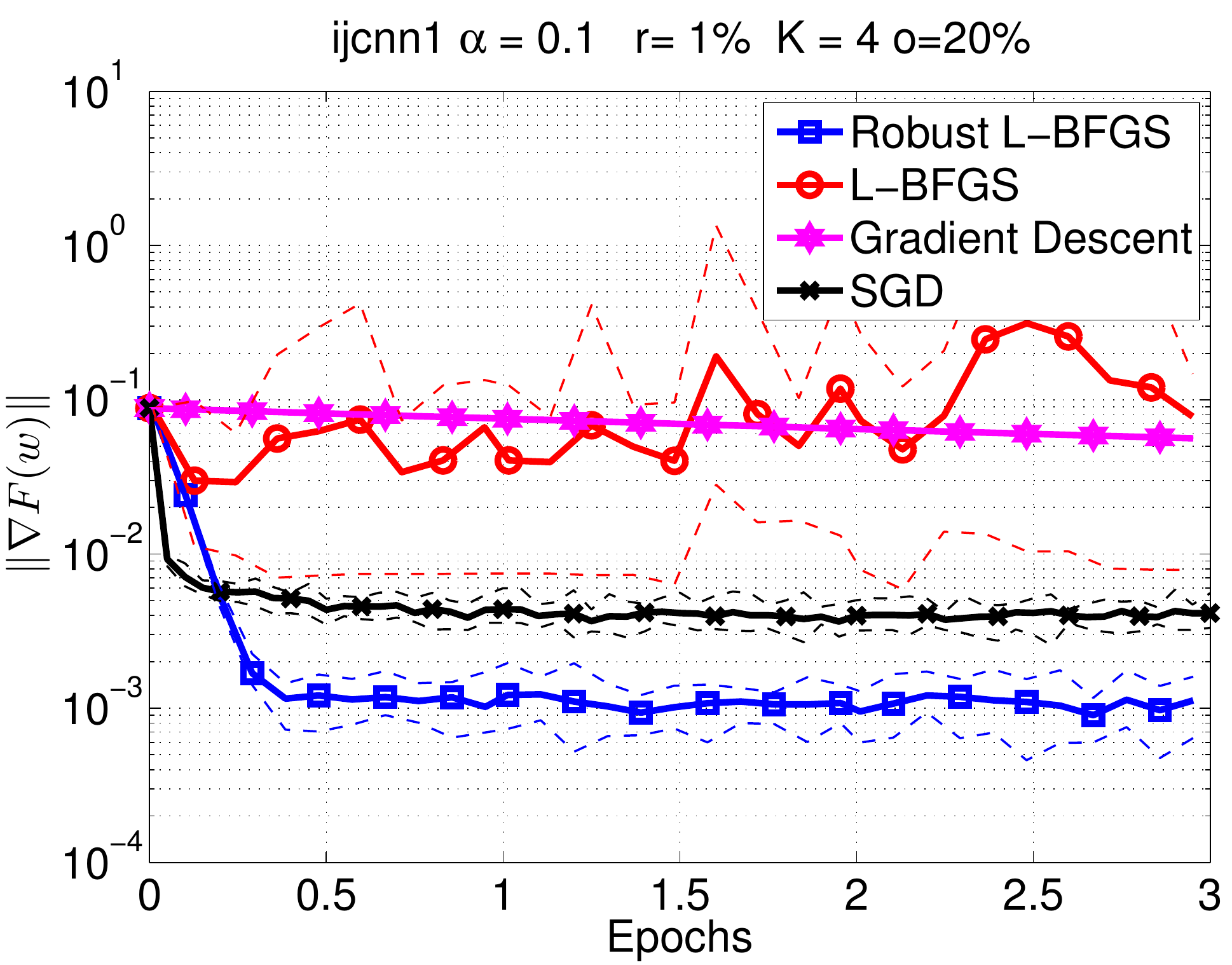}
\includegraphics[width=4.6cm]{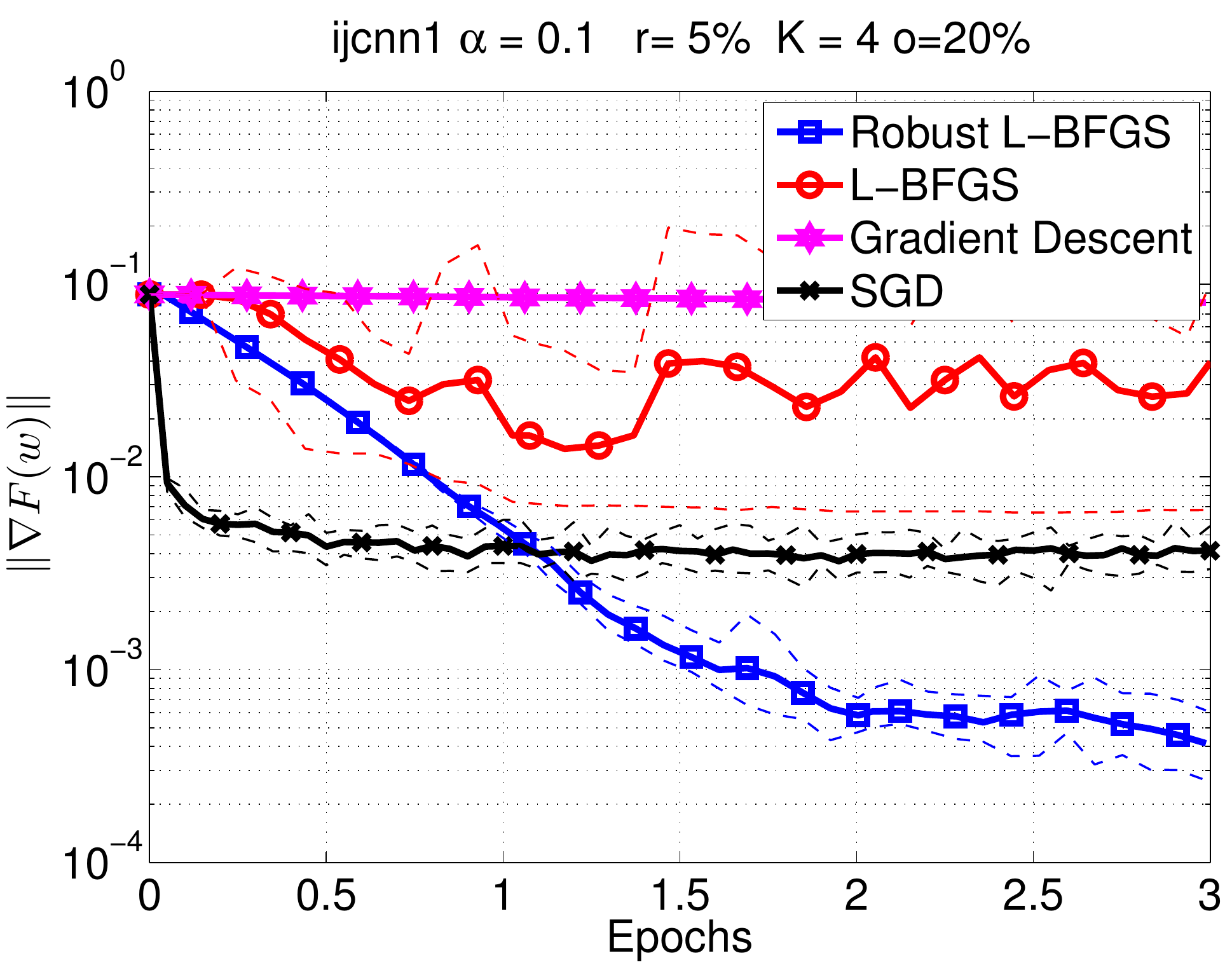}
\includegraphics[width=4.6cm]{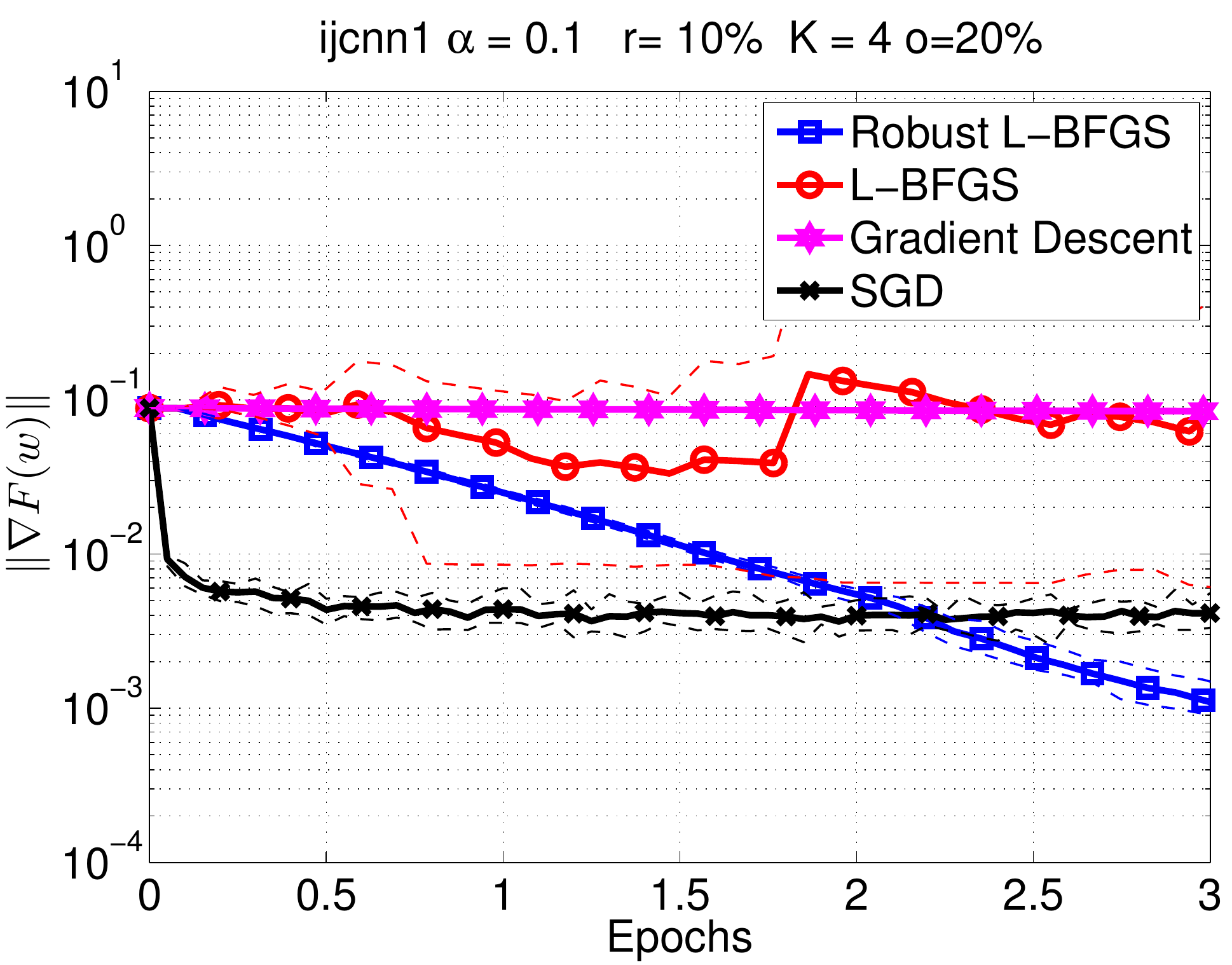}

\hrule 

\includegraphics[width=4.6cm]{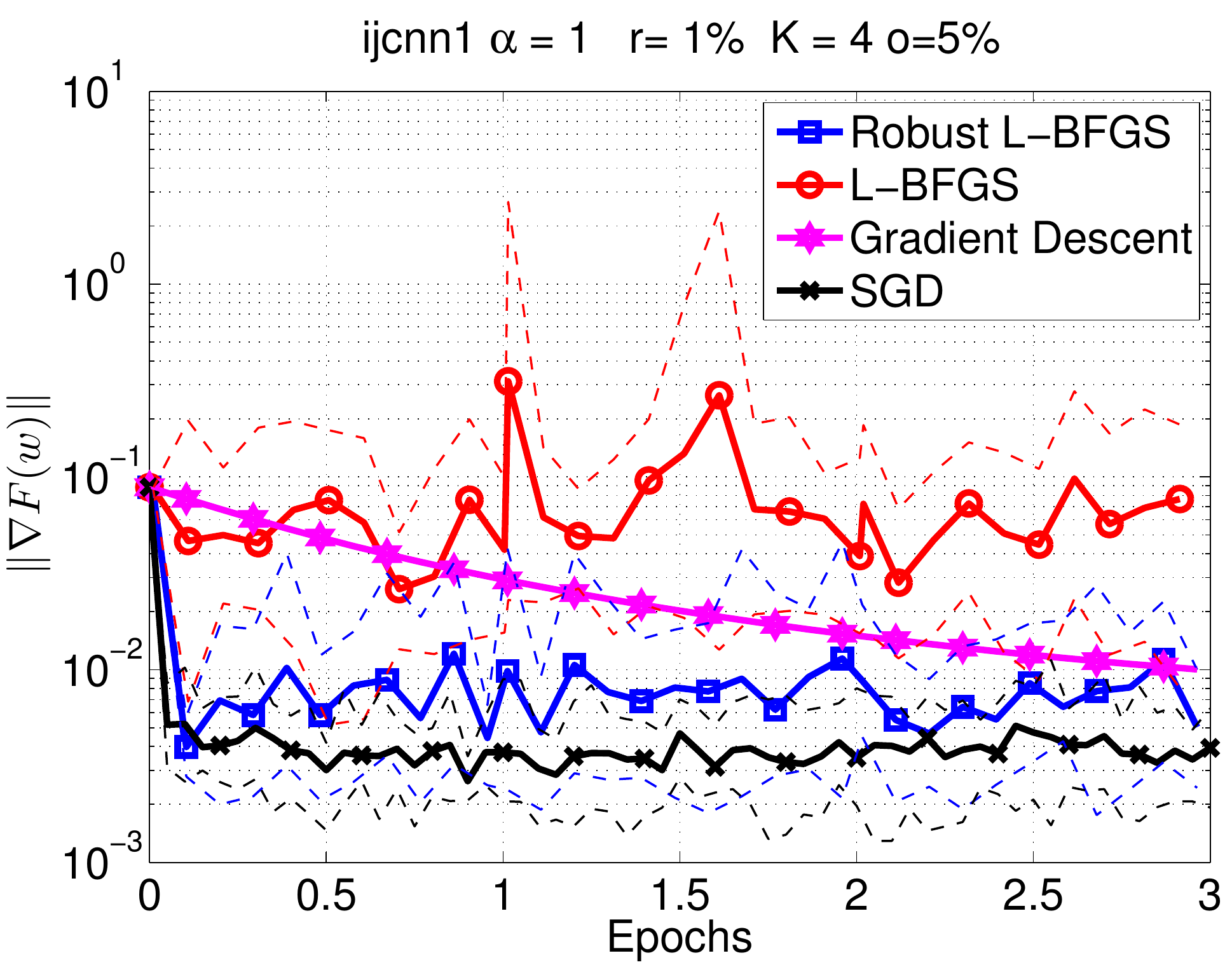}
\includegraphics[width=4.6cm]{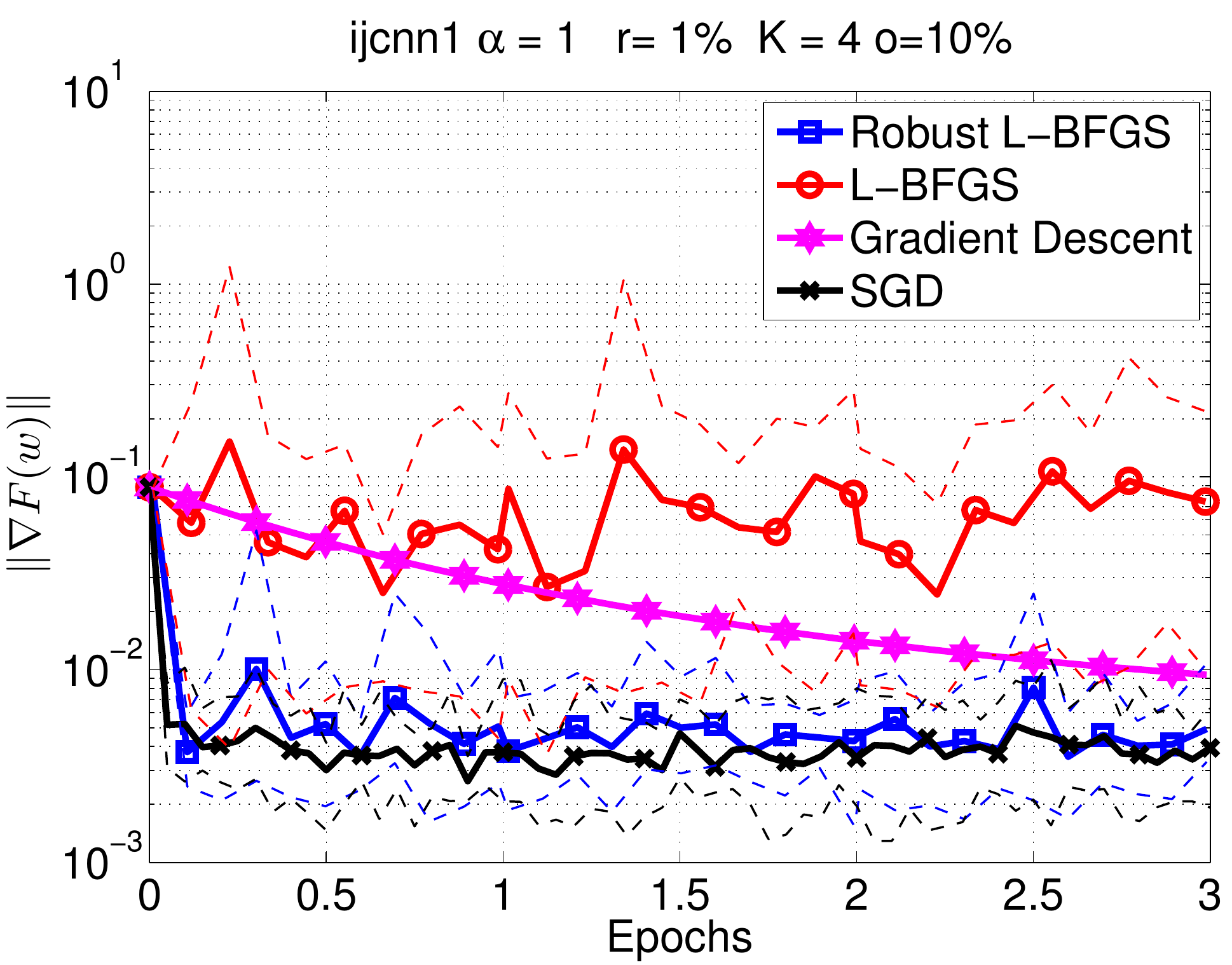}
\includegraphics[width=4.6cm]{ijcnn1_mb_1_0_01_4_0_2-eps-converted-to.pdf}
\includegraphics[width=4.6cm]{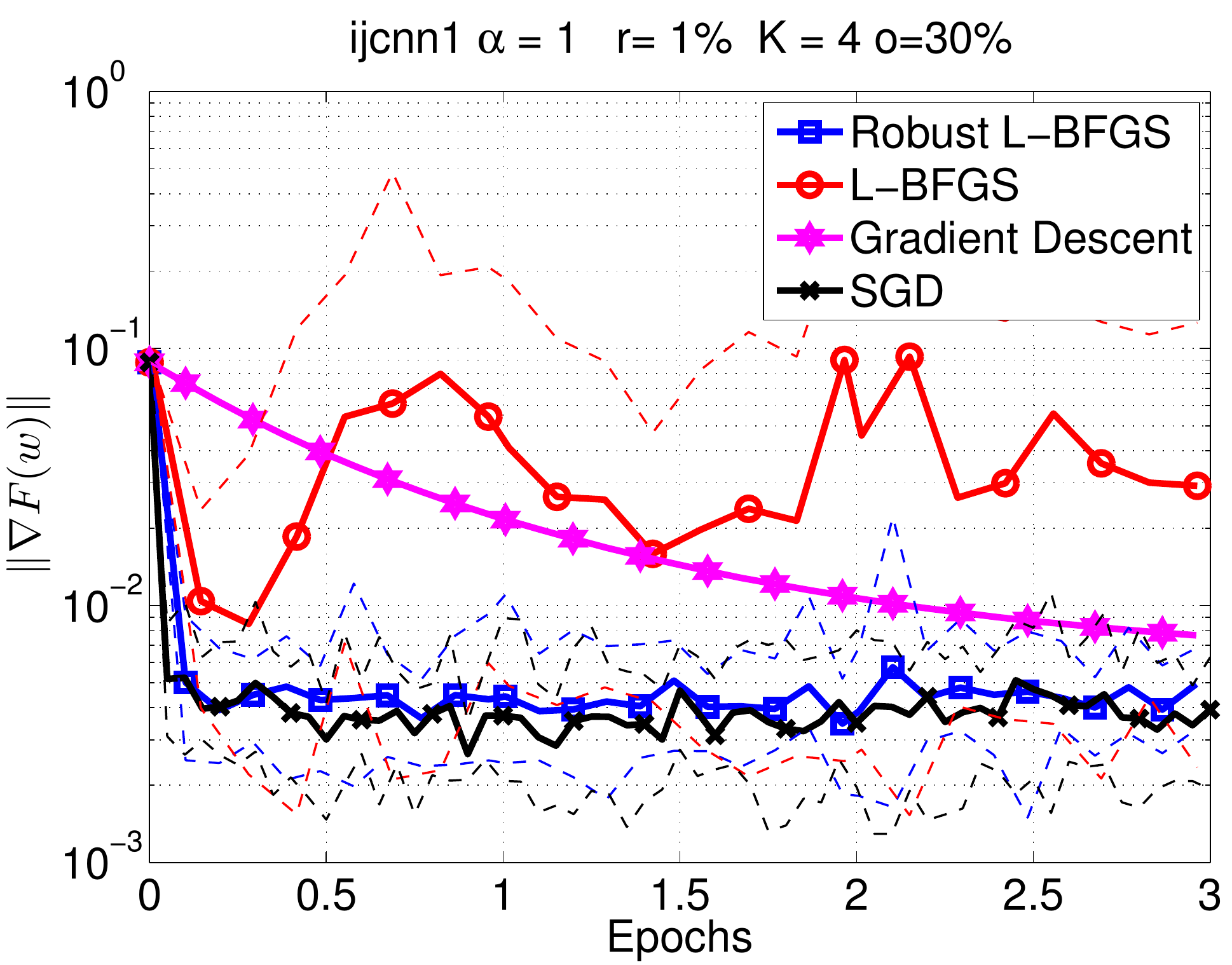}

\caption{\textbf{ijcnn1 dataset}. Comparison of Robust L-BFGS, L-BFGS (multi-batch L-BFGS without enforcing sample consistency), Gradient Descent (multi-batch Gradient method) and SGD. Top part:
we used $\alpha \in \{1, 0.1\}$,
$r\in \{1\%,  5\%,  10\%\}$ and $o=20\%$.
Bottom part: we used $\alpha=1$, $r=1\%$ and
$o\in \{5\%,  10\%, 20\%, 30\%\}$. Solid lines show average performance, and dashed lines show worst and best performance, over 10 runs (per algorithm). $K=4$ MPI processes.}
\label{fig:ijcnn}
\end{figure}

\begin{figure}
\centering
\includegraphics[width=4.6cm]{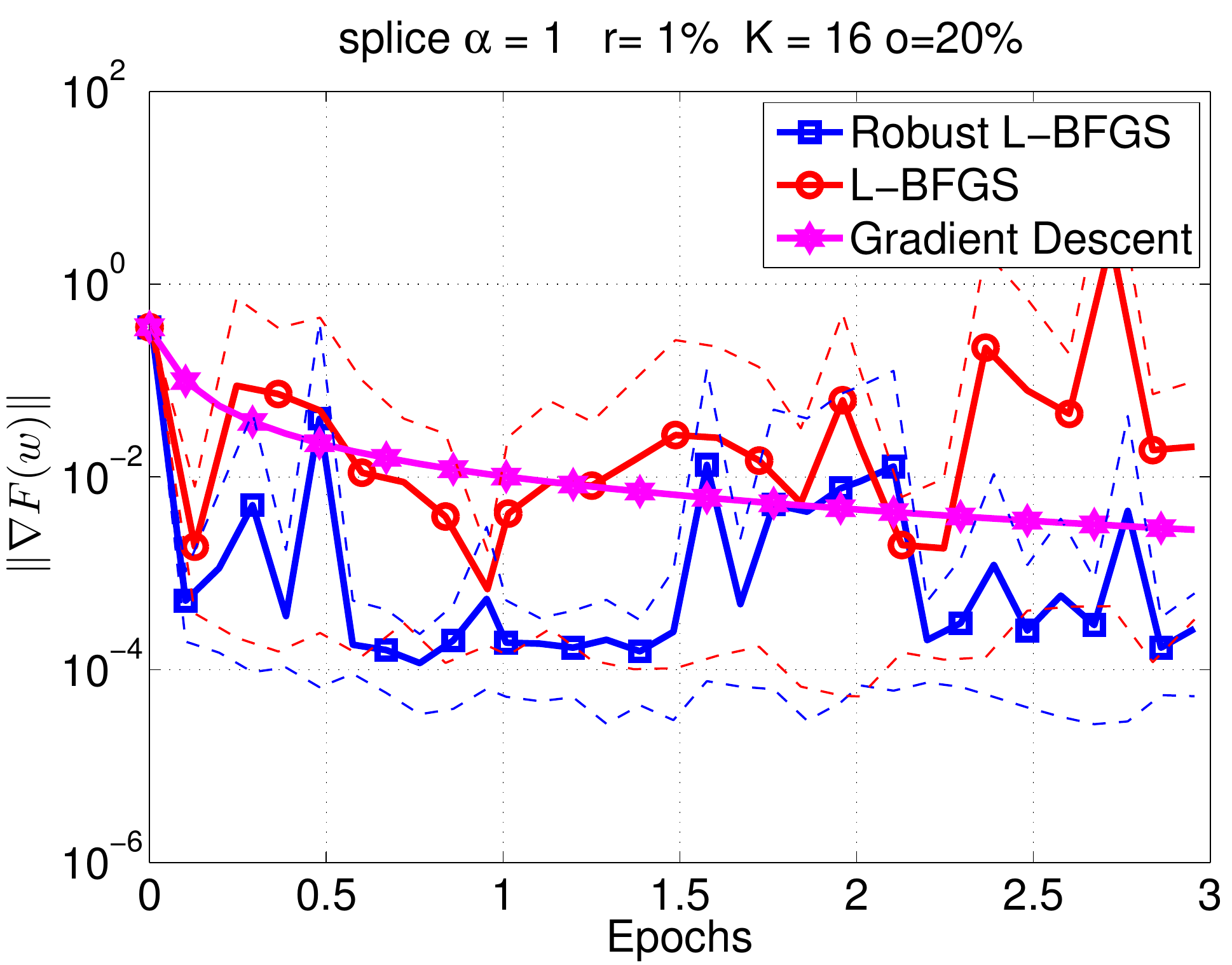}
\includegraphics[width=4.6cm]{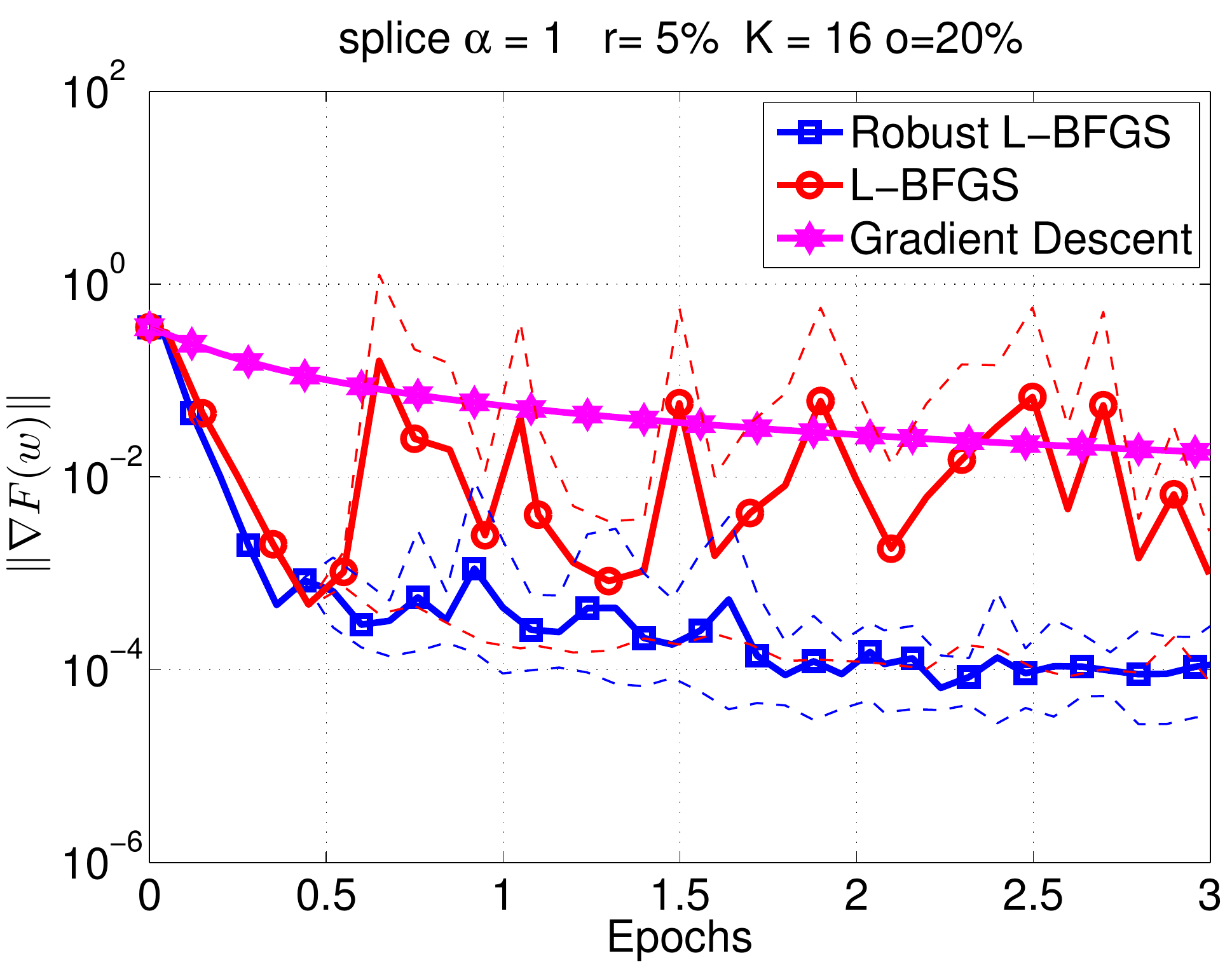}
\includegraphics[width=4.6cm]{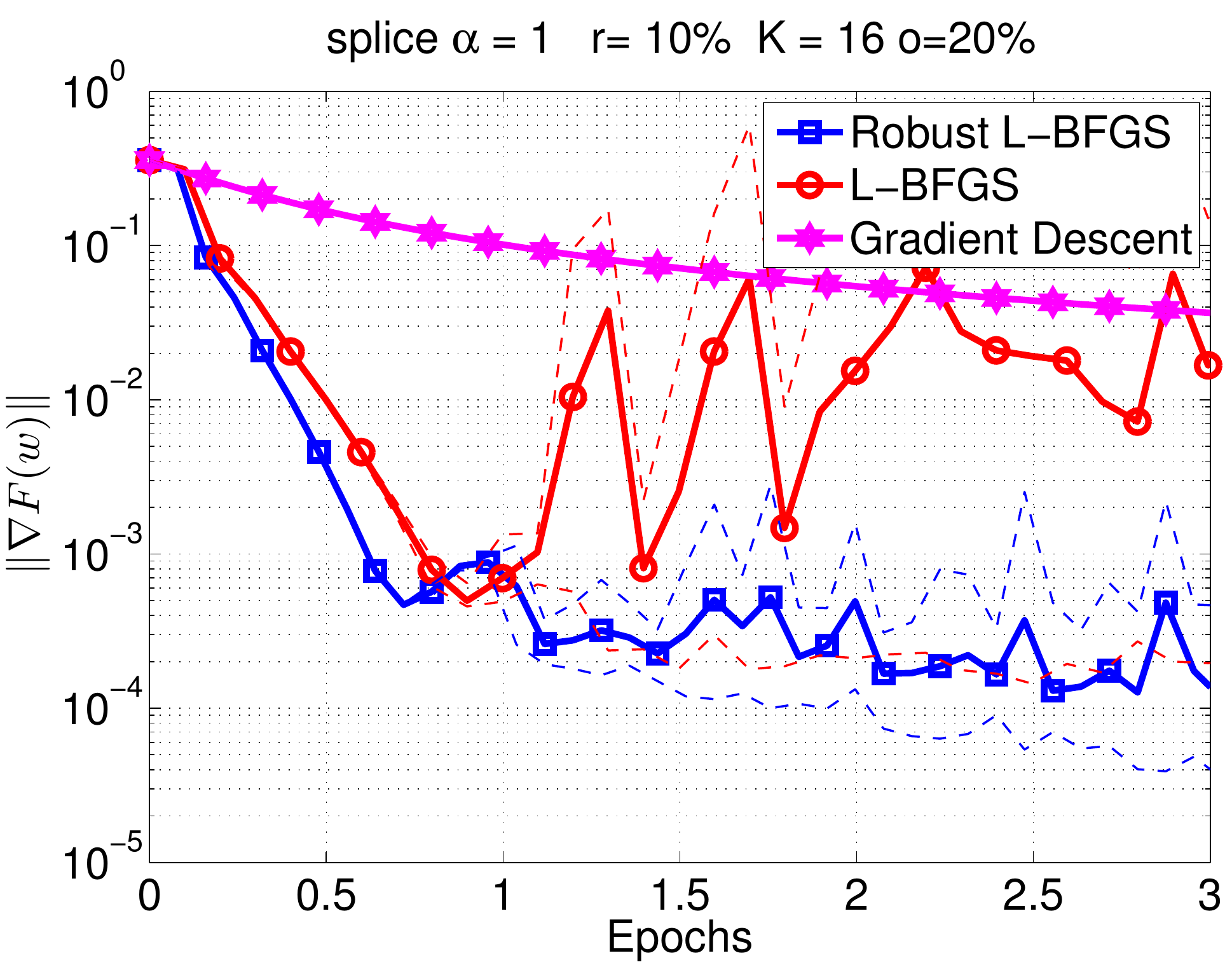}

\includegraphics[width=4.6cm]{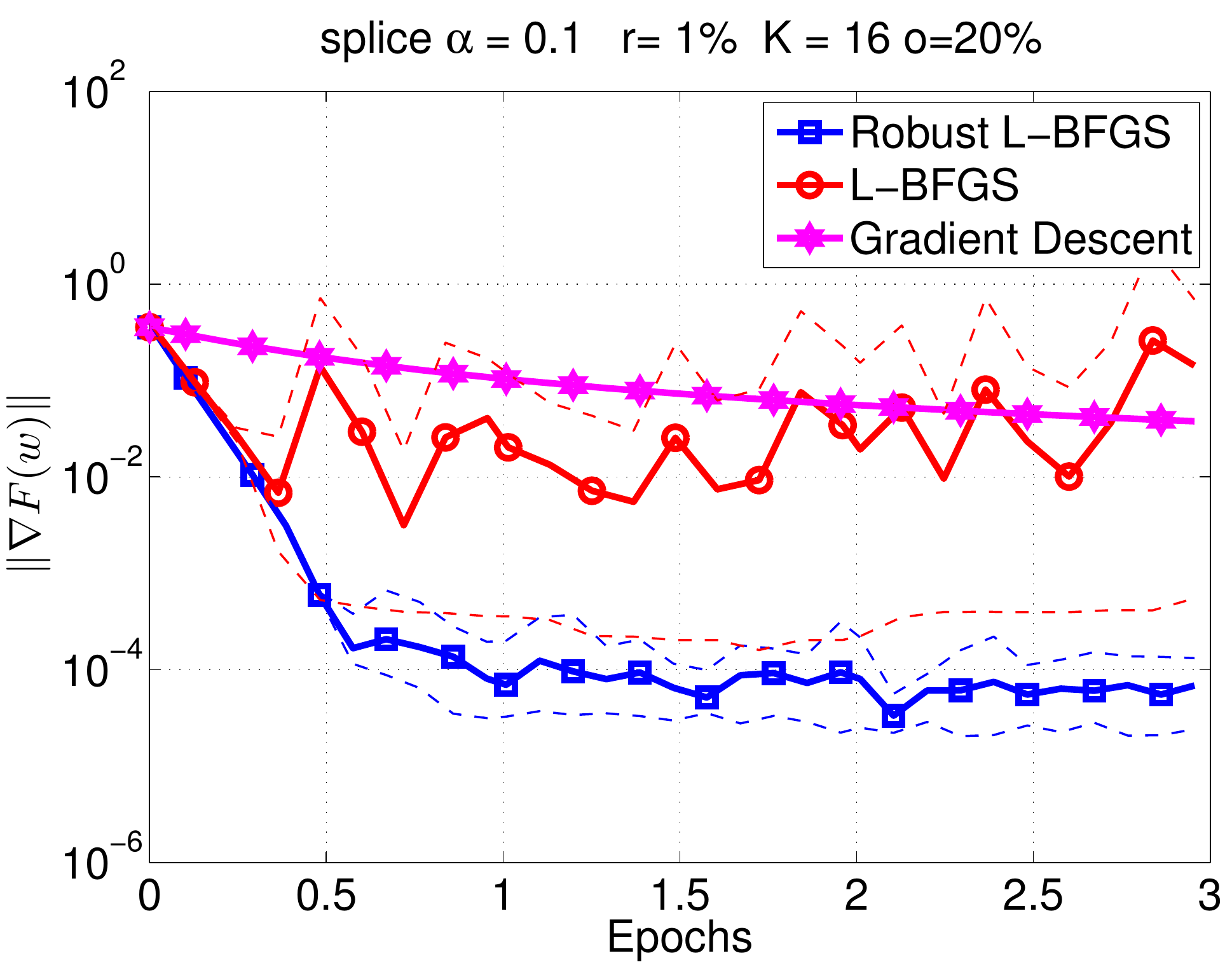}
\includegraphics[width=4.6cm]{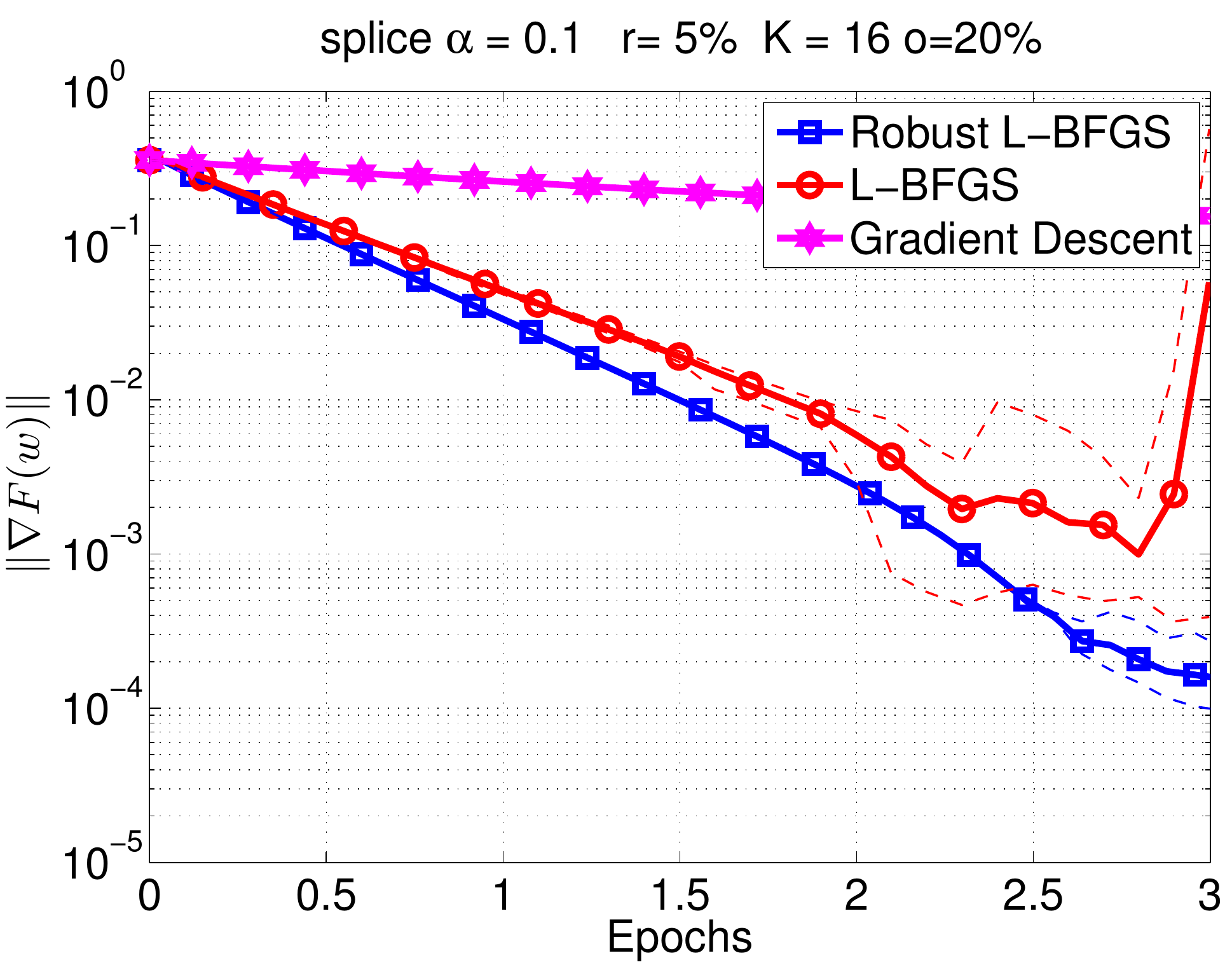}
\includegraphics[width=4.6cm]{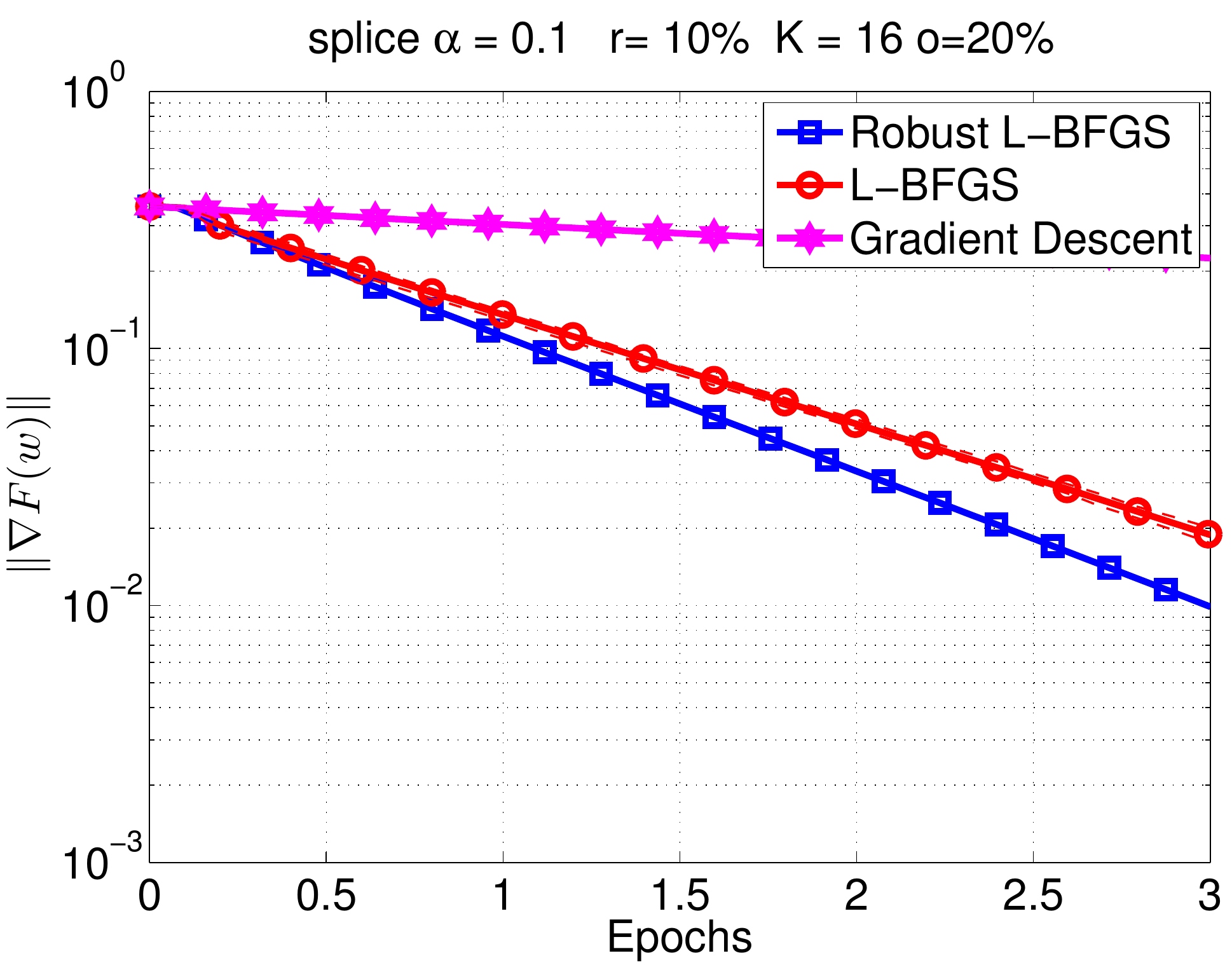}

\hrule 

\includegraphics[width=4.6cm]{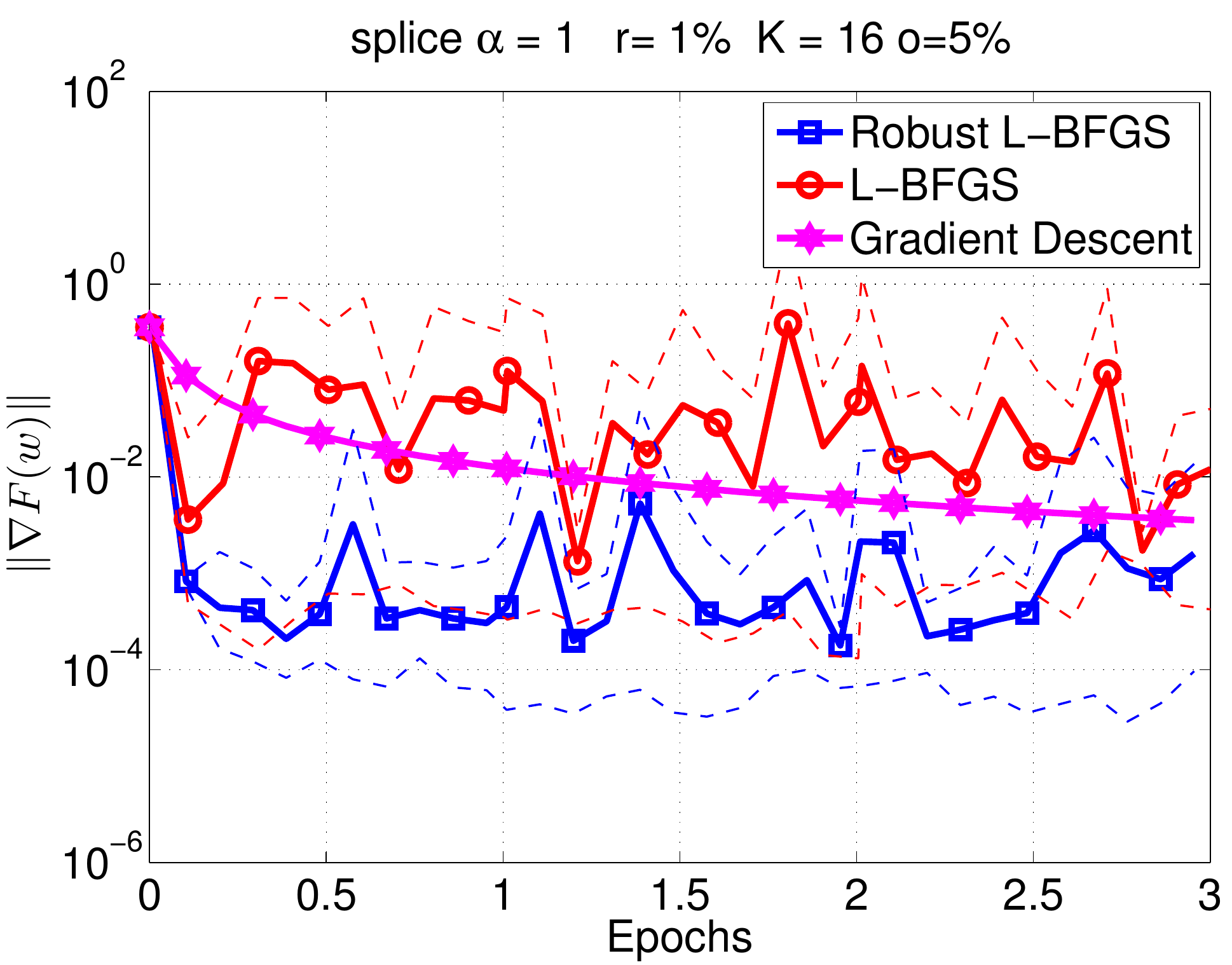}
\includegraphics[width=4.6cm]{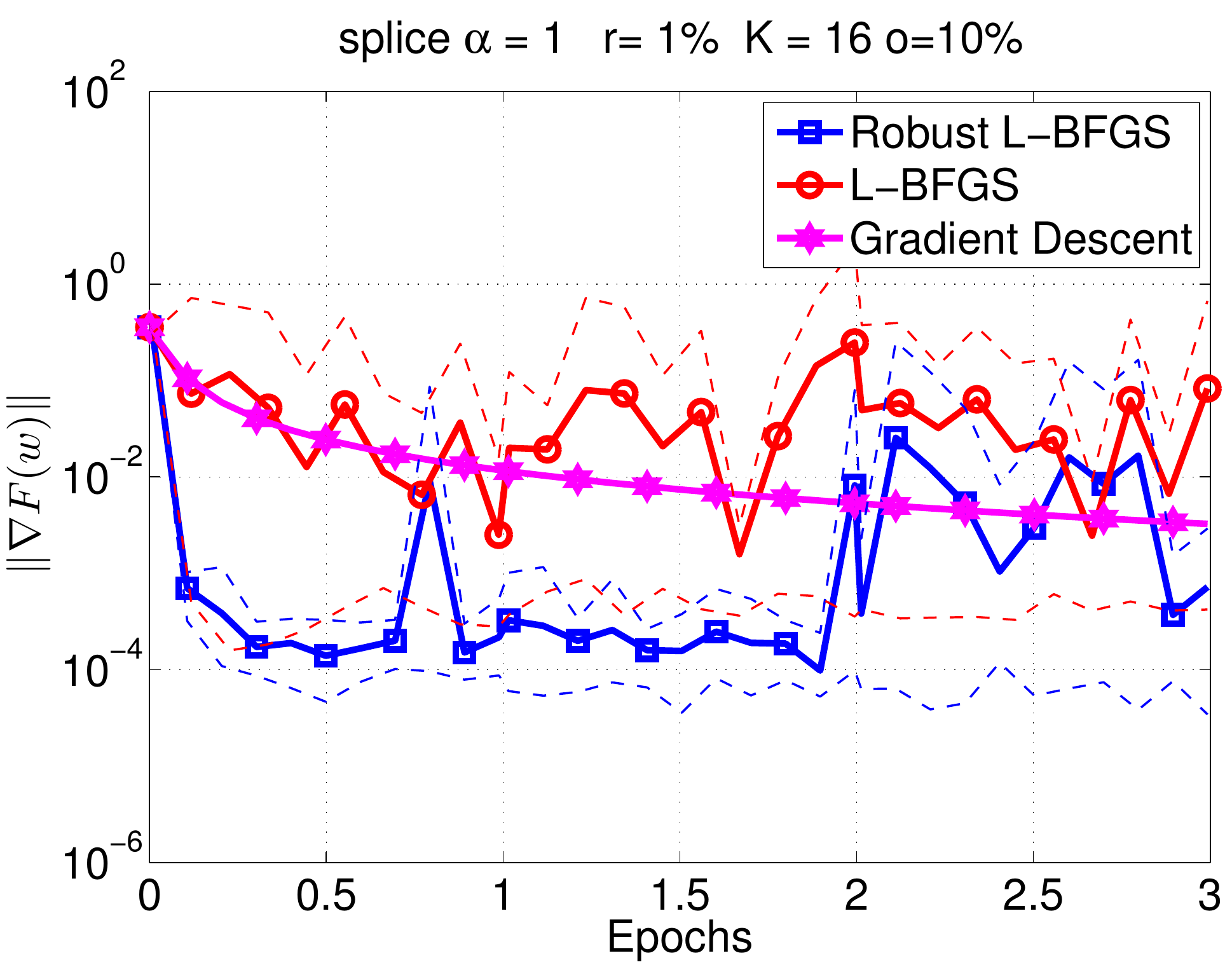}
\includegraphics[width=4.6cm]{splice_mb_1_0_01_16_0_2-eps-converted-to.pdf}
\includegraphics[width=4.6cm]{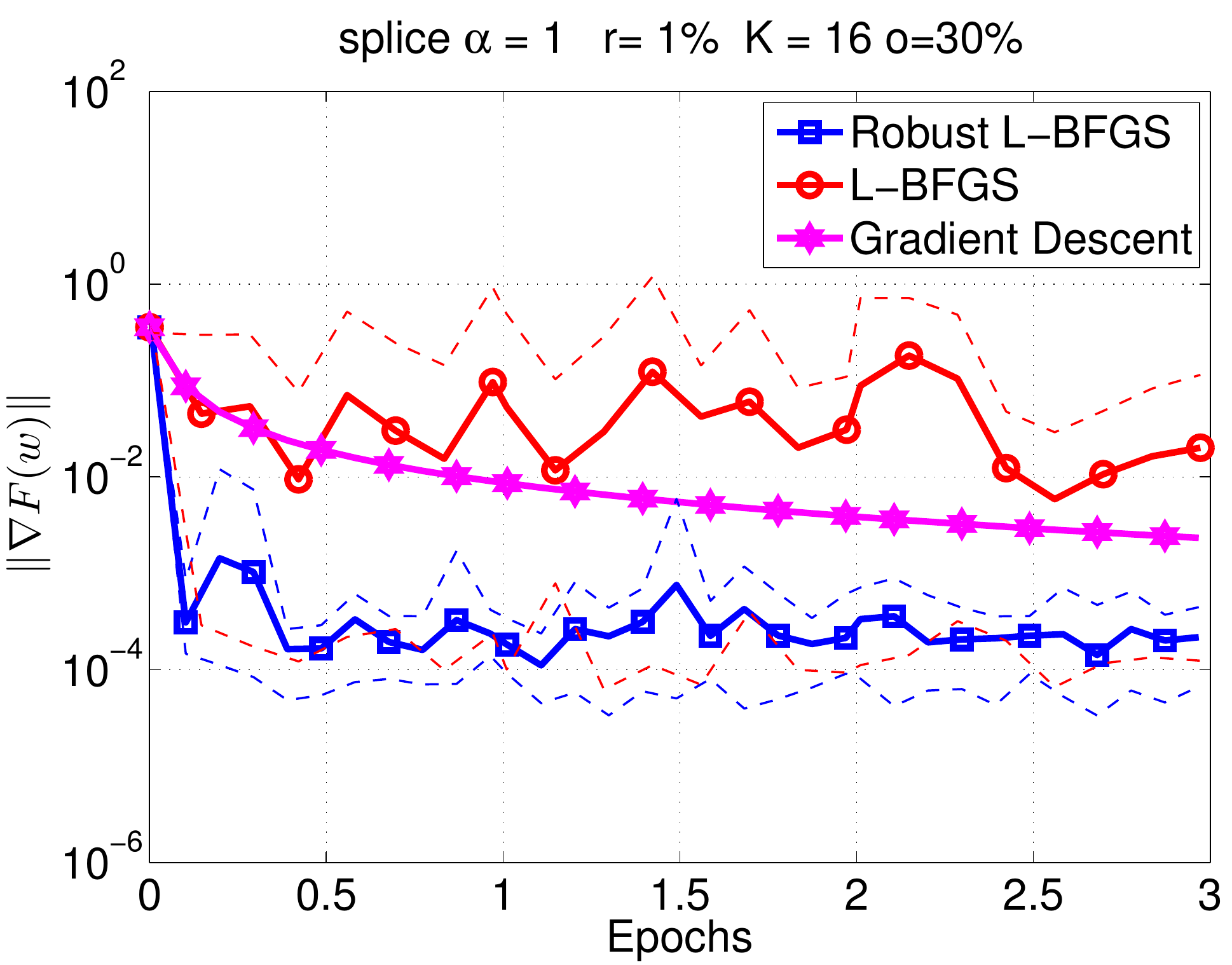}

\caption{\textbf{splice-cite dataset}. Comparison of Robust L-BFGS, L-BFGS (multi-batch L-BFGS without enforcing sample consistency), Gradient Descent (multi-batch Gradient method) and SGD. Top part:
we used $\alpha \in \{1, 0.1\}$,
$r\in \{1\%,  5\%,  10\%\}$ and $o=20\%$.
Bottom part: we used $\alpha=1$, $r=1\%$ and
$o\in \{5\%,  10\%, 20\%, 30\%\}$. Solid lines show average performance, and dashed lines show worst and best performance, over 10 runs (per algorithm). $K=16$ MPI processes. (No Serial SGD experiments due to memory limitations of our cluster.)
}
\label{fig:splice}
\end{figure}

For the experiments in this section (Figures \ref{fig:ijcnn}-\ref{fig:splice}), we run four methods:
\begin{itemize}
	\item (Robust L-BFGS) robust multi-batch L-BFGS (Algorithm \ref{alg:multi}),
	\item (L-BFGS) multi-batch L-BFGS without enforcing sample consistency; gradient differences are computed using different samples, i.e., $y_k = g_{k+1}^{S_{k+1}}-g_{k}^{S_{k}}$,
	\item (Gradient Descent) multi-batch gradient descent; obtained by setting $H_k = I$ in Algorithm \ref{alg:multi},
	\item (SGD) serial SGD; at every iteration one sample is used to compute the gradient.
\end{itemize}
In Figures \ref{fig:ijcnn}-\ref{fig:splice} we 
show the evolution of $\|\nabla F(w)\|$ for different step lengths $\alpha$, and for various batch ($\left| S\right| = r\cdot n $) and overlap ($\left| O\right| = o \cdot \left| S\right| $) sizes. Each Figure (\ref{fig:ijcnn}-\ref{fig:splice}) consists of 10 plots that illustrate the performance of the methods with the following parameters:
\begin{itemize}
	\item Top 3 plots: $\alpha=1$, $o=20\%$ and $r=1\%,5\%,10\%$
	\item Middle 3 plots: $\alpha=0.1$, $o=20\%$ and $r=1\%,5\%,10\%$
	\item Bottom 4 plots: $\alpha=1$, $r=1\%$ and $o=5\%,10\%,20\%,30\%$
\end{itemize}
 As is expected for quasi-Newton methods, robust L-BFGS performs best with a step-size $\alpha=1$, for the most part.


\begin{figure}
\centering
\includegraphics[width=4.6cm]{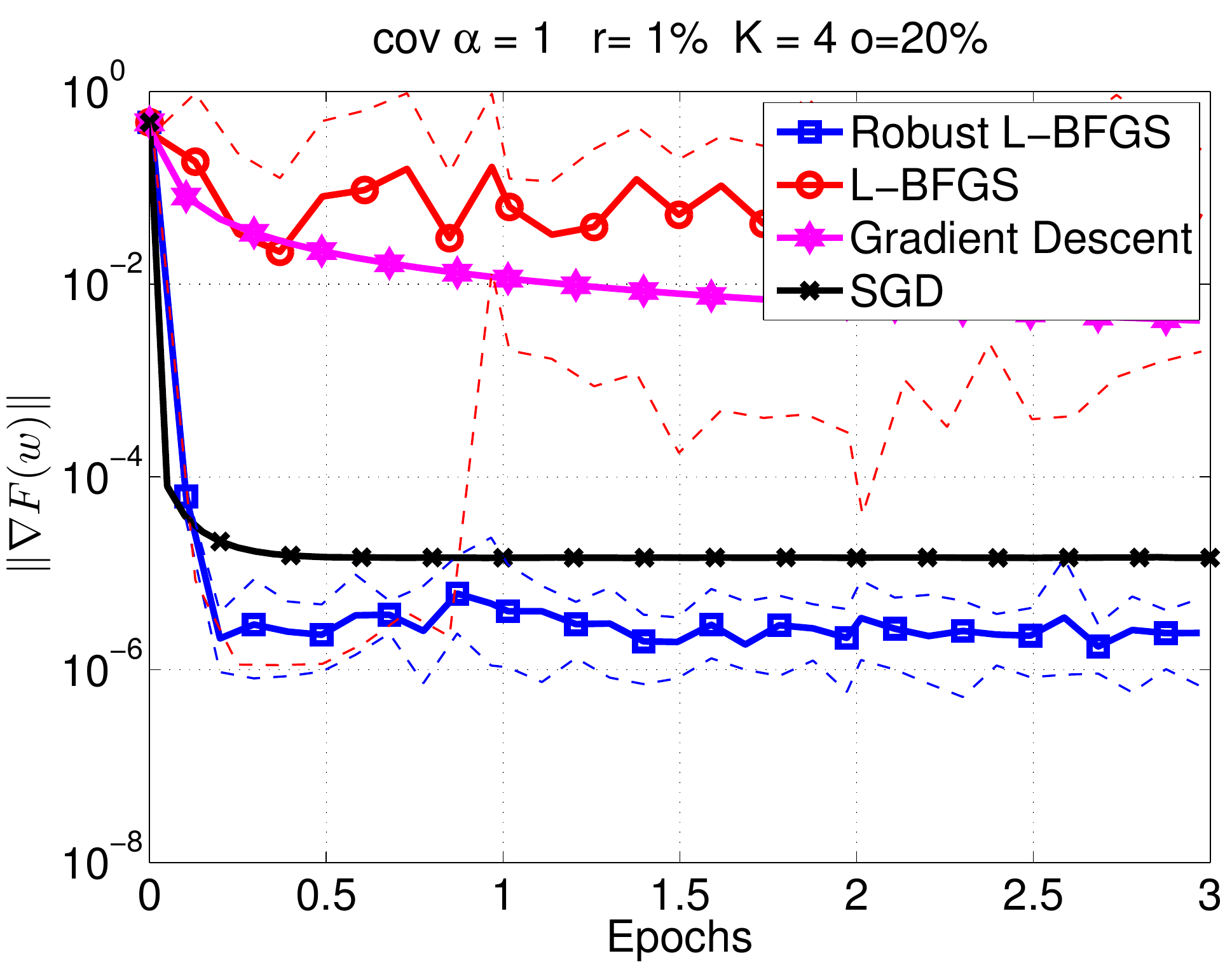}
\includegraphics[width=4.6cm]{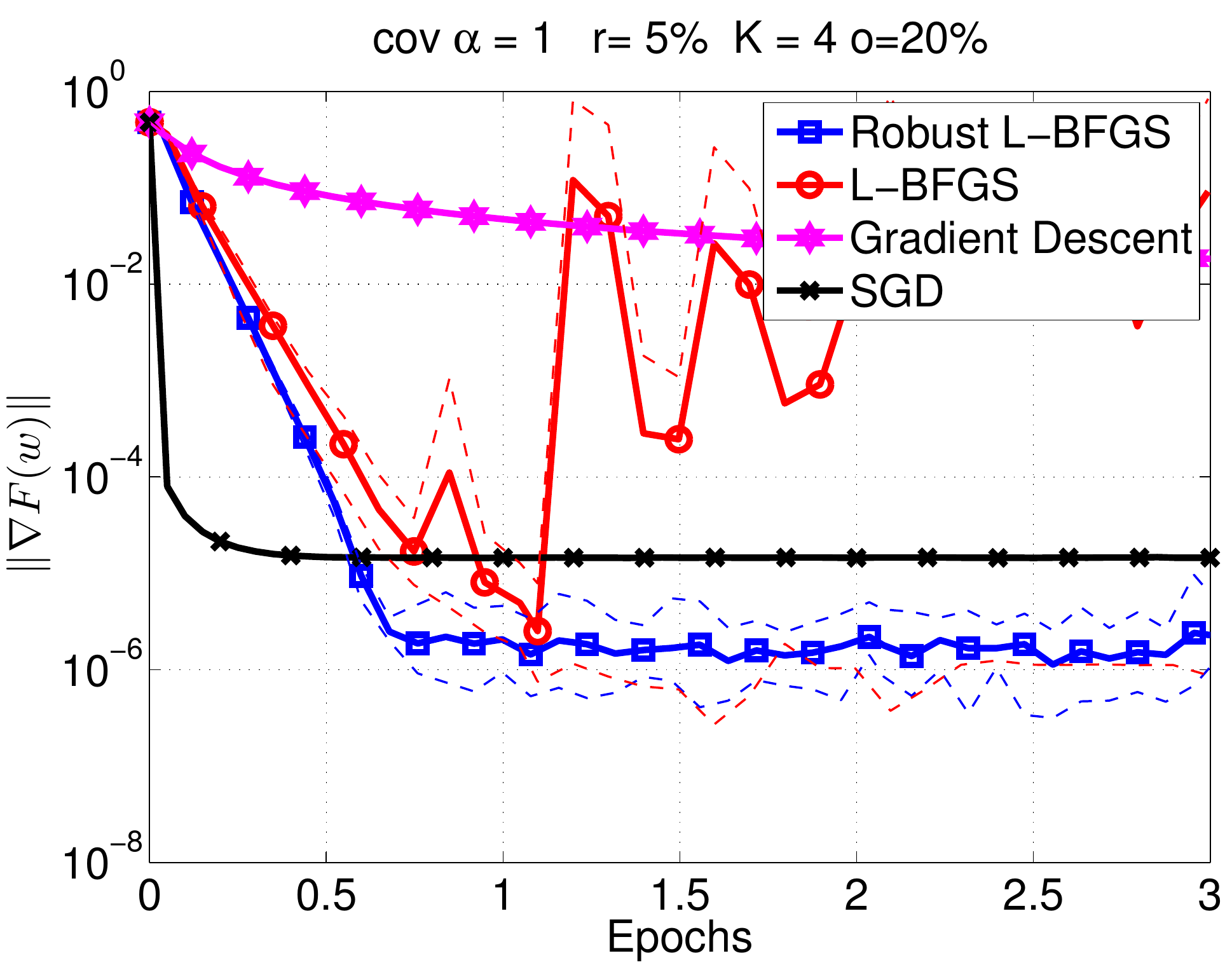}
\includegraphics[width=4.6cm]{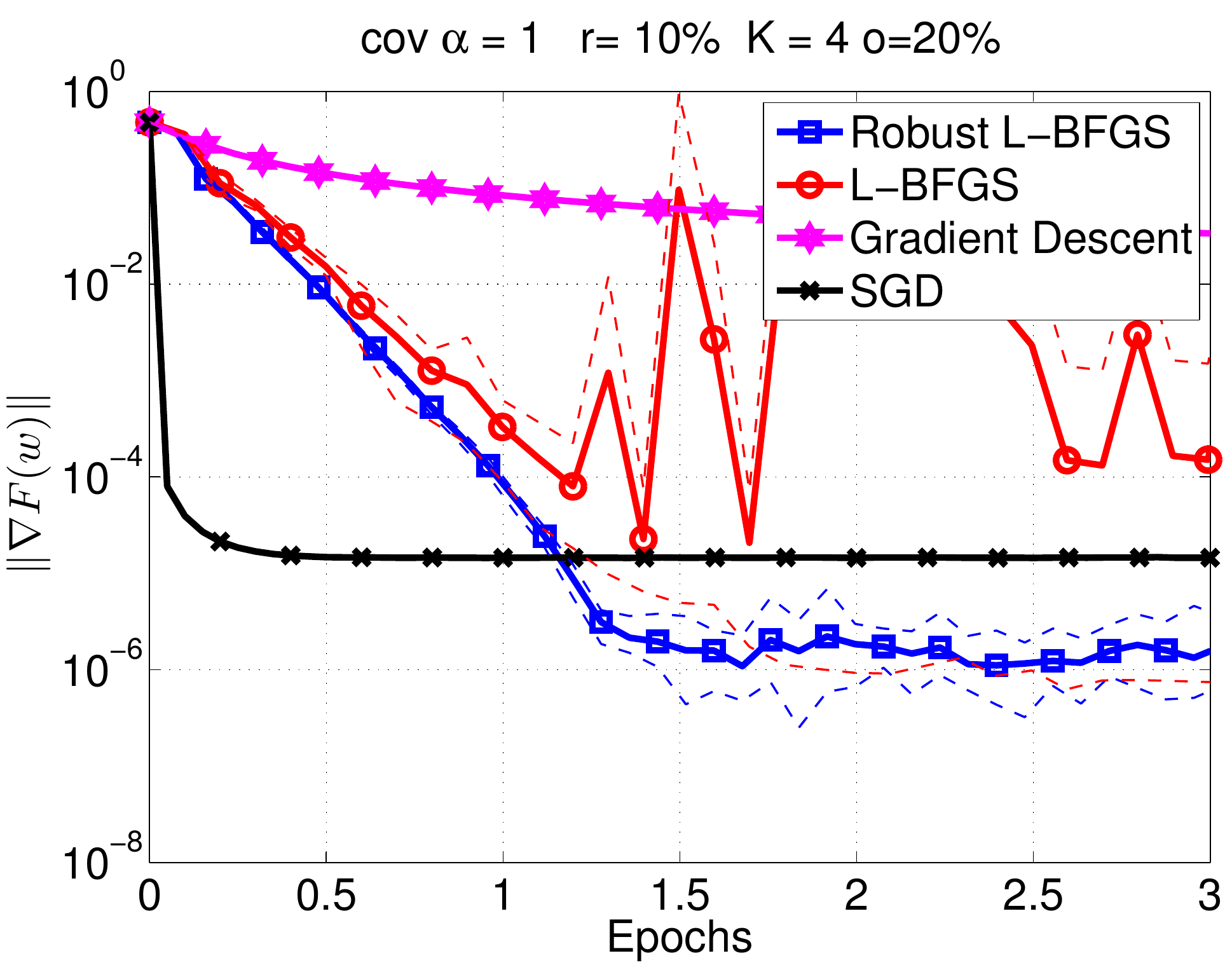}

\includegraphics[width=4.6cm]{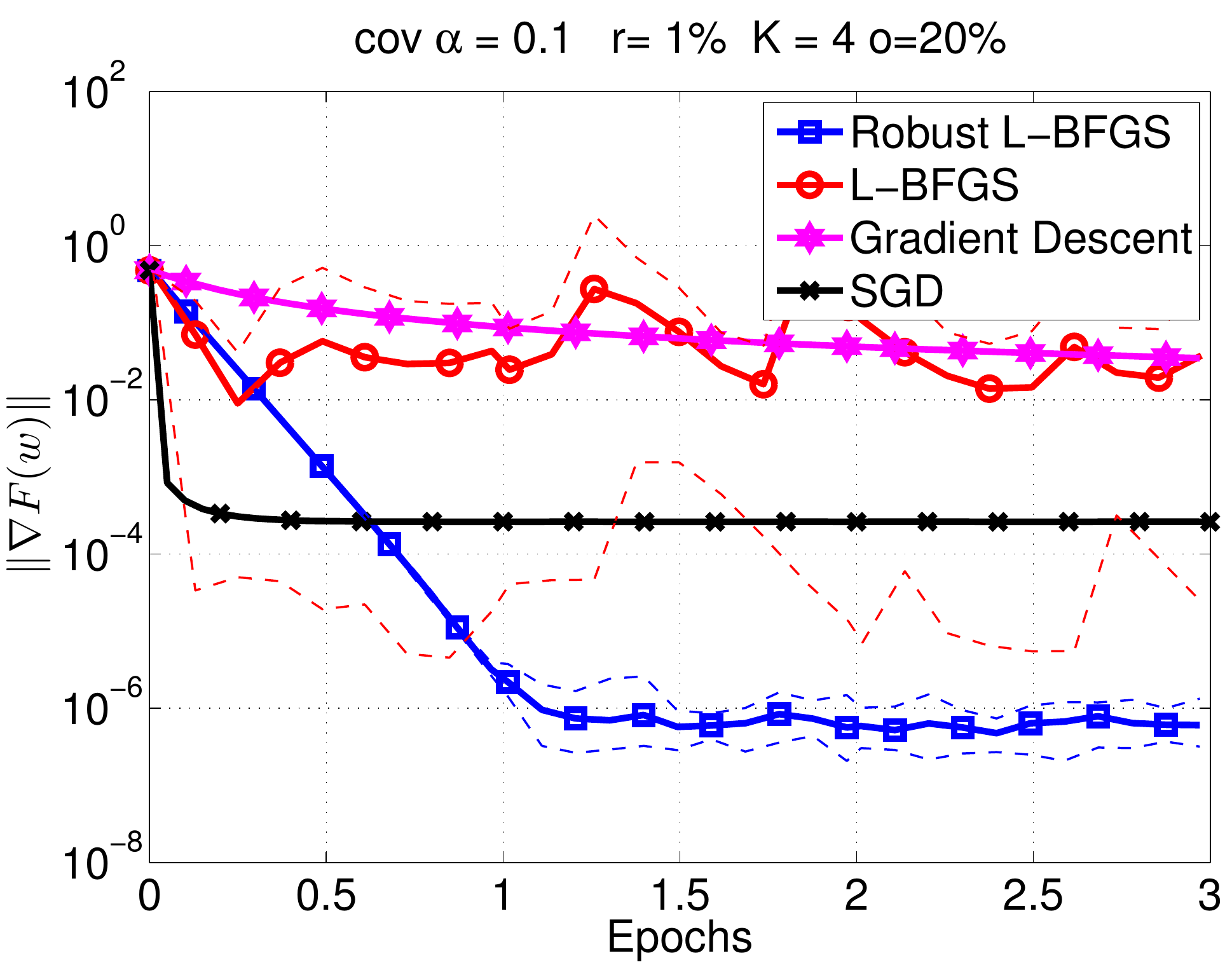}
\includegraphics[width=4.6cm]{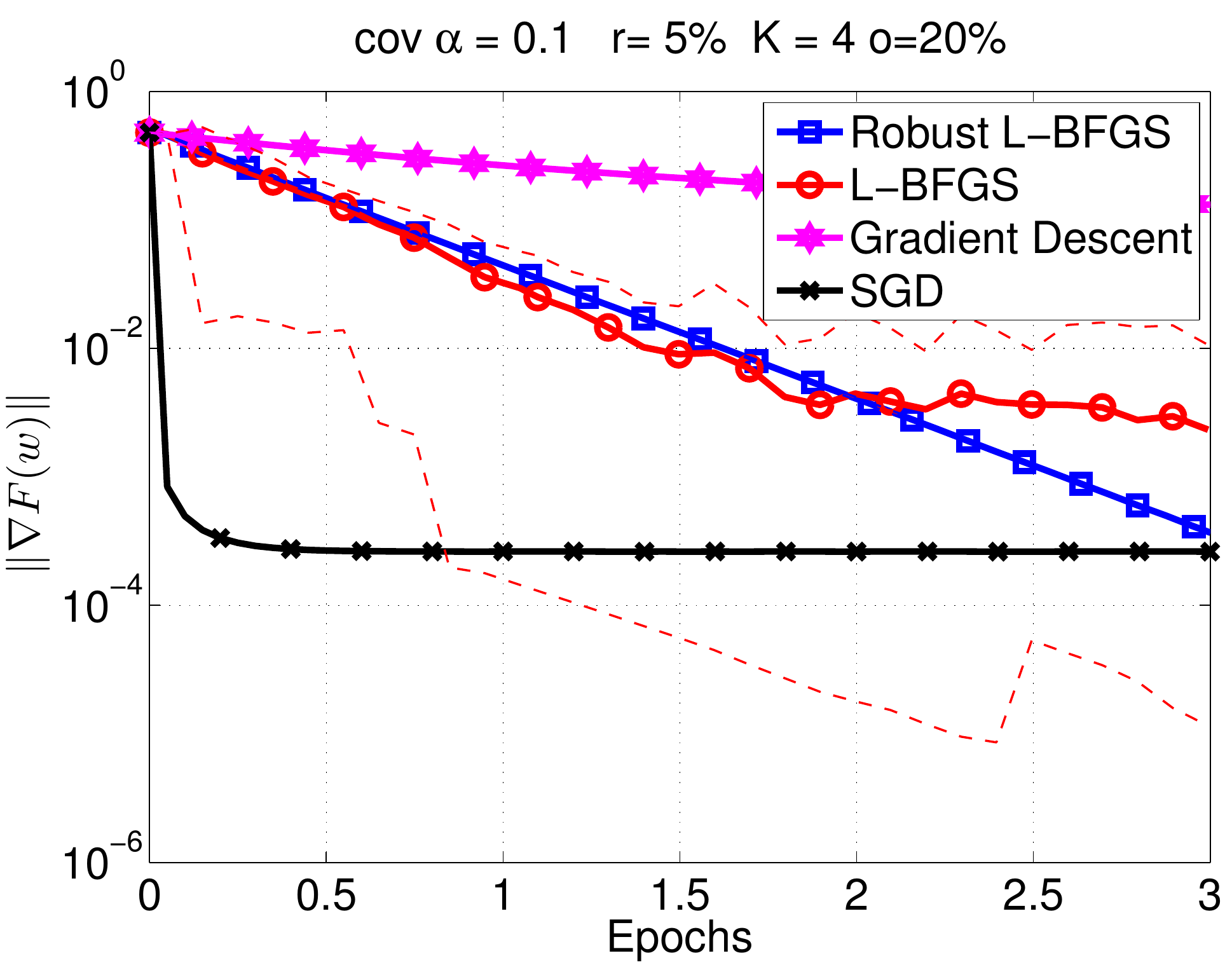}
\includegraphics[width=4.6cm]{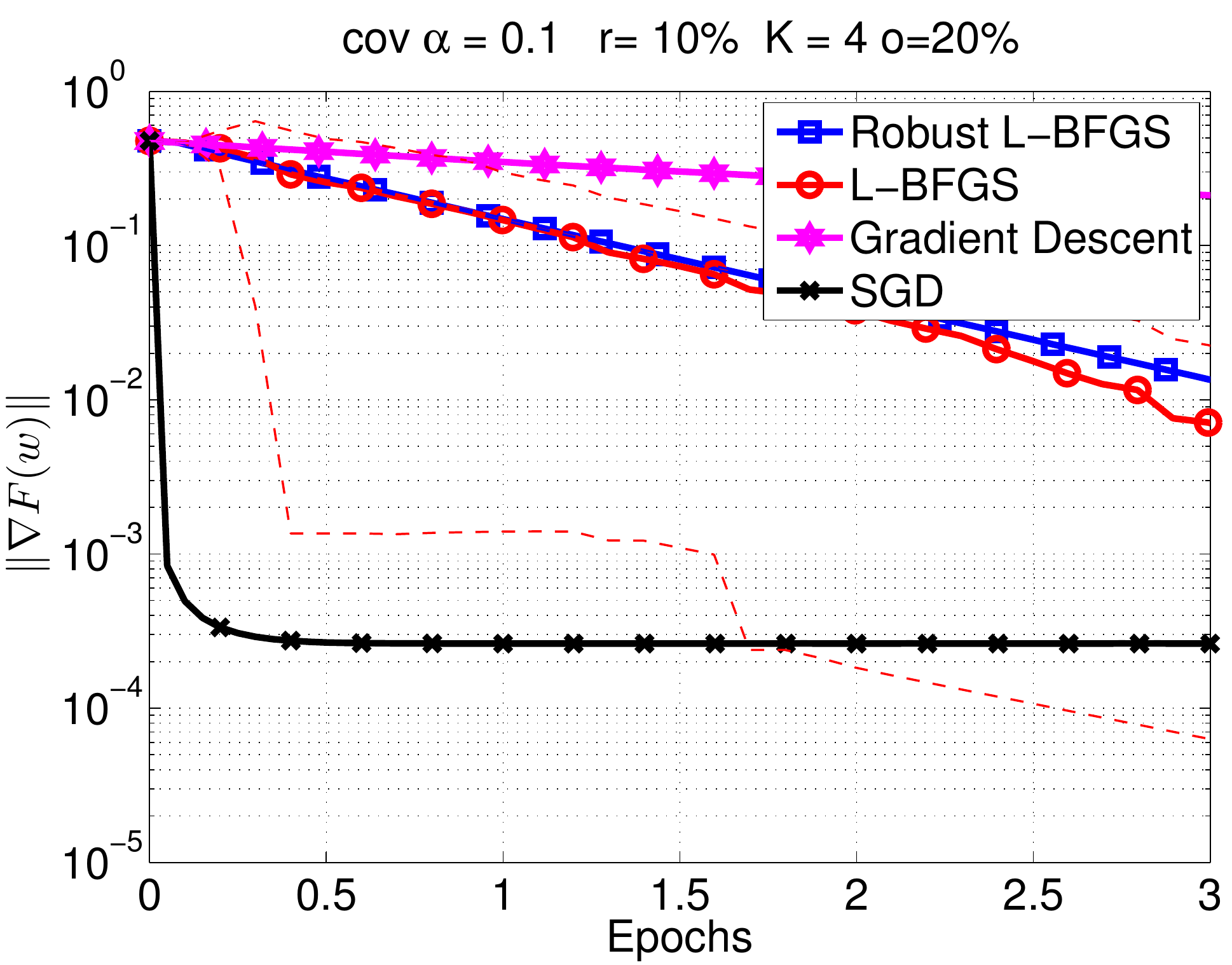}

 \hrule 

\includegraphics[width=4.6cm]{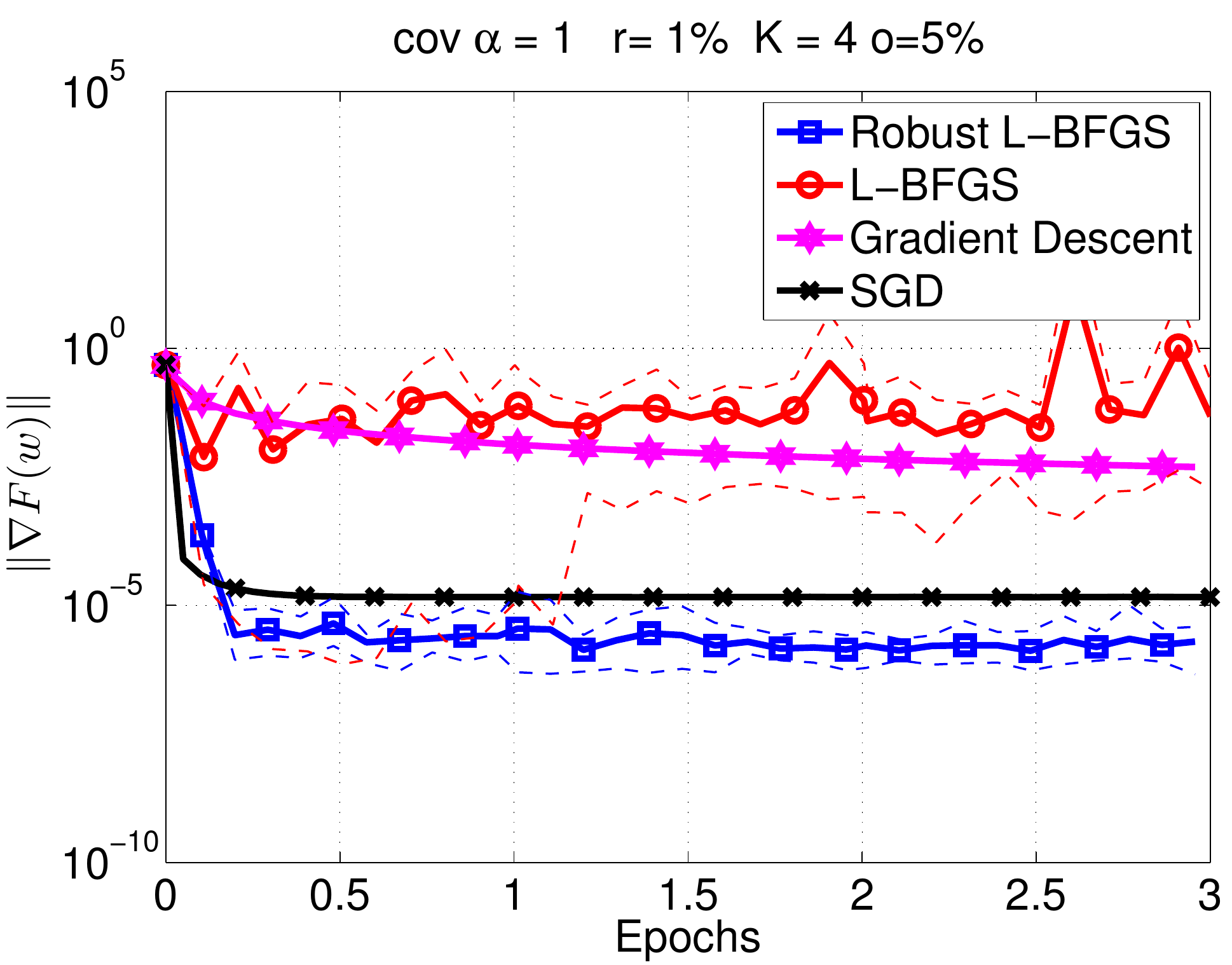}
\includegraphics[width=4.6cm]{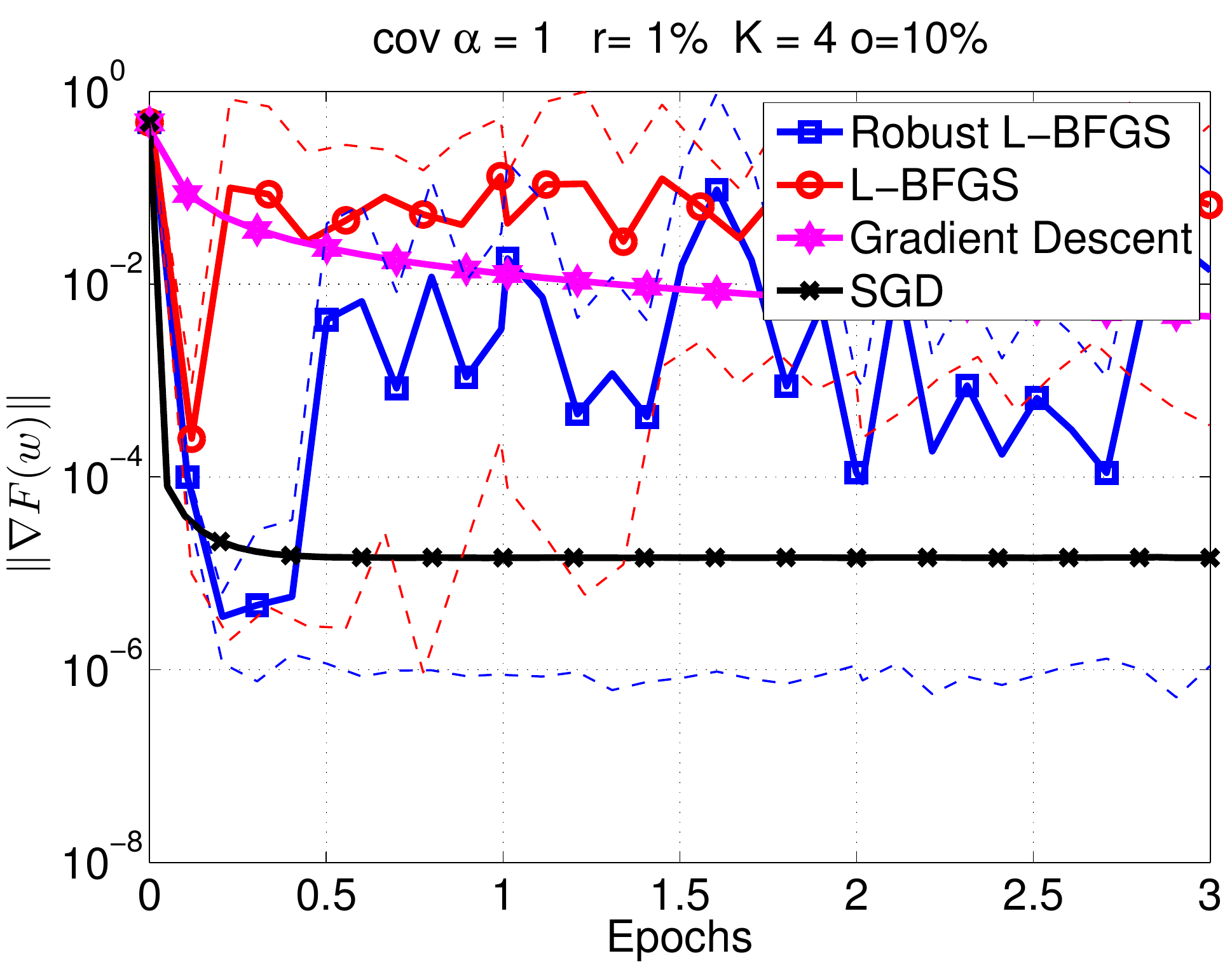}
\includegraphics[width=4.6cm]{cov_mb_1_0_01_4_0_2-eps-converted-to.pdf}
\includegraphics[width=4.6cm]{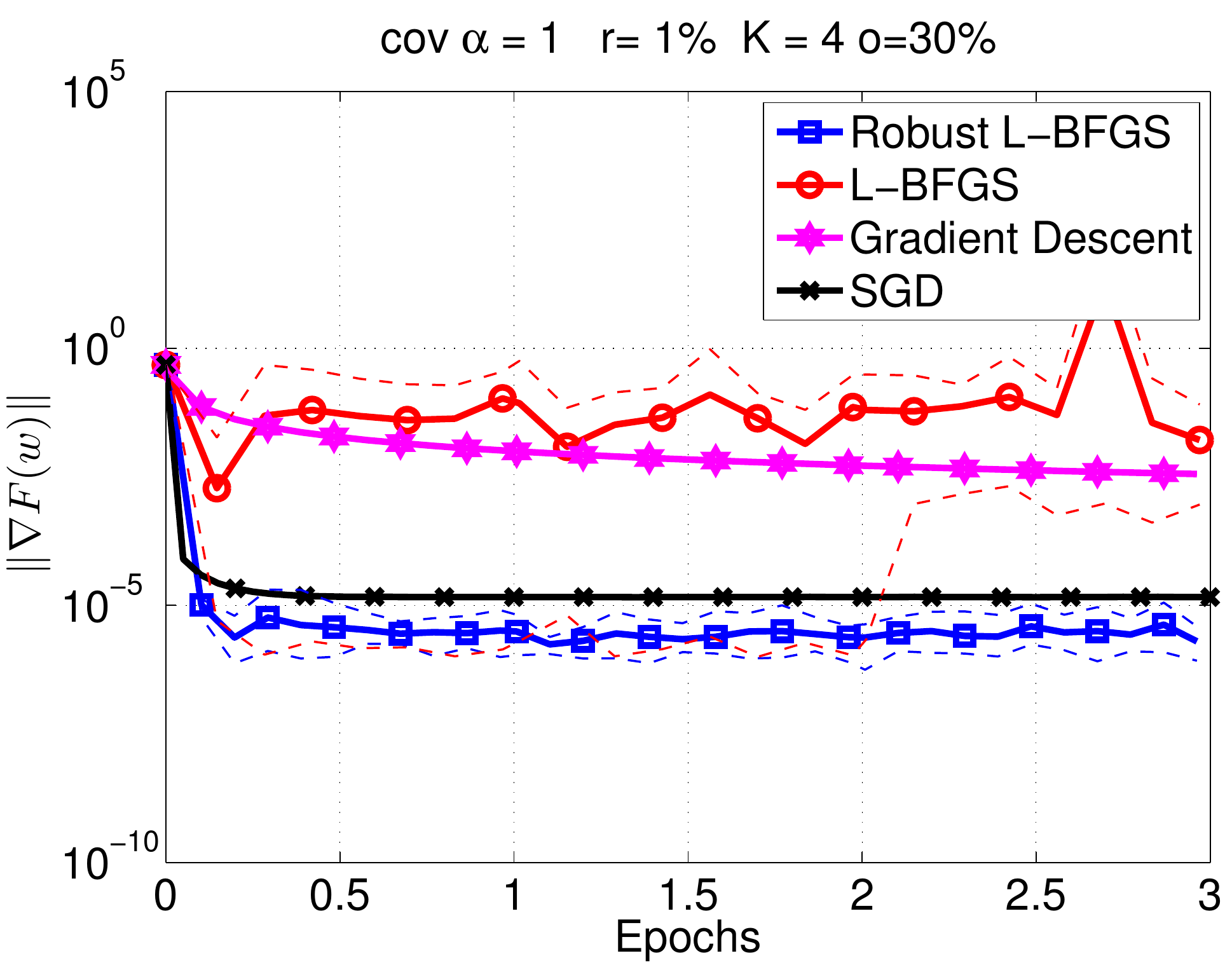}

\caption{\textbf{cov dataset}. Comparison of Robust L-BFGS, L-BFGS (multi-batch L-BFGS without enforcing sample consistency), Gradient Descent (multi-batch Gradient method) and SGD. Top part:
we used $\alpha \in \{1, 0.1\}$,
$r\in \{1\%,  5\%,  10\%\}$ and $o=20\%$.
Bottom part: we used $\alpha=1$, $r=1\%$ and
$o\in \{5\%,  10\%, 20\%, 30\%\}$. Solid lines show average performance, and dashed lines show worst and best performance, over 10 runs (per algorithm). $K=4$ MPI processes.}
\end{figure}

\begin{figure}
\centering
\includegraphics[width=4.6cm]{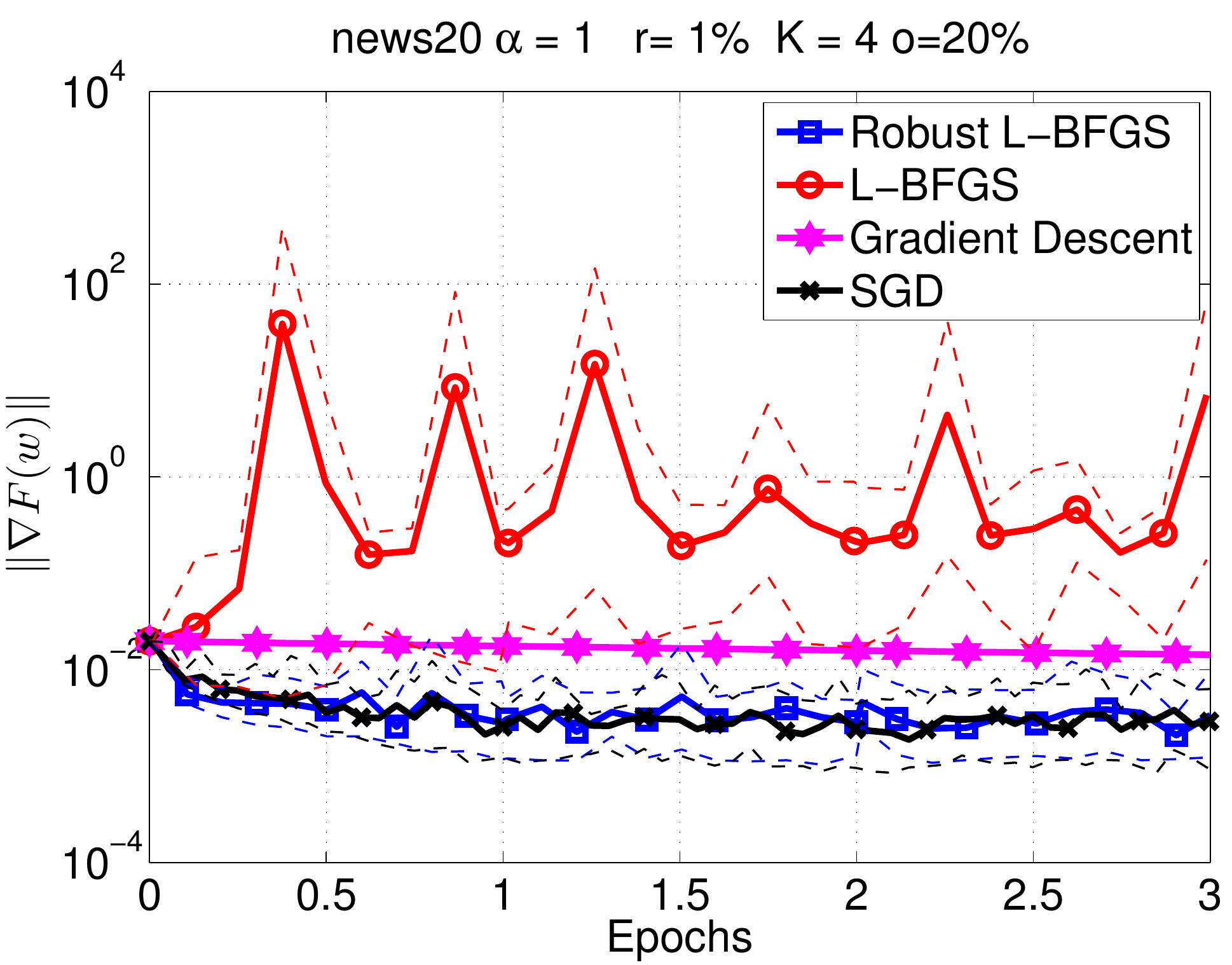}
\includegraphics[width=4.6cm]{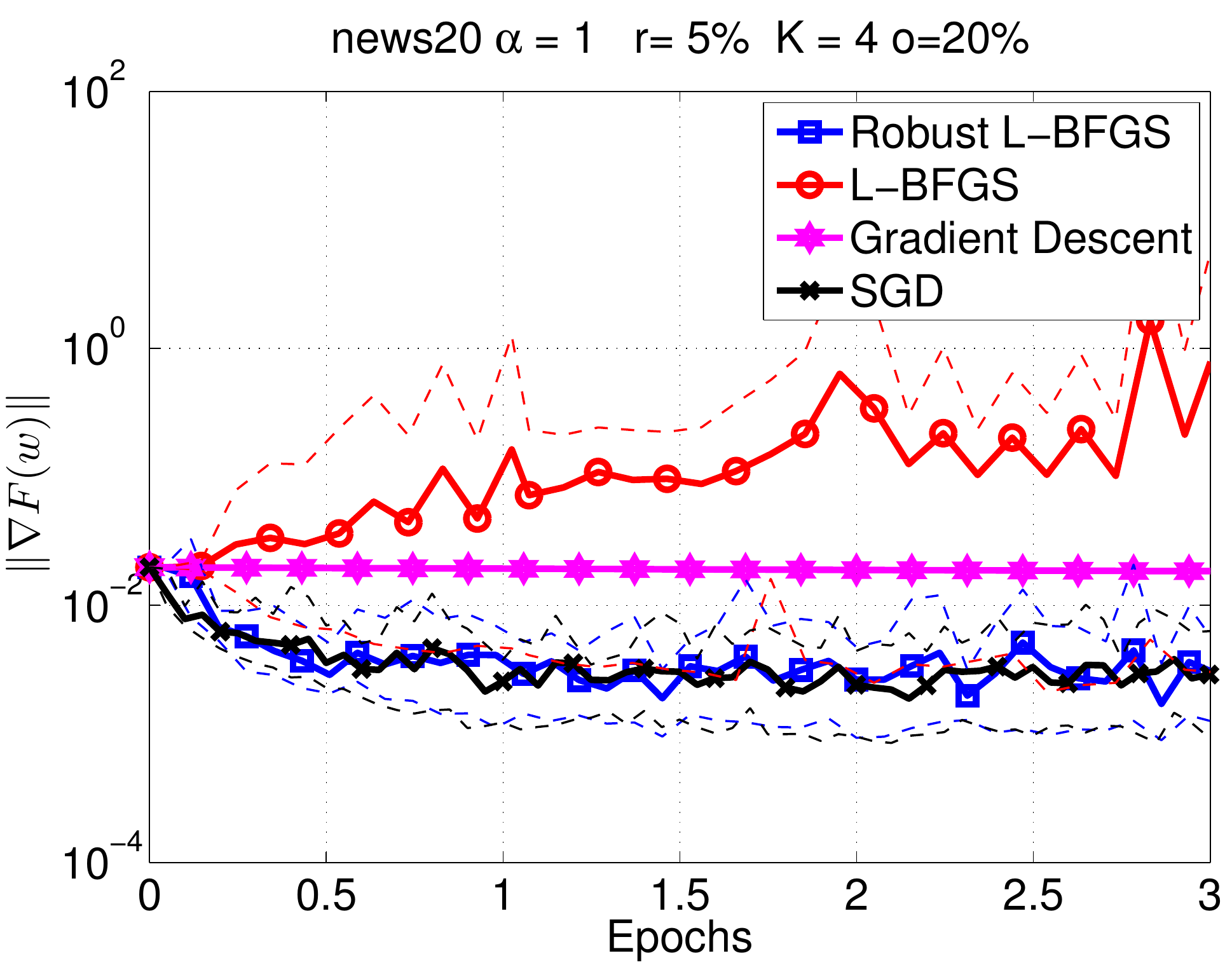}
\includegraphics[width=4.6cm]{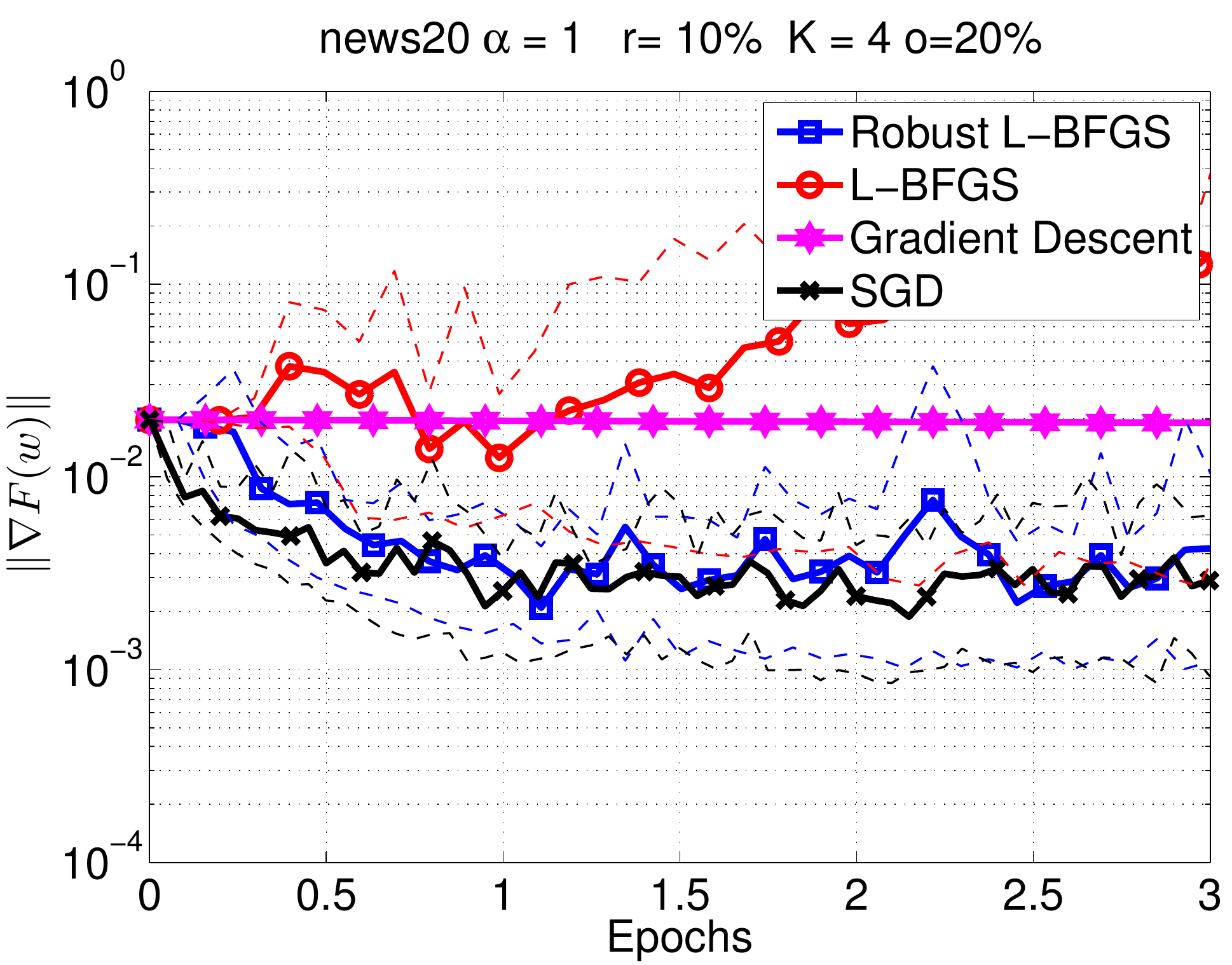}

\includegraphics[width=4.6cm]{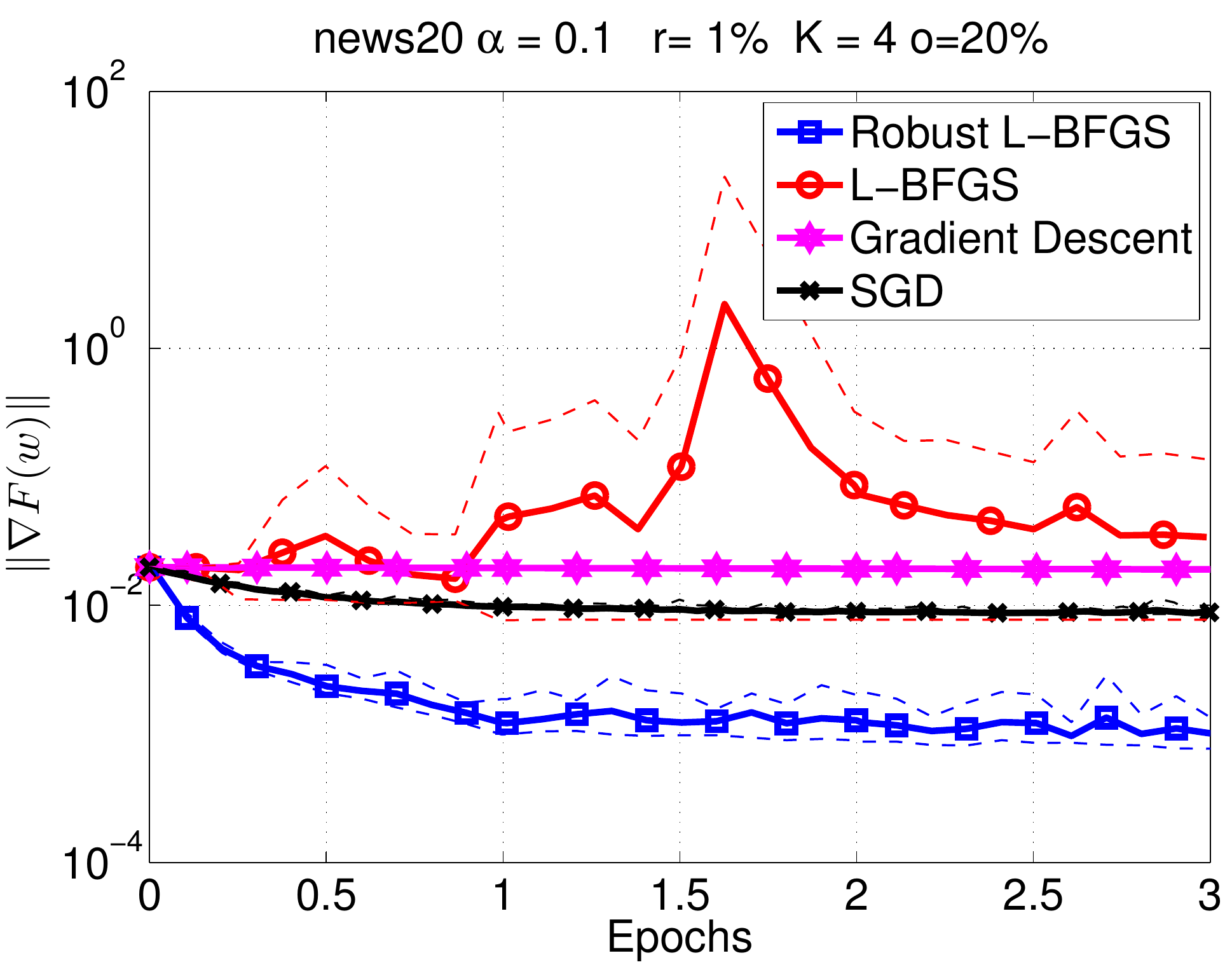}
\includegraphics[width=4.6cm]{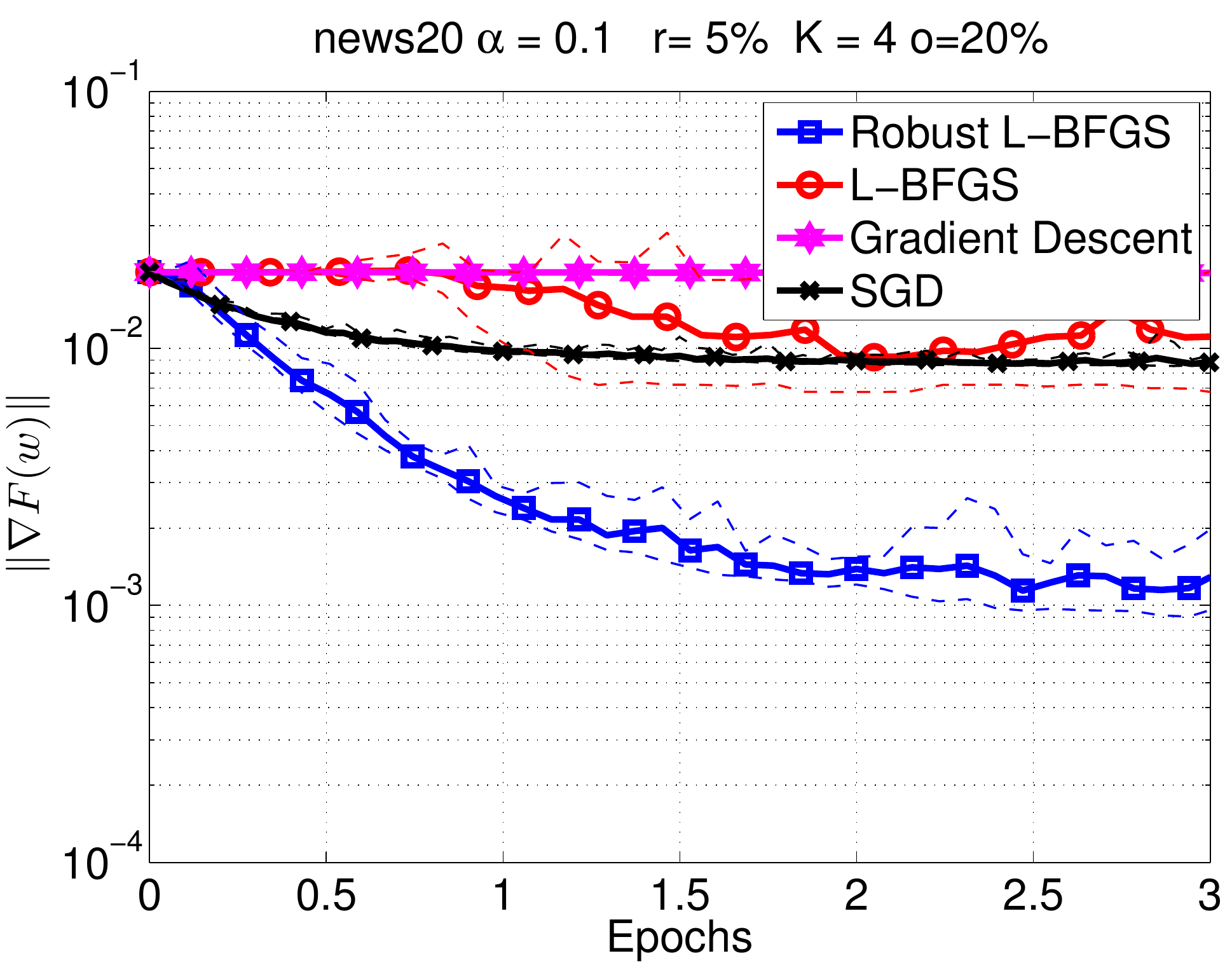}
\includegraphics[width=4.6cm]{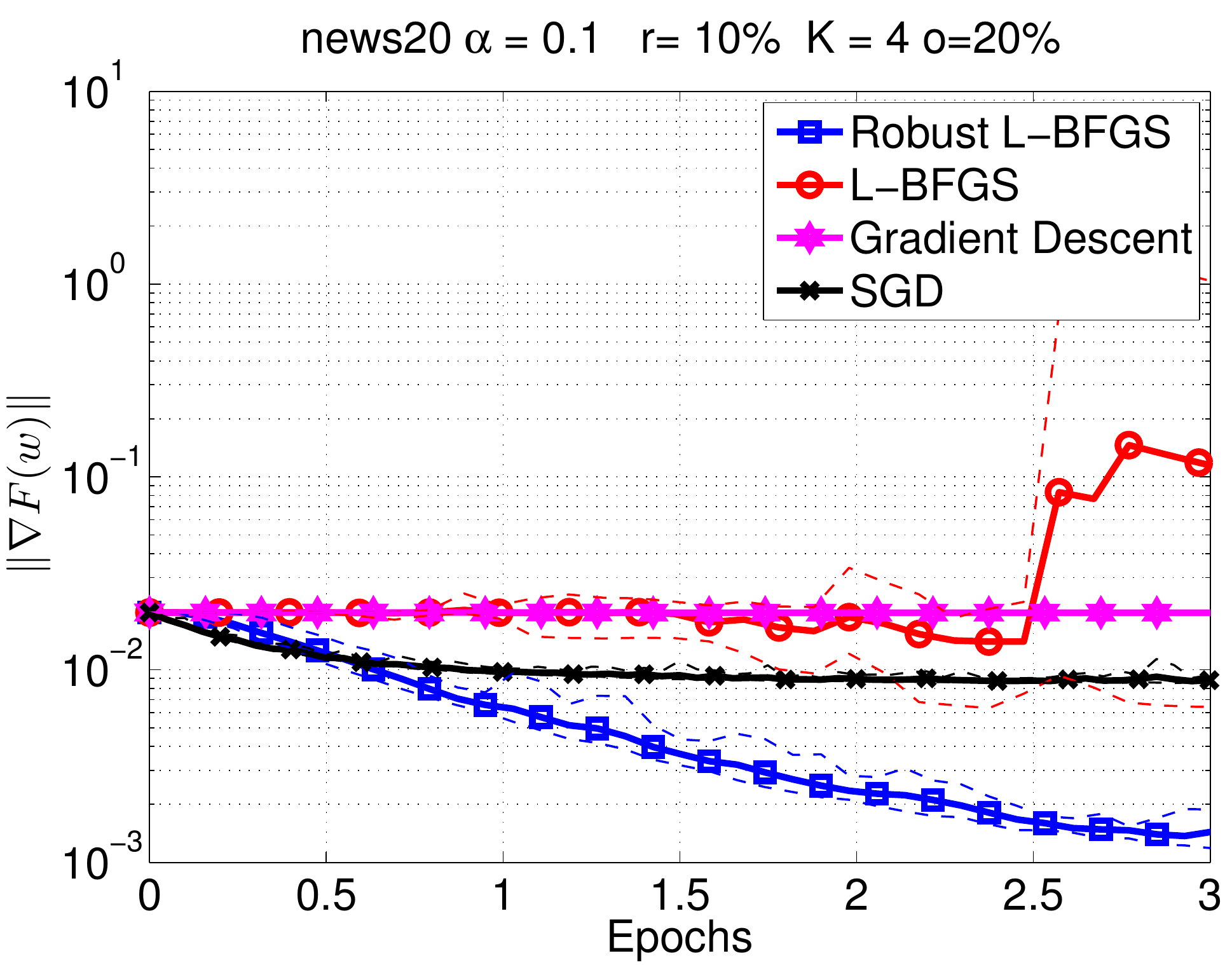}

\hrule 

\includegraphics[width=4.6cm]{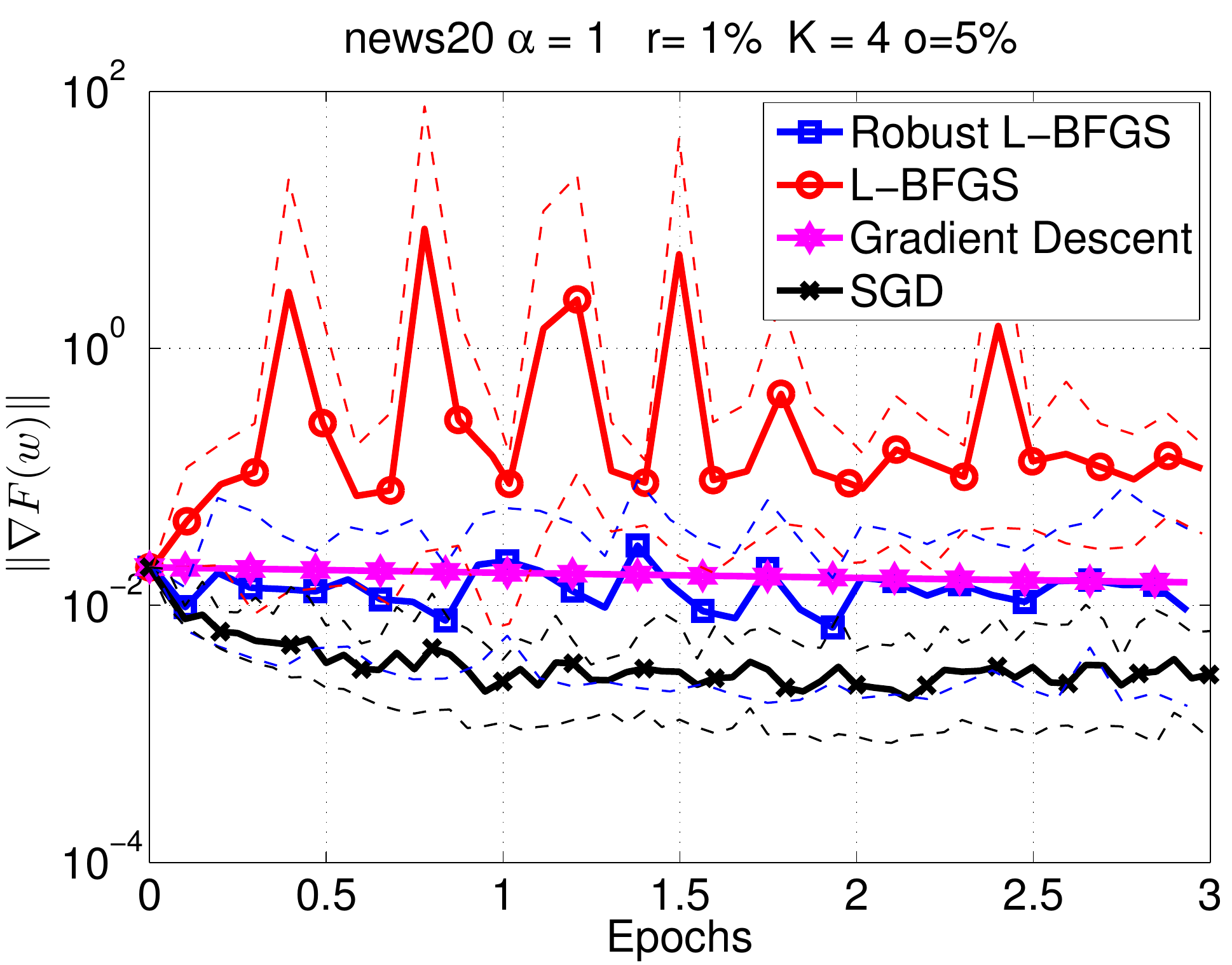}
\includegraphics[width=4.6cm]{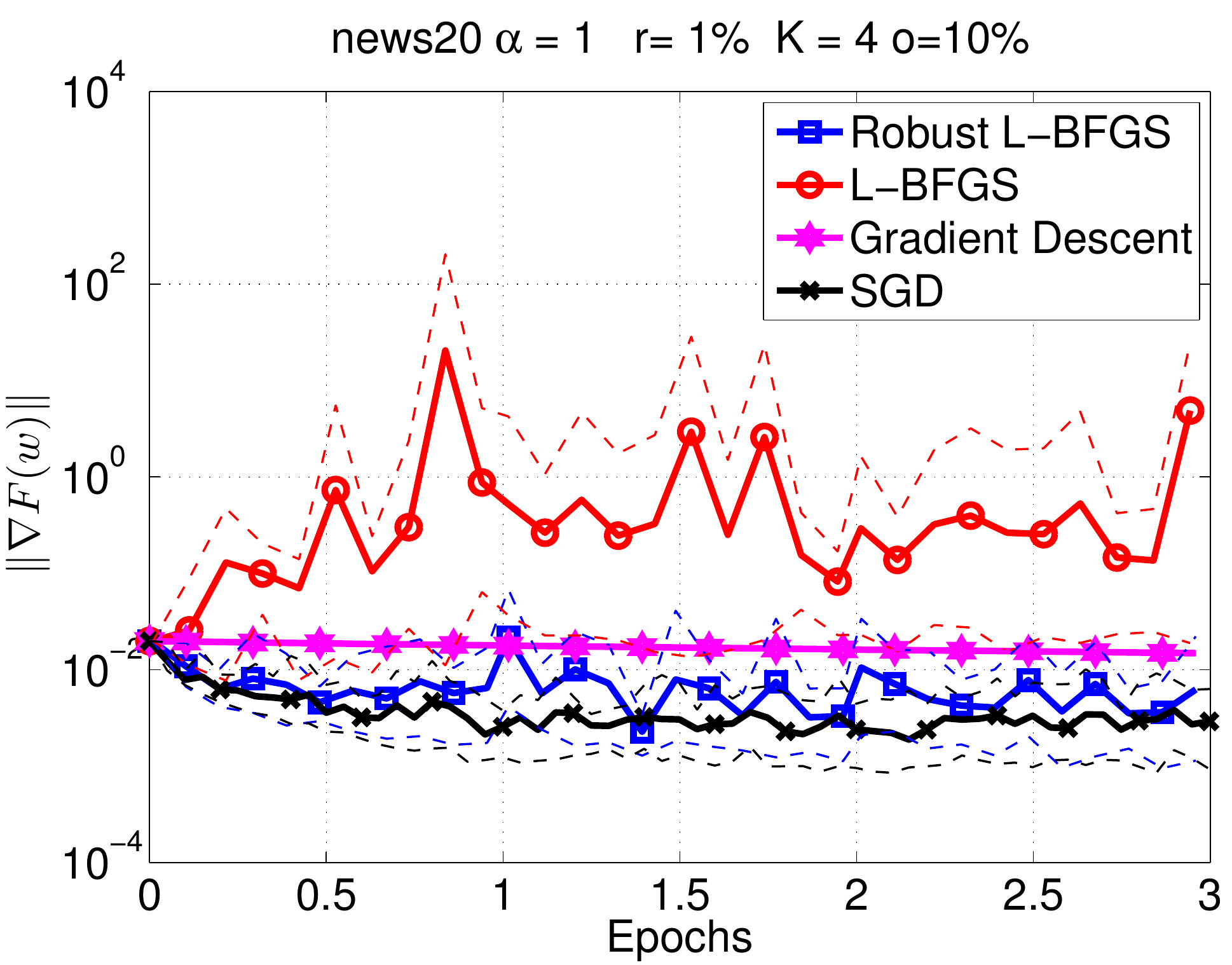}
\includegraphics[width=4.6cm]{news20_mb_1_0_01_4_0_2-eps-converted-to.pdf}
\includegraphics[width=4.6cm]{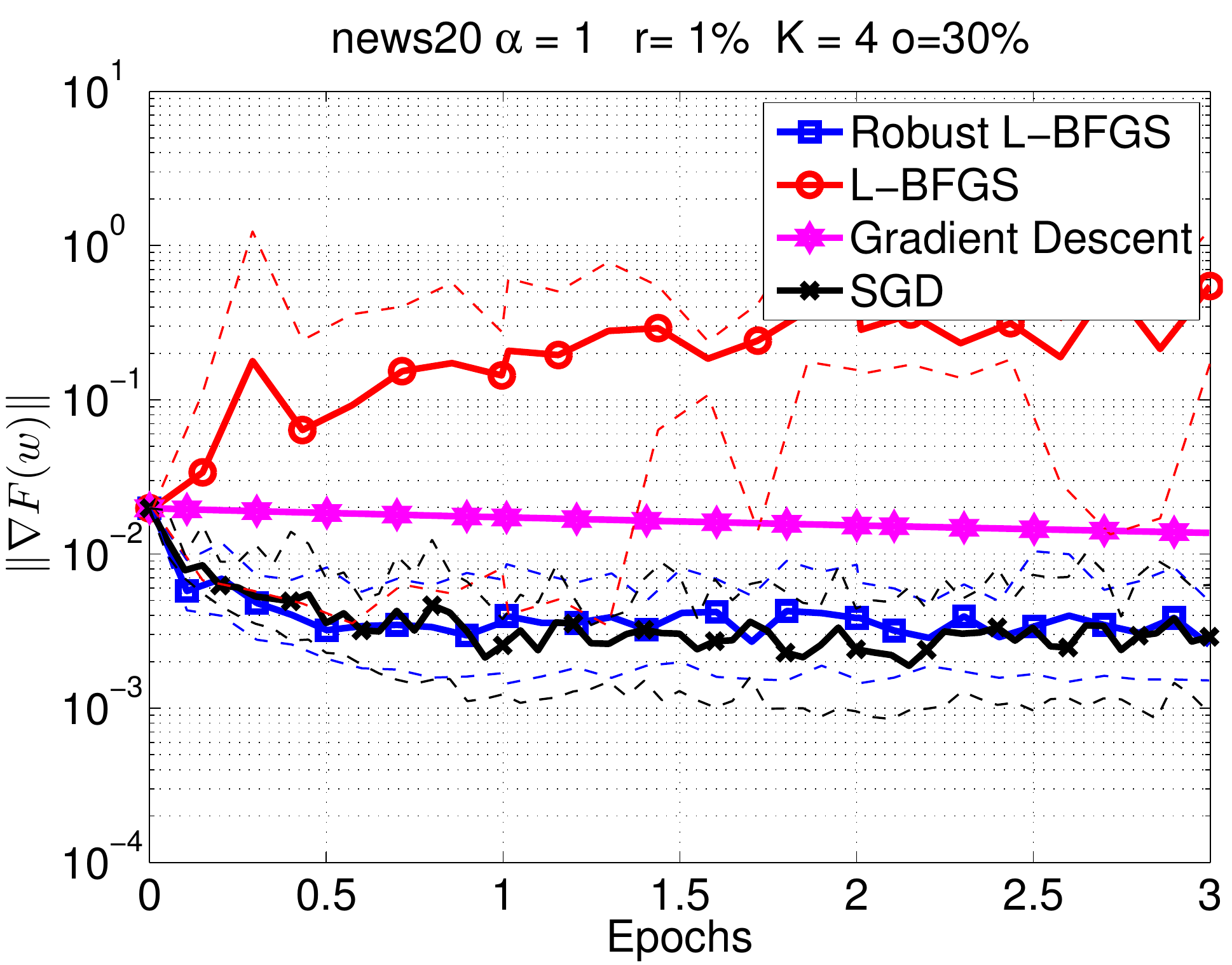}

\caption{\textbf{news20 dataset}. Comparison of Robust L-BFGS, L-BFGS (multi-batch L-BFGS without enforcing sample consistency), Gradient Descent (multi-batch Gradient method) and SGD. Top part:
we used $\alpha \in \{1, 0.1\}$,
$r\in \{1\%,  5\%,  10\%\}$ and $o=20\%$.
Bottom part: we used $\alpha=1$, $r=1\%$ and
$o\in \{5\%,  10\%, 20\%, 30\%\}$. Solid lines show average performance, and dashed lines show worst and best performance, over 10 runs (per algorithm). $K=4$ MPI processes.}
\end{figure}

\begin{figure}
\centering
\includegraphics[width=4.6cm]{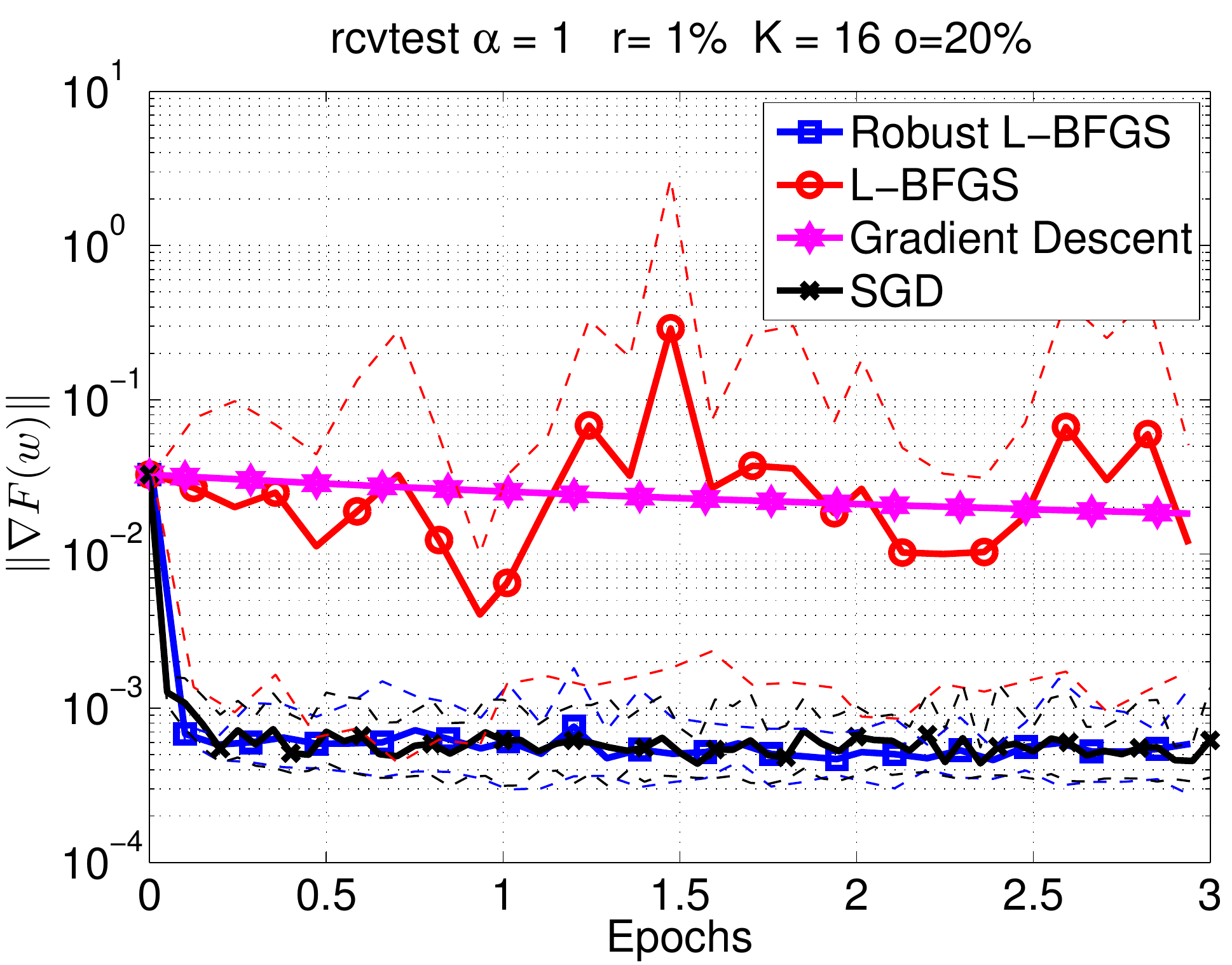}
\includegraphics[width=4.6cm]{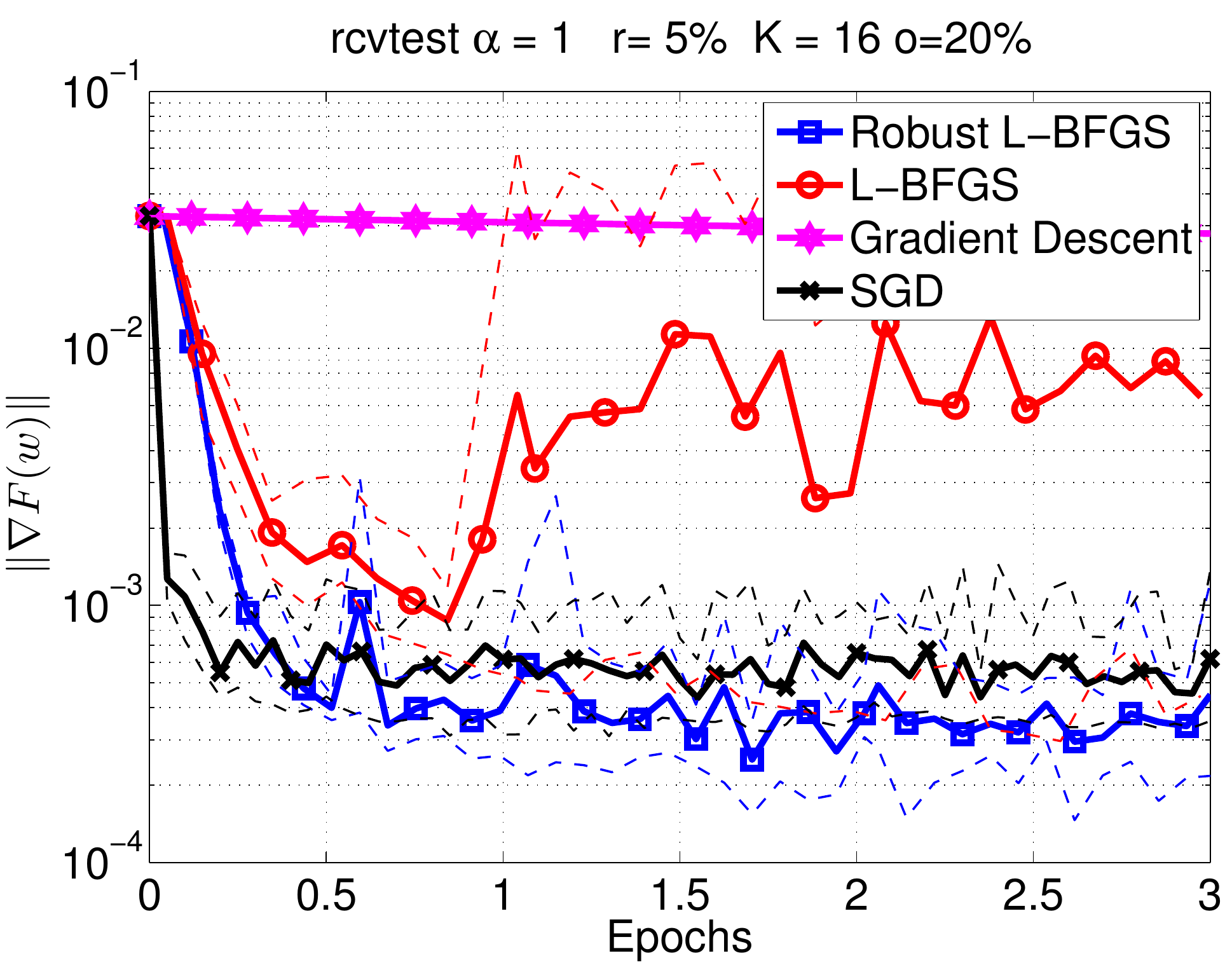}
\includegraphics[width=4.6cm]{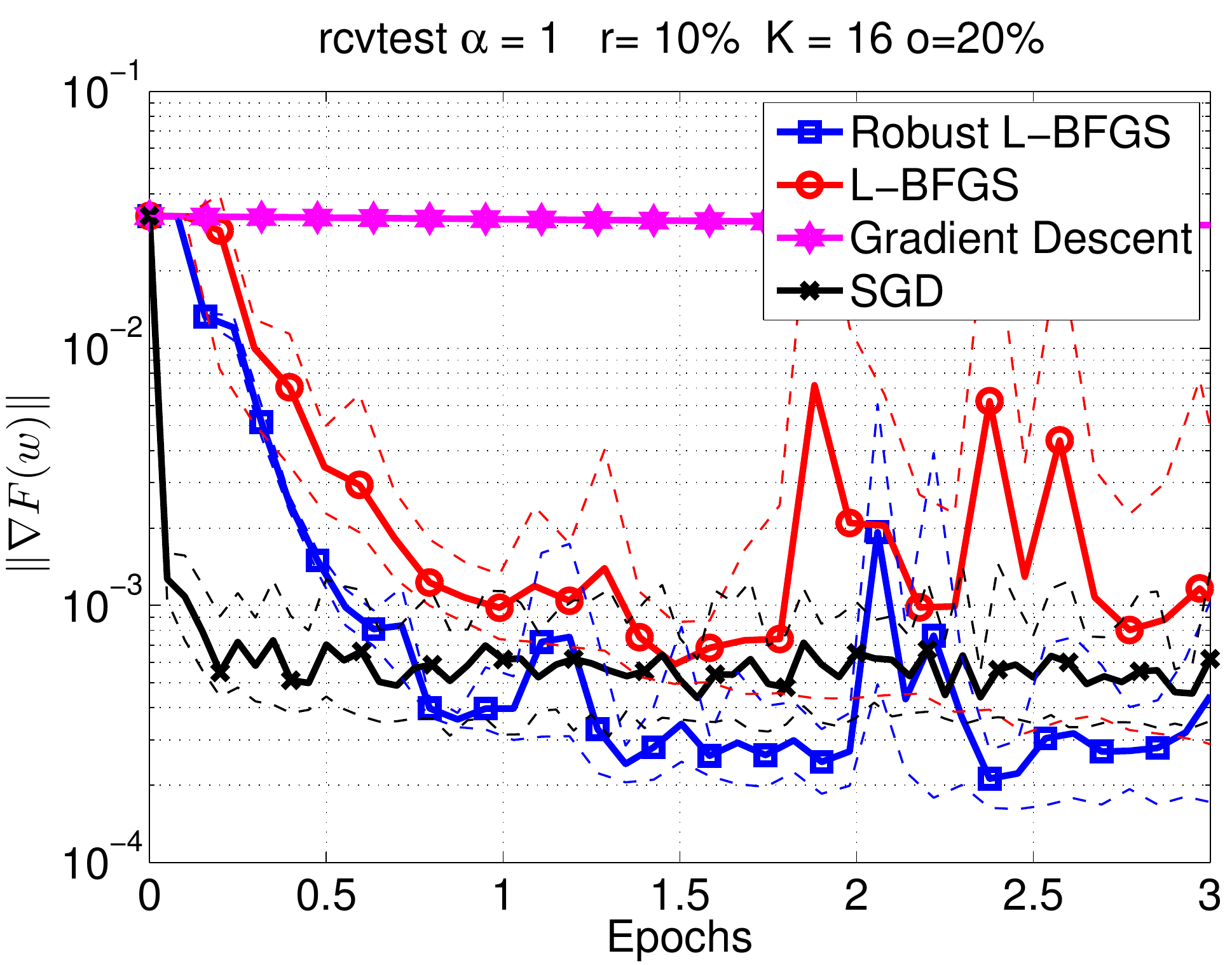}

 \includegraphics[width=4.6cm]{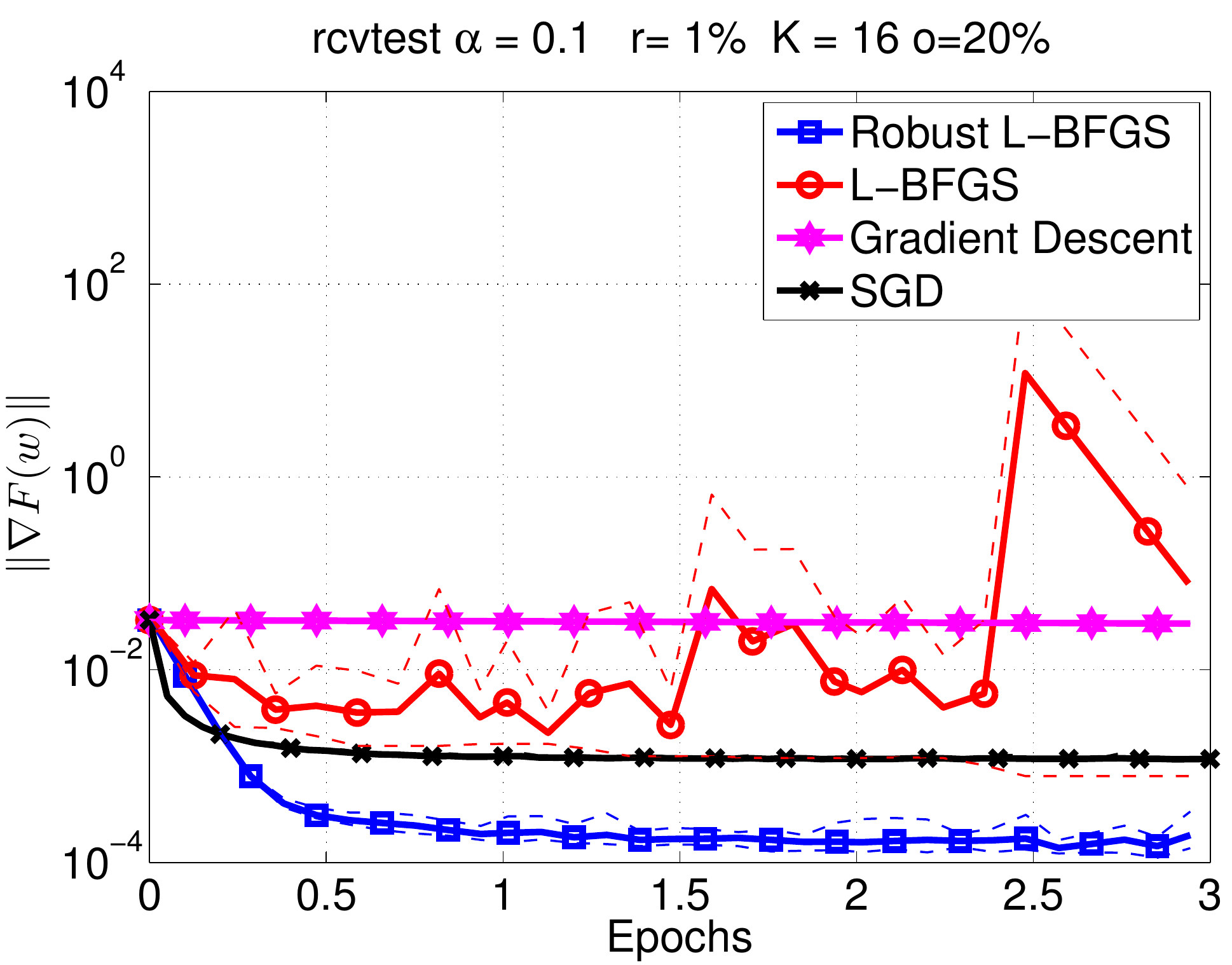}
\includegraphics[width=4.6cm]{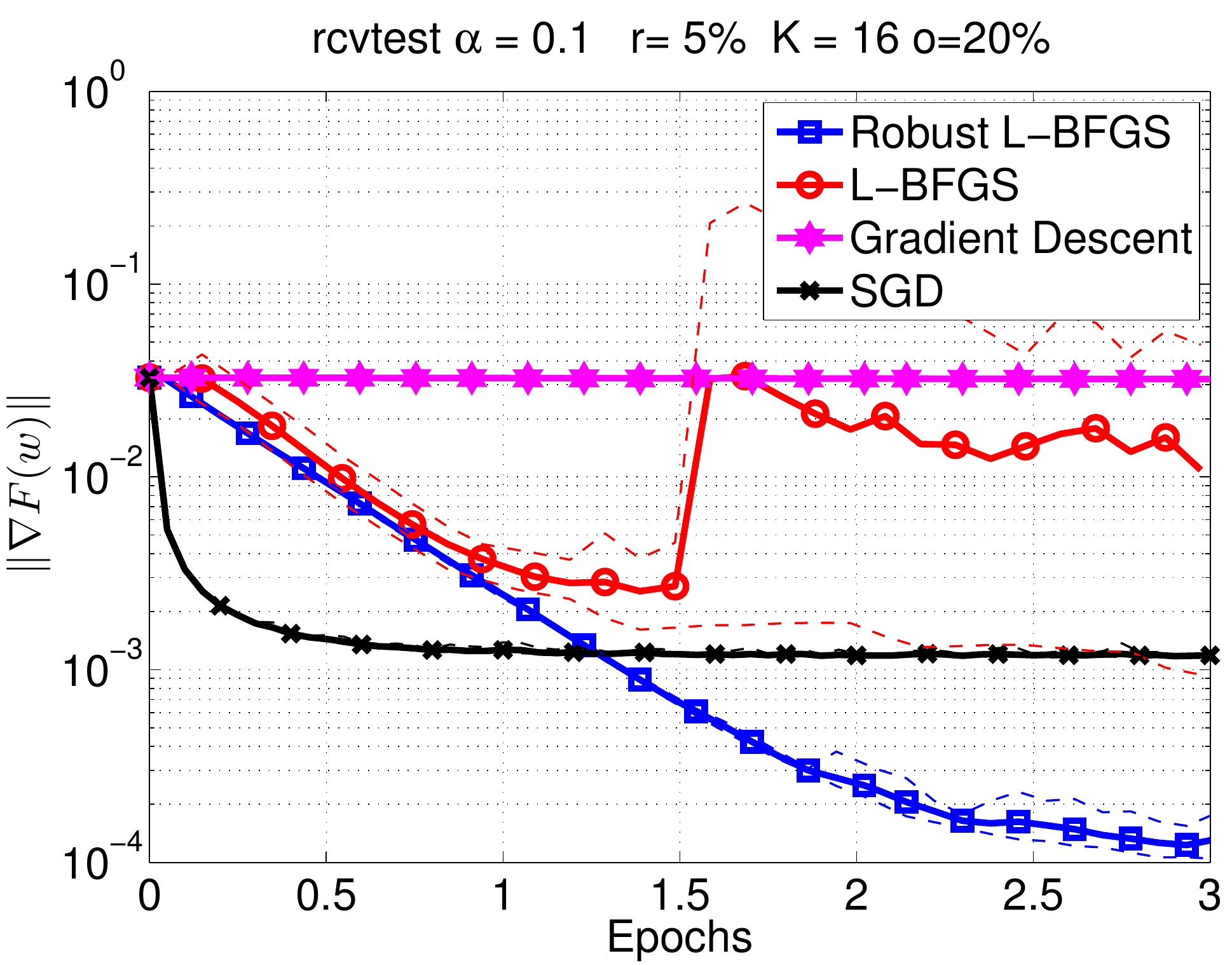}
\includegraphics[width=4.6cm]{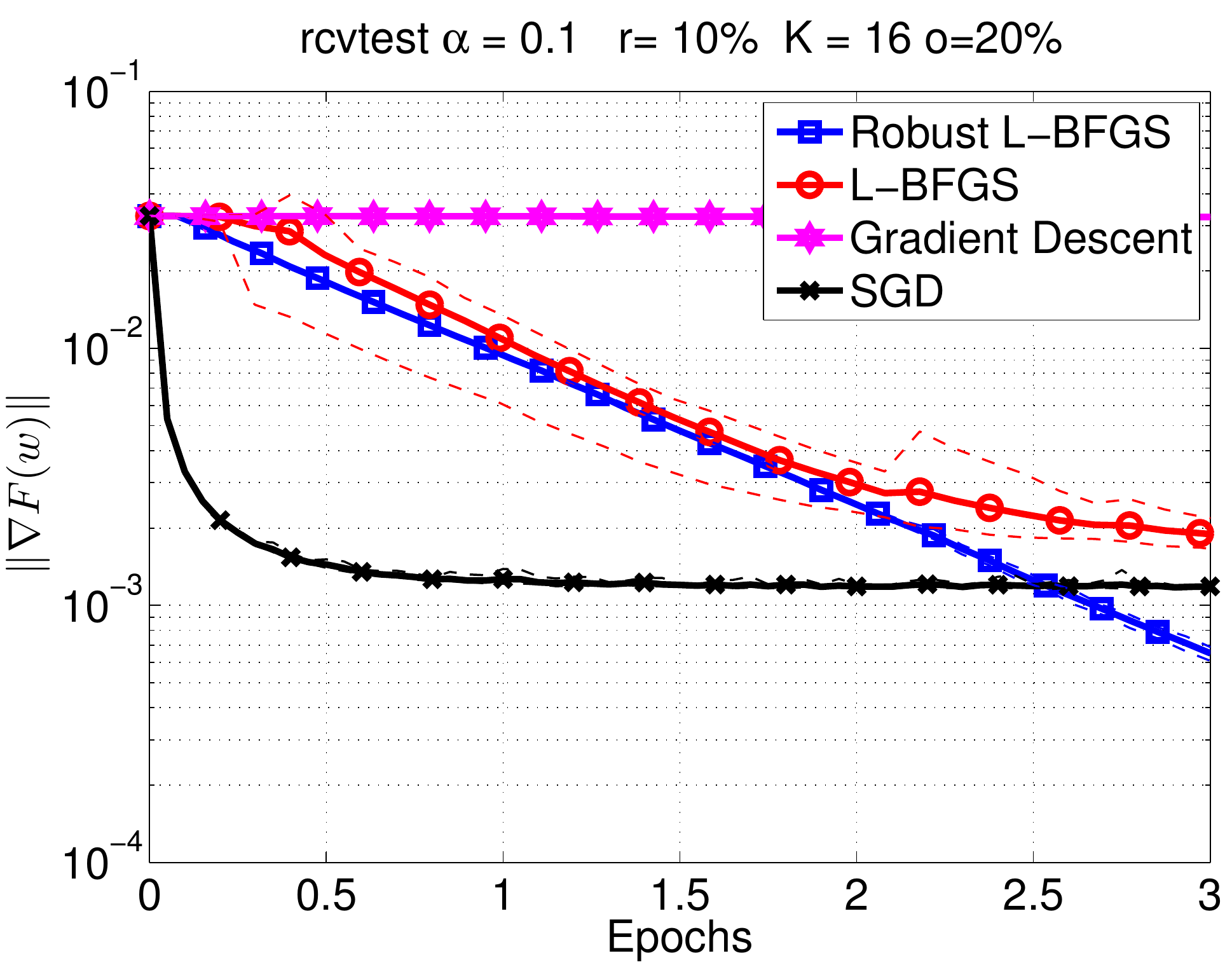}

\hrule 

\includegraphics[width=4.6cm]{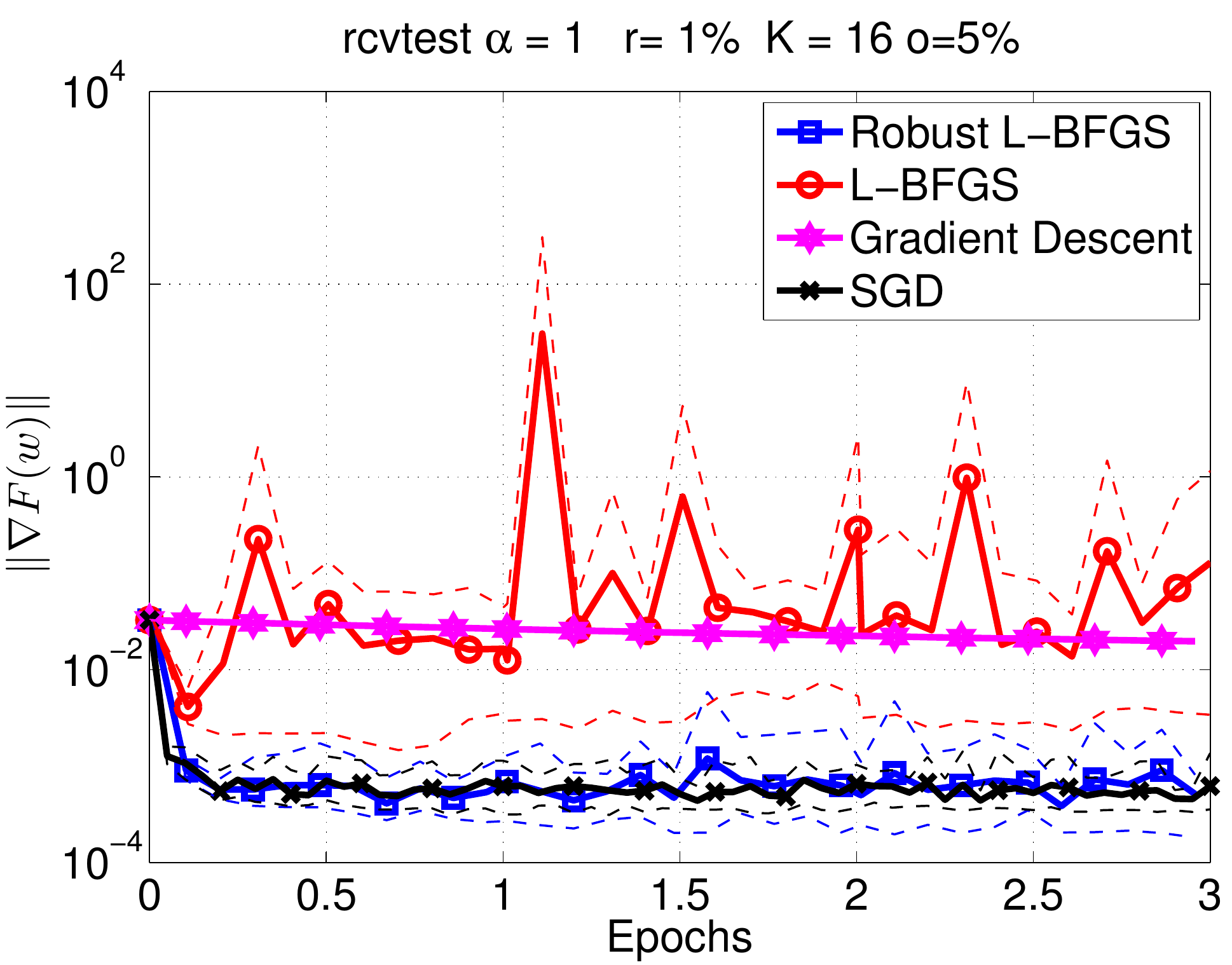}
\includegraphics[width=4.6cm]{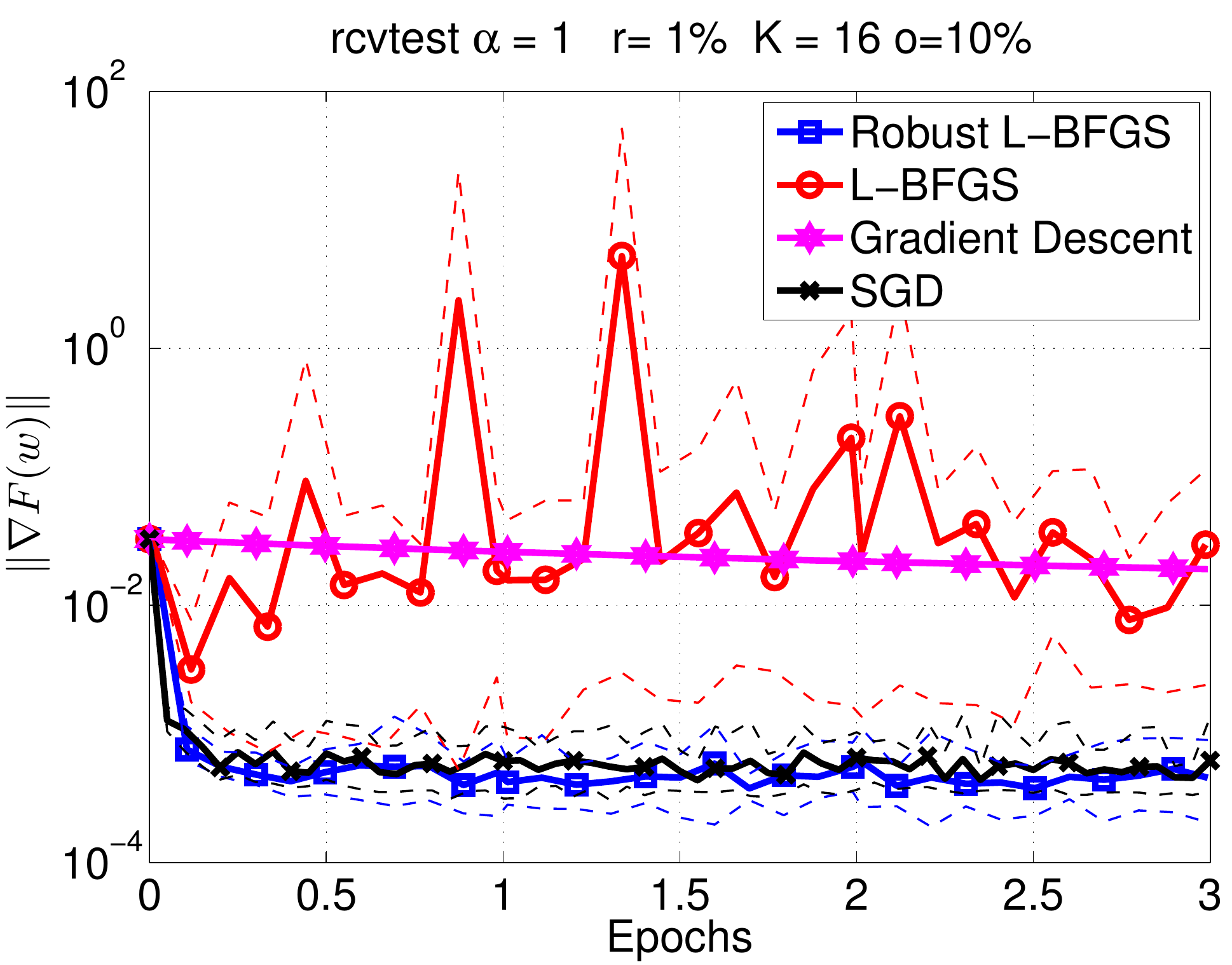}
\includegraphics[width=4.6cm]{rcvtest_mb_1_0_01_16_0_2-eps-converted-to.pdf}
\includegraphics[width=4.6cm]{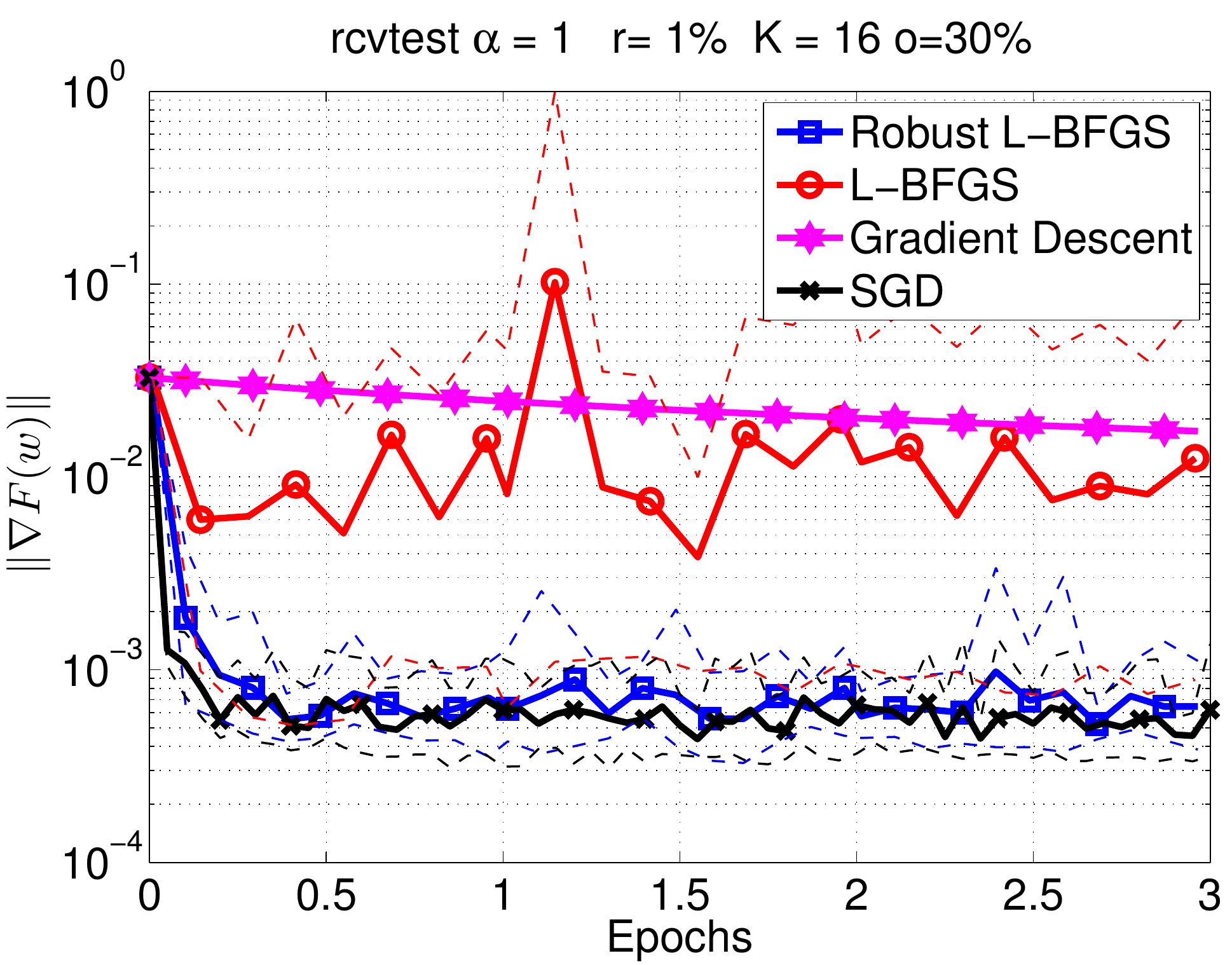}

\caption{\textbf{rcvtest dataset}. Comparison of Robust L-BFGS, L-BFGS (multi-batch L-BFGS without enforcing sample consistency), Gradient Descent (multi-batch Gradient method) and SGD. Top part:
we used $\alpha \in \{1, 0.1\}$,
$r\in \{1\%,  5\%,  10\%\}$ and $o=20\%$.
Bottom part: we used $\alpha=1$, $r=1\%$ and
$o\in \{5\%,  10\%, 20\%, 30\%\}$. Solid lines show average performance, and dashed lines show worst and best performance, over 10 runs (per algorithm). $K=16$ MPI processes.}
\end{figure}

\begin{figure}
\centering
\includegraphics[width=4.6cm]{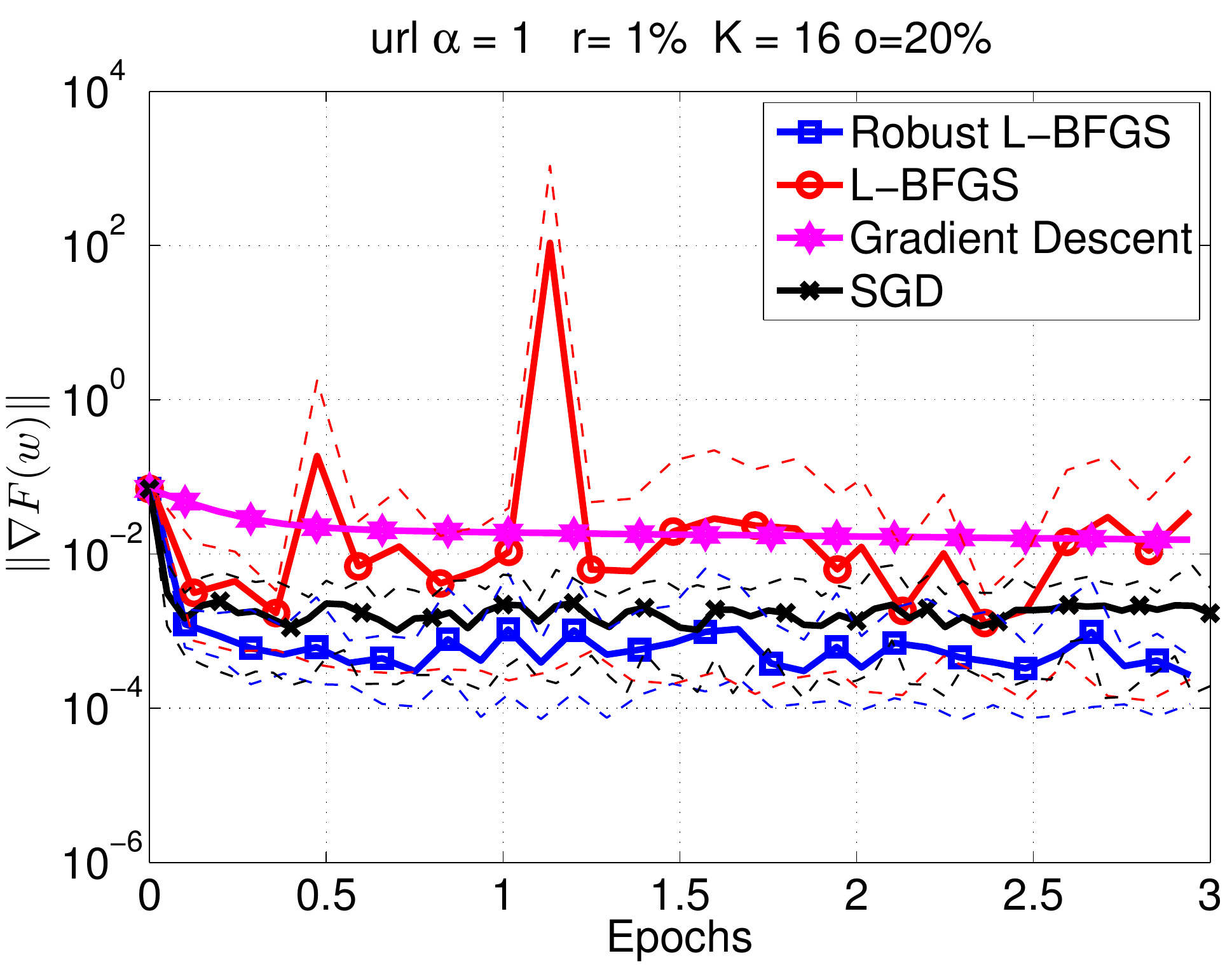}
\includegraphics[width=4.6cm]{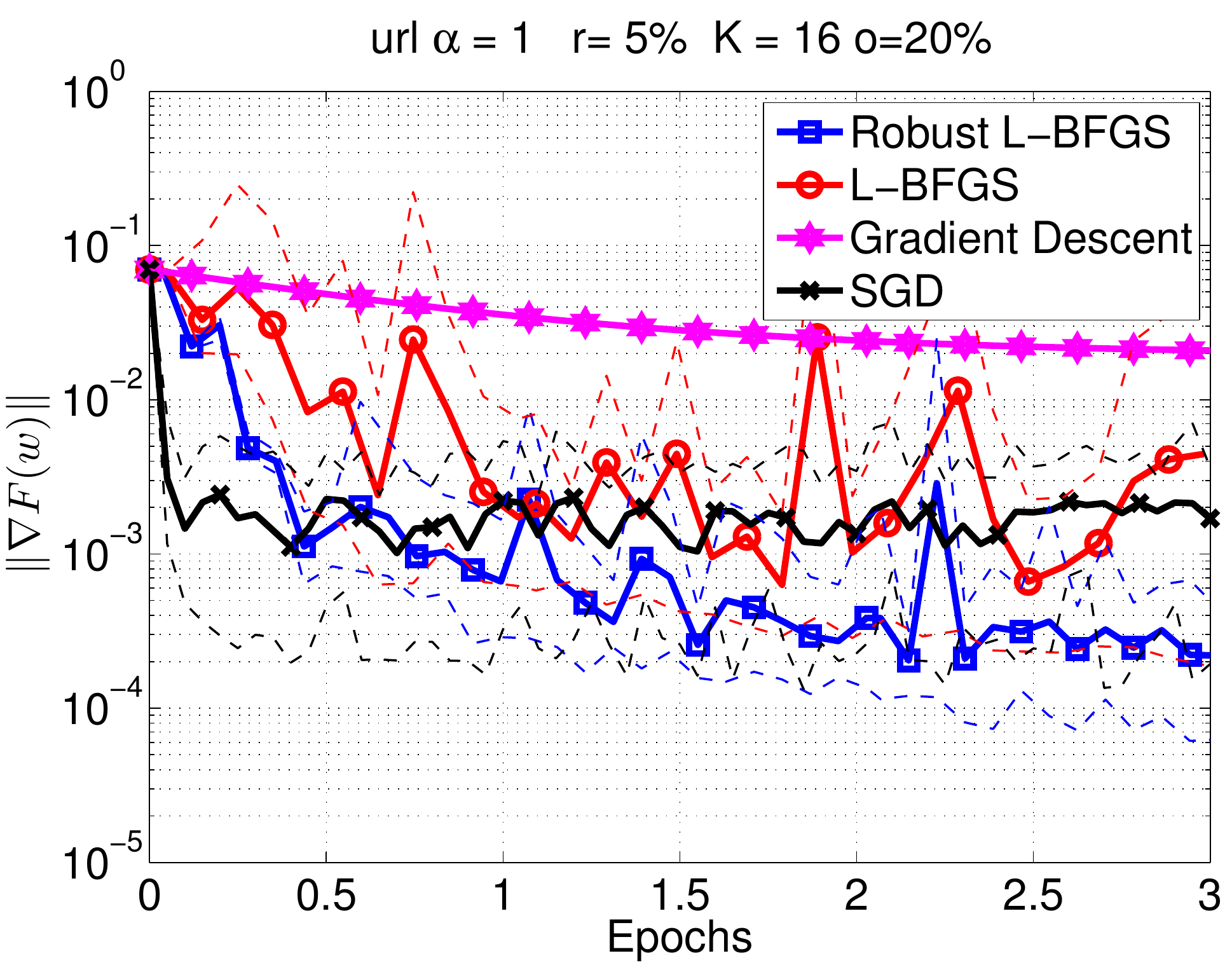}
\includegraphics[width=4.6cm]{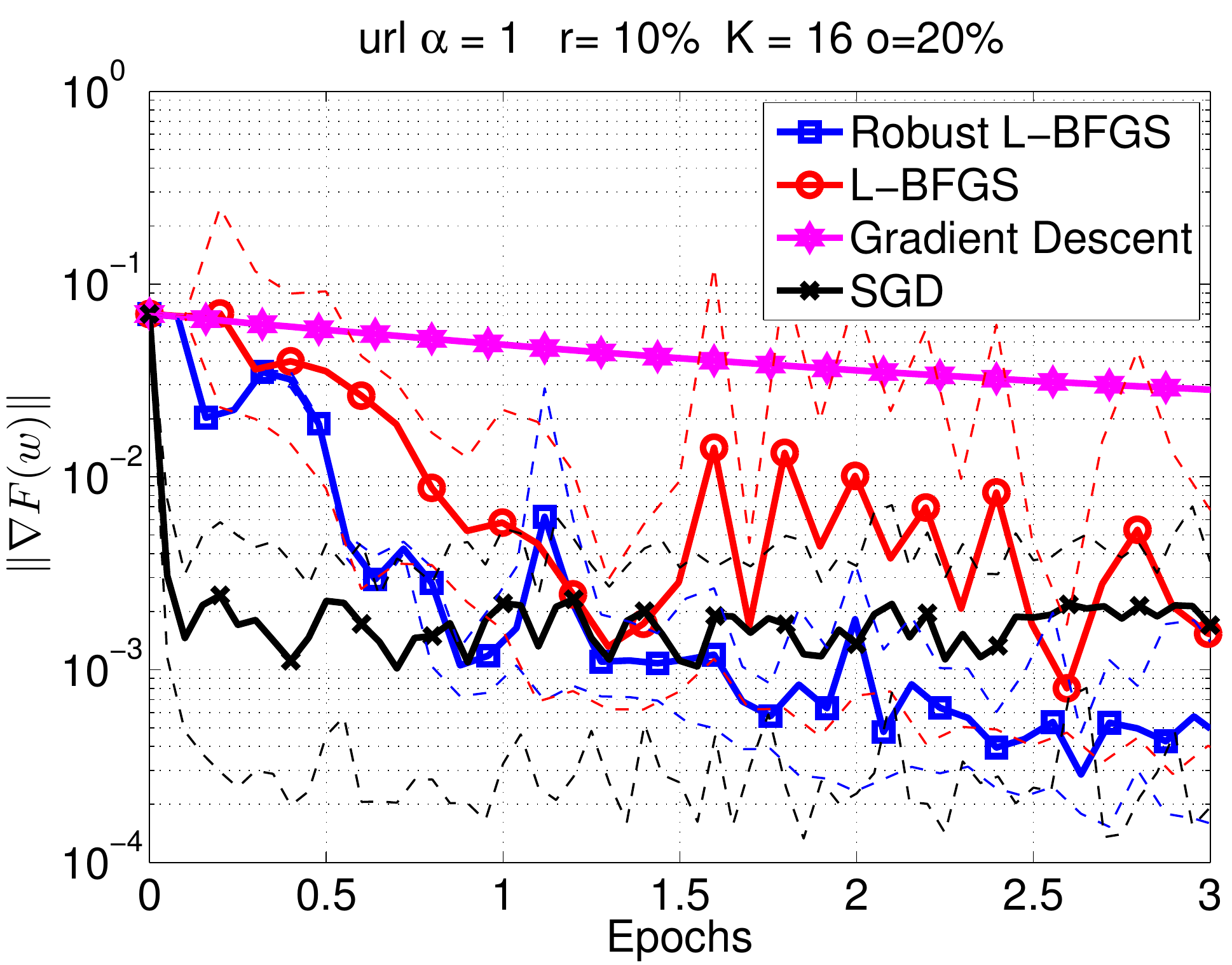}

\includegraphics[width=4.6cm]{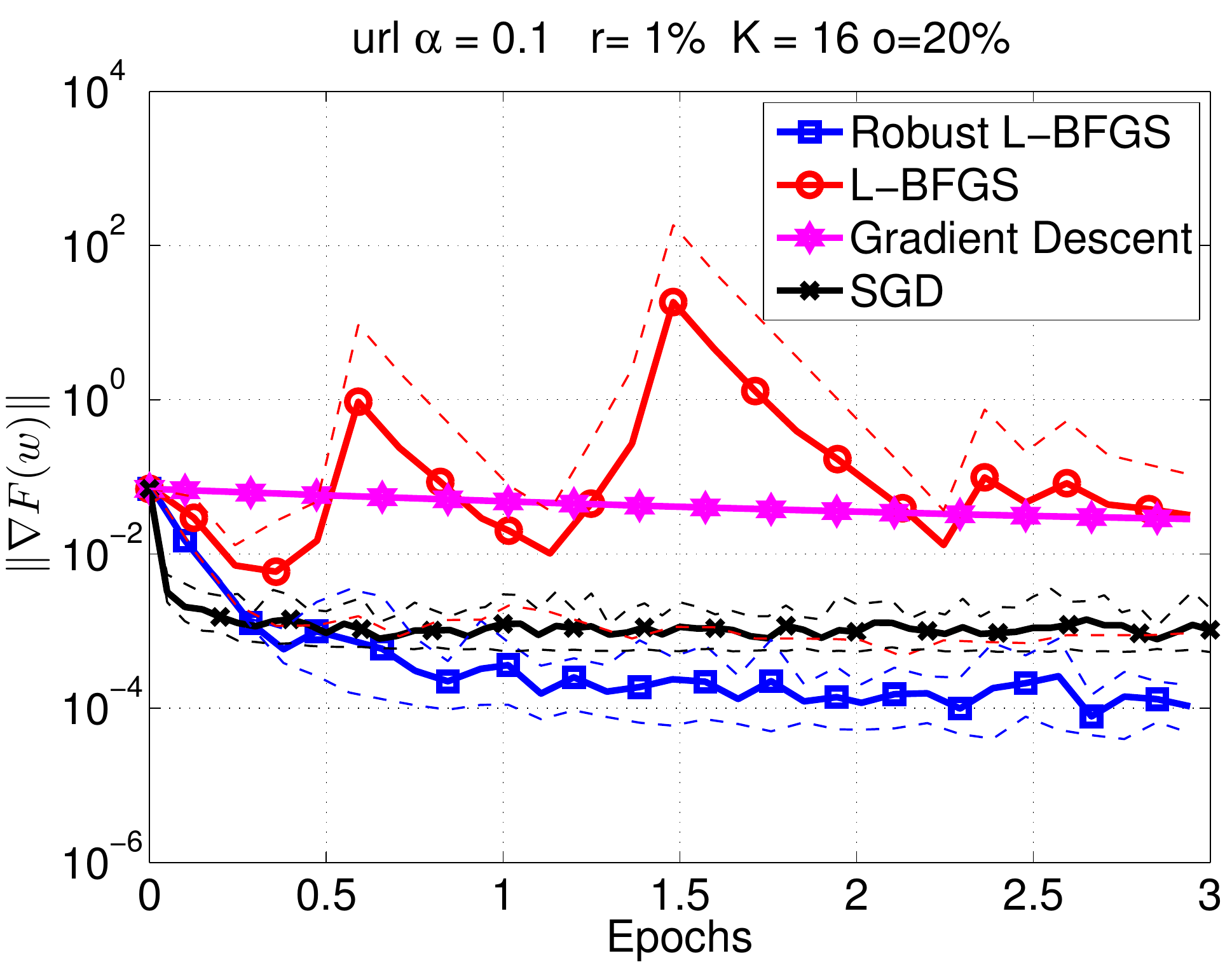}
\includegraphics[width=4.6cm]{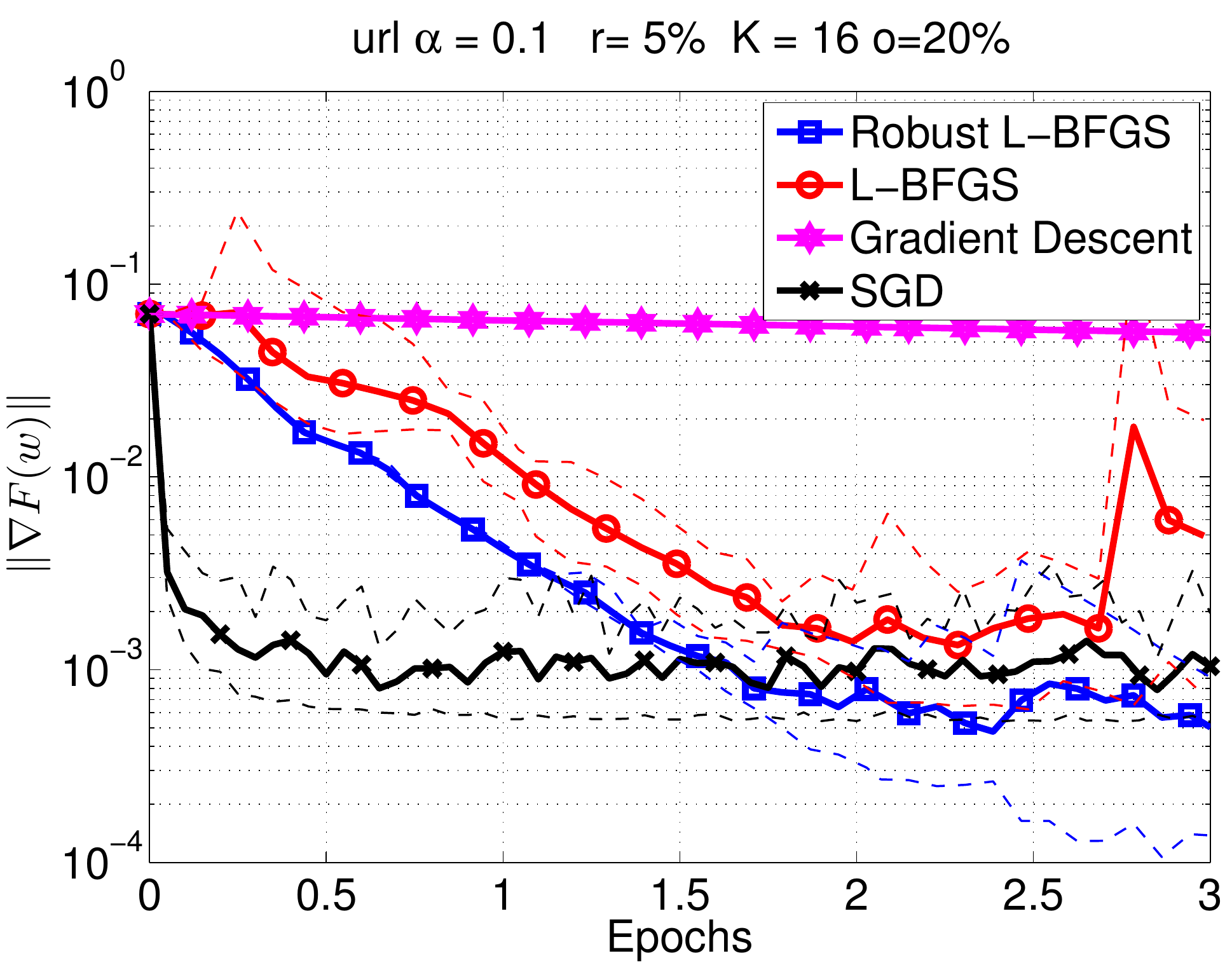}
\includegraphics[width=4.6cm]{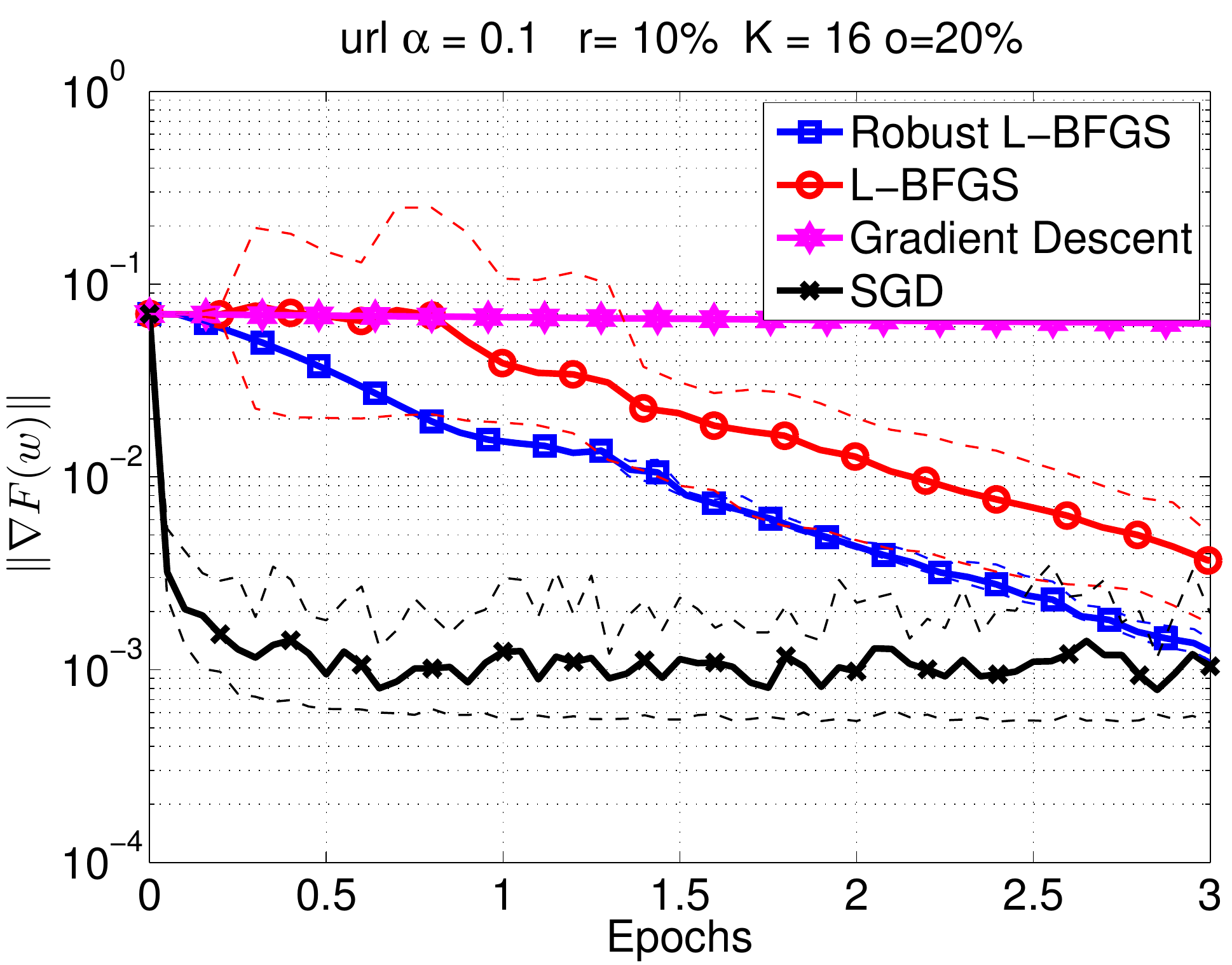}

\hrule 

\includegraphics[width=4.6cm]{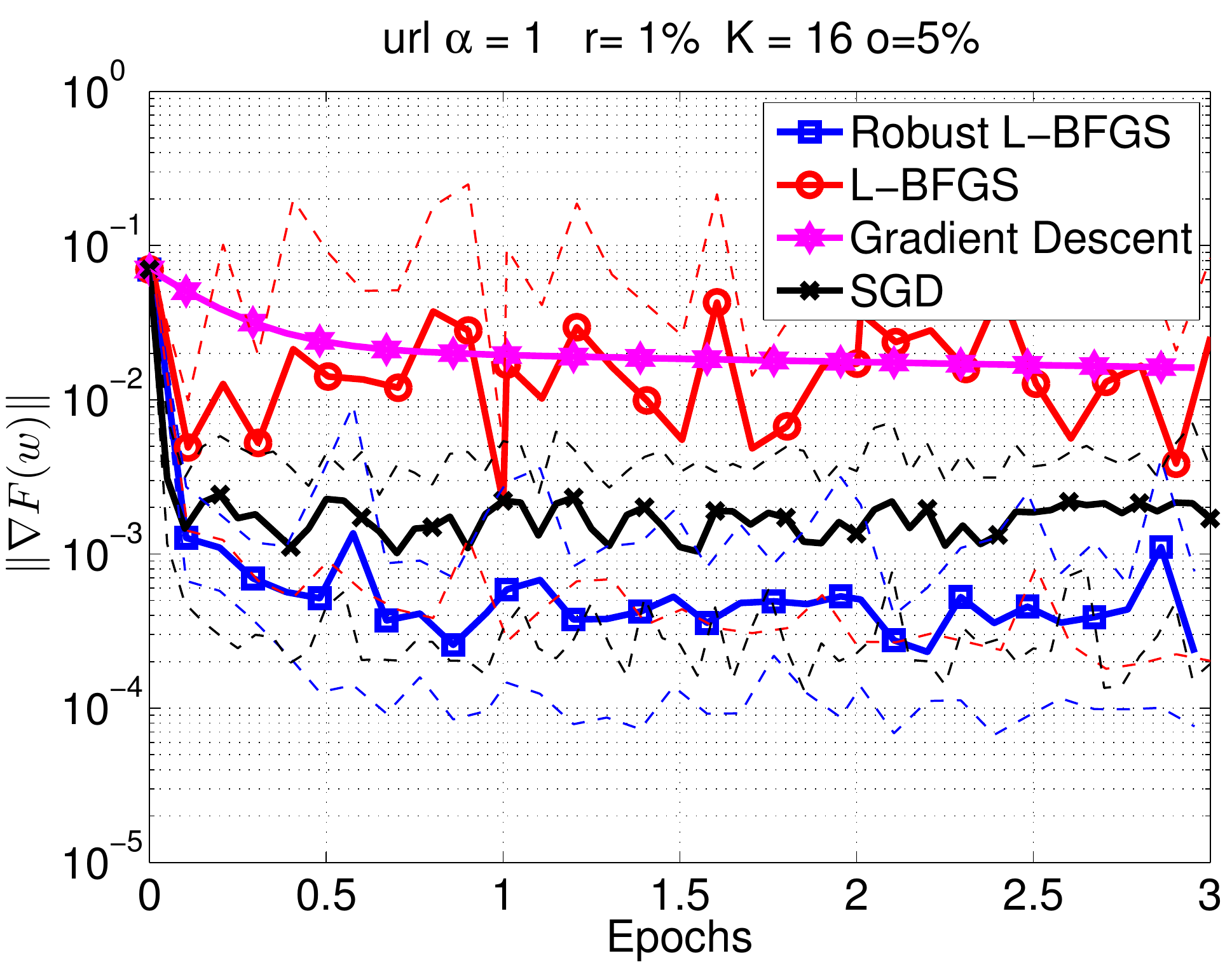}
\includegraphics[width=4.6cm]{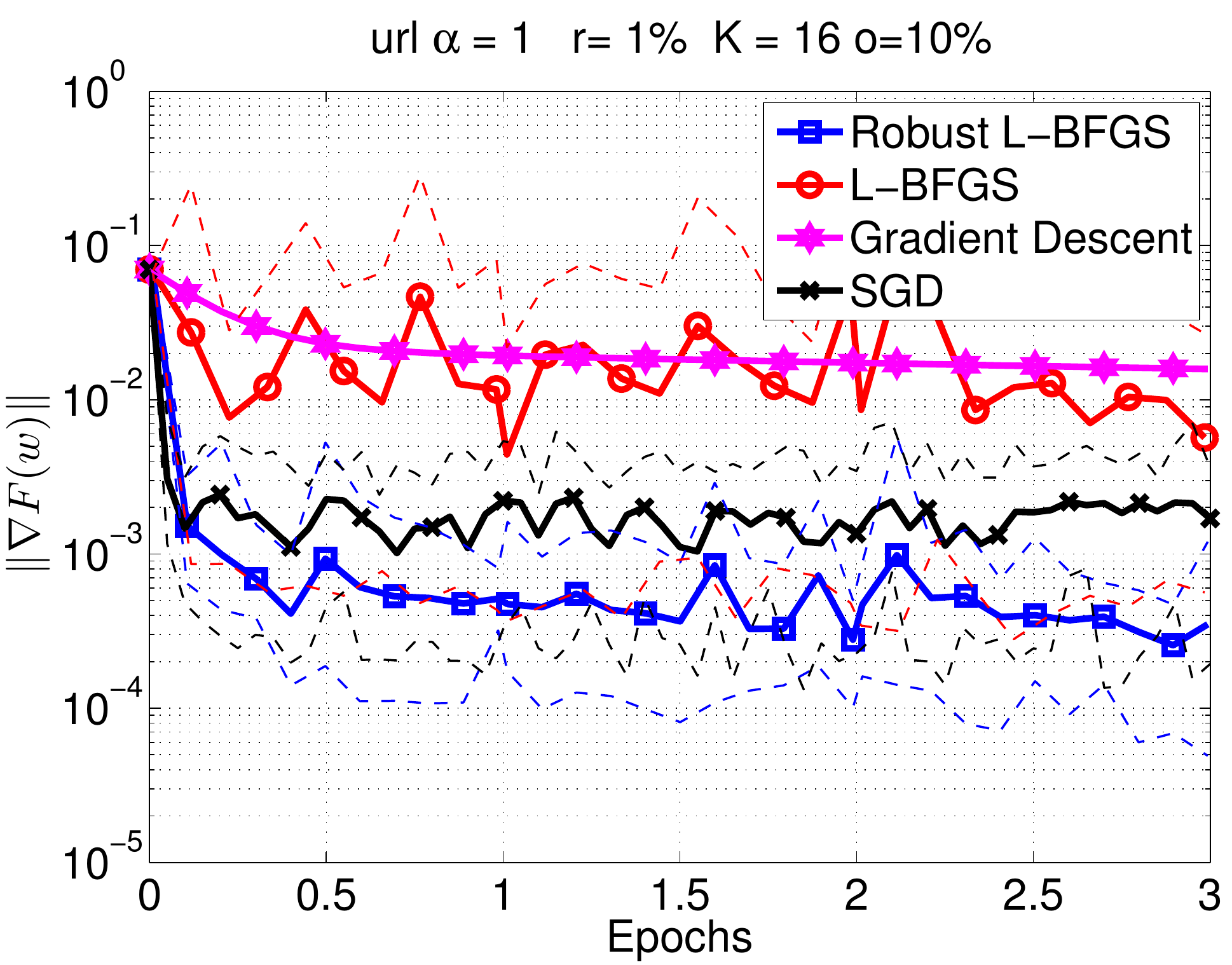}
\includegraphics[width=4.6cm]{url_mb_1_0_01_16_0_2-eps-converted-to.pdf}
\includegraphics[width=4.6cm]{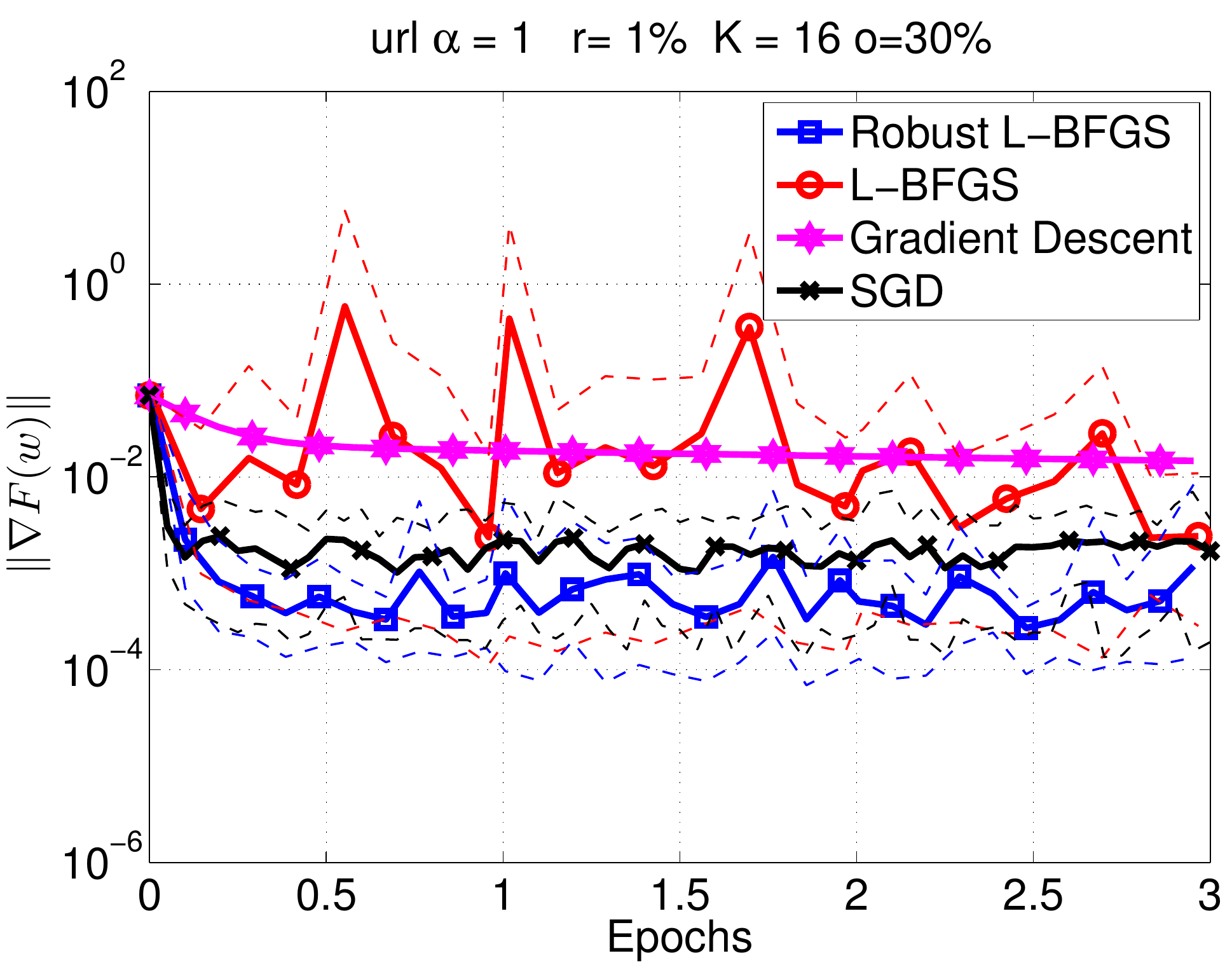}

\caption{\textbf{url dataset}. Comparison of Robust L-BFGS, L-BFGS (multi-batch L-BFGS without enforcing sample consistency), Gradient Descent (multi-batch Gradient method) and SGD. Top part:
we used $\alpha \in \{1, 0.1\}$,
$r\in \{1\%,  5\%,  10\%\}$ and $o=20\%$.
Bottom part: we used $\alpha=1$, $r=1\%$ and
$o\in \{5\%,  10\%, 20\%, 30\%\}$. Solid lines show average performance, and dashed lines show worst and best performance, over 10 runs (per algorithm). $K=16$ MPI processes.}
\end{figure}

\begin{figure}
\centering
\includegraphics[width=4.6cm]{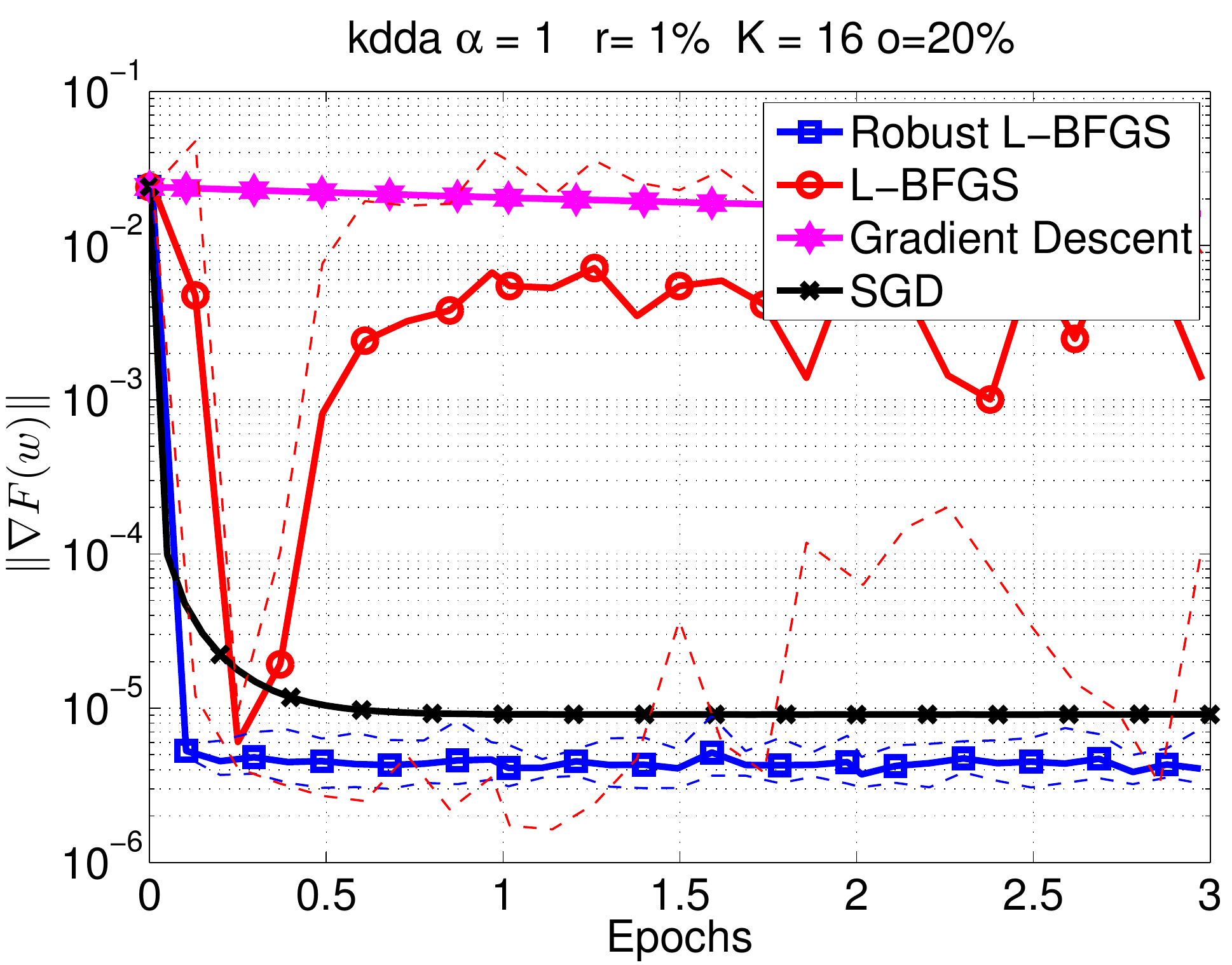}
\includegraphics[width=4.6cm]{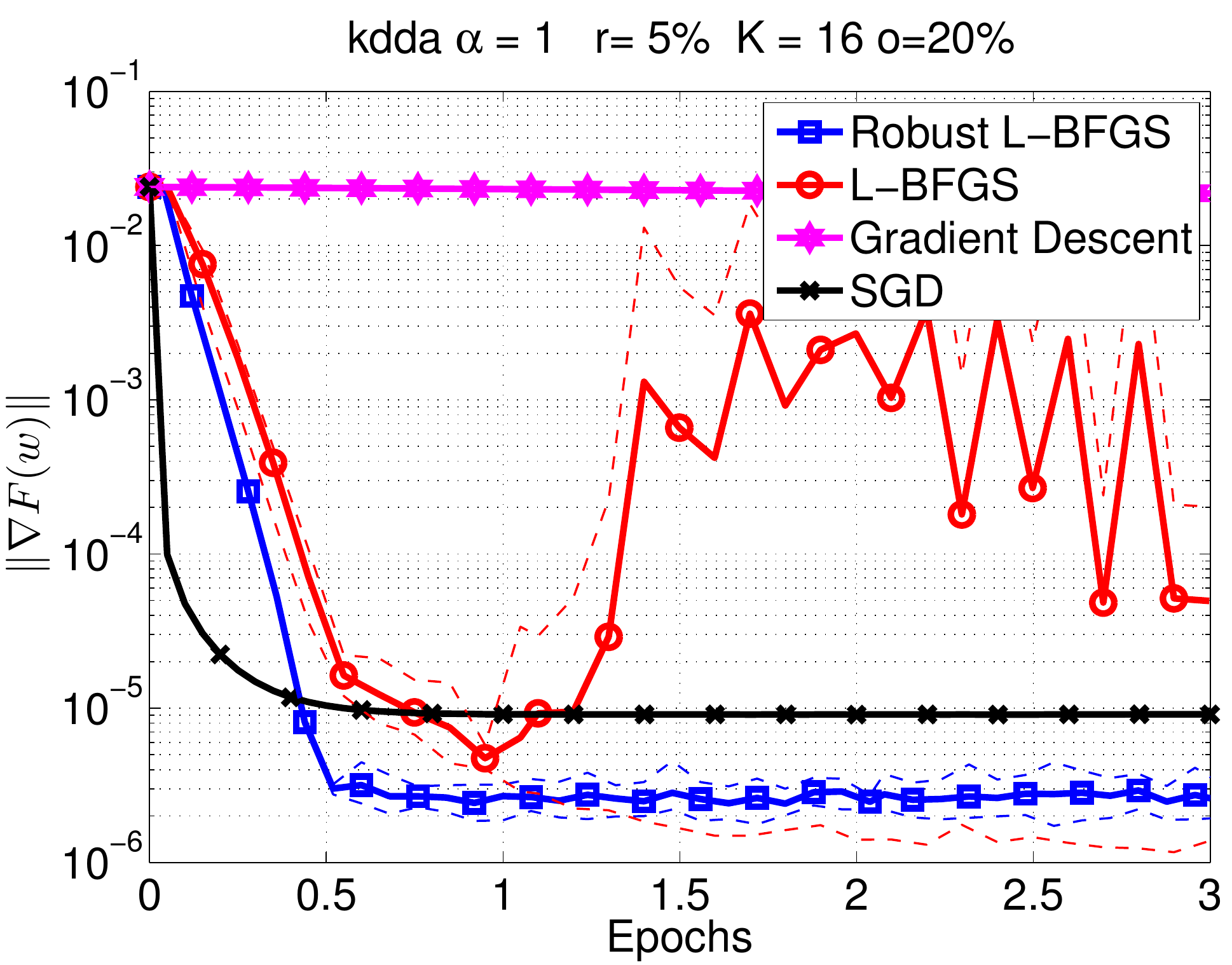}
\includegraphics[width=4.6cm]{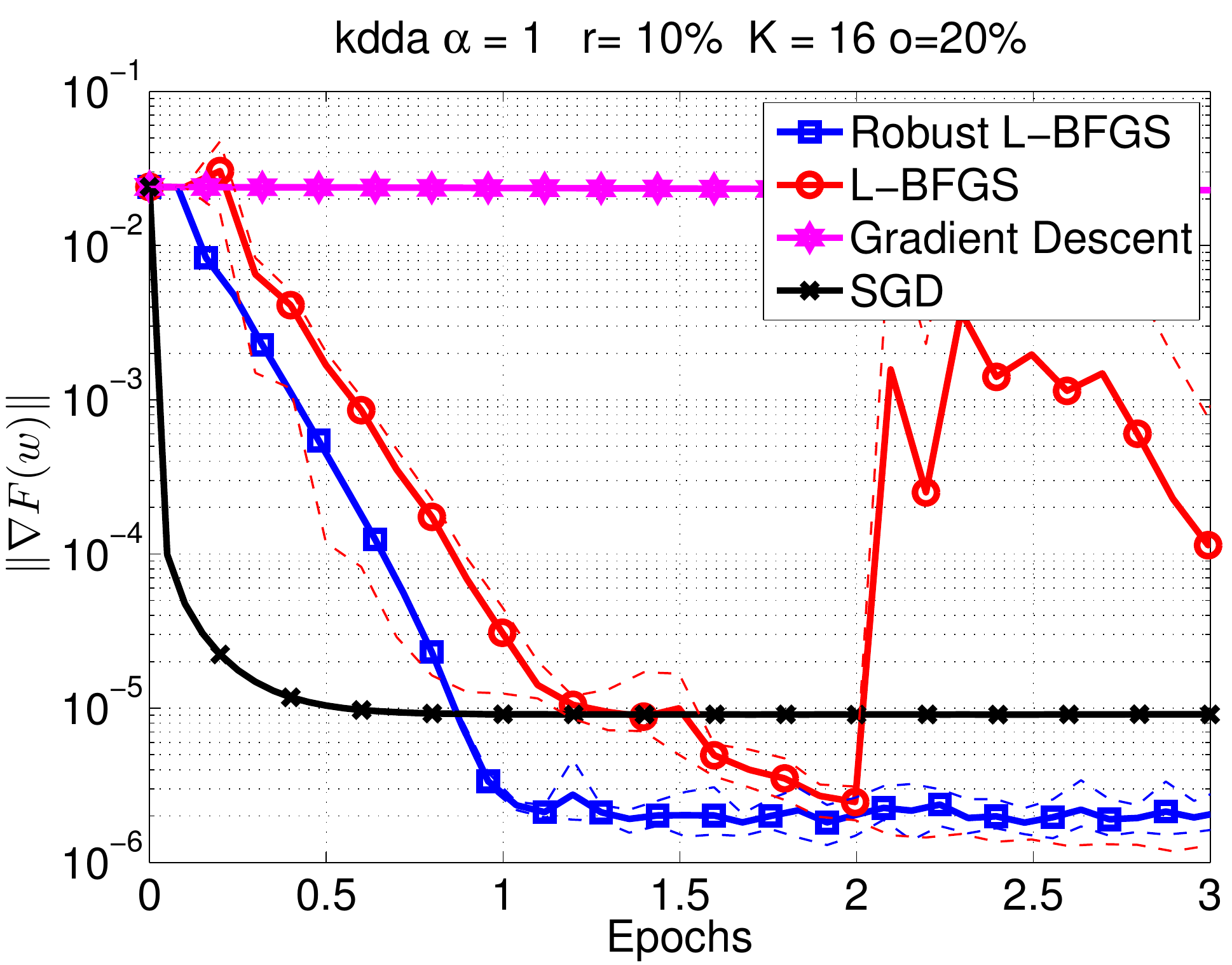}

 \includegraphics[width=4.6cm]{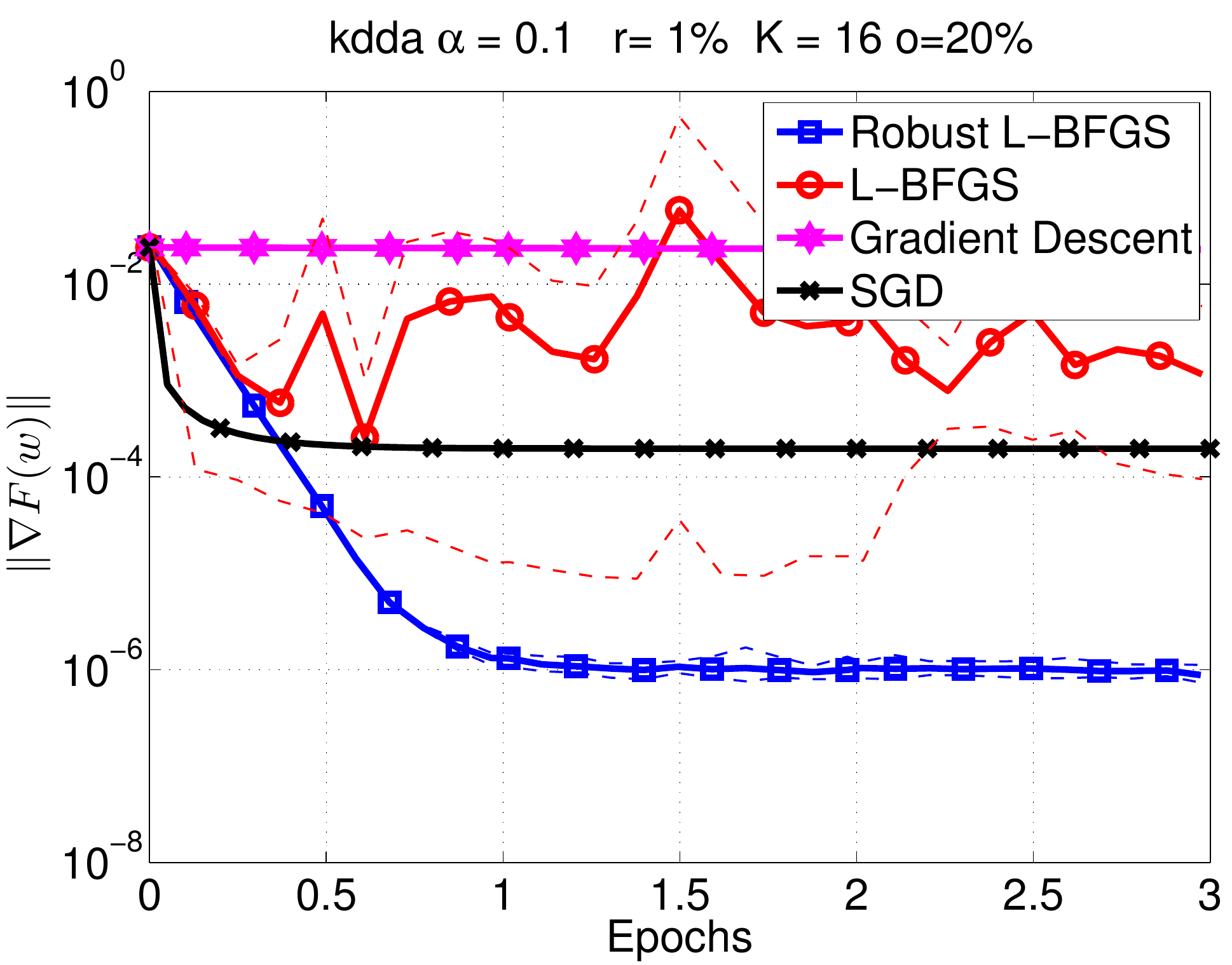}
\includegraphics[width=4.6cm]{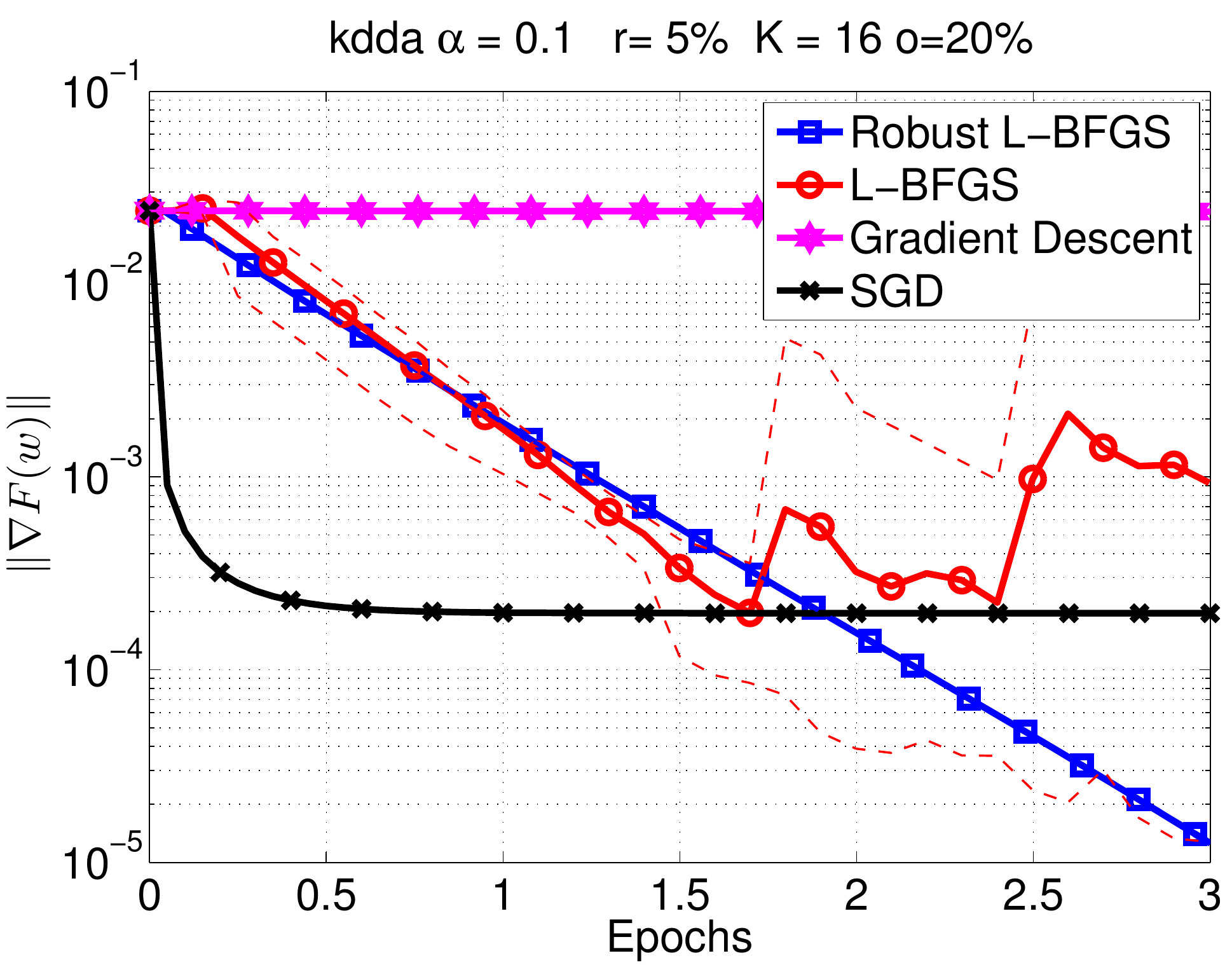}
\includegraphics[width=4.6cm]{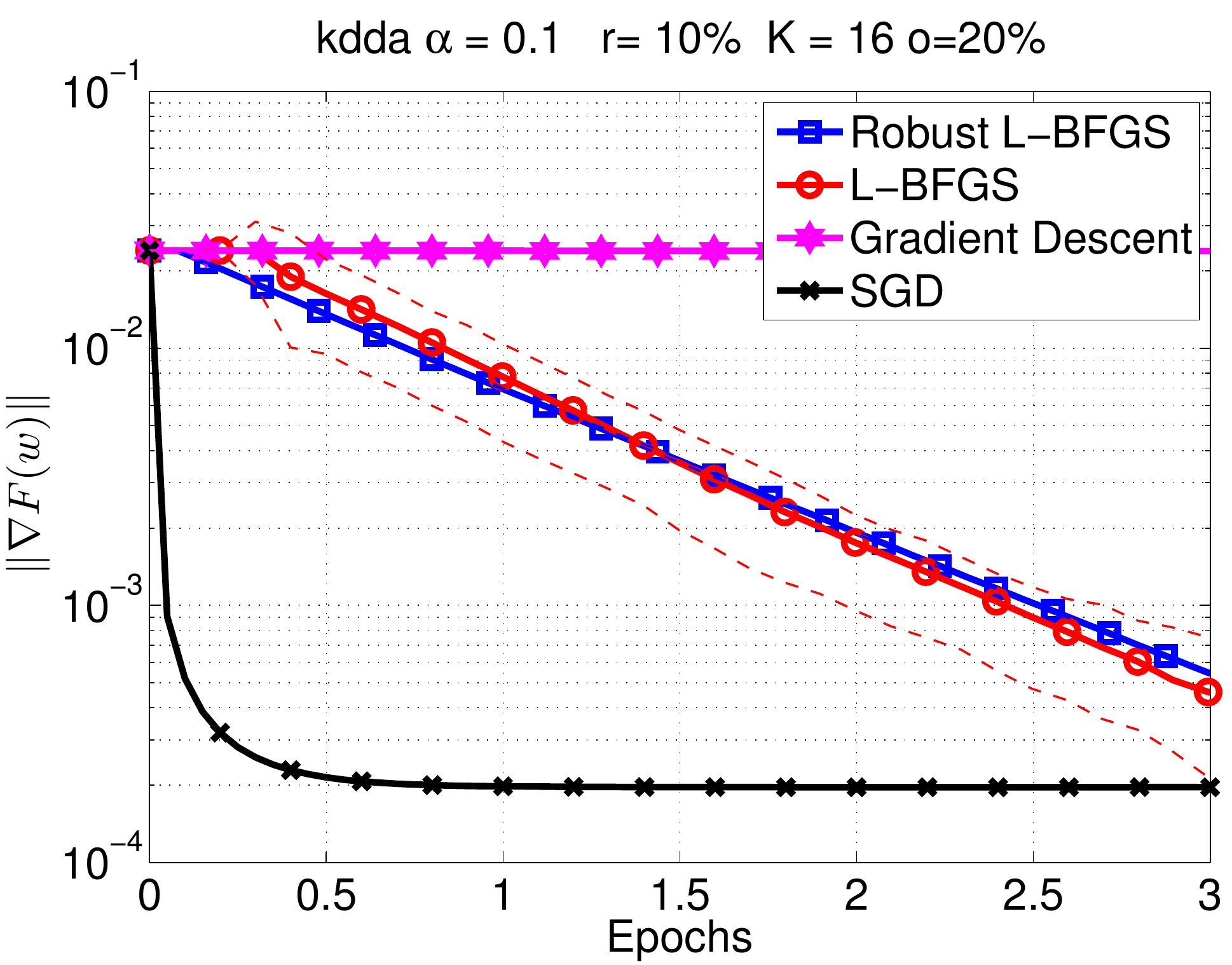}

\hrule 

\includegraphics[width=4.6cm]{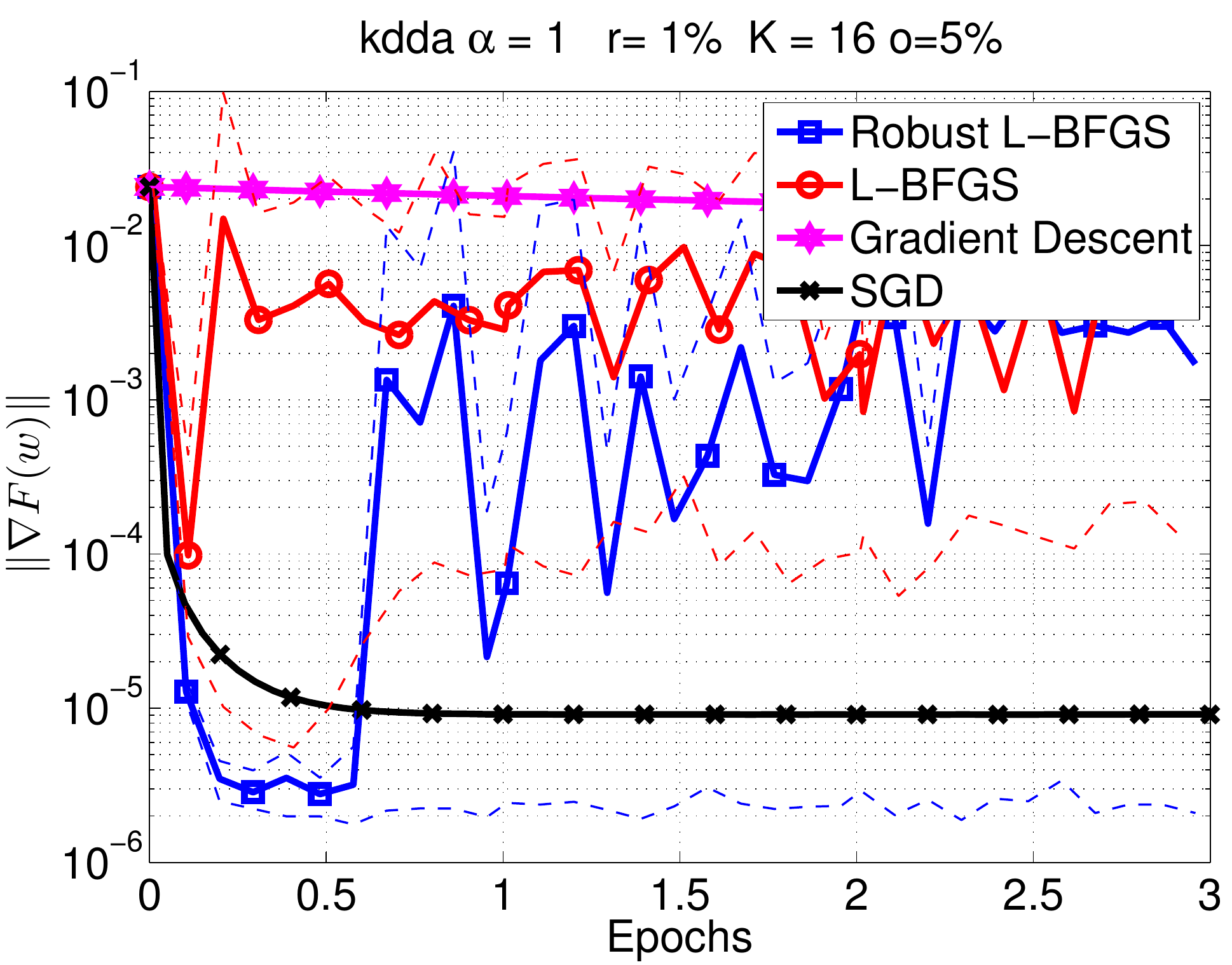}
\includegraphics[width=4.6cm]{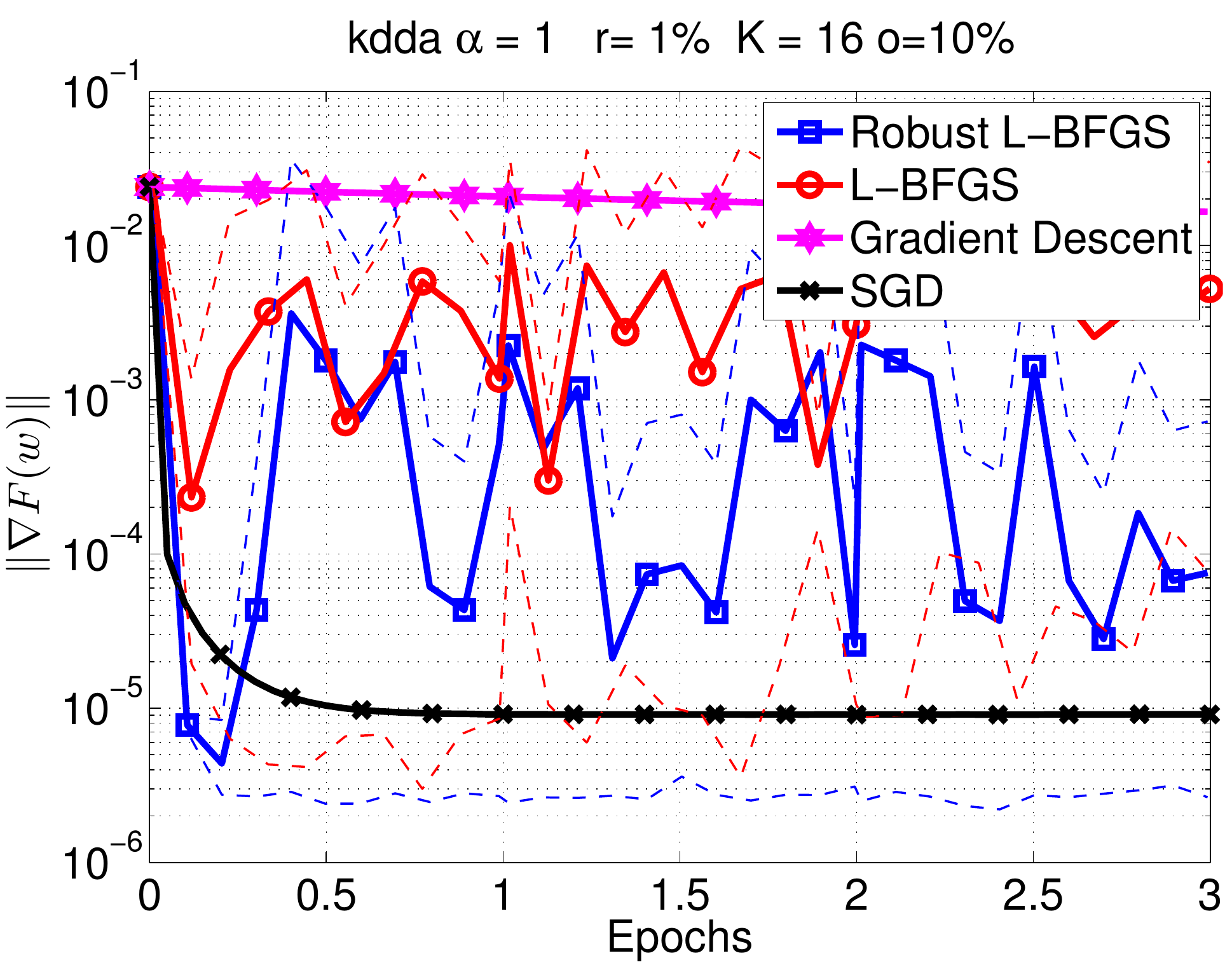}
\includegraphics[width=4.6cm]{kdda_mb_1_0_01_16_0_2-eps-converted-to.pdf}
\includegraphics[width=4.6cm]{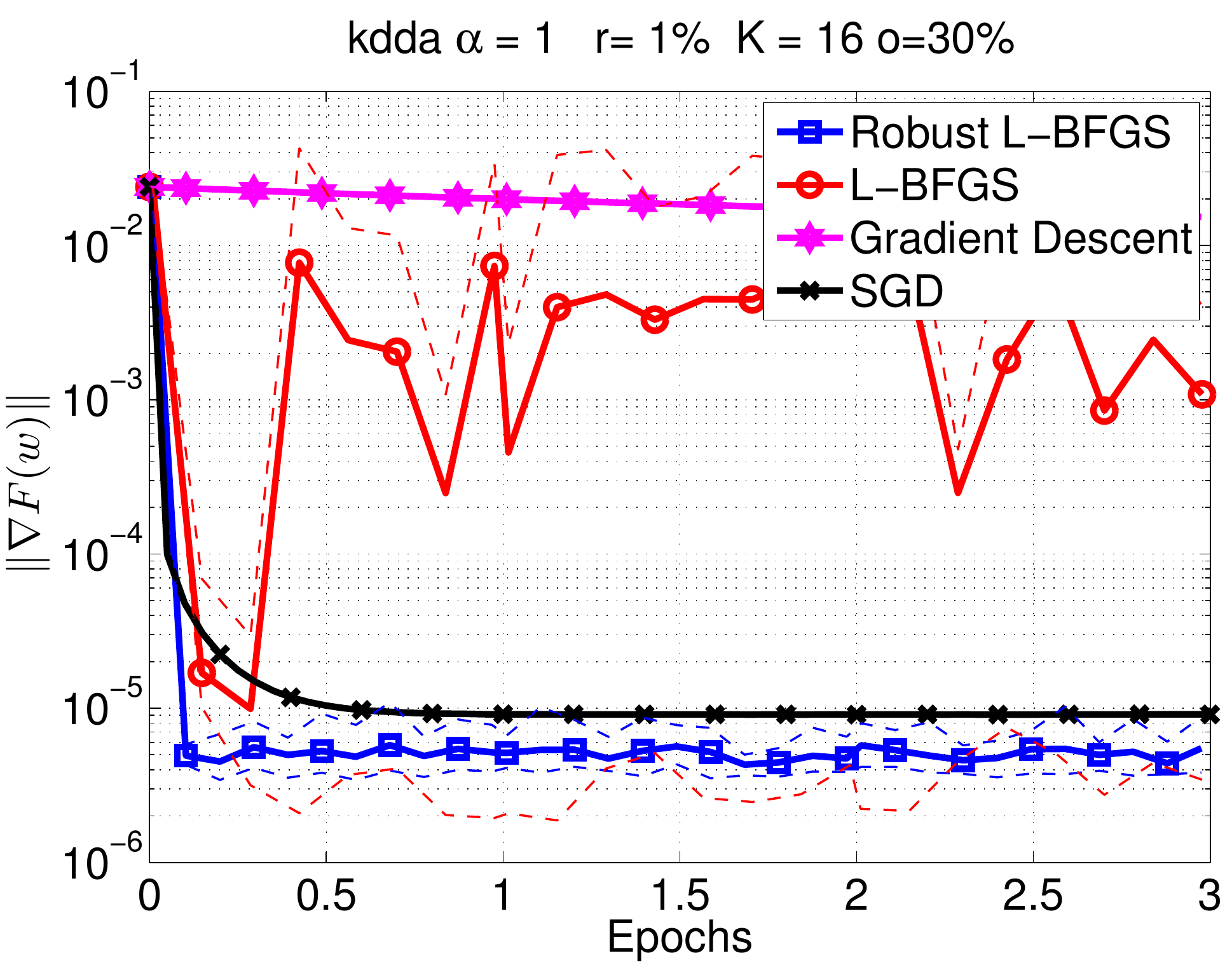}

\caption{\textbf{kdda dataset}. Comparison of Robust L-BFGS, L-BFGS (multi-batch L-BFGS without enforcing sample consistency), Gradient Descent (multi-batch Gradient method) and SGD. Top part:
we used $\alpha \in \{1, 0.1\}$,
$r\in \{1\%,  5\%,  10\%\}$ and $o=20\%$.
Bottom part: we used $\alpha=1$, $r=1\%$ and
$o\in \{5\%,  10\%, 20\%, 30\%\}$. Solid lines show average performance, and dashed lines show worst and best performance, over 10 runs (per algorithm). $K=16$ MPI processes.}
\end{figure}

\begin{figure}
\centering
\includegraphics[width=4.6cm]{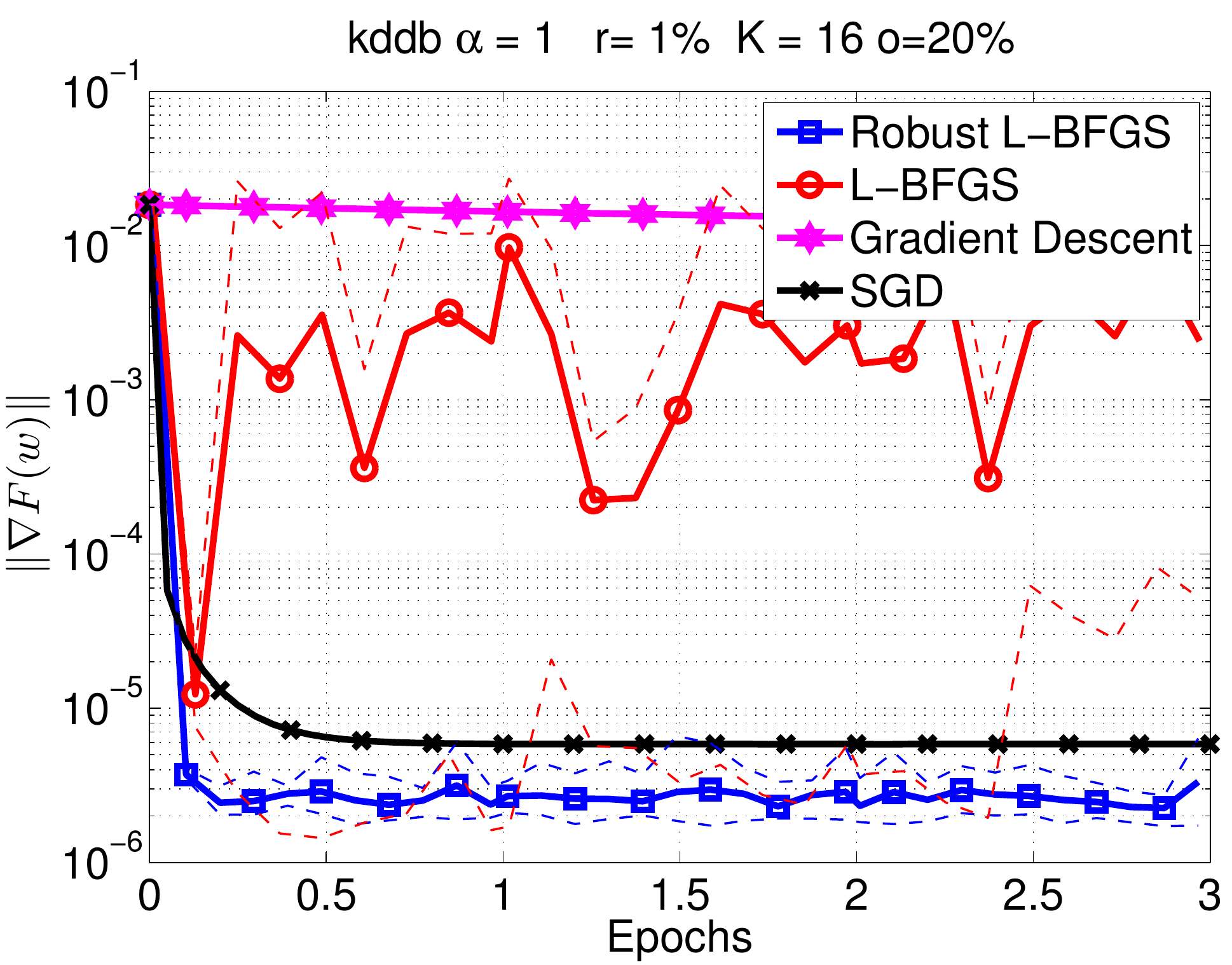}
\includegraphics[width=4.6cm]{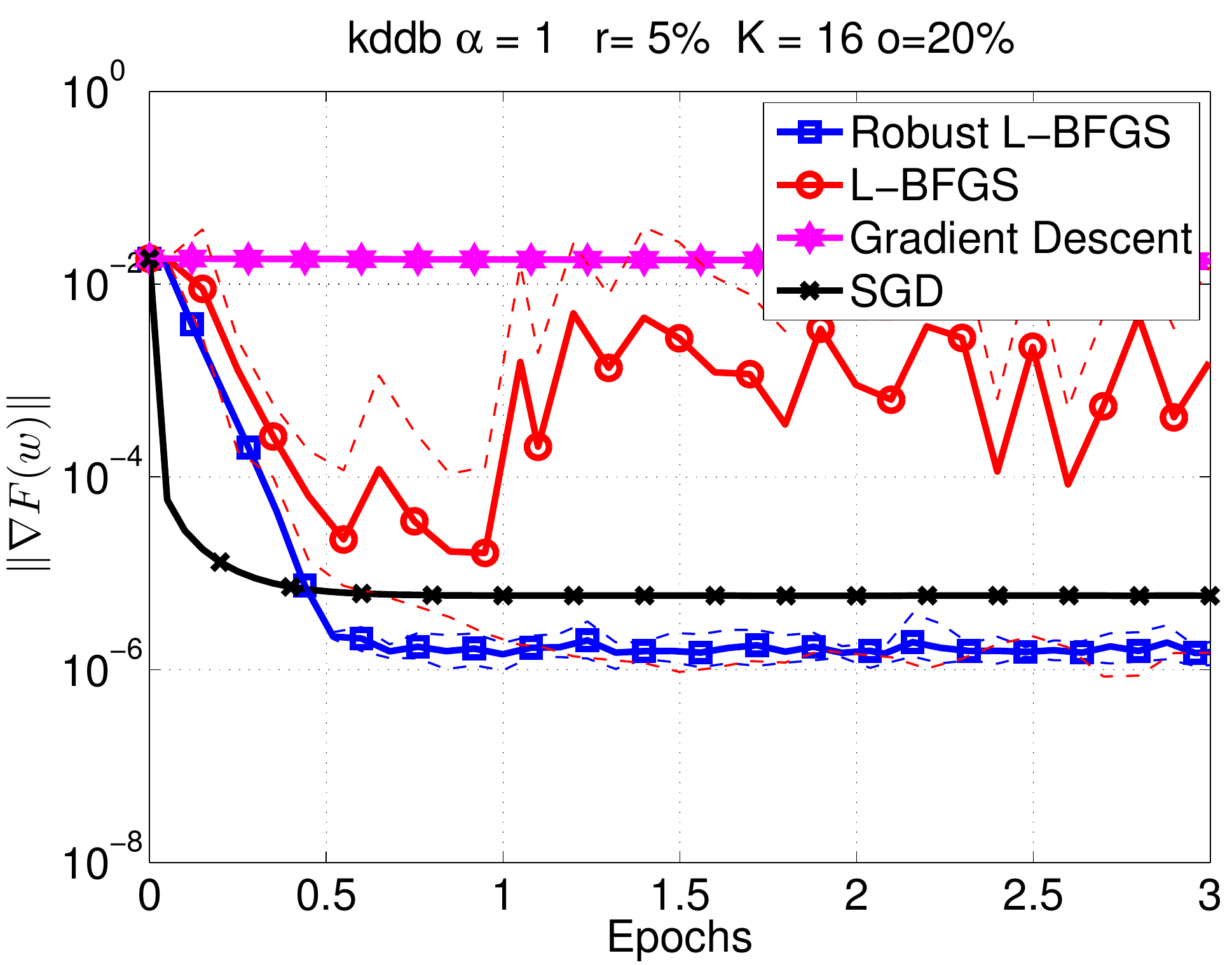}
\includegraphics[width=4.6cm]{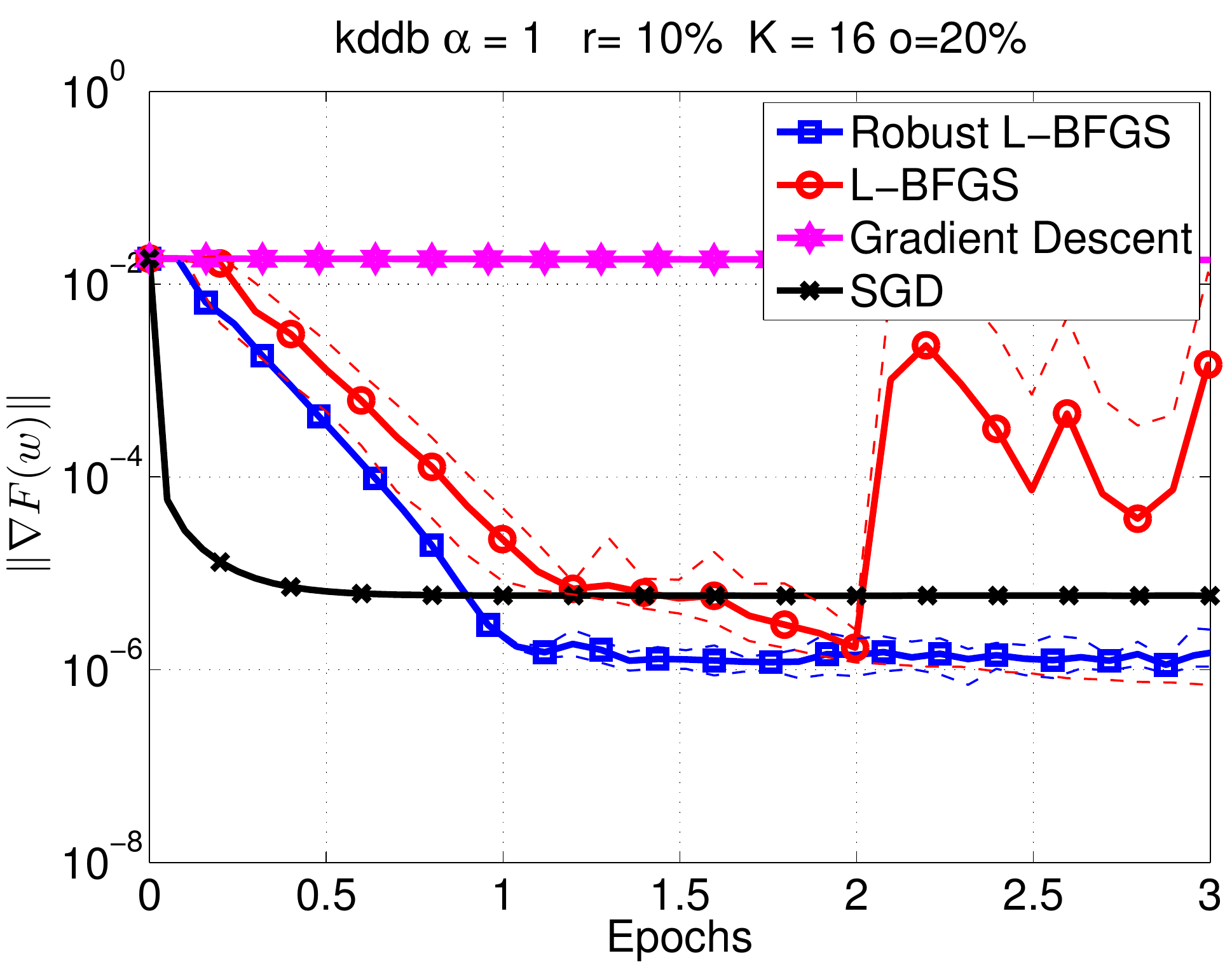}

\includegraphics[width=4.6cm]{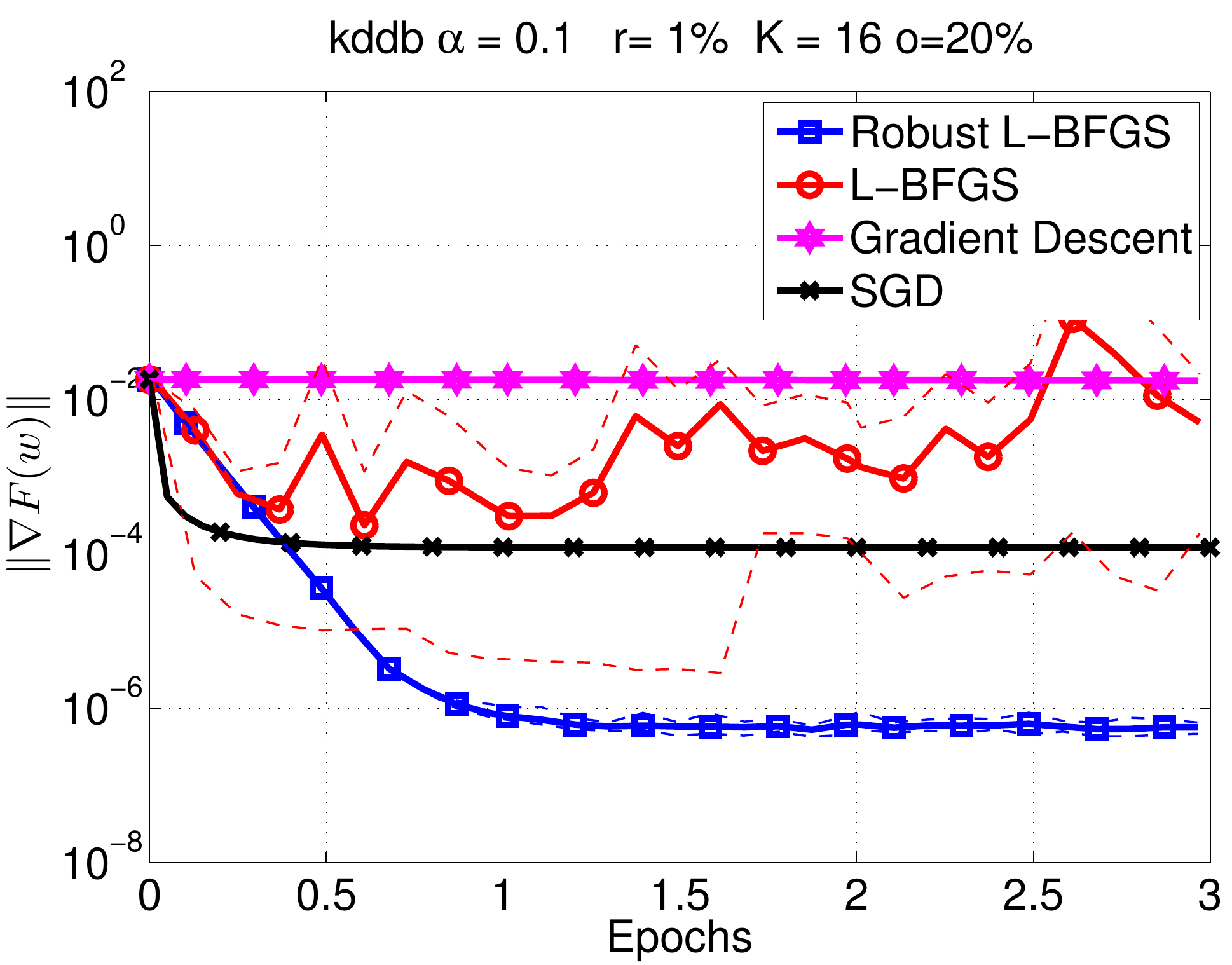}
\includegraphics[width=4.6cm]{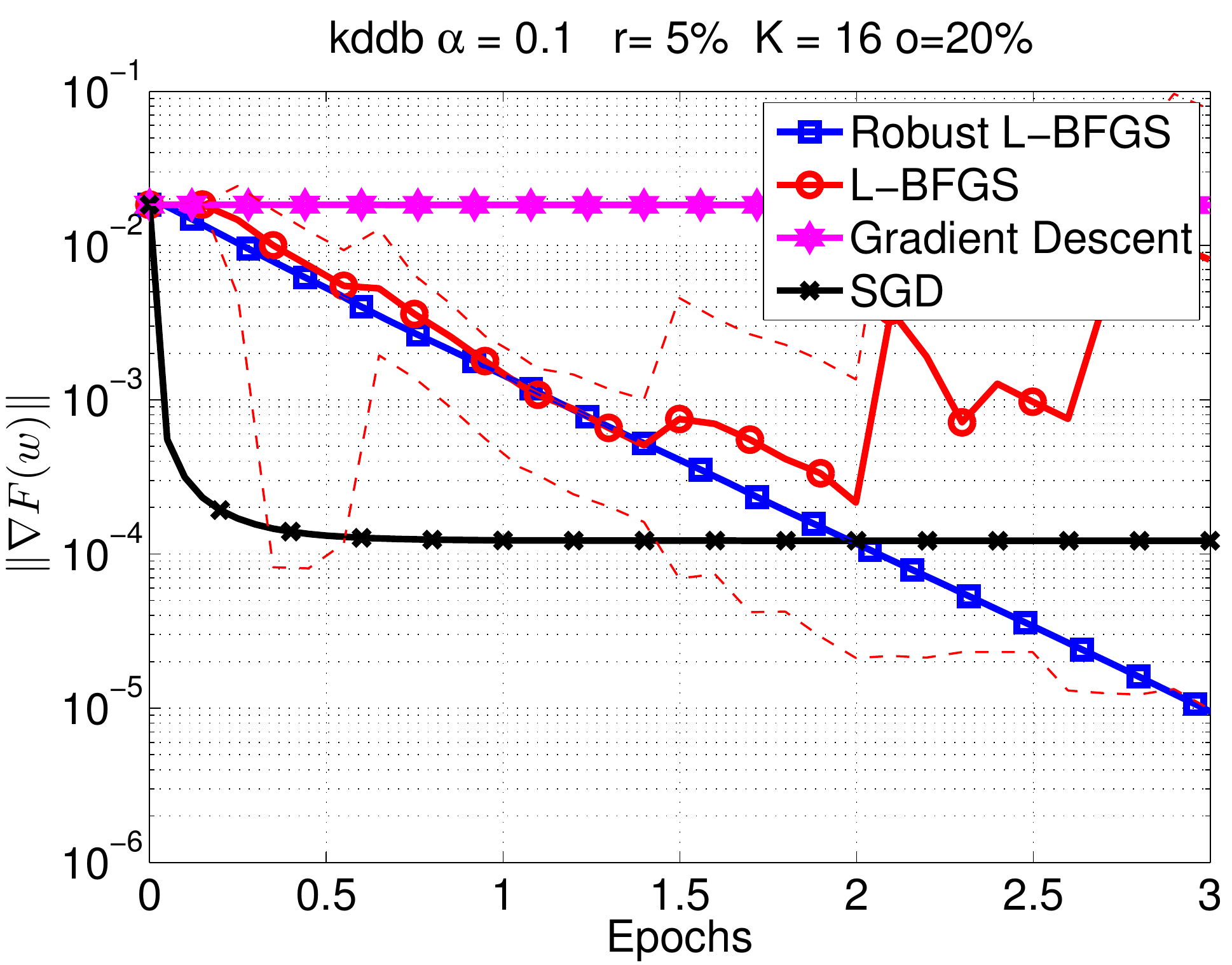}
\includegraphics[width=4.6cm]{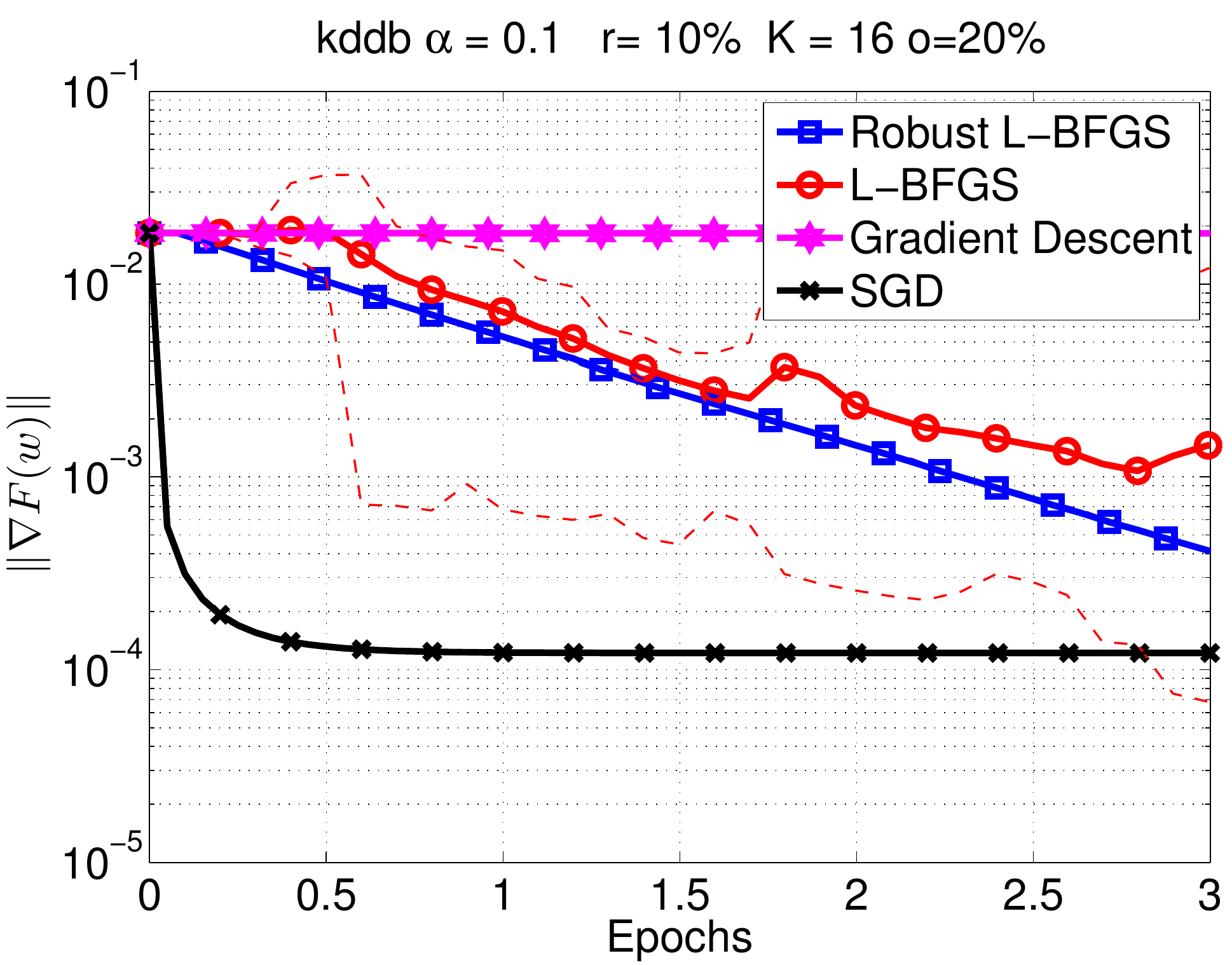}

\hrule 

\includegraphics[width=4.6cm]{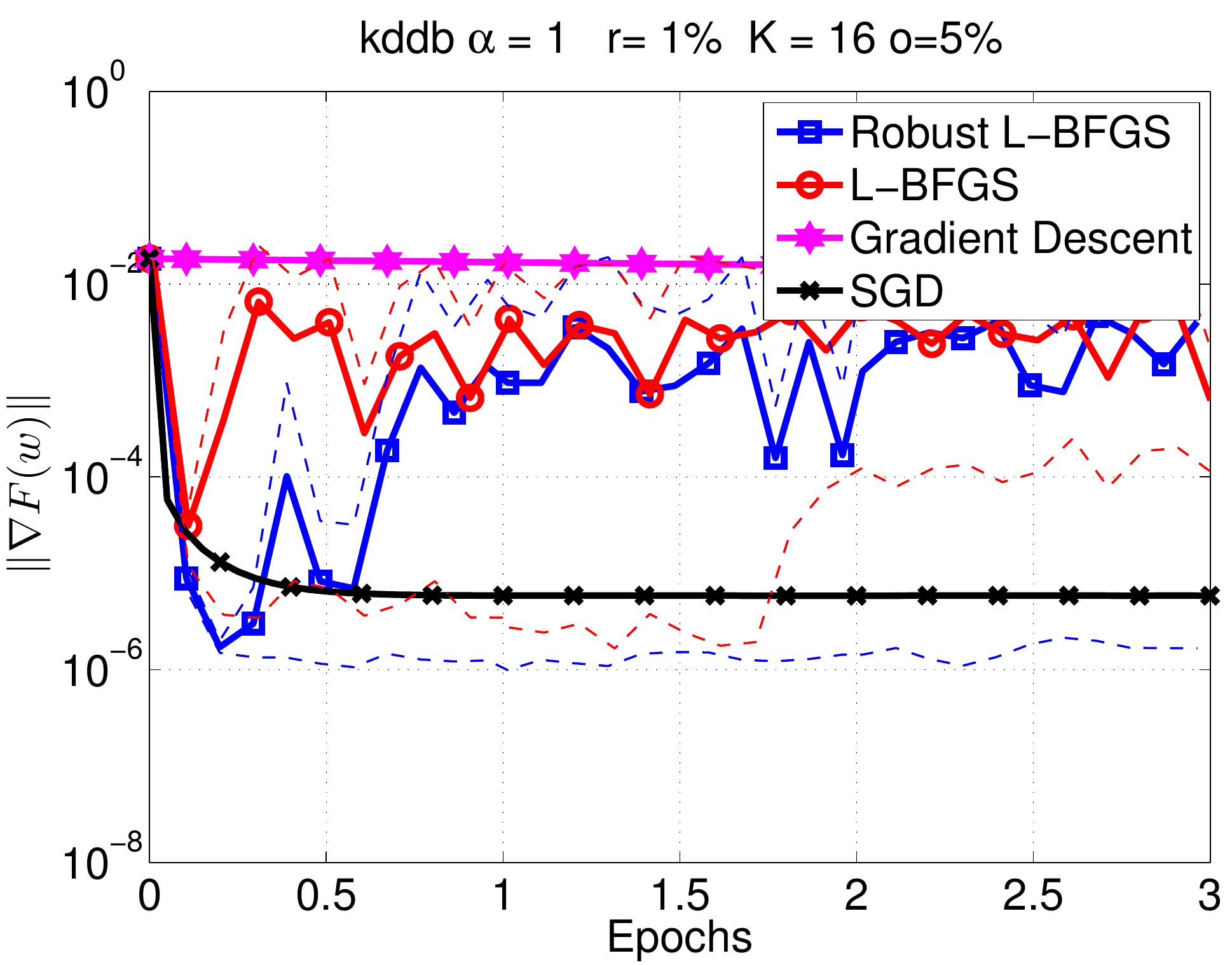}
\includegraphics[width=4.6cm]{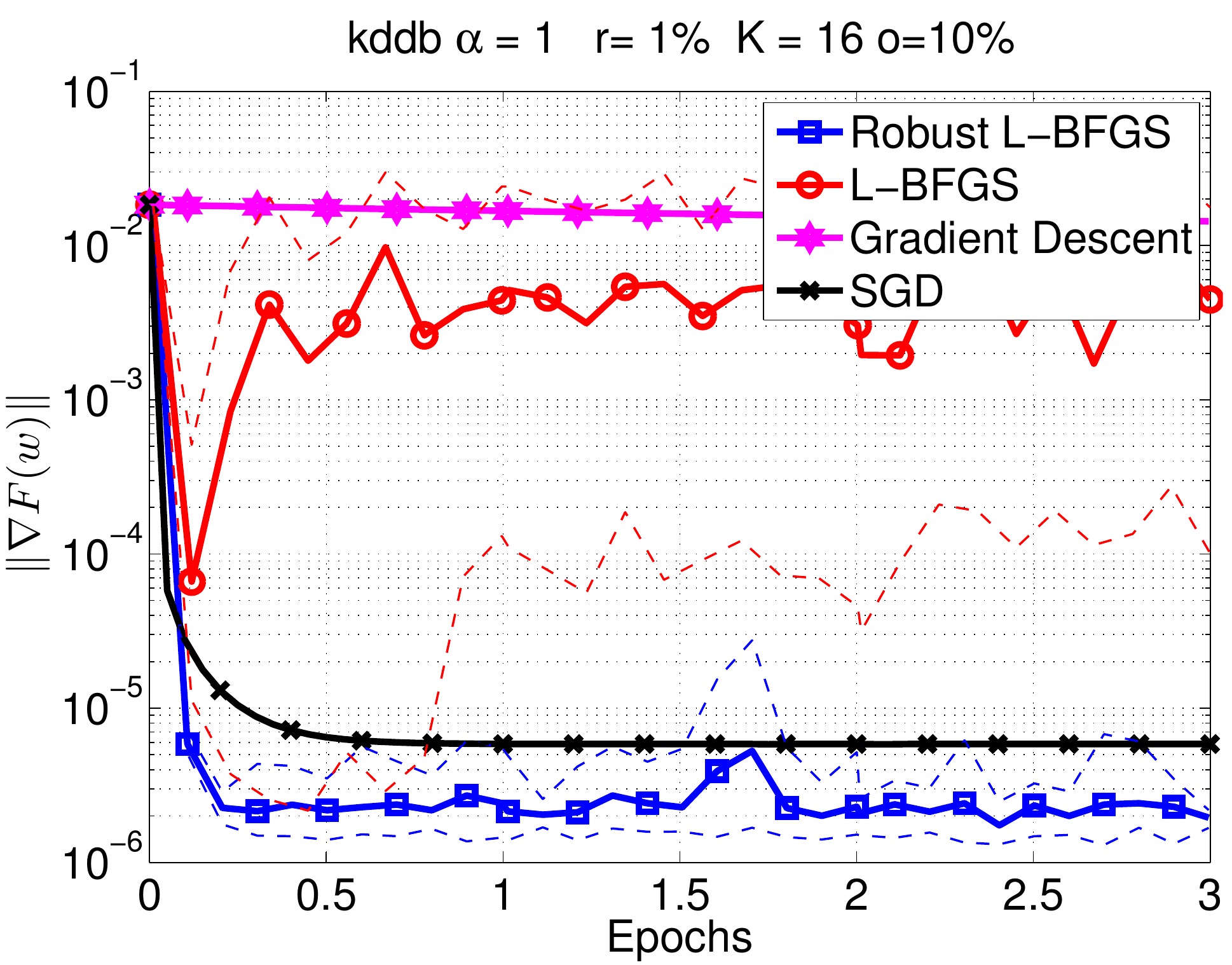}
\includegraphics[width=4.6cm]{kddb_mb_1_0_01_16_0_2-eps-converted-to.pdf}
\includegraphics[width=4.6cm]{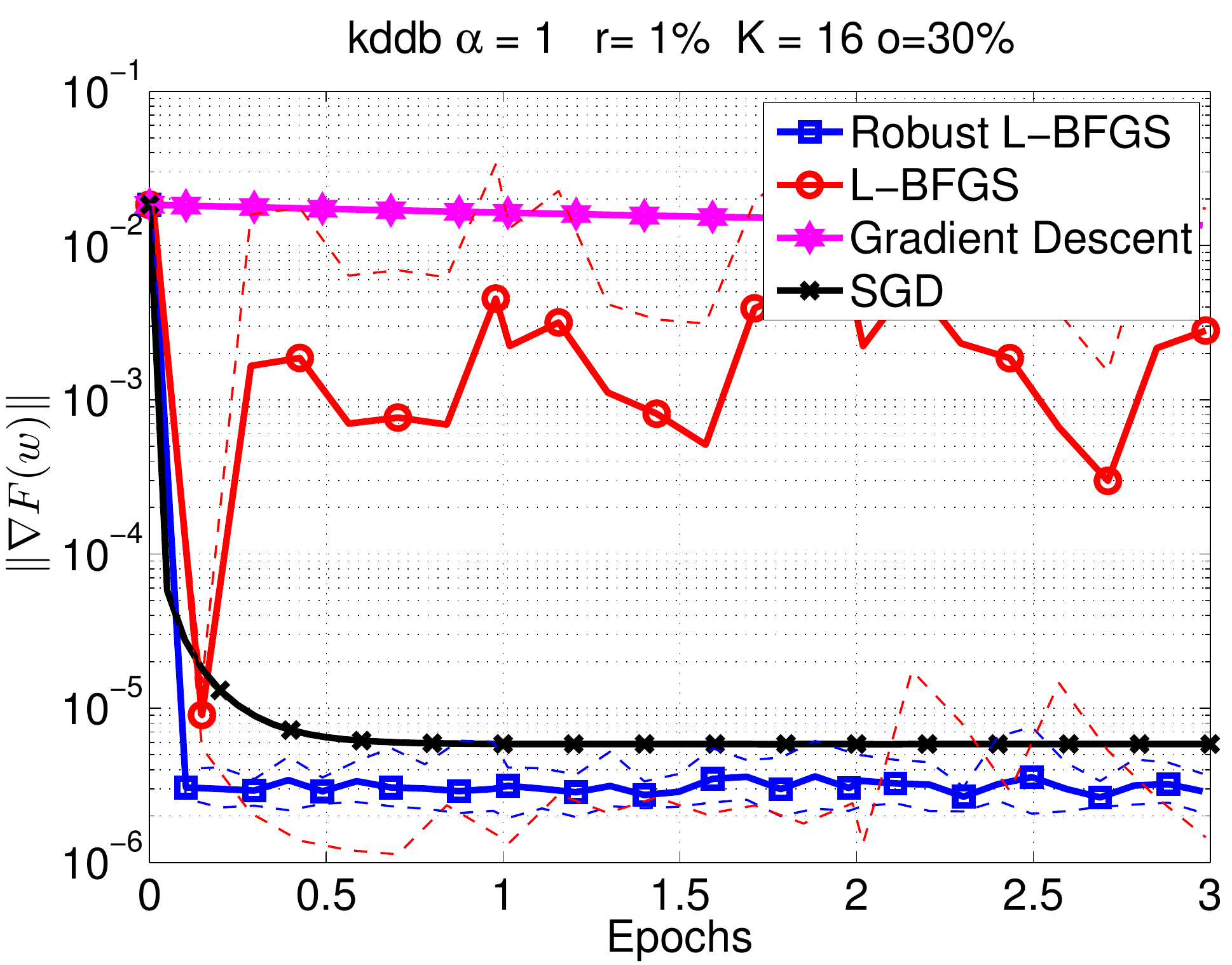}

\caption{\textbf{kddb dataset}. Comparison of Robust L-BFGS, L-BFGS (multi-batch L-BFGS without enforcing sample consistency), Gradient Descent (multi-batch Gradient method) and SGD. Top part:
we used $\alpha \in \{1, 0.1\}$,
$r\in \{1\%,  5\%,  10\%\}$ and $o=20\%$.
Bottom part: we used $\alpha=1$, $r=1\%$ and
$o\in \{5\%,  10\%, 20\%, 30\%\}$. Solid lines show average performance, and dashed lines show worst and best performance, over 10 runs (per algorithm). $K=16$ MPI processes.}
\end{figure}

\begin{figure}
\centering
\includegraphics[width=4.6cm]{webspam_mb_1_0_01_16_0_2-eps-converted-to.pdf}
\includegraphics[width=4.6cm]{webspam_mb_1_0_05_16_0_2-eps-converted-to.pdf}
\includegraphics[width=4.6cm]{webspam_mb_1_0_1_16_0_2-eps-converted-to.pdf}

\includegraphics[width=4.6cm]{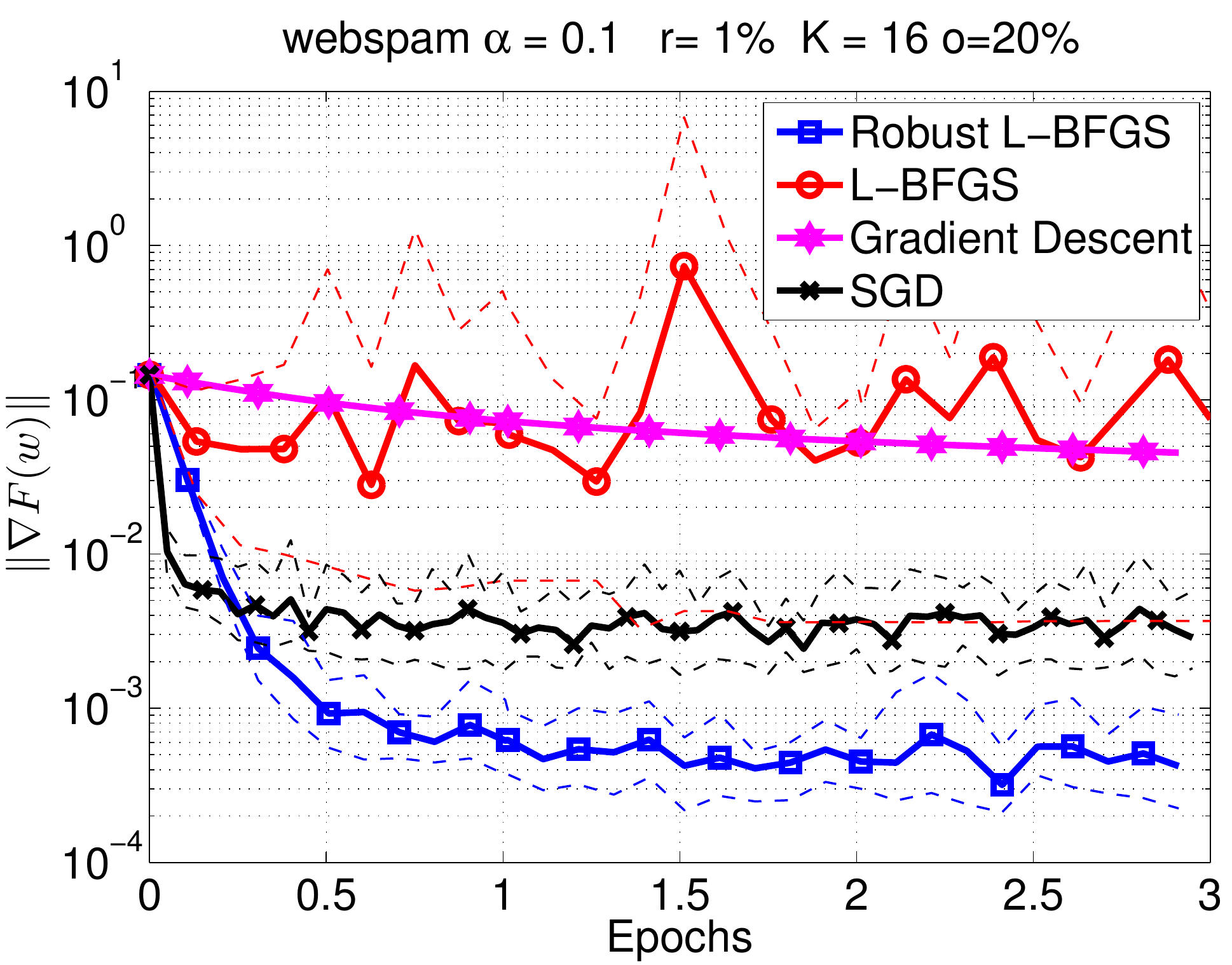}
\includegraphics[width=4.6cm]{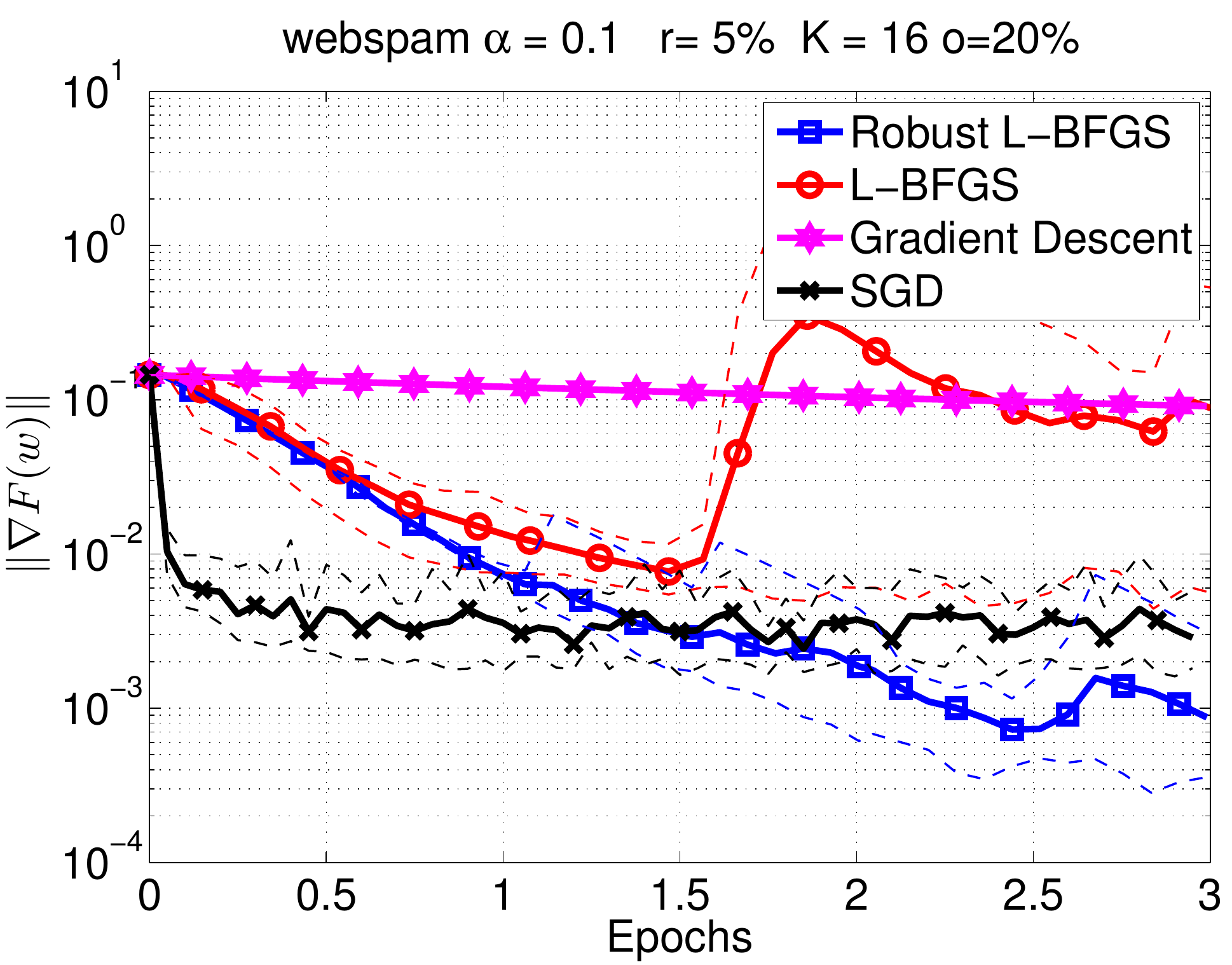}
\includegraphics[width=4.6cm]{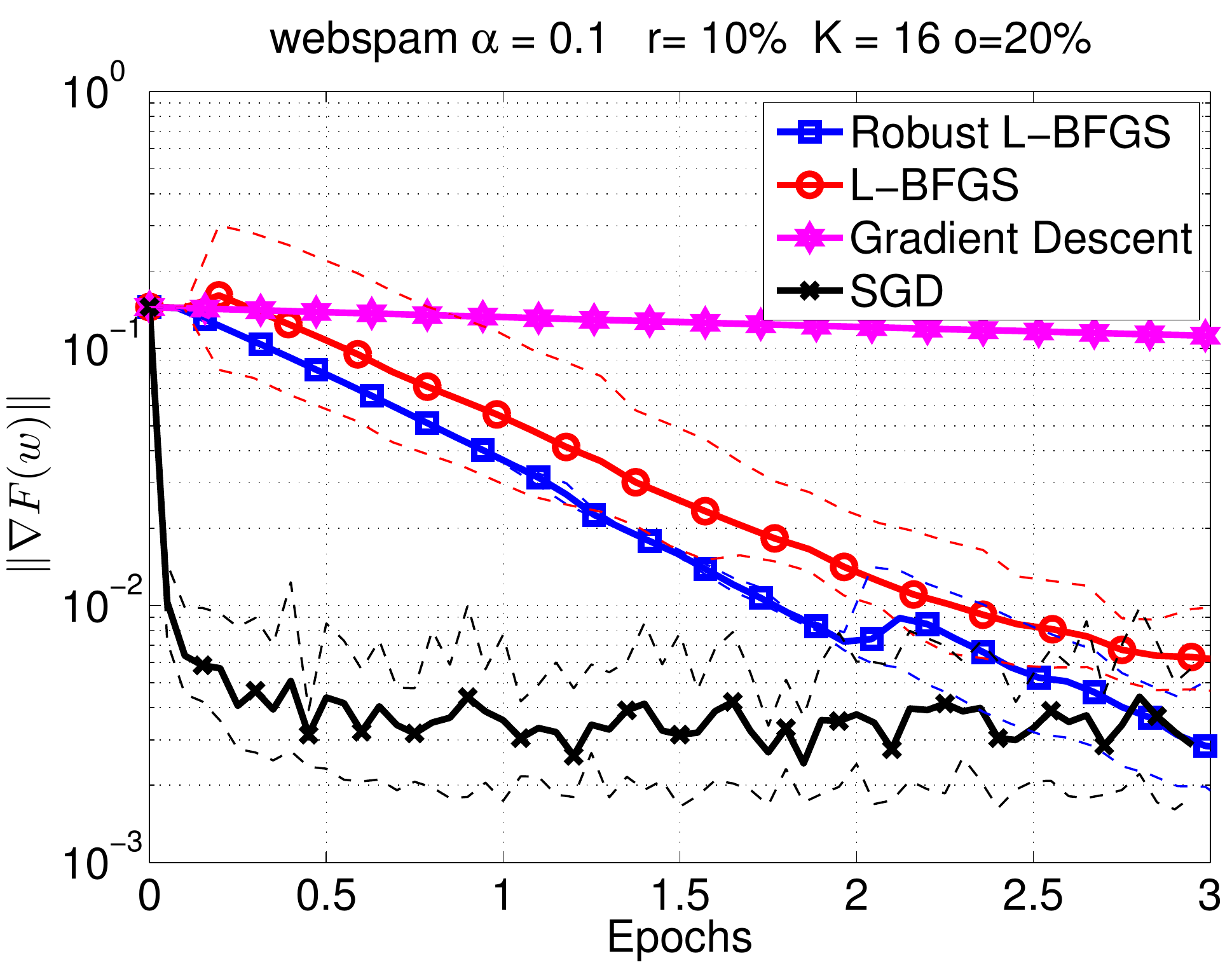}

\hrule 

\includegraphics[width=4.6cm]{webspam_mb_1_0_01_16_0_05-eps-converted-to.pdf}
\includegraphics[width=4.6cm]{webspam_mb_1_0_01_16_0_1-eps-converted-to.pdf}
\includegraphics[width=4.6cm]{webspam_mb_1_0_01_16_0_2-eps-converted-to.pdf}
\includegraphics[width=4.6cm]{webspam_mb_1_0_01_16_0_3-eps-converted-to.pdf}

\caption{\textbf{webspam dataset}. Comparison of Robust L-BFGS, L-BFGS (multi-batch L-BFGS without enforcing sample consistency), Gradient Descent (multi-batch Gradient method) and SGD. Top part:
we used $\alpha \in \{1, 0.1\}$,
$r\in \{1\%,  5\%,  10\%\}$ and $o=20\%$.
Bottom part: we used $\alpha=1$, $r=1\%$ and
$o\in \{5\%,  10\%, 20\%, 30\%\}$. Solid lines show average performance, and dashed lines show worst and best performance, over 10 runs (per algorithm). $K=16$ MPI processes.}
\end{figure}

\clearpage
\subsection{Fault-tolerant L-BFGS Implementation}
\label{sec:ext_numerical_fault}

If we run a distributed algorithm, for example on a shared computer cluster, then we may experience  delays. Such delays can be caused by other processes running on the same compute node, node failures and for other reasons. As a result, given a computational (time) budget, these delays may cause nodes to fail to return a value. To illustrate this behavior, and to motivate the robust fault-tolerant L-BFGS method, we run a simple benchmark MPI code on two different environments:
\begin{itemize}
\item {\bf Amazon EC2} -- Amazon EC2 is a cloud system provided by Amazon. It is expected that if load balancing is done properly, the execution time will have small noise; however, the network and communication can still be an issue. (4 MPI processes)
\item {\bf Shared Cluster} -- In our shared cluster, multiple jobs run on each node, with some jobs being more demanding than others. Even though each node has 16 cores, the amount of resources each job can utilize changes over time. In terms of communication, we have a GigaBit network. (11 MPI processes, running on 11 nodes)
\end{itemize}

We run a simple code on the cloud/cluster, with MPI communication. We generate two matrices $A,B \in  R^{n \times n}$, then  synchronize all MPI processes and compute $C=A\cdot B$ using the GSL C-BLAS library. The time is measured and recorded as computational time. After the matrix product is computed, the result is sent to the master/root node using asynchronous communication, and the time required is recorded. The process is repeated 3000 times.



\begin{figure}[h!]

\centering

\includegraphics[width=4.2cm]{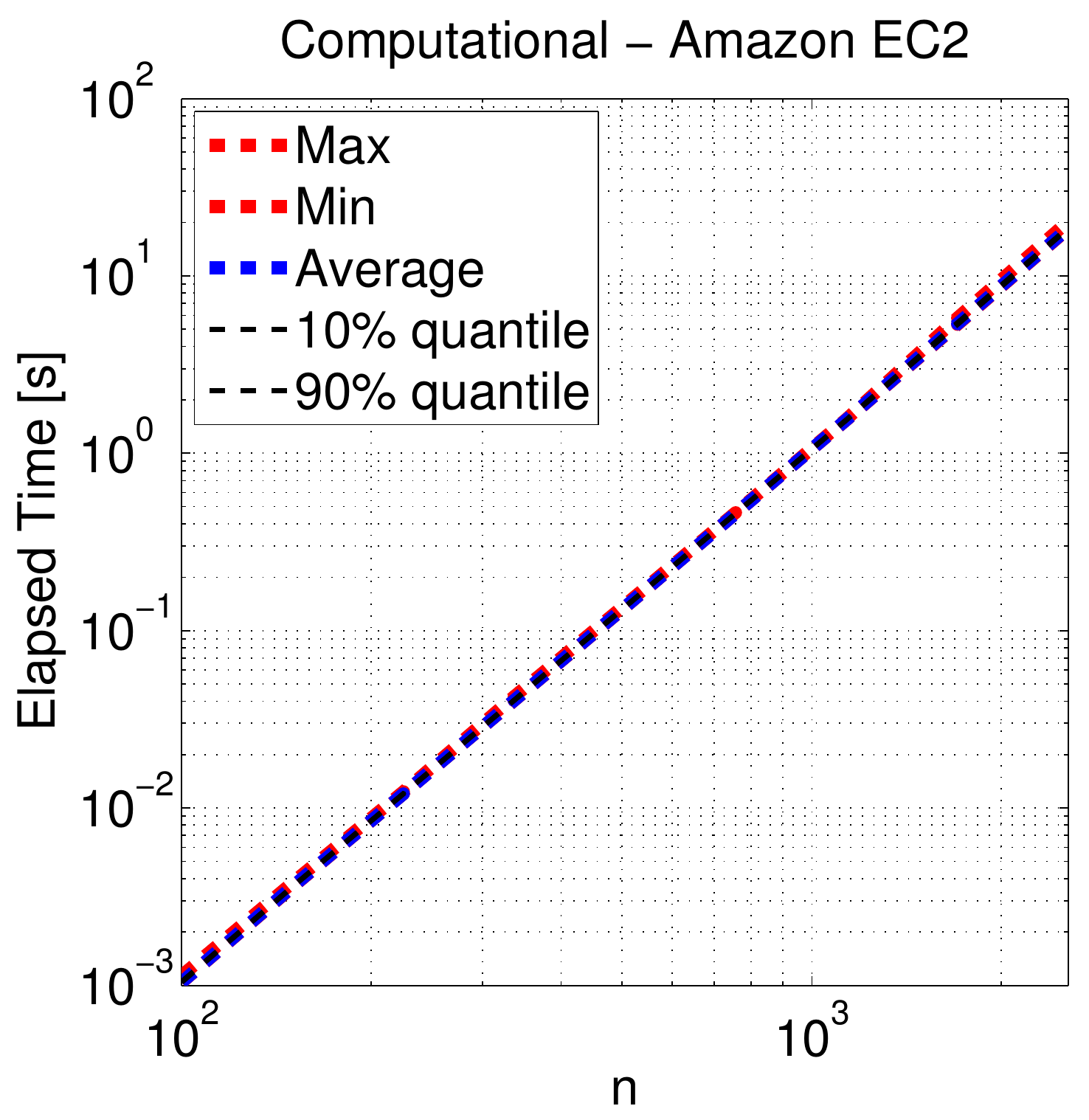}
\includegraphics[width=4.2cm]{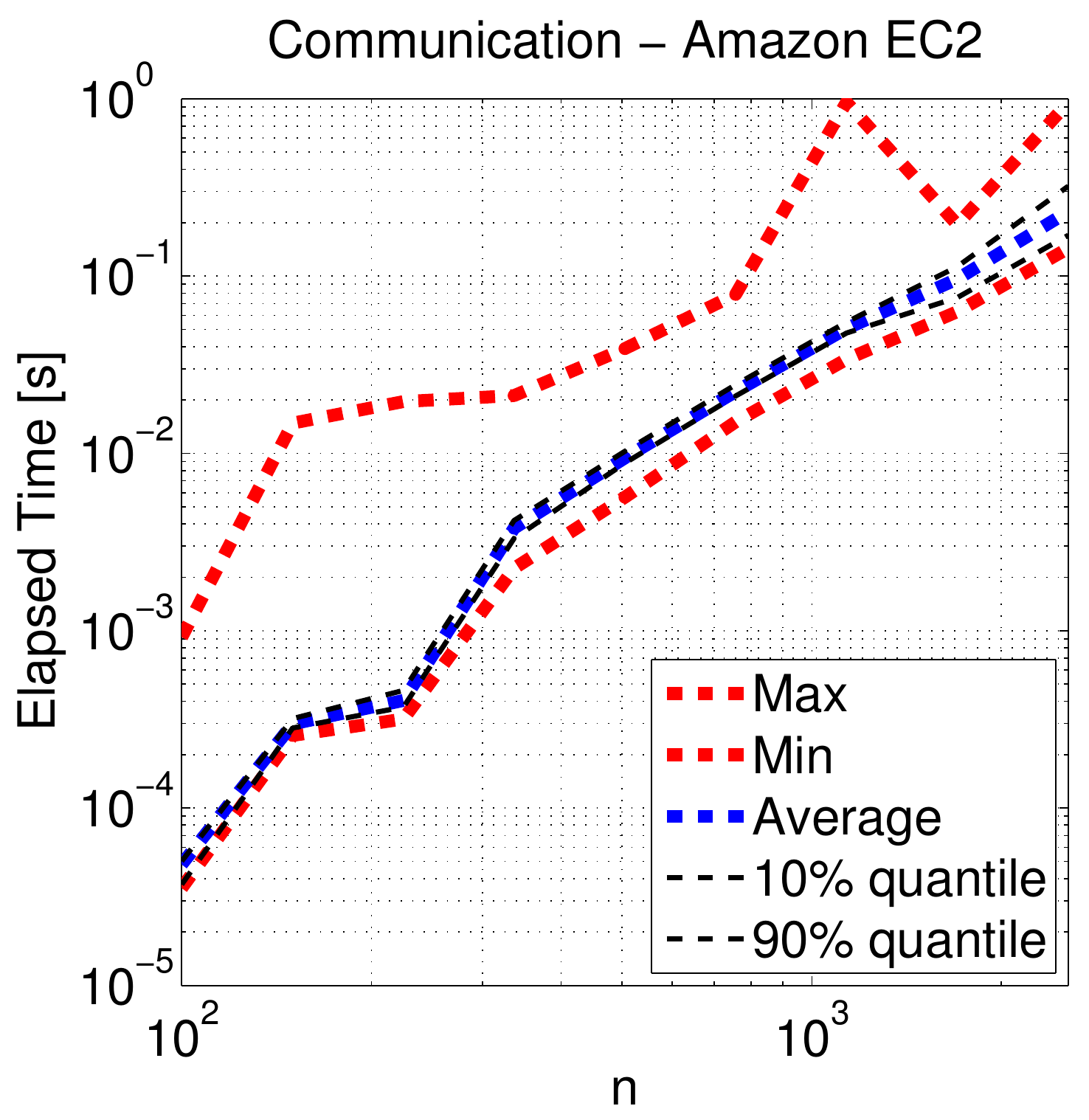}

\includegraphics[width=4.2cm]{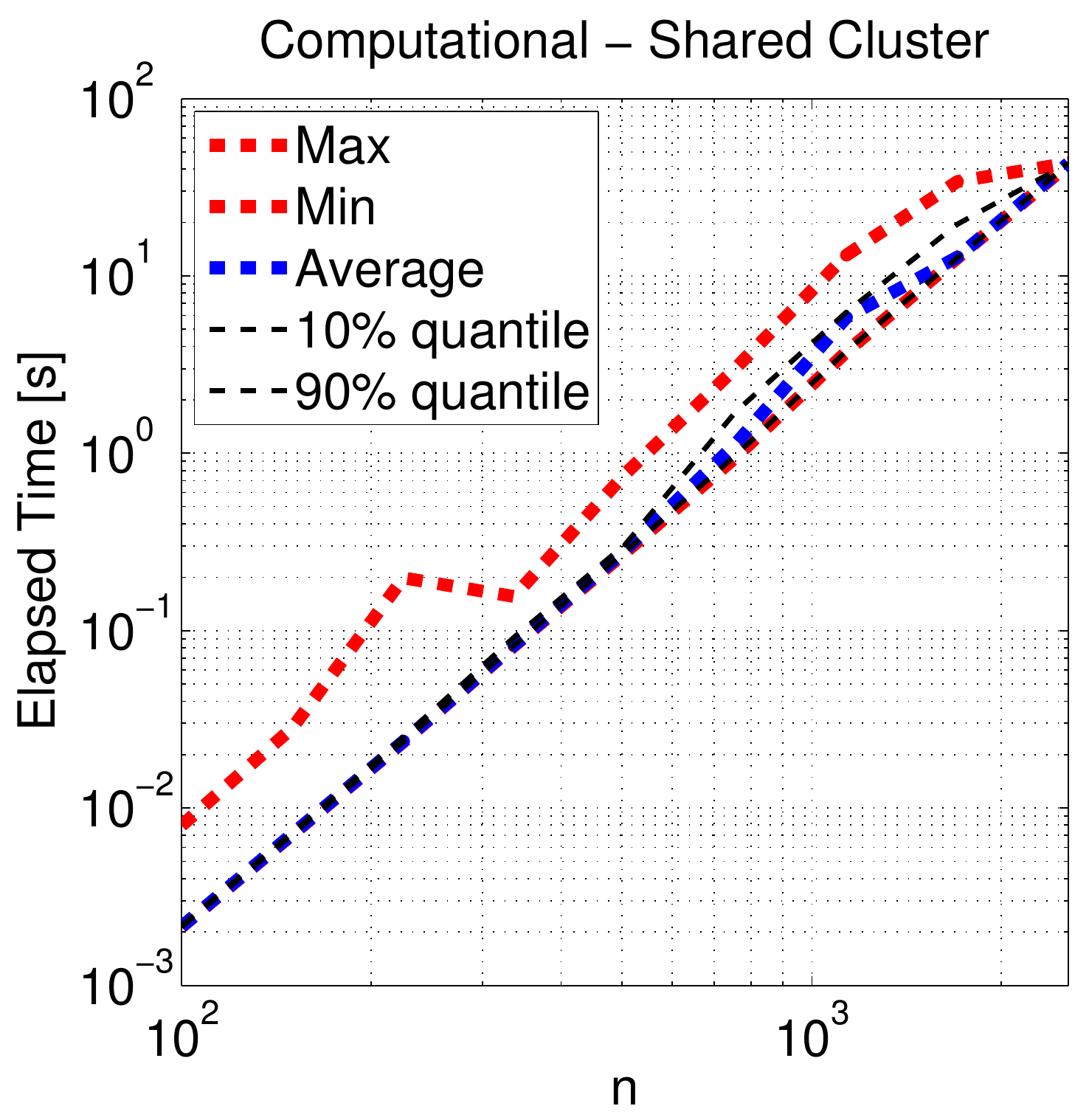}
\includegraphics[width=4.2cm]{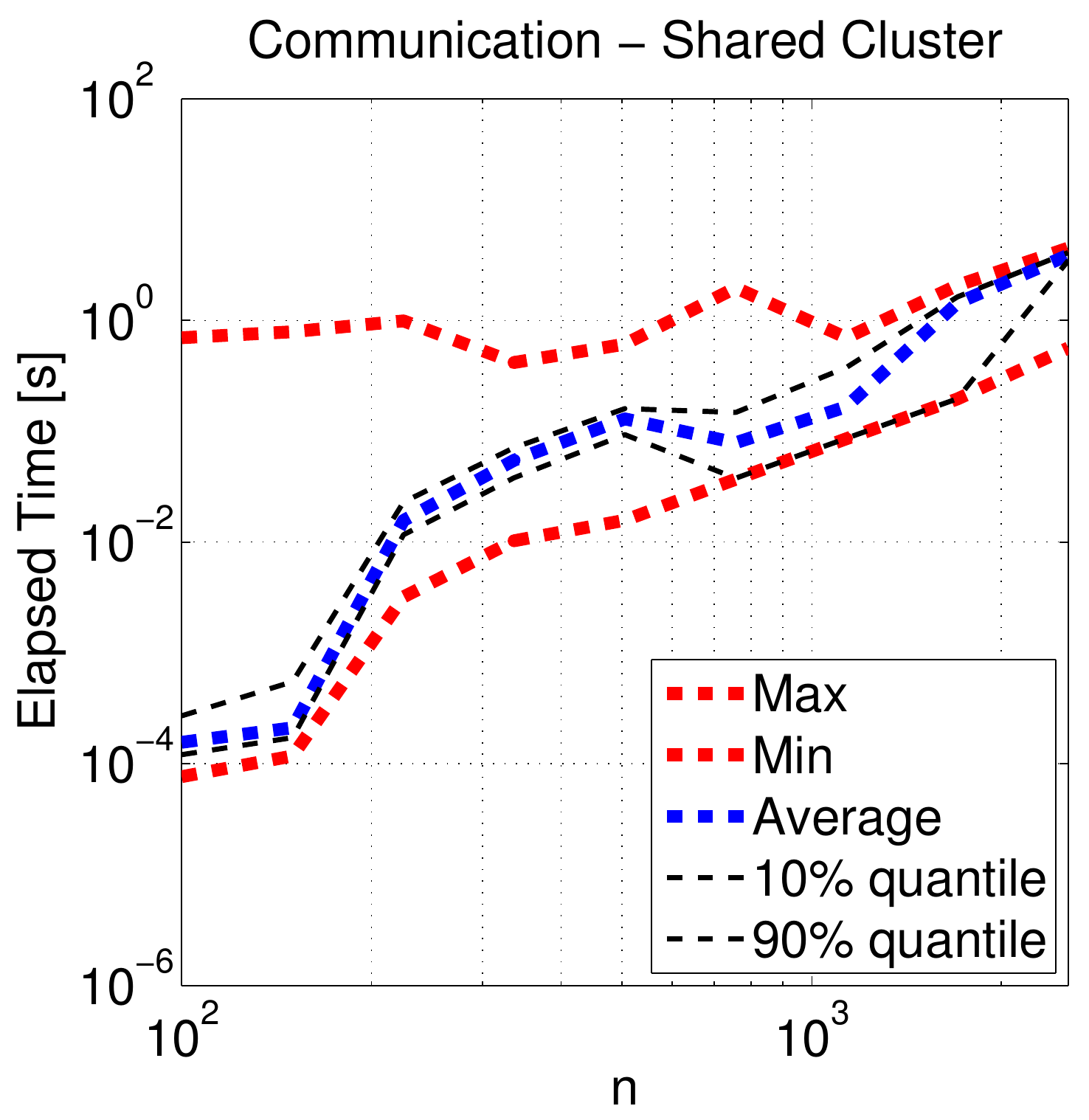}

\caption{Distribution of Computation and Communication Time for  Amazon EC2 and Shared Cluster.
Figures show worst and best time, average time and 10\% and 90\% quantiles.
Amazon Cloud EC: In the experiment: 4 MPI processes; Shared Cluster: 11 MPI processes.}
\label{fig.motivation}

\end{figure}

The results of the experiment described above are captured in Figure \ref{fig.motivation}. As expected, on the Amazon EC2 cloud, the matrix-matrix multiplication takes roughly the same time for all replications and the noise in communication is relatively small. In this example the cost of communication is negligible when compared to the cost of computation. On our shared cluster, one cannot guarantee that all resources are exclusively used for a specific process, and thus, the computation and communication time is considerably more stochastic and unbalanced. For some cases the difference between the minimum and maximum computation (communication) time varies by an order of magnitude or more. Hence, on such a platform a fault-tolerant algorithm that only uses information from nodes that return an update within a preallocated budget is a natural choice.



\begin{figure}[h!]
\centering

\includegraphics[width=4.6cm]{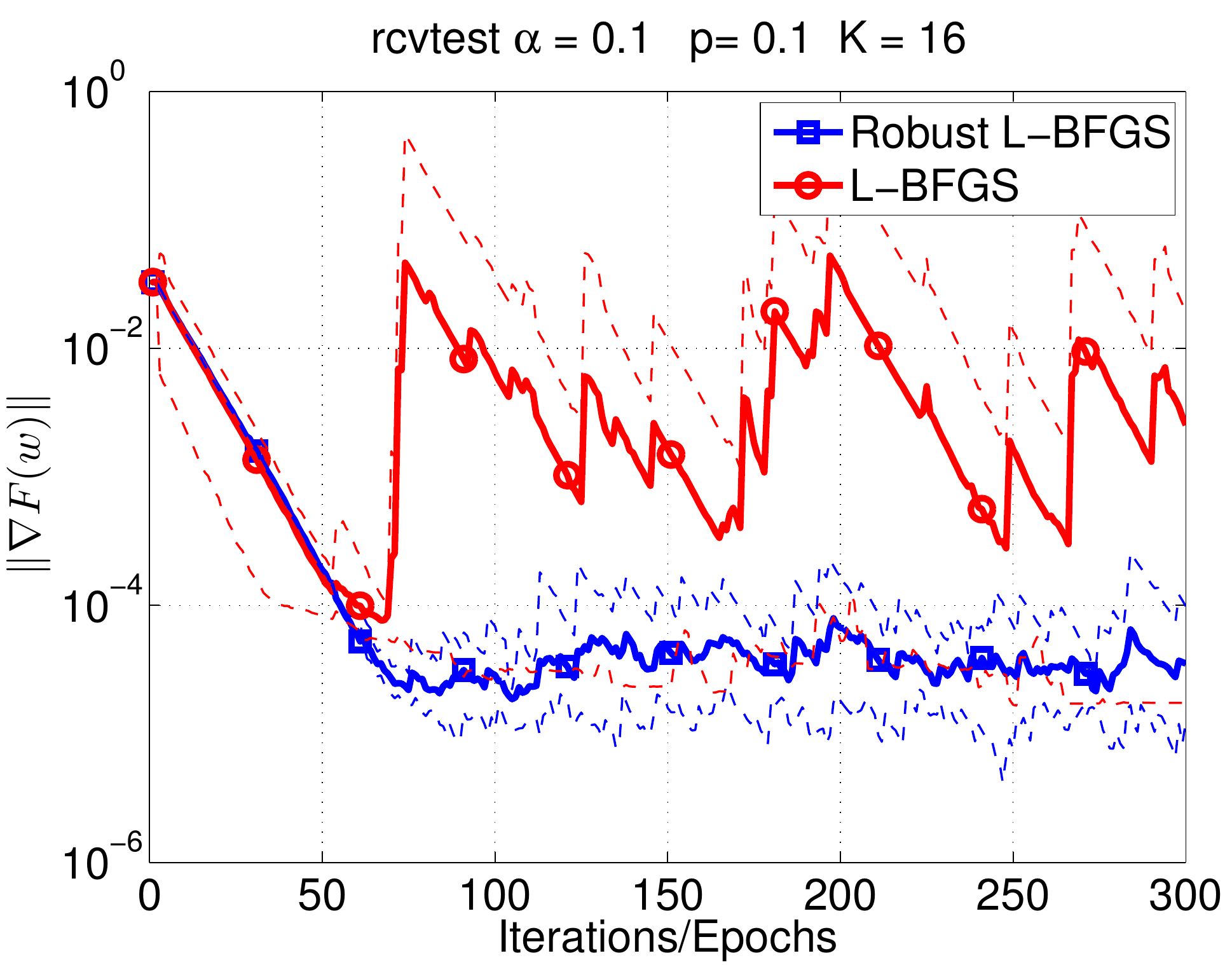}
\includegraphics[width=4.6cm]{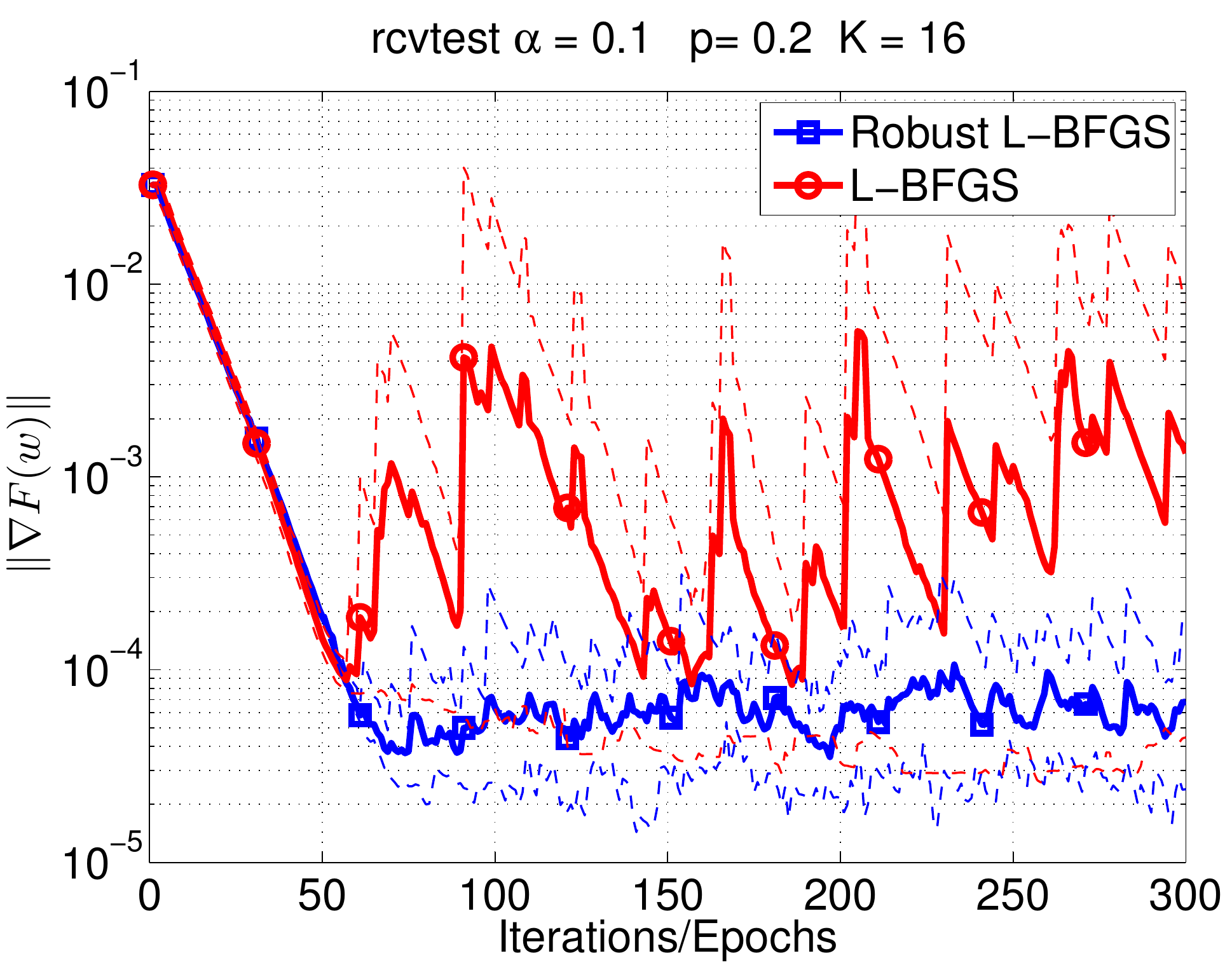}
\includegraphics[width=4.6cm]{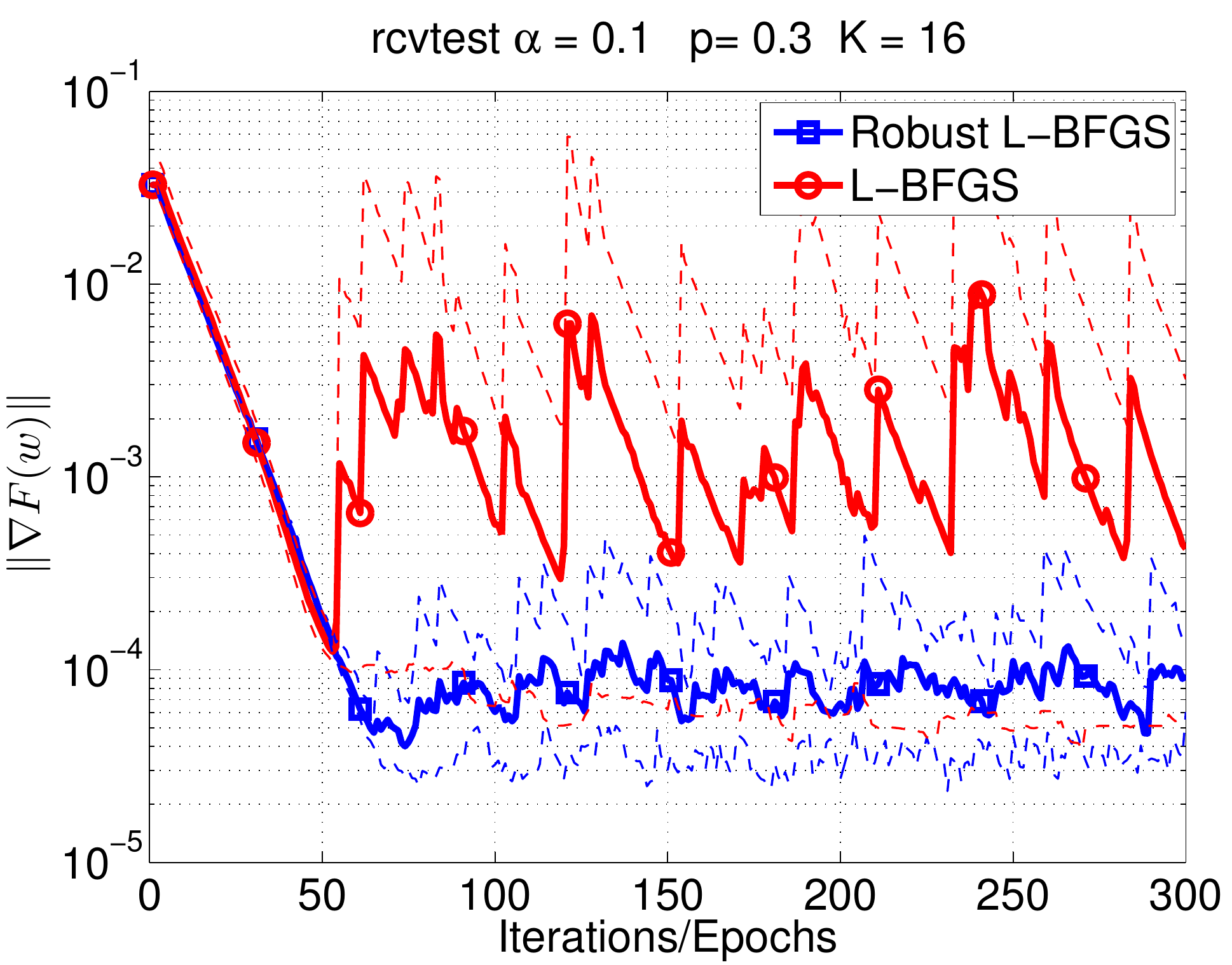}

\includegraphics[width=4.6cm]{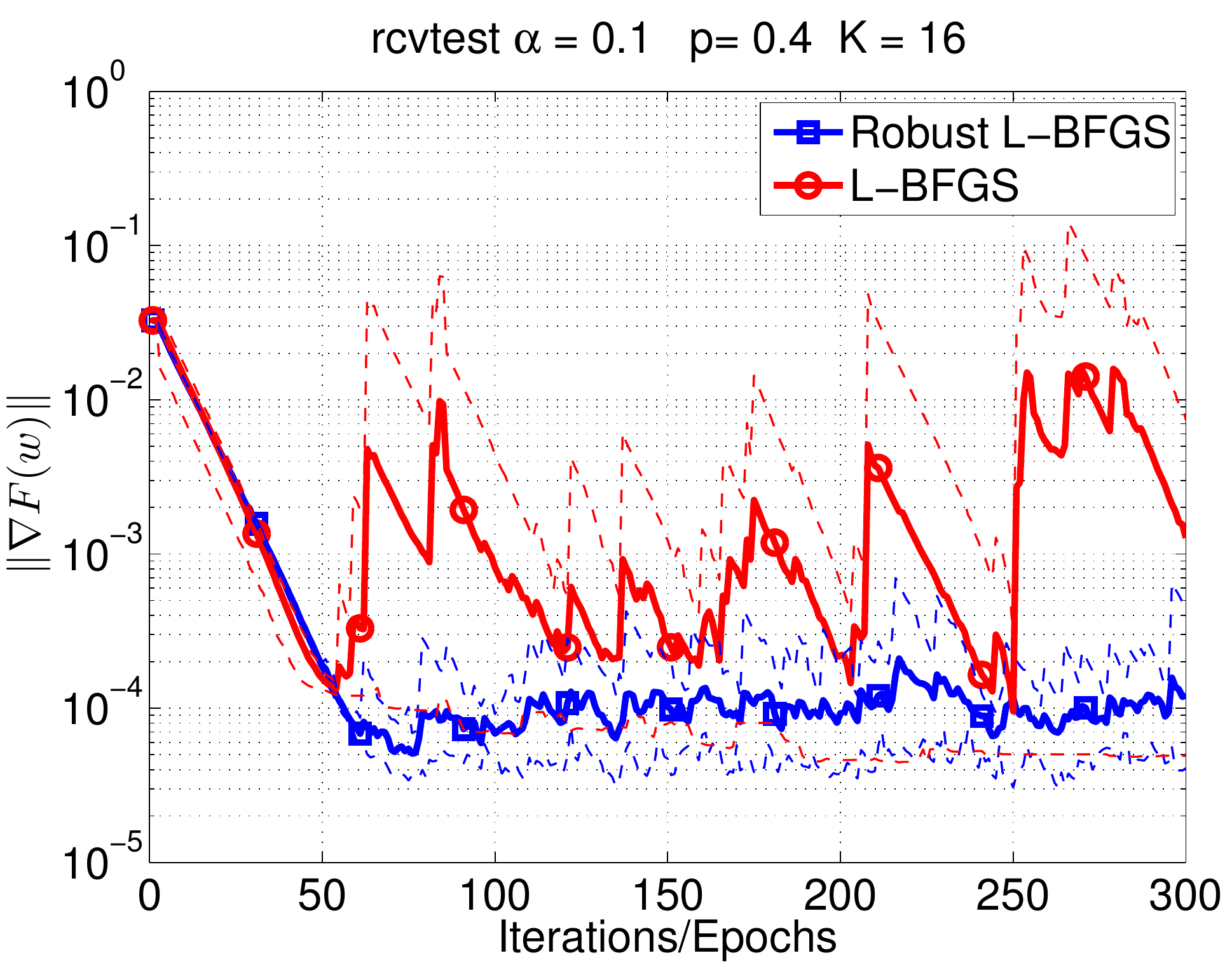}
\includegraphics[width=4.6cm]{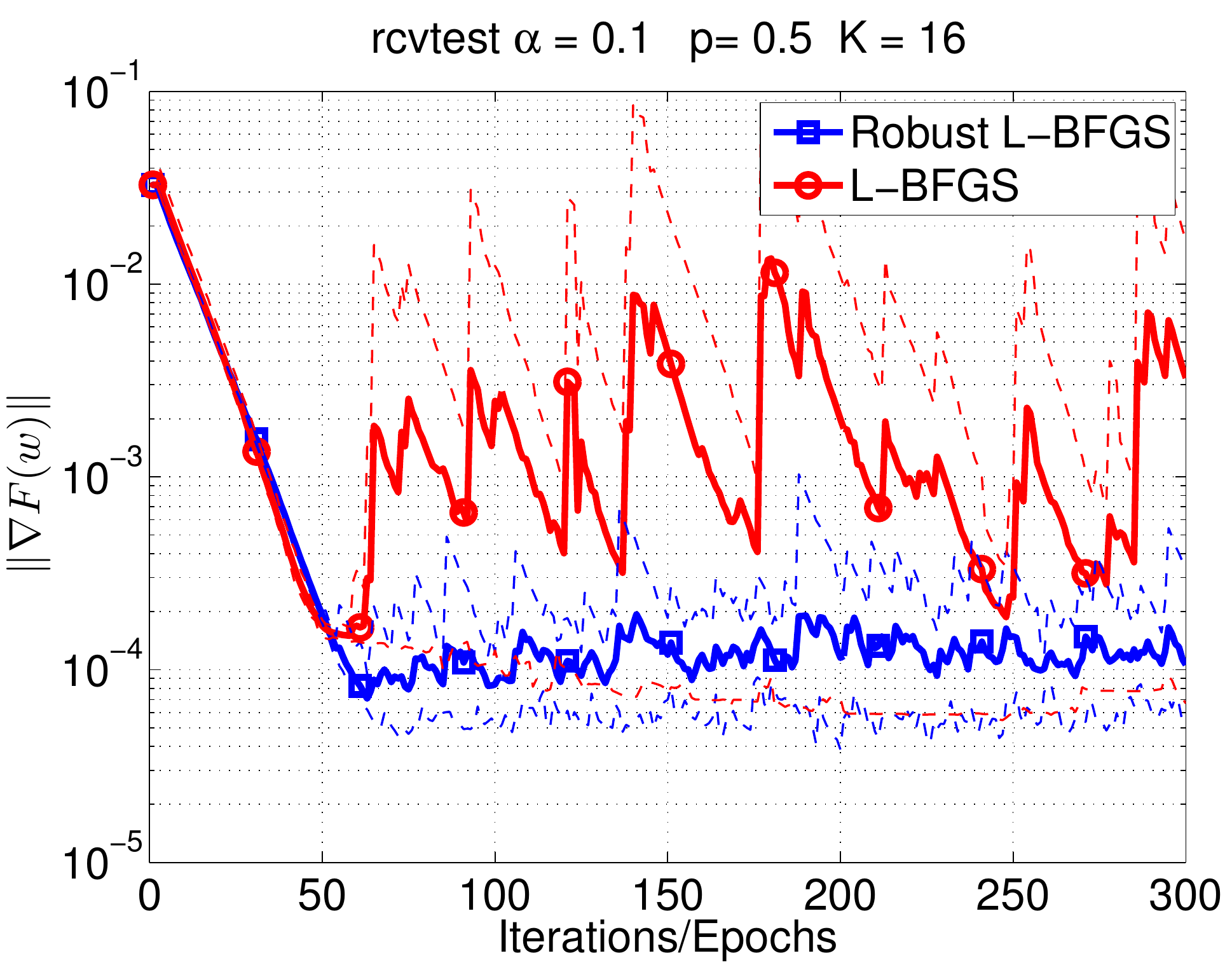}

\caption{\textbf{rcvtest dataset}. Comparison of Robust L-BFGS and L-BFGS in the presence of faults.
We used $\alpha=0.1$ and $p\in \{0.1, 0.2, 0.3, 0.4, 0.5\}$. Solid lines show average performance, and dashed lines show worst and best performance, over 10 runs (per algorithm). $K=16$ MPI processes.}
\label{fig:ft:rcv}
\end{figure}
\begin{figure}[h!]
\centering

\includegraphics[width=4.6cm]{webspam_faulttolerant_0_1_0_1_16-eps-converted-to.pdf}
\includegraphics[width=4.6cm]{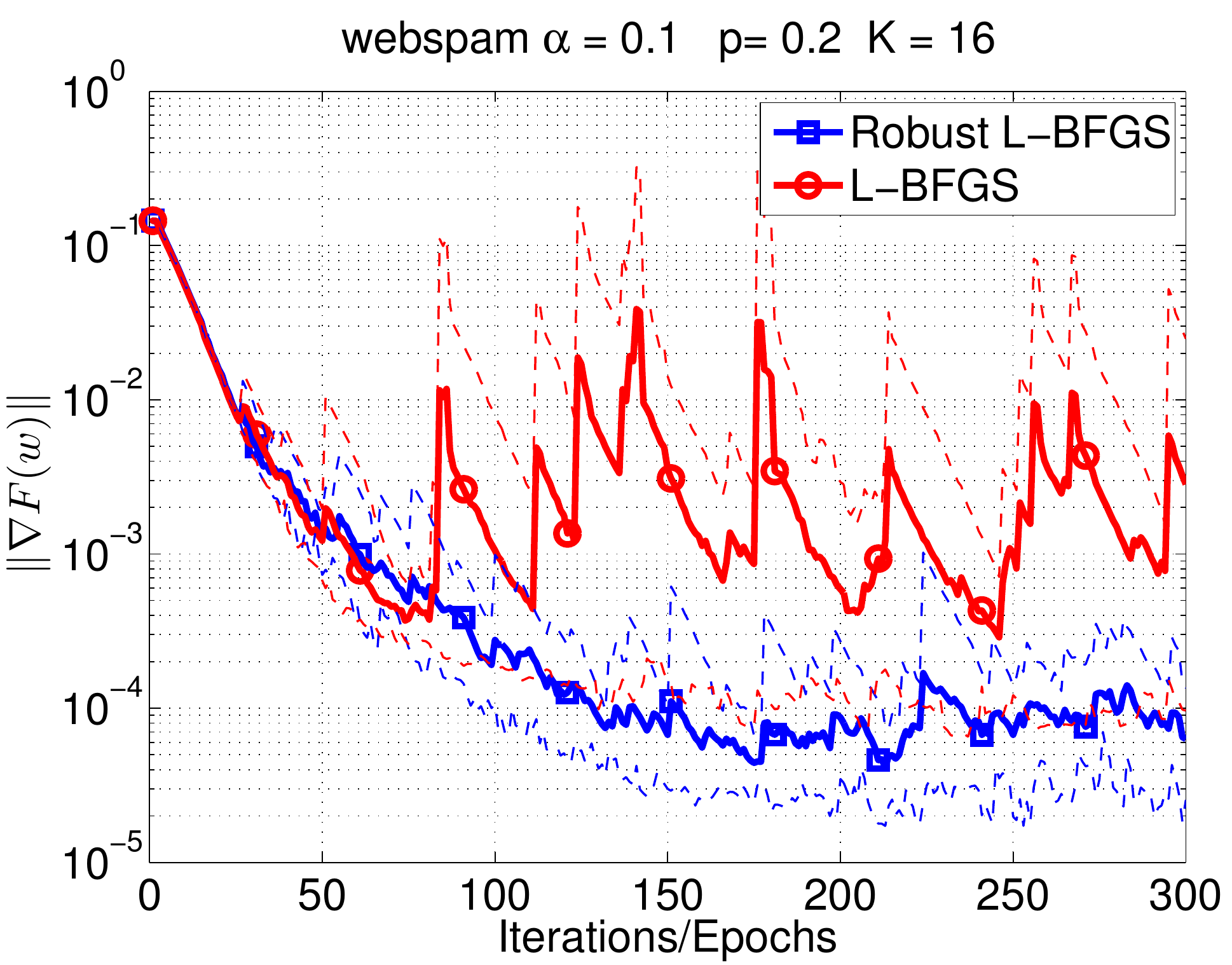}
\includegraphics[width=4.6cm]{webspam_faulttolerant_0_1_0_3_16-eps-converted-to.pdf}

\includegraphics[width=4.6cm]{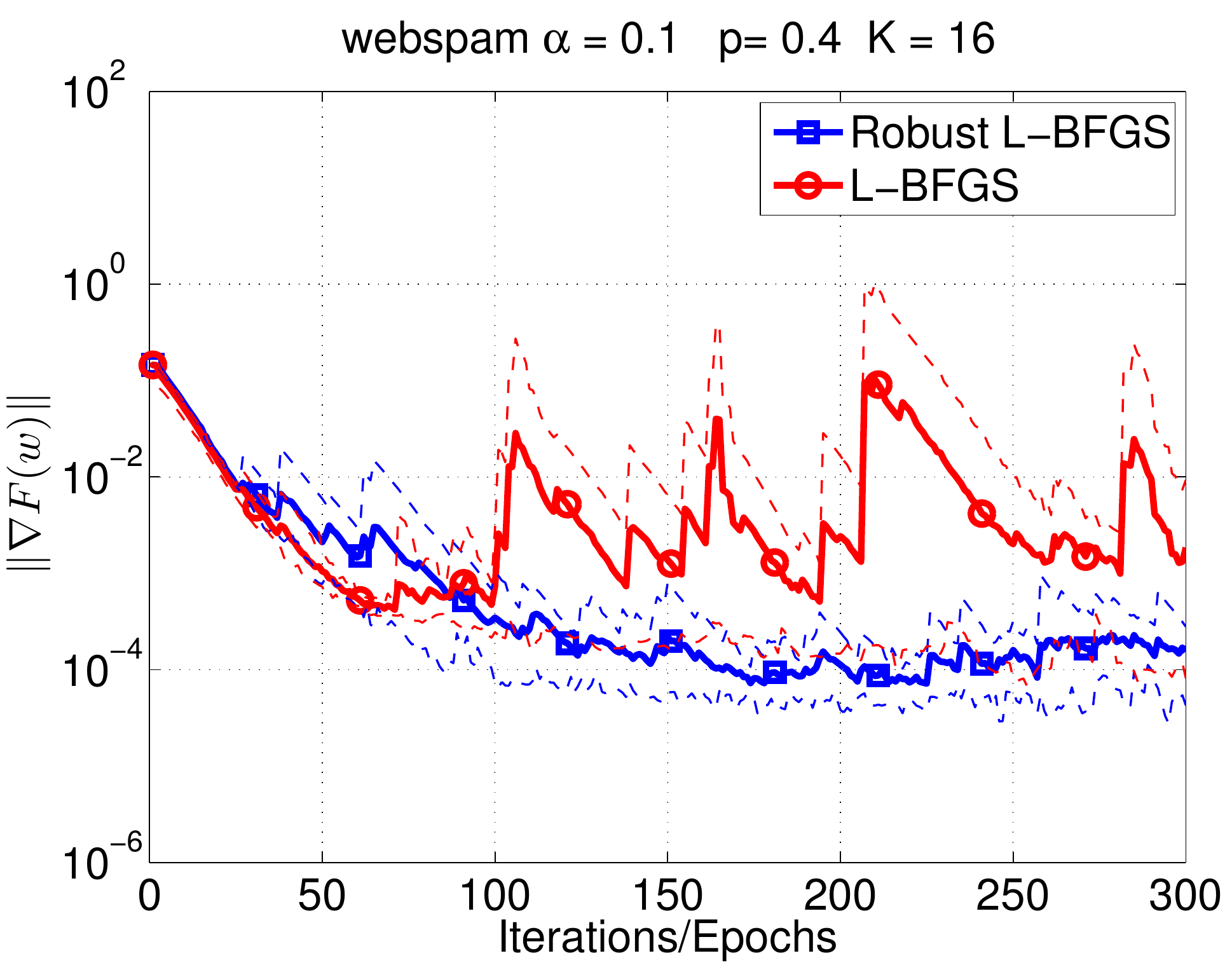}
\includegraphics[width=4.6cm]{webspam_faulttolerant_0_1_0_5_16-eps-converted-to.pdf}
 
\caption{\textbf{webspam dataset}. Comparison of Robust L-BFGS and L-BFGS in the presence of faults.
We used $\alpha=0.1$ and $p\in \{0.1, 0.2, 0.3, 0.4, 0.5\}$. Solid lines show average performance, and dashed lines show worst and best performance, over 10 runs (per algorithm). $K=16$ MPI processes.
}\label{fig:ft:webspam}
\end{figure}
\begin{figure}[h!]
\centering

\includegraphics[width=4.6cm]{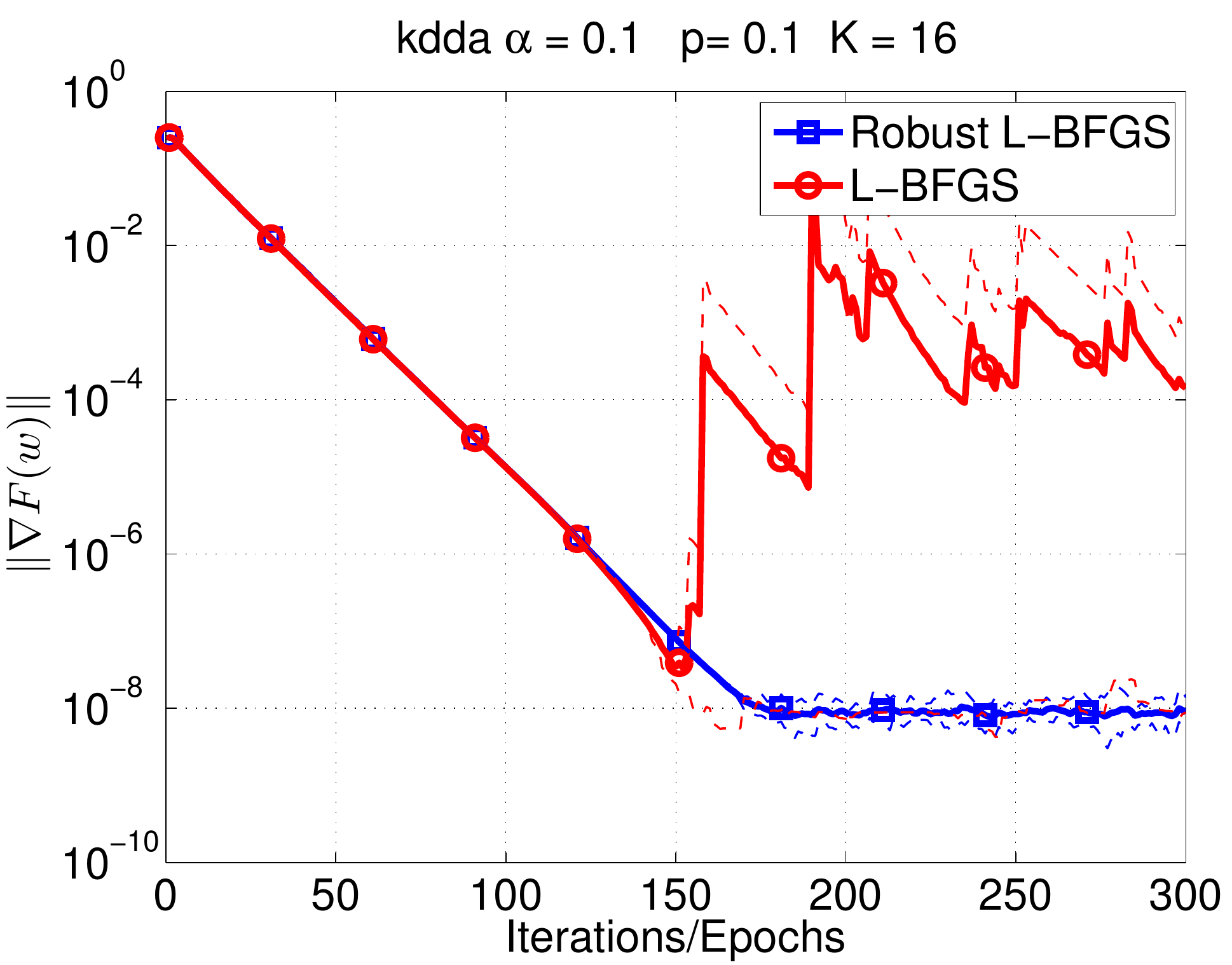}
\includegraphics[width=4.6cm]{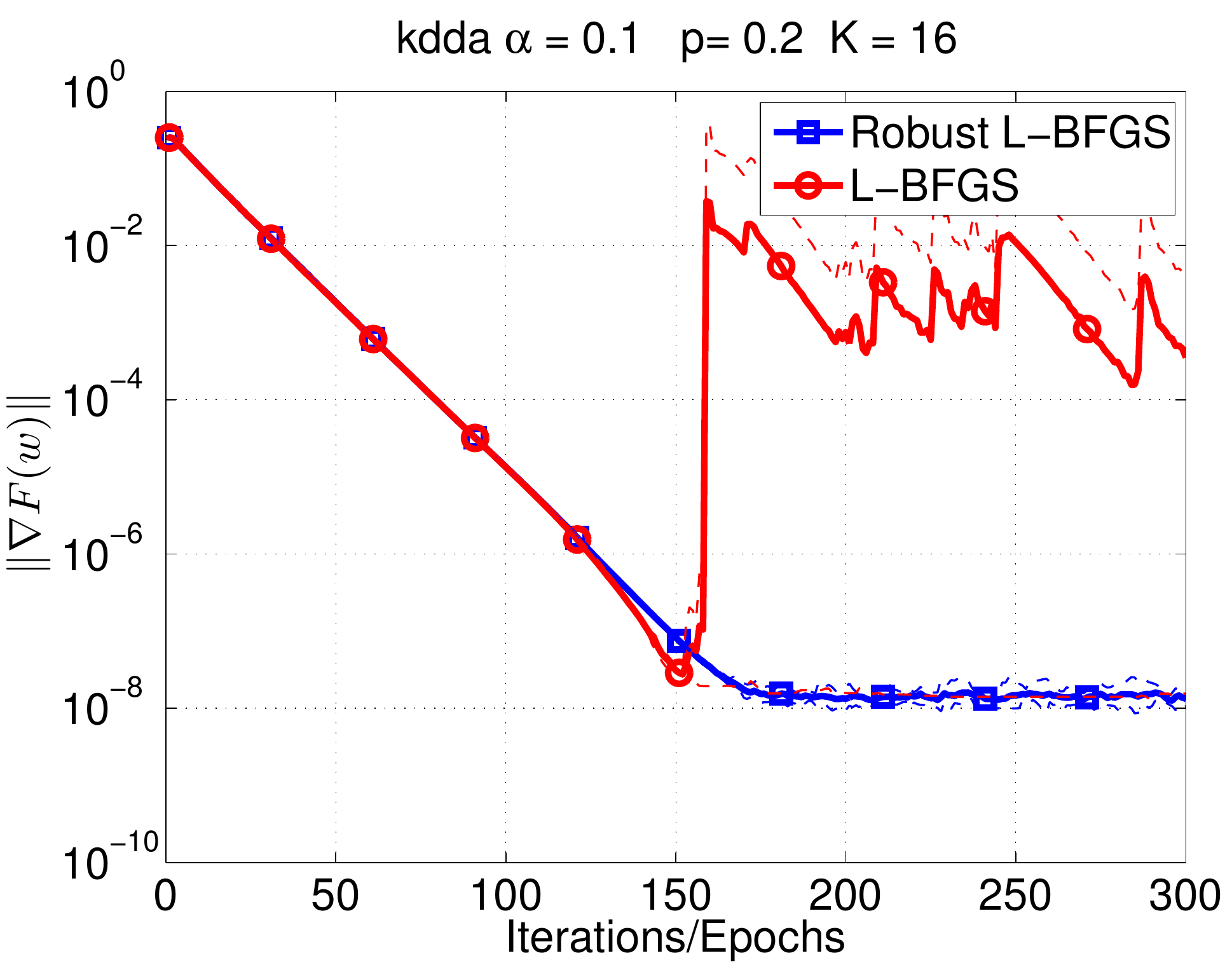}
\includegraphics[width=4.6cm]{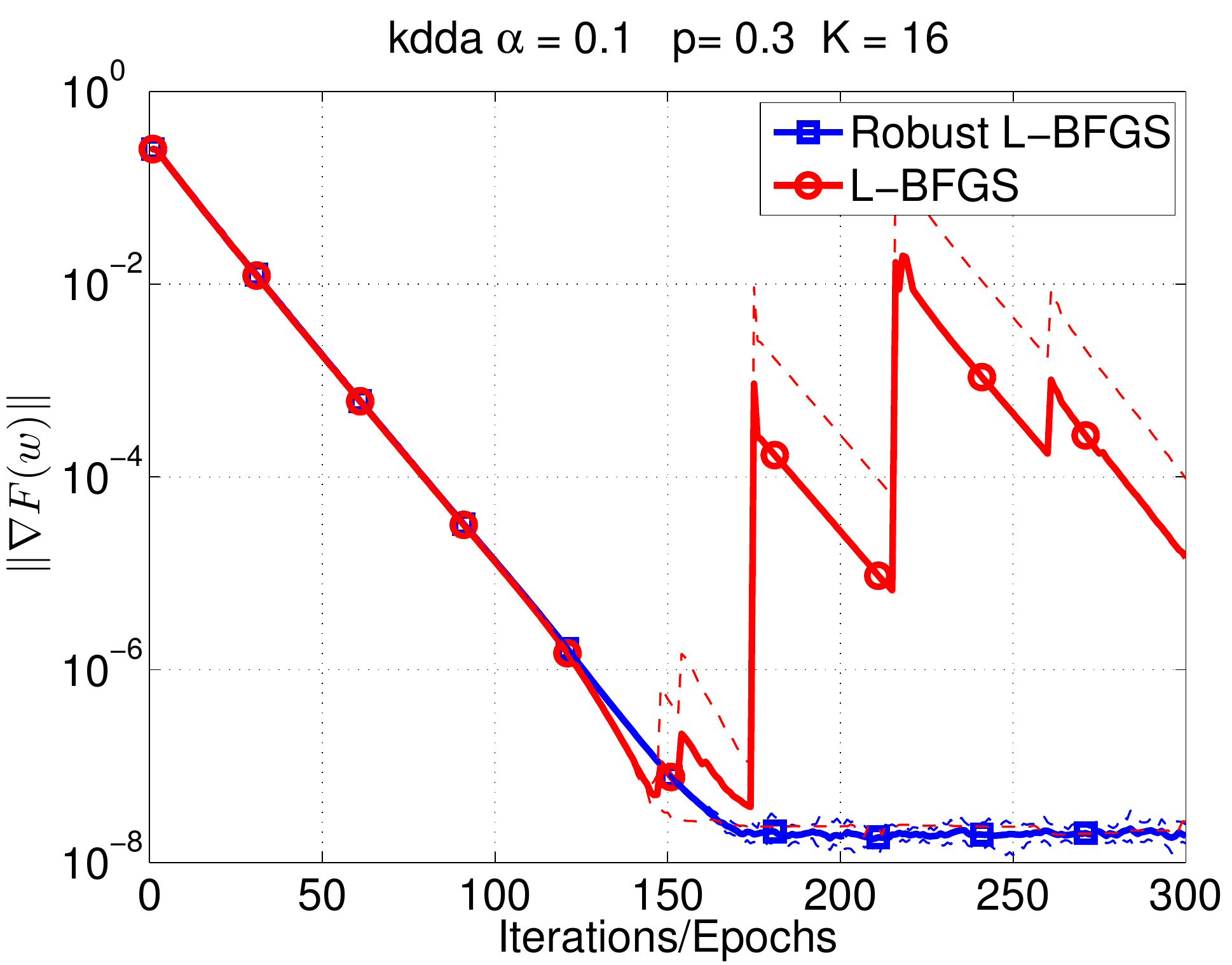}

\includegraphics[width=4.6cm]{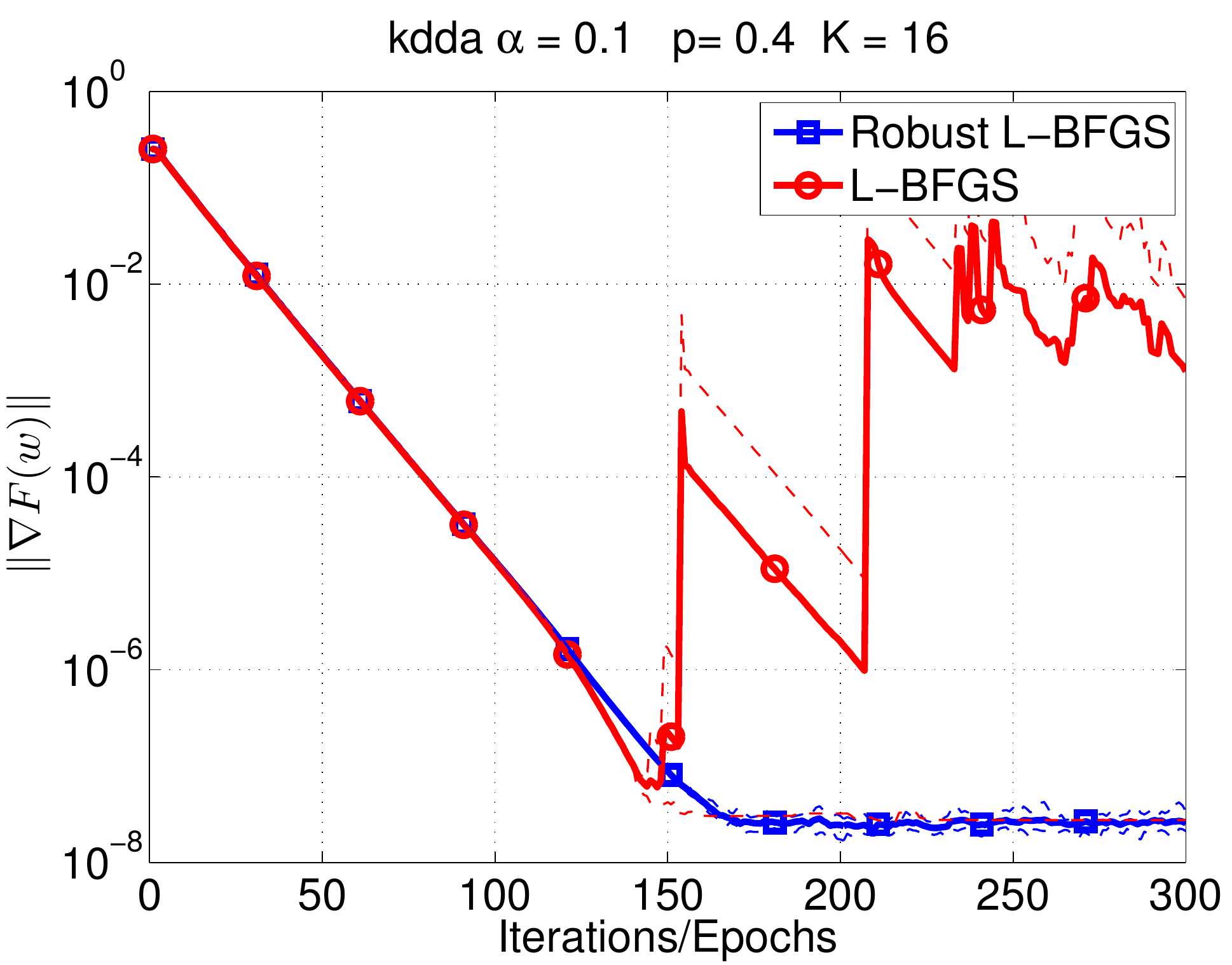}
\includegraphics[width=4.6cm]{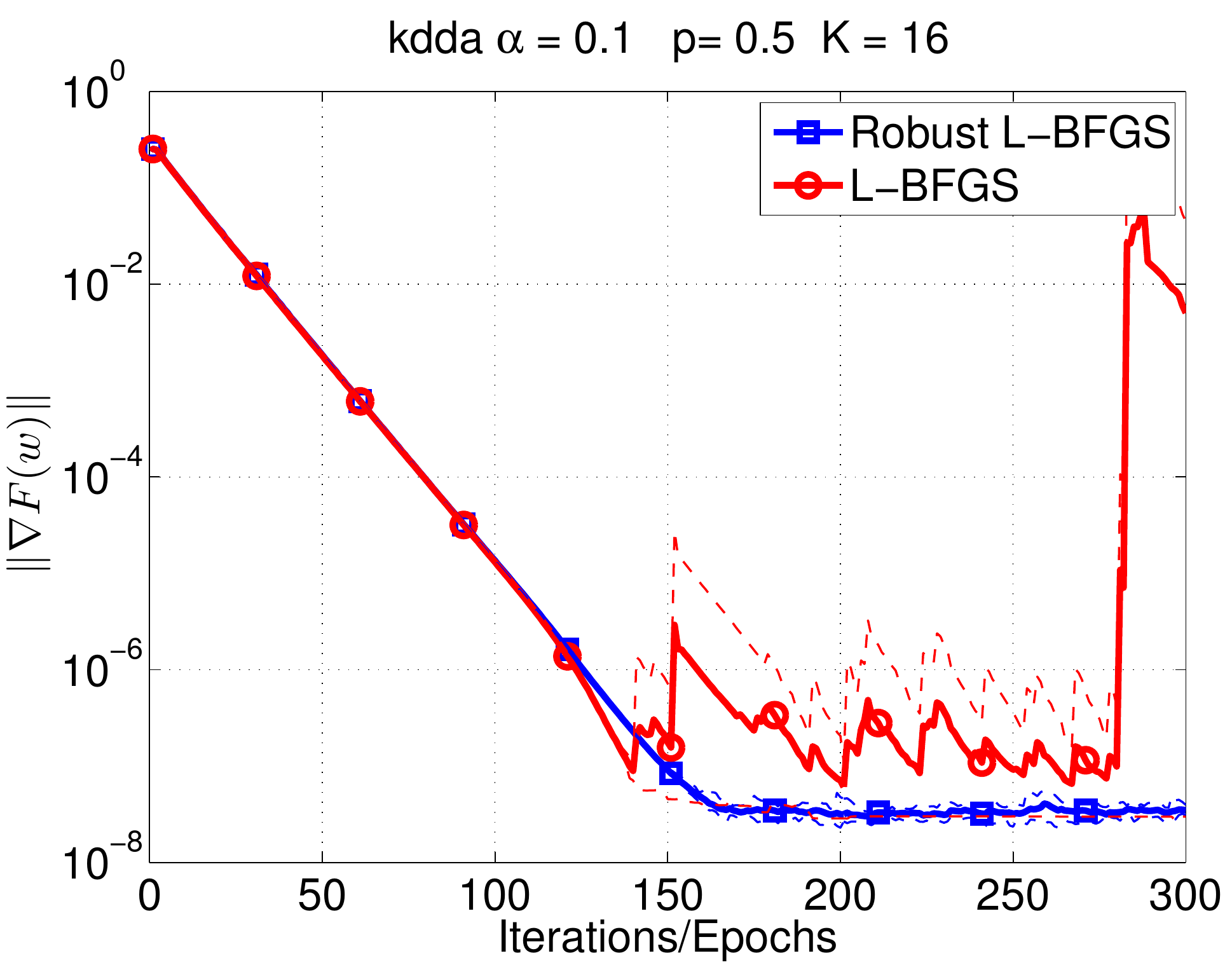}
 
\caption{\textbf{kdda dataset}. Comparison of Robust L-BFGS and L-BFGS in the presence of faults.
We used $\alpha=0.1$ and $p\in \{0.1, 0.2, 0.3, 0.4, 0.5\}$. Solid lines show average performance, and dashed lines show worst and best performance, over 10 runs (per algorithm). $K=16$ MPI processes.
}\label{fig:ft:kdda}
\end{figure}
\begin{figure}[h!]
\centering

\includegraphics[width=4.6cm]{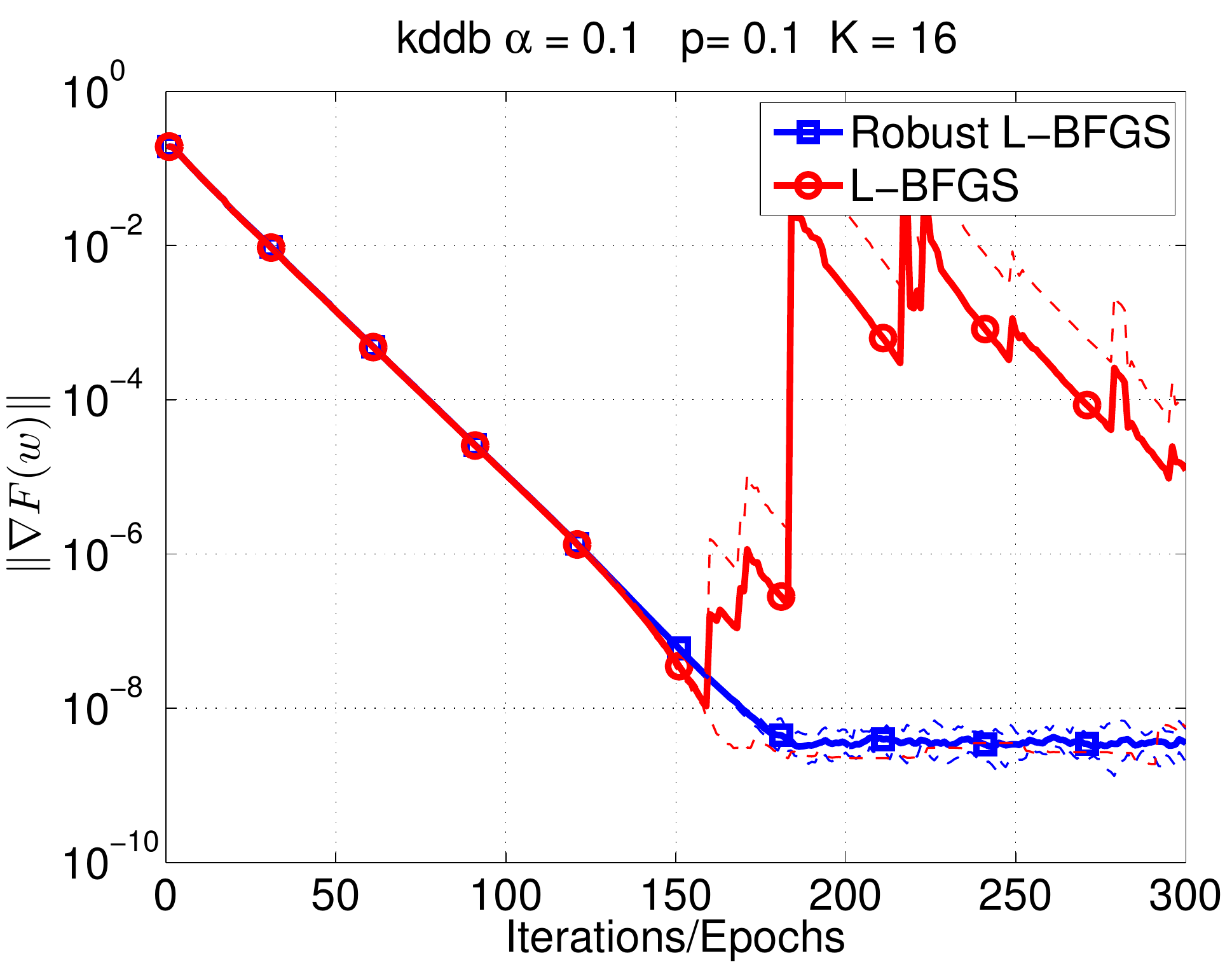}
\includegraphics[width=4.6cm]{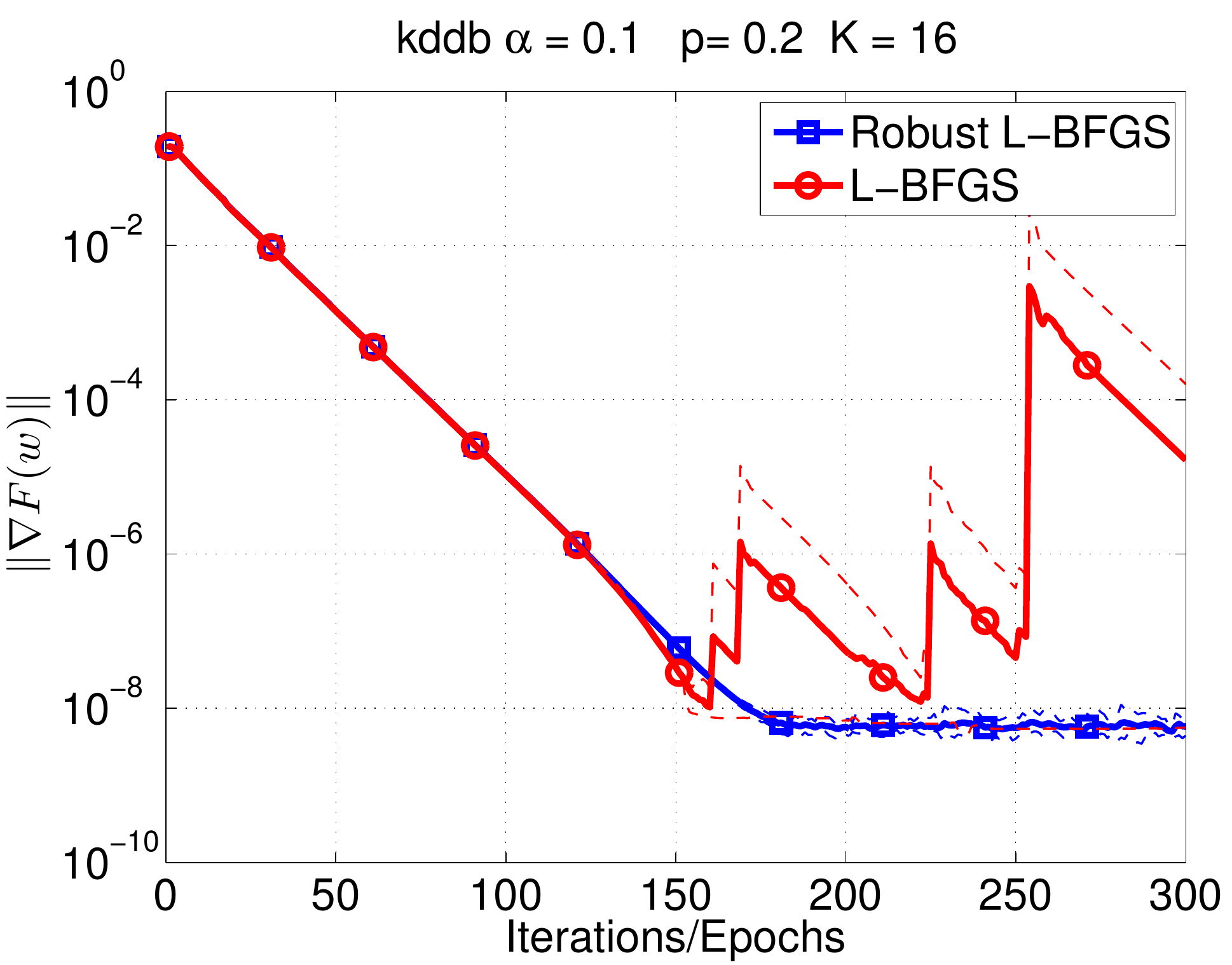}
\includegraphics[width=4.6cm]{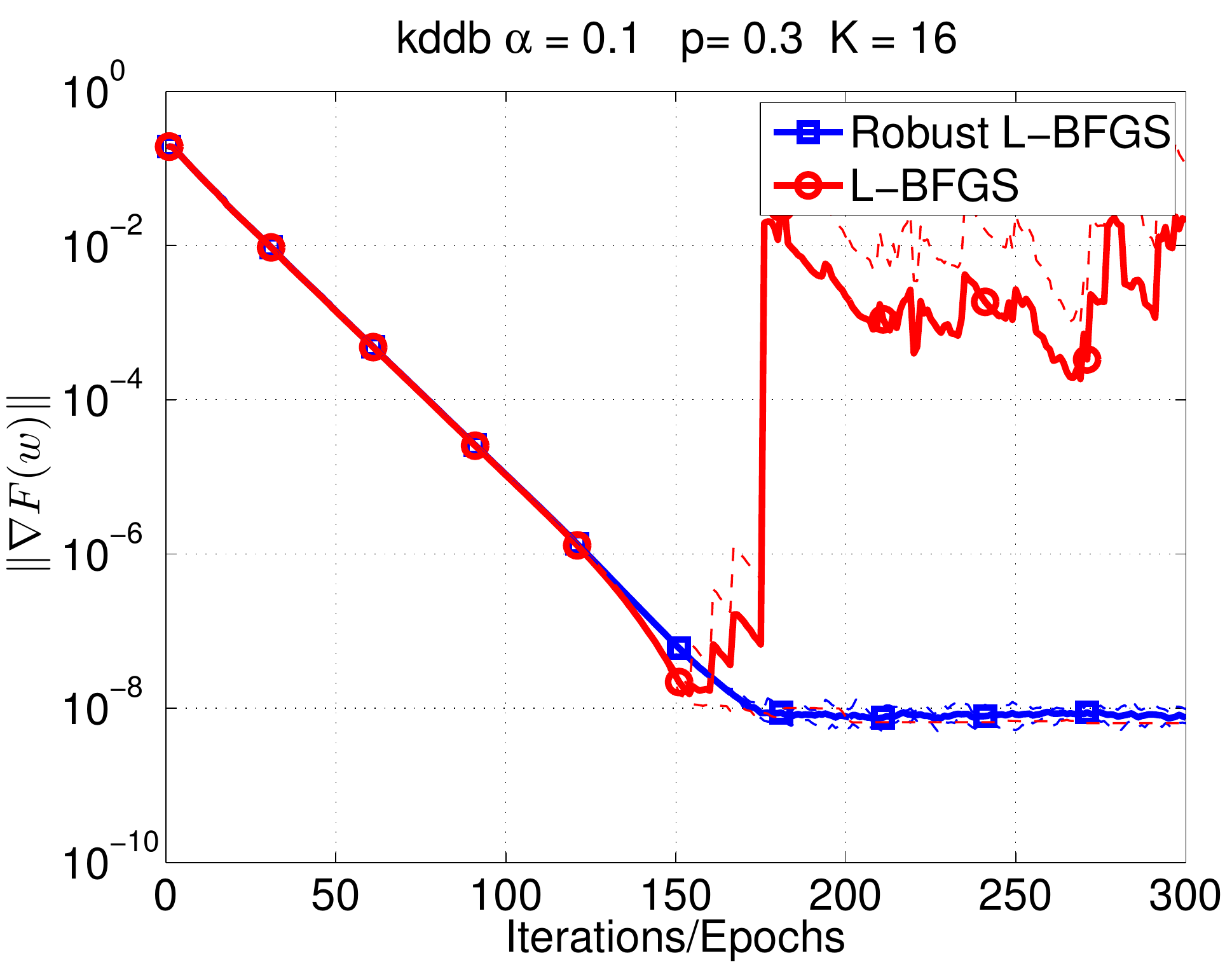}

\includegraphics[width=4.6cm]{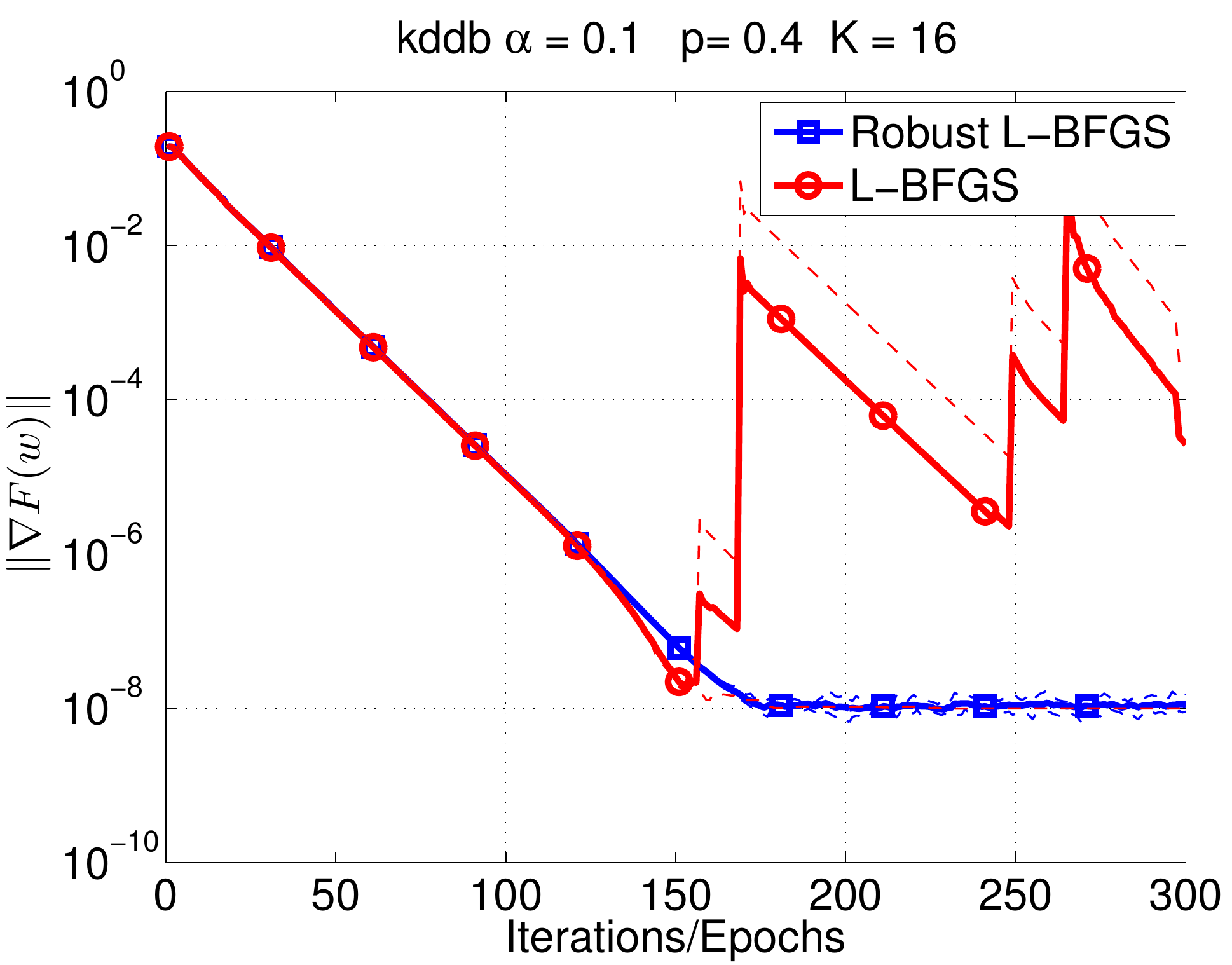}
\includegraphics[width=4.6cm]{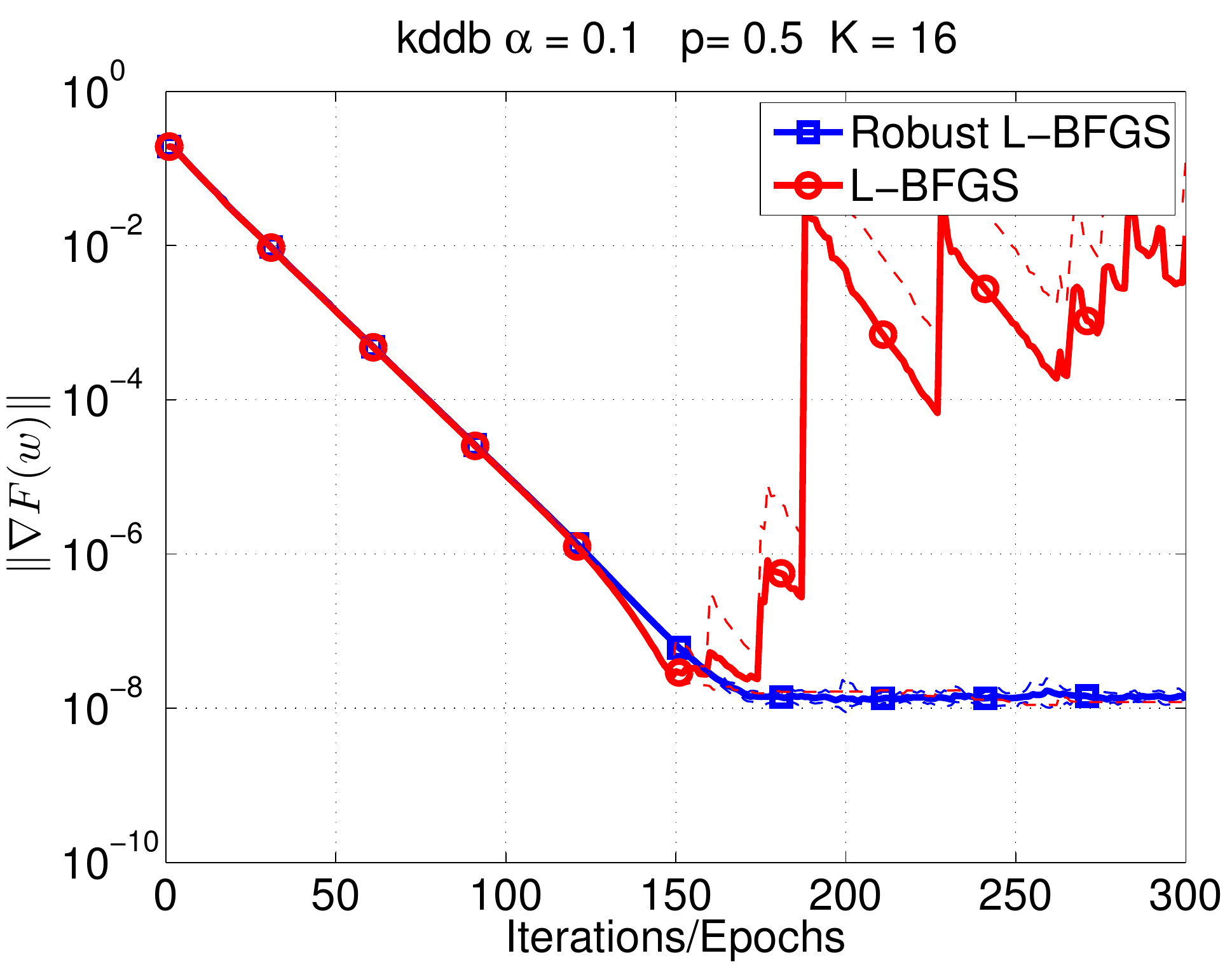}
 
\caption{\textbf{kddb dataset}. Comparison of Robust L-BFGS and L-BFGS in the presence of faults.
We used $\alpha=0.1$ and $p\in \{0.1, 0.2, 0.3, 0.4, 0.5\}$. Solid lines show average performance, and dashed lines show worst and best performance, over 10 runs (per algorithm). $K=16$ MPI processes.
}\label{fig:ft:kddb}
\end{figure}
\begin{figure}[h!]
\centering

\includegraphics[width=4.6cm]{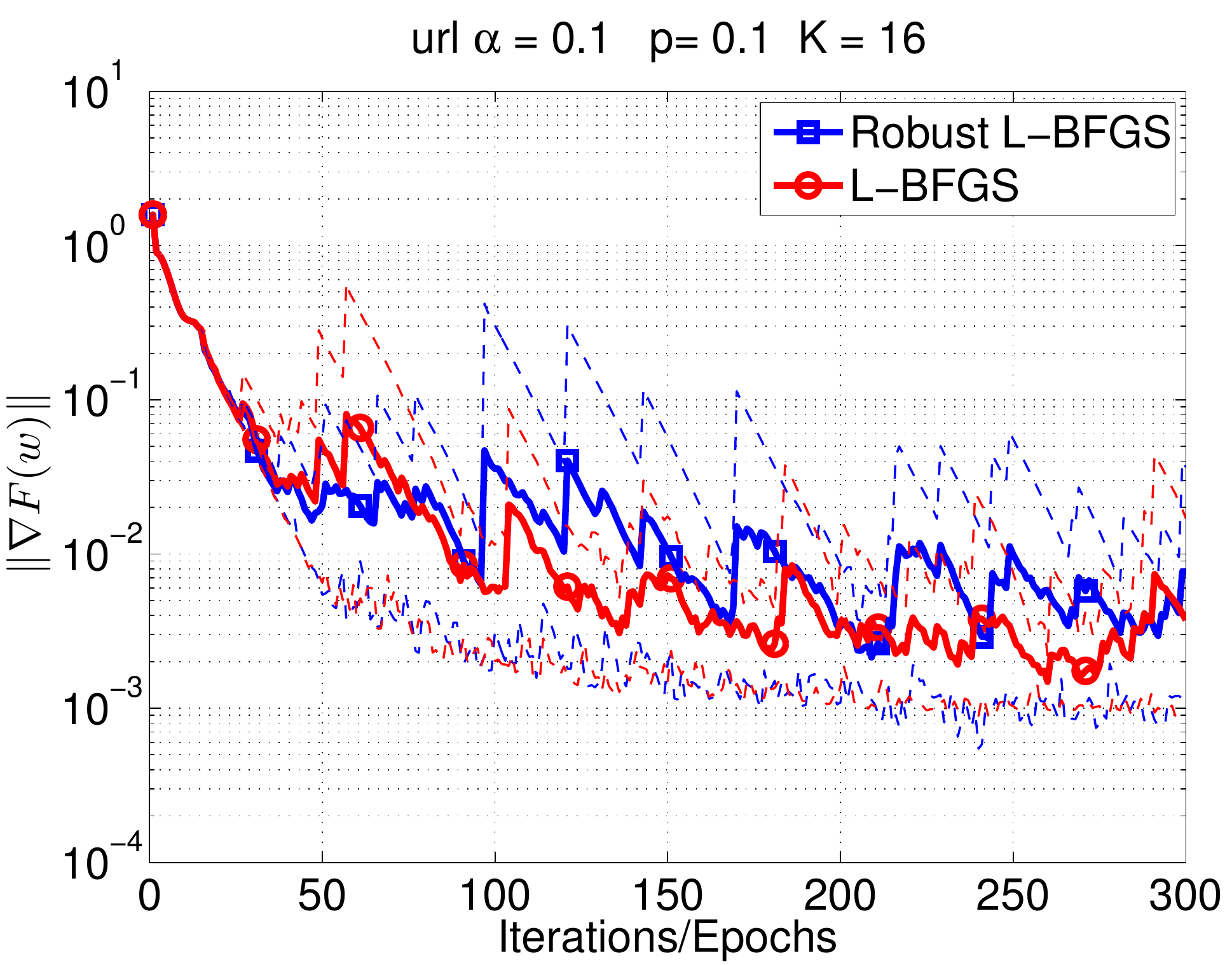}
\includegraphics[width=4.6cm]{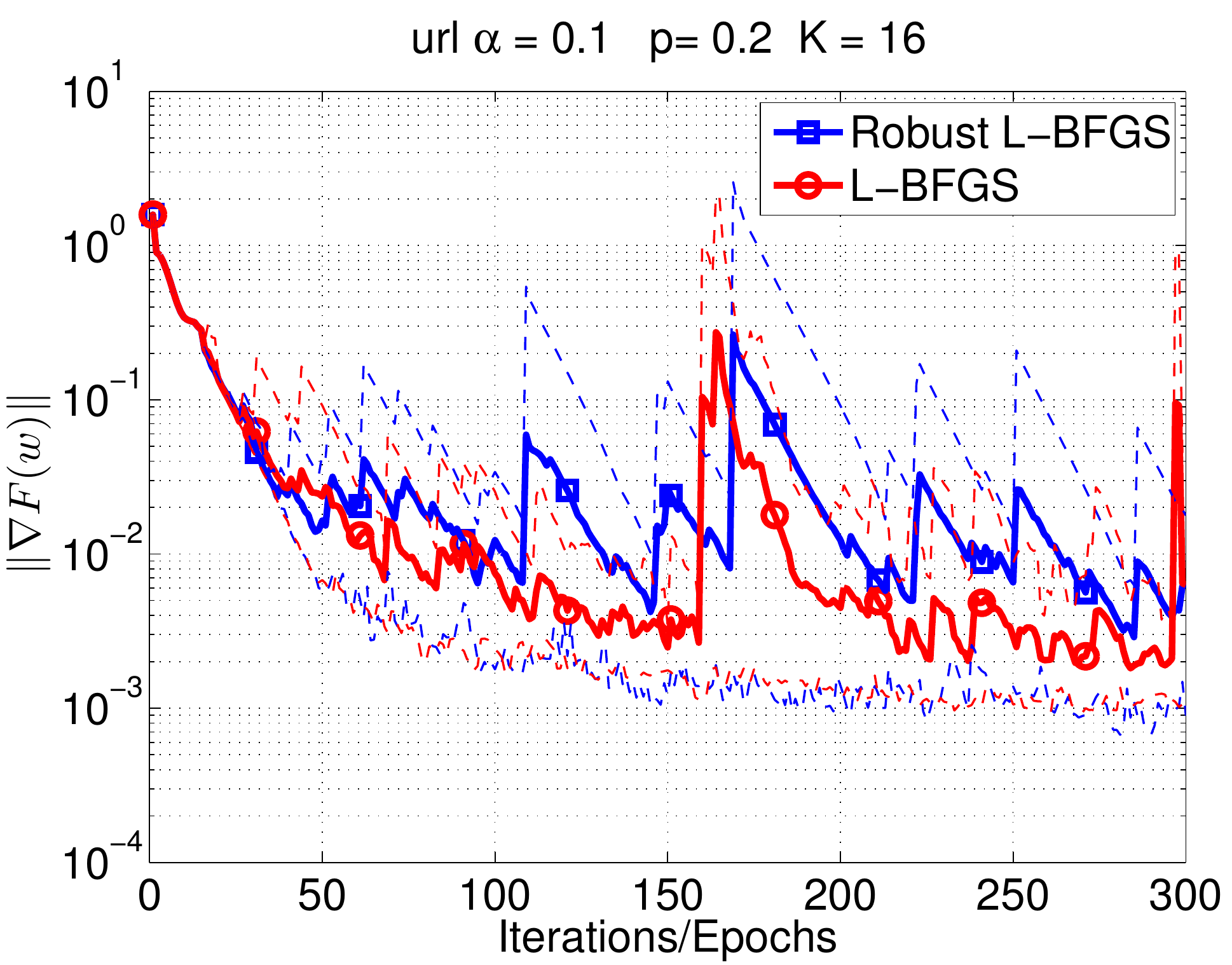}
\includegraphics[width=4.6cm]{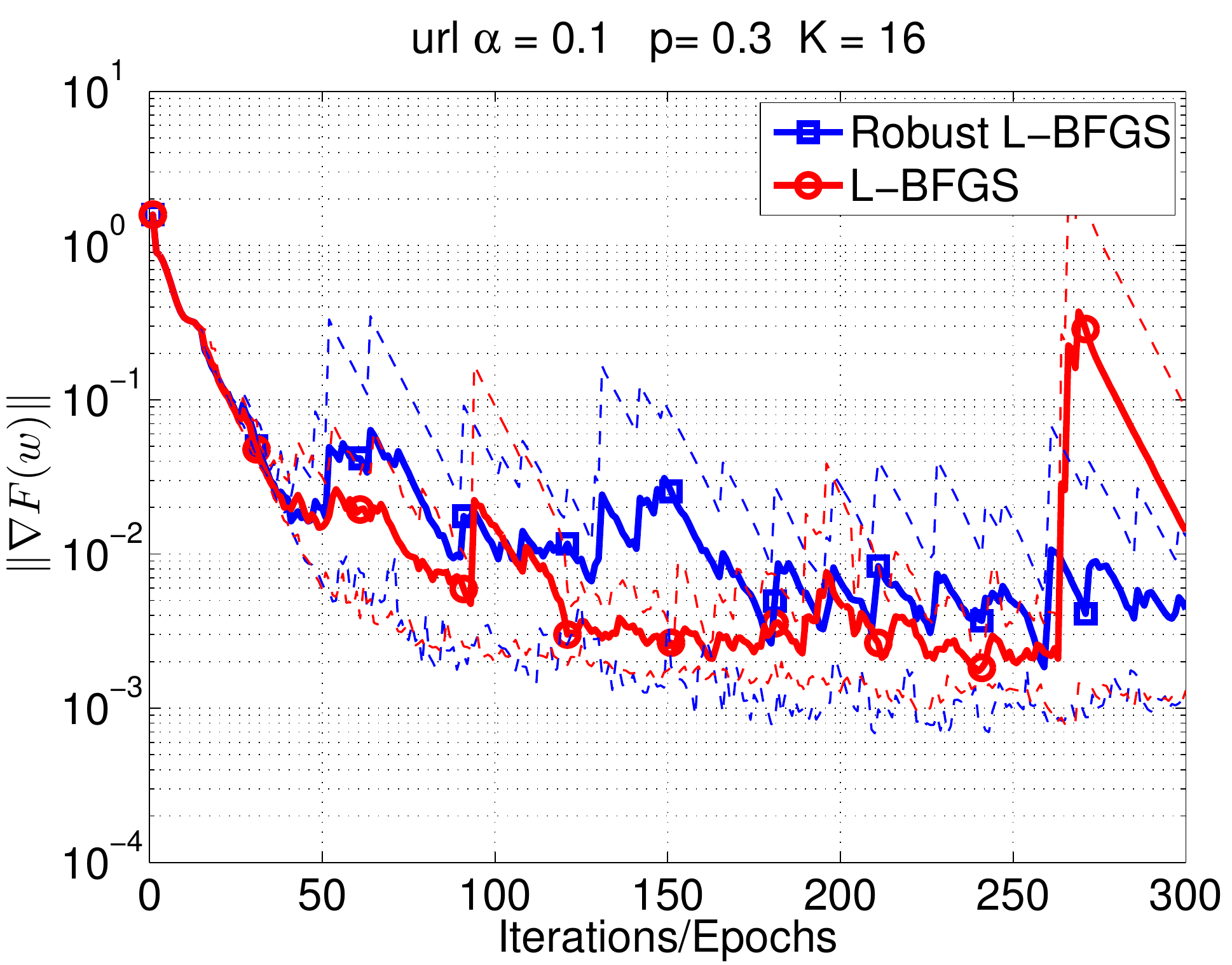}

\includegraphics[width=4.6cm]{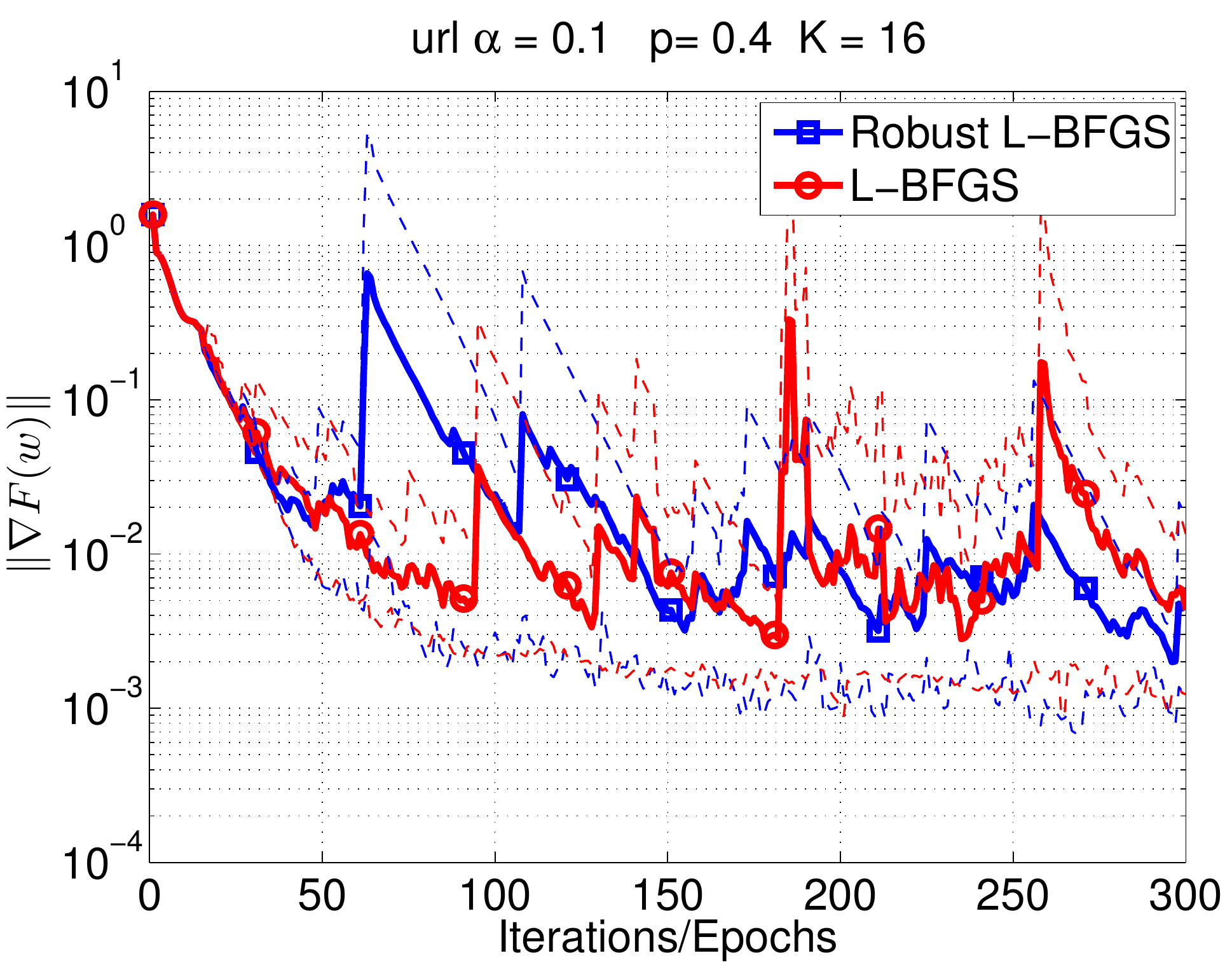}
\includegraphics[width=4.6cm]{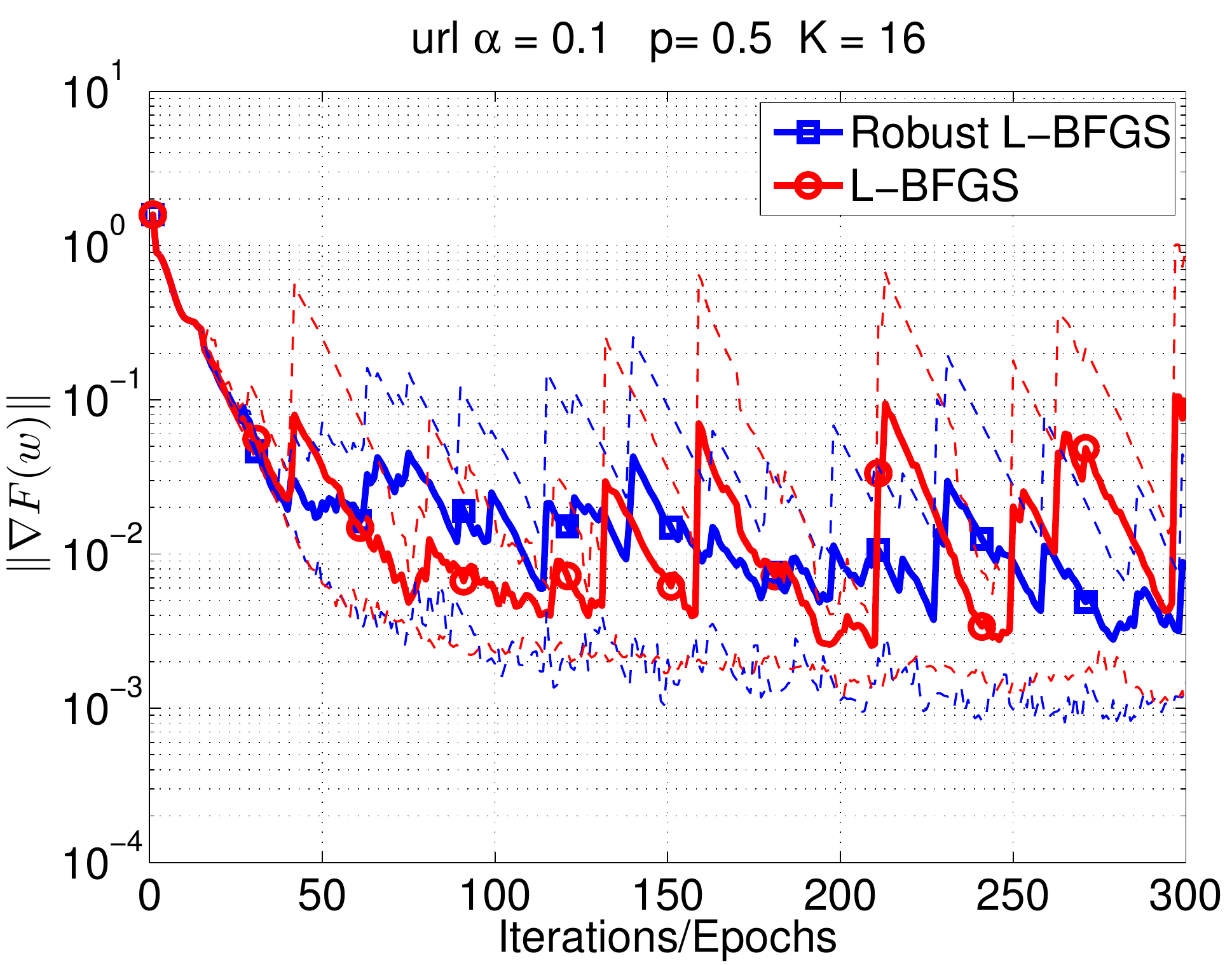}
 
\caption{\textbf{url dataset}. Comparison of Robust L-BFGS and L-BFGS in the presence of faults.
We used $\alpha=0.1$ and $p\in \{0.1, 0.2, 0.3, 0.4, 0.5\}$. Solid lines show average performance, and dashed lines show worst and best performance, over 10 runs (per algorithm). $K=16$ MPI processes.
}\label{fig:ft:url}
\end{figure}

In Figures  \ref{fig:ft:rcv}-\ref{fig:ft:url} 
we show a comparison of the proposed robust multi-batch L-BFGS method and the multi-batch L-BFGS method that does not enforce sample consistency (L-BFGS). In these experiments, $p$ denotes the probability that a single node (MPI process) will not return a gradient evaluated on local data within a given time budget. We illustrate the performance of the methods for $\alpha=0.1$ and $p\in \{0.1, 0.2, 0.3, 0.4, 0.5\}$. We observe that the robust implementation is not affected much by the failure probability $p$.

\clearpage
\section{Scaling of Robust Multi-Batch L-BFGS Implementation}
\label{sec:scaling_multi}

In this Section, we study the strong and weak scaling properties of the robust multi-batch L-BFGS method on an artificial dataset.
For various values of $r$ and $K$, we measure the time needed to compute a gradient (Gradient) and the time needed to compute and communicate the gradient (Gradient+C), as well as, the time needed to compute the L-BFGS direction (L-BFGS) and the associated communication overhead (L-BFGS+C).

\subsection{Strong Scaling}

\begin{figure}[h!]
  \centering

\includegraphics[width=8cm]{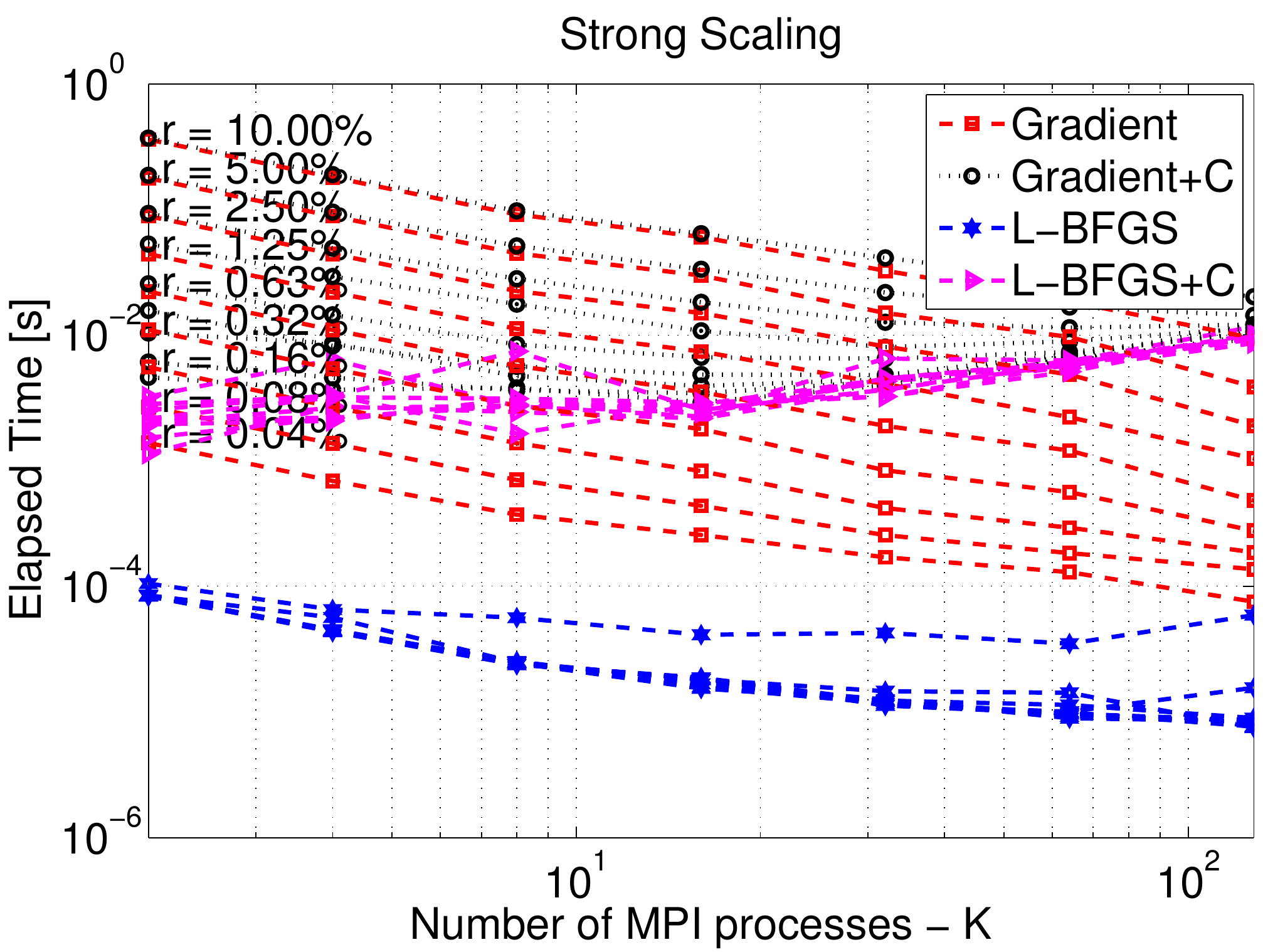}

  \caption{ Strong scaling of robust multi-batch L-BFGS on a problem with artificial data; 
  $n=10^7$ and $d=10^4$. Each sample has $160$ non-zero elements. $+C$ indicates that we include communication time to the gradient computation and L-BFGS update computation.
}  
\label{strongscaling}  
   \end{figure}

Figure \ref{strongscaling} depicts the strong scaling properties of our proposed algorithm. We generate a dataset with $n=10^7$ samples and $d=10^4$ dimensions, where each sample has 160 randomly chosen non-zero elements (dataset size 24GB). We run our code for different values of $r$ (different batch sizes $S_k$), with $K= 1, 2, \dots, 128$ number of MPI processes. 

One can observe that the compute time for the gradient and the L-BFGS direction decreases as $K$ is increased. However, when communication time is considered, the combined cost increases slightly as $K$ is increased. Notice that for large $K$, even when $r=10\%$ (i.e., $10\%$ of all samples processed in one iteration, $\sim$18MB of data), the amount of local work is not sufficient to overcome the communication cost.

\subsection{Weak Scaling - Fixed Problem Dimension, Increasing Data Size}

In order to illustrate the weak scaling properties of the algorithm, we generate a data-matrix $X \in R^{10^7	\times 10^4}$, and run it on a shared cluster with $K=1,2,4,8,\dots,128$ MPI processes. For a given number of MPI processes ($K$), each sample contains $10\cdot K$ non-zero elements. Effectively, the dimension of the problem is fixed, but sparsity of the data is decreased as more MPI processes are used. The size of the input data is 1.5 $\cdot K$ GB (i.e., 1.5GB per MPI process).

\begin{figure}[h!]
  \centering

\includegraphics[width=8cm]{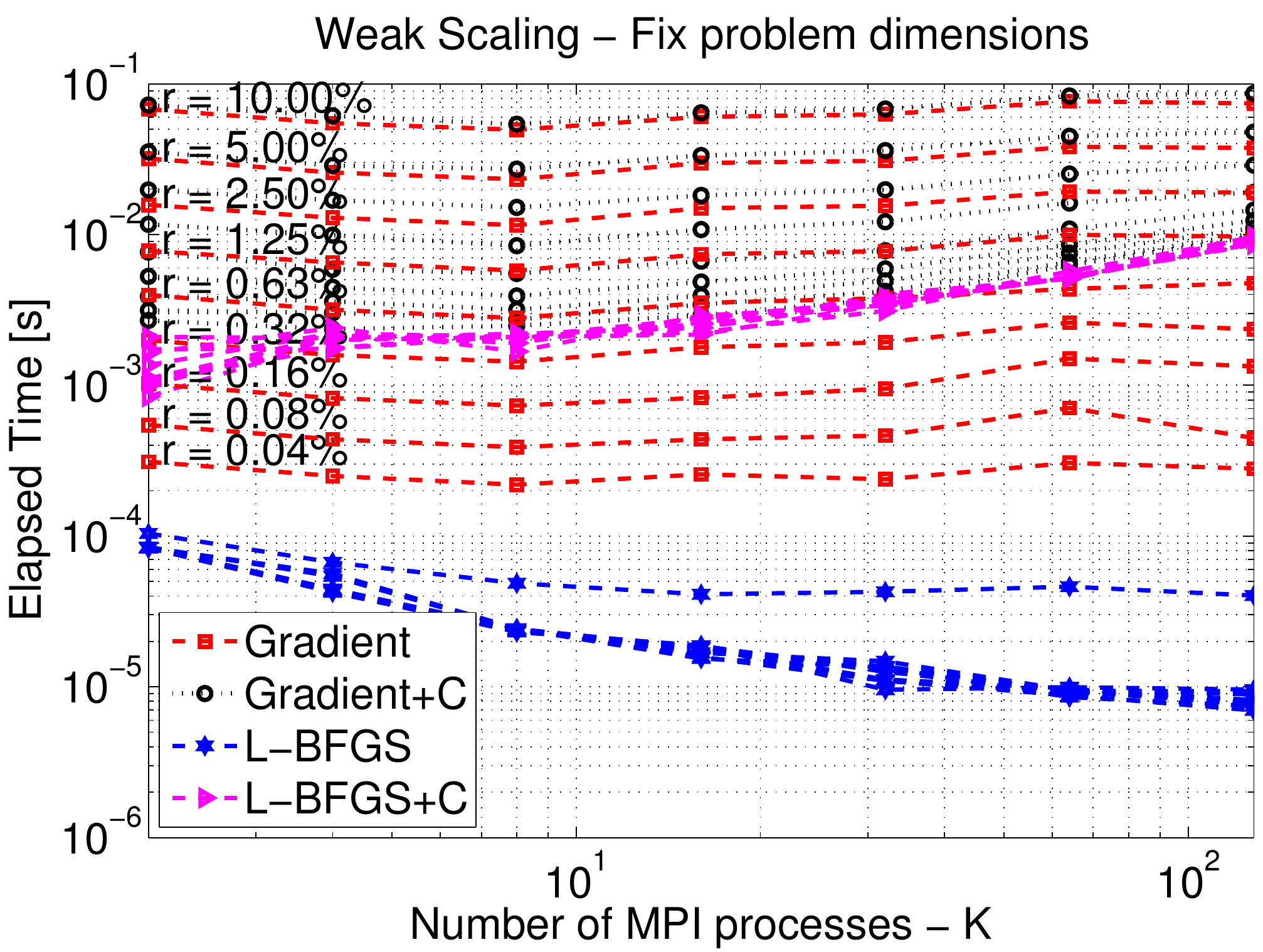}

  \caption{Weak scaling of robust multi-batch L-BFGS on a problem with artificial data;  
  $n=10^7$ and $d=10^4$. Each sample has $10\cdot K$ non-zero elements. $+C$ indicates that we also include communication time to the gradient computation and L-BFGS update computation.}  
\label{weakscaling2}  
  \end{figure}

The compute time for the gradient is almost constant, this is because the amount of work per MPI process (rank) is almost identical; see Figure \ref{weakscaling2}. On the other hand, because we are using a Vector-Free L-BFGS implementation \cite{chen2014large} for computing the L-BFGS direction, the amount of time needed for each node to compute the L-BFGS direction is decreasing as $K$ is increased. However, increasing $K$ does lead to larger communication overhead, which can be observed in Figure \ref{weakscaling2}. For $K=128$ (192GB of data) and $r=10\%$, almost 20GB of data are processed per iteration in less than 0.1 seconds, which implies that one epoch would take around 1 second. 

\subsection{Increasing Problem Dimension, Fixed Data Size and $K$}

\begin{figure}[h!]
  \centering

\includegraphics[width=8cm]{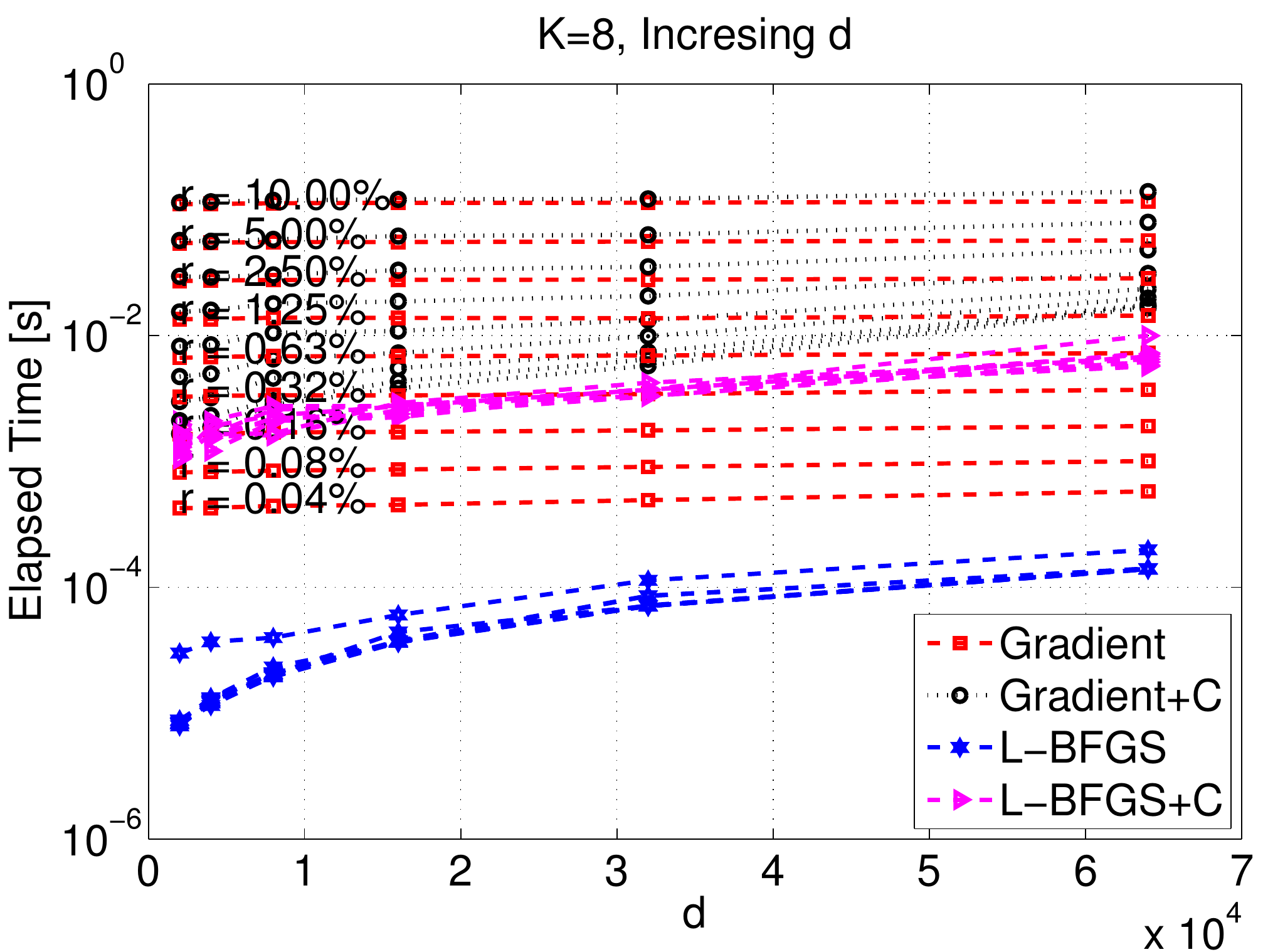}

  \caption{Scaling of robust multi-batch L-BFGS on a problem with artificial data; $n=10^7$ samples, with increasing $d$ and $K=8$ MPI processes. Each sample had 200 non-zero elements. $+C$ indicates that we also include communication time to the gradient computation and L-BFGS update computation.}  
\label{weakscaling}  
  
  \end{figure}

In this experiment, we investigate the effect of a change in the dimension $d$ of the problem on the performance of the algorithm. We fix the size of data ($29GB$) and the number of MPI processes ($K=8$). We generate data with $n=10^7$ samples, where each sample has 200 non-zero elements. Figure \ref{weakscaling} shows that increasing the dimension $d$ has a mild effect on the computation time of the gradient, while the effect on the time needed to compute the L-BFGS direction is more apparent. However, if communication time is taken into consideration, the time required for the gradient computation and the L-BFGS direction computation increase as $d$ is increased.

\clearpage

\end{document}